\documentclass[12pt, english, intlimits, fleqn, twoside, bibliography=totoc,numbers=noenddot]{scrbook}
\makeatletter
\usepackage[headings]{fullpage}
\usepackage{dsfont}
\usepackage{url}
\usepackage{mathtools}
\usepackage[T1]{fontenc}

\renewcommand{\geq}{\geqslant}
\renewcommand{\leq}{\leqslant}

\usepackage{amsmath}
\usepackage{amssymb}
\usepackage{amsfonts}
\usepackage{remreset}
\usepackage[english]{babel}

\usepackage{typearea}
\usepackage{enumerate}

\usepackage{xparse}

\usepackage[scale=0.95]{tgpagella}
\usepackage[scale=0.92]{tgheros}
\usepackage[scaled=0.83]{beramono}
\usepackage{mathpazo}
\usepackage{eucal}

\usepackage[thmmarks,amsmath]{ntheorem}

\usepackage[all,matrix,arrow,arc,curve,color,frame]{xy}
\usepackage{enumitem}
\setlist{noitemsep}

\newtheorem{Theorem}{Theorem}[section]

\newtheorem{Corollary}[Theorem]{Corollary}
\newtheorem{Lemma}[Theorem]{Lemma}
\newtheorem{Proposition}[Theorem]{Proposition}

\theoremstyle{plain}
\theoremsymbol{\ensuremath{\triangle}}

\newtheorem{Definition}[Theorem]{Definition}
\newtheorem{Hypothesis}[Theorem]{Hypothesis}

\theoremstyle{plain}
\theoremsymbol{\scalebox{1.5}[1.5]{$\diamond$}}
\theorembodyfont{\normalfont\upshape}
\newtheorem{Remark}[Theorem]{Remark}
\newtheorem{Example}[Theorem]{Example}

\theoremstyle{nonumberplain}
\theoremheaderfont{\normalfont\itshape}
\theorembodyfont{\normalfont\upshape}
\theoremsymbol{\ensuremath{\square}}
\newtheorem{Proof}{Proof}

\DeclareMathOperator{\spt}{spt}

\DeclareMathOperator{\TextRe}{Re}

\renewcommand{\Re}{\TextRe}

\DeclareMathOperator{\lin}{lin}

\DeclareMathOperator{\dom}{dom}
\DeclareMathOperator{\kar}{ker}
\DeclareMathOperator{\ran}{ran}

\DeclareMathOperator{\id}{id}

\newcommand{\N}{\mathbb{N}}

\newcommand{\R}{\mathbb{R}}

\newcommand{\s}{\mathcal}

\newcommand{\eps}{\varepsilon}

\DeclareMathAccent{\Circ}{\mathalpha}{operators}{"17}

\DeclareMathOperator{\esssup}{ess\ \! sup}

\renewcommand{\1}{\mathds{1}}
\renewcommand{\epsilon}{\varepsilon}

\newcommand{\dd}{\ \mathrm{d}}

\DeclareMathOperator{\grad}{grad}

\DeclareMathOperator{\dive}{div}

\DeclareMathOperator{\curl}{\operatorname{curl}}

\newcommand{\abs}[1]{\left\lvert#1\right\rvert}
\newcommand{\Abs}[1]{\left\lVert#1\right\rVert}

\newcommand{\ben}{\begin{enumerate}[(i)]}
\newcommand{\een}{\end{enumerate}}

\newcommand{\mysection}[3]{%
    \section[#1]{#2{\newline\normalfont\small\textit{#3}}}
}

\usepackage{tikz}

\usepackage{url}
\usepackage[intoc]{nomencl}
\usepackage{babel}
\usepackage[bf]{caption}
\usepackage{subfigure}
\makeatother

\typearea[1cm]{15}

\begin{document}

\thispagestyle{empty}
\begin{center}
\LARGE{

\textsc{On the continuous dependence on the coefficients of evolutionary equations}}

\medskip

\Large{

\vfill

Habilitationsschrift 

\medskip
}

\vfill

{\large

\emph{vorgelegt der}}\\ 
\medskip
\textbf{Fakult\"at der Mathematik und Naturwissenschaften \\ der Technischen Universit\"at Dresden}
\\\medskip

 \emph{von}
 \\
{\Large
Dr.\ rer.\ nat.\ Marcus Waurick
}\\
{\large \emph{geboren am}}\\
{\Large
10.09.1984
}\\
{\large
\emph{in}
}\\
{\Large
Zittau
}\\
\medskip\medskip
\medskip\medskip \medskip\medskip \medskip\medskip\medskip\medskip
{\large \emph{Submitted:} January 29, 2016}\\
{\Large \medskip\medskip\medskip\medskip\medskip}
{\large \emph{Defended:} June 10, 2016}\\
{\Large \medskip\medskip\medskip\medskip\medskip}

\end{center}






\chapter*{Preface}

When mathematically modeling physical phenomena, the mathematical models almost always include certain unknown parameters. These parameters are usually to be determined via experiments. The experiments are subject to variability and chance. Having modeled physical phenomena mathematically with appropriate parameters determined from the experiments, one tries to predict the outcome of certain physical processes. For this prediction a computer may be used. Since machine precision is of limited accuracy, the results of these computations are subject to chance and variability due to rounding errors. 

With these difficulties in mind is the mathematical model derived correct? Do small variations of the set of parameters only result in small variations of the outcome?  Are the computations carried out by the computer meaningful? In summary, what is the effect of small variations in the (parameters of the mathematical) model to the solution?

The aim of this thesis is to provide a mathematical framework for addressing the latter question. Moreover, discussing several notions of ``smallness'', we provide results of the type: ``small variations in the parameters lead to small variations of the solution.'' The focus is on evolutionary equations, that is, ordinary or partial differential equations involving the time derivative.

\newpage
{
\chapter*{Acknowledgments} 

The thesis in its present form would not have been possible without the groundbreaking work of my `Doktorvater', my co-author and friend Prof Rainer Picard. It is always great fun chatting with him and to discuss anything in both mathematical and non-mathematical areas. He is a role model in kindness, honesty and integrity. 

I want to thank Prof Ralph Chill for giving me the opportunity to do research in many places of the world. He provided a beautiful atmosphere  and it was and is always possible to discuss anything with him. In fact, I want to thank a great deal of researchers. Among them my friend Dr Sascha Trostorff, who read a preliminary version of this thesis and spotted flaws, inaccuracies, misprints and, thus, has helped to improve this text in various places. But there are plenty more maths scholars I'd like to thank for making the work at the TU Dresden such a wonderful time: My first analysis teacher Prof J\"urgen Voigt and his ``real analogue'' Dr Hendrik Vogt, Prof Werner Timmermann, Prof Stefan Siegmund, Jun-Prof Friedrich Martin Schneider, Dr Anke Kalauch, Dr Eva Fasangova, Dr Norbert Koksch, Dr Daniel Borchmann, Tom Hanika, Jeremias Epperlein, and Sebastian Mildner. Thanks also go to the administrative staff Monika Gaede-Samat, Christine Heinke and Karola Schreiter for solving plenty of non-mathematical problems.

I also like to thank Dr Kirill Cherednichenko at the University of Bath for useful discussions and comments on the present manuscript.

Ich m\"ochte mich auch bei meinen Eltern, Peter und Gabriele, als auch bei meinem Bruder Steffen, seiner Frau Bernadette, deren Kindern Hannah and Julius bedanken. Au\ss erdem danke (Warum gibt es in der deutschen Sprachen so wenig Alternativen für das Wort 'Danke'?) ich allen weiteren Mitgliedern meiner ganzen Familie, meinen Freunden, insbesondere, TonKat, Frantje, JoIsi, Martin, Daniel und Sophia, Dario und Judith, Hauke, Julio, Schek, Mille, Schnappo, Herrn Handschrift, Kleopatra Maximowna, dem Engler (ich habe hier sicher viel zu viele vergessen; Sorry! Die Nennungsreihenfolge ist zufällig.) und Susn. Ihr alle habt mich zu dem gemacht, was ich bin; das habt ihr nun davon. Mathematik macht mir zwar Spa\ss, ihr habt aber alle dazu beigetragen, dass mir das Leben Spa\ss\ macht! Damit niemand, den ich pers\"onlich kenne, in dieser Aufz\"ahlung zuletzt genannt werden muss, danke ich schlie\ss lich noch George Lucas und Gene Roddenberry f\"ur die Schaffung zweier unverwechselbarer(!) Science-Fiction Abenteuer.

Marcus `moppi' Waurick,

Bath, January 2016.
}

%
%
%
 \tableofcontents


%



%
%
%
%
%

\cleardoublestandardpage

\addchap{Introduction}

In 2009 Rainer Picard \cite{PicPhy} developed a breakthrough functional analytic approach, taking advantage of a common structural property of many linear partial differential equations modeling dynamical processes of mathematical physics. It is well-known that many equations can be written in the form of a certain `balance law'
\begin{equation}\label{eq:inm1}
   \partial_t V + \s A U =F,
\end{equation}
which relates the time derivative of an unknown quantity $V$ and an operator $\s A$ comprising the spatial derivatives of another unknown $U$ to a given external source term $F$. Equation \eqref{eq:inm1} is usually complemented by a `material law' relating $U$ and $V$, which in the present thesis we assume to be a linear operator $\s M$, as
\[
   \s M U = V.
\]
The resulting equation, where we substituted the material law into \eqref{eq:inm1}, reads
\begin{equation}\label{eq:in0}
   \partial_t \s M U+ \s A U= F.
\end{equation}
In this introduction, we shall refer to equations of the form \eqref{eq:in0} as \emph{evolutionary equations} and to $\s M$ as the \emph{material law}. Various aspects of equations of the form \eqref{eq:in0} have been studied extensively from the mathematical perspective using a large variety of techniques. The pioneering technique of \cite{PicPhy} is based on a new structural observation about equations \eqref{eq:in0} and an associated functional analytic framework in a special Hilbert space setting. It allows one to prove existence and uniqueness of solutions $U$ as well as their continuous dependence on $F$ for a wide class of $\s M$ and $\s A$, which was not possible within the previously existing approaches. More precisely, Picard provided a space-time Hilbert space in \cite{PicPhy} such that the operator $\s S(\s M)\coloneqq (\partial_t\s M+\s A)^{-1}$ is well-defined and continuous. Moreover, in the particular setting in \cite{PicPhy}, $\s S(\s M)$ is also shown to be causal, that is, $\s S(\s M)F$ vanishes up to time $t$ if $F$ does. The Hilbert space setting developed in combination with the possibly provocative philosophy that any ``reasonable'' linear time-dependent problem in mathematical physics should be an \emph{evolutionary equation}, that is, an equation as described in \cite{PicPhy}, leads to important structural insights and well-posedness theorems for problems that emerge in diverse applied fields. The standard evolution equations in mathematical physics, that is, Maxwell's equations, the wave equation, the elasticity equations, and the heat equation fit into the aforementioned class (\cite[Section 6]{Picard}). Possibly degenerate cases like the eddy-current approximation for Maxwell's equations or problems with change of type ranging from elliptic to parabolic to hyperbolic on different space-time regions are evolutionary as well, see \cite[Examples 2.7 and Examples 2.43]{Waurick2014PROC_Survey}. The framework has found to be useful for problems in control theory with unbounded control and observation operators \cite{Waurick2014IMAJMCI_ComprehContr,Waurick2014OAM_BdyCtrSys,Waurick2013REISSIG_Control}. As it was already pointed out in \cite[Section 3.5]{PicPhy}, the structural perspective developed is particularly useful, when discussing coupled phenomena in the light of so-called `multi-physics' problems, see also \cite[Sections 4.3 and 4.4]{PicPhy} for a brief account on thermoelasticity and piezo-electro-magnetism (cf.~also \cite{Mulholland2015}). The coupling of several equations of elastic type leads to well-posedness results in problems with micromorphic media, \cite[Section 3]{Waurick2014ZAMM_Micromorph}. A combination of the heat equation and the elastic equations form the description of an elastic material, that changes its elastic properties upon thermal excitation. A class of models of thermoelasticity is discussed in \cite{Waurick2014MMS_Thermo,Mukhopadyay2015}; several equations of Maxwell type lead to a description of linearized versions of Maxwell-Dirac systems or the equations of gravito-electro-magnetism, \cite{Waurick2014MMAS_MaxGrav}.

The framework of evolutionary equations has also natural applications to problems with memory involving integral expressions in $\s M$. Typical examples are fractional derivatives in time (\cite{Waurick2014MMAS_Frac,Waurick2014SIAM_HomFrac}) or integro-differential equations involving (other) con\-vo\-lution-in-time type operators (\cite{Trostorff2015}).

Once the unique existence of $U$ in \eqref{eq:in0} is established, it is natural to address questions like energy conservation or (exponential) stability. For the former we refer to \cite{Waurick2013REISSIG_Control} and for the latter we refer to \cite[Theorem 3.2]{Trostorff2014PAMM}, see also \cite{Trostorff2013a,Trostorff2015a} where criteria for the respective properties are given.

Of course, in applications, one is interested in actually computing the solution to a given problem numerically. In order to reduce computational costs when treating heterogeneous media there is a need for simplification of the constitutive relations. In particular, if the coefficients describing the material properties are highly oscillatory, computational costs for computing the respective solutions might even exceed the capabilities of modern computers. In the late 1960's, Spagnolo \cite{Spagnolo1967,Spagnolo1968a} mathematically approached this problem. After that many other researchers devoted a great deal of effort into the development of the newly founded mathematical theory of homogenization. We refer to some standard references \cite{CioDon,BenLiPap,TarIntro,Zhikov1983,Gcon1} and the references given there for a more detailed account on homogenization theory. With the structural perspective of evolutionary equations as in \eqref{eq:in0} for the particular case of time-shift invariant operators in mind, the author of the present manuscript considered homogenization problems for mathematical physics in \cite{Waurick2011}. Further development was achieved for ordinary (delay) differential equations in \cite{Waurick2012MMAS_ODES,Waurick2014JAA_G}. A contribution for partial differential equations can be found in \cite{Waurick2012_HowFar,Waurick2014SIAM_HomFrac,Waurick2013AA_Hom}. In the studies mentioned, homogenization theory is viewed as a certain property of the solution operator $\s S(\s M)$ associated to \eqref{eq:in0}: Is $\s M\mapsto \s S(\s M)$ continuous under the \emph{weak operator topology}? In essence, the author's work in this line of research reduces to proving criteria, where $\s M\mapsto \s S(\s M)$ indeed satisfies the continuity property mentioned. 

The present thesis addresses the continuous dependence of $\s M\mapsto \s S(\s M)$ under various topologies for a particular class of material laws. We focus on three possible topologies the set of material laws may be endowed with: The norm topology as well as the strong and weak operator topologies. 

We comment on the set of material laws as follows. In non-autonomous problems, the material law $\s M$ does not satisfy time-shift invariance. In \cite{Waurick2013JEE_Non} with minor extensions in \cite{Waurick2014MMAS_Non} we developed an adapted solution theory for problems as in \eqref{eq:in0}, which are not time-shift invariant. The class of admissible $\s M$ that lead to a solution theory for \eqref{eq:in0} has been introduced in \cite{Waurick2014JAA_G,Waurick2014MMAS_Non}: \emph{The class of evolutionary mappings}. This class will be the central object to study in the bulk of this manuscript. Inspired by a similar notion introduced in \cite[Definition 3.1.14]{Picard}, we develop a theory of evolutionary mappings. To the best of the author knowledge this has not been done before.

Chapter \ref{ch:TSO} is devoted to introducing the time derivative in exponentially weighted $L^2$ spaces. Moreover, we will present a well-known representation theorem, which characterizes operator-valued functions of the time derivative introduced. This representation theorem will render the relationship of time-shift invariant operators to evolutionary mappings. In Chapter \ref{ch:EMC} we introduce the central concept of this exposition, the notion of evolutionary mappings. This chapter also contains some preliminary results particularly useful in the forthcoming chapters. We provide a solution theory for linear ordinary differential equations (in infinite-dimensional state spaces) and abstract partial differential equations of the form \eqref{eq:in0} in Chapter \ref{ch:ST}. The method of proof for the solution theory just mentioned is similar to the one used in \cite{Waurick2013JEE_Non} or \cite{Waurick2014MMAS_Non}. However, as the focus is on (causal) evolutionary mappings, we will be able to derive more properties of the solution operator to the extent that we show that the solution operator $\s S$ is causal and evolutionary itself. We introduce the topologies on evolutionary mappings in Chapter \ref{Chapter:ODE}. In a first step towards applications, we provide continuous dependence results of $\s M\mapsto \s S(\s M)=(\partial_t\s M+\s A)^{-1}$ under the topologies introduced for $\s A=0$. A natural example is the Drude--Born--Fedorov model in the theory of electro-magnetism, which we shall treat as an application of the continuous dependence results developed in Chapter \ref{Chapter:ODE}. The corresponding results for partial differential equations are provided in Chapter \ref{ch:pde}. We will apply the respective results to the eddy-current approximation in electro-magnetic theory, where a hyperbolic equation is approximated by a parabolic one. Further, we provide an application to non-autonomous thermodynamics. The last application concerns homogenization theory and relates our findings to a homogenization problem for the equations modeling acoustic wave propagation.

All the chapters are accompanied by a comments section, which puts the results obtained into perspective of other research. Moreover, almost all sections are added a subtitle. In this subtitle we give some keywords for this section and/or mention the most important results of the respective section. We emphasize that the strengths of the results developed lie in its general applicability. Furthermore, the main focus is on the discussion of properties of $\s M$ rather than $\s A$ in \eqref{eq:in0}. 

Given the structural point of view of evolutionary equations, the results developed in this thesis are original and shed some light on the influence of the choice of $\s M$ for modeling aspects and computational treatments. Moreover, as to the best of the author knowledge, the question of continuous dependence on the coefficients has not been addressed so far to this extent of generality. Regarding known results obtained by the author of this manuscript in earlier works, we highlight that the methods developed here do not use the concept of abstract distribution spaces such as Sobolev lattices/towers.

There is some research for specific equations under particular topologies, which complement the result of the thesis. We refer to the concluding section of Chapter \ref{ch:pde} for a more detailed account of available results.

\section*{Some Remarks on the Notation}

When writing this manuscript, we made an effort in using the most reasonable notation in the whole of the present thesis. However, we are aware of the fact that some readers (if not all) may find the notation used counterintuitive at some point. We mention some of the idiosyncrasies as follows. The time derivative will be denoted by \begin{equation}\tag{time derivatives}\label{eq:rnot}
  \partial_t,\ \partial_{t,\nu},\ \partial_{t,\mu},\ \text{or }\partial_{t,\eta}                                                                                                                                                                                                                                                                                                                                                                                           \end{equation}
 (the additional subscript refers to the domain of definition). All vector-analytic operators such as 
 \begin{equation}\tag{spatial operators}\label{eq:rnot2}
  \curl,\ \dive,\  \text{or }\grad
  \end{equation} are thought of acting on the spatial variables only. Operators acting in space-time are denoted with capital calligraphic letters, that is, 
  \begin{equation}\tag{space time operators}\label{eq:rnot35}
  \s M,\ \s N,\ \s S,\ \s T,\ \text{etc.}
  \end{equation}  
   General abstract Hilbert spaces are denoted by 
  \begin{equation}\tag{abstract Hilbert spaces}\label{eq:rnot3}
  \s X,\ \s Y,\ \s Z,\ \text{or } \s W
  \end{equation}  
   and operators acting in these spaces are denoted by capital Latin letters, such as, $M$, $N$, etc. We use the letters 
   \begin{equation}\tag{superscripts}\label{eq:rnot4}
   \textnormal{n},\ \textnormal{s},\ \text{and }\textnormal{w}
   \end{equation}
 as superscripts of certain sets to indicate the underlying topology being the norm, the strong or the weak operator topology, respectively. The space of continuous linear operators from a Hilbert space $\s X$ to a Hilbert space $\s Y$ is denoted by 
   \begin{equation}\tag{continuous linear operators}\label{eq:rnot5}
      L(\s X,\s Y),\text{ or }L(\s X), \text{if $\s X=\s Y$.}
   \end{equation}
 Multiplication operators induced by certain functions $f$ are denoted by 
 \begin{equation}
   T_f \tag{multiplication operators}\label{eq:rnot6}
 \end{equation}
or, to avoid unnecessarily cluttered notation, just by $f$ again. The operator of multiplication by $\1_{(-\infty,t]}$, the characteristic function of the real interval $(-\infty,t]$, will be written as 
\[
  Q_t = T_{\1_{(-\infty,t]}}.
\]
 Furthermore, we write
 \[
  P_t\coloneqq 1-Q_t,\ t\in \mathbb{R}.
 \]
We denote the identity operator just by 
\begin{equation}\tag{identity}\label{eq:rnot7}
 1 \text{ or } \id.
\end{equation}
If we consider a mapping acting as $x\mapsto x$, we also write 
\begin{equation}\tag{canonical embedding}\label{eq:rnot8}
 \hookrightarrow
\end{equation}
and if this mapping happens to be compact, we may use 
\begin{equation}\tag{compact embedding}\label{eq:rnot9}
 \hookrightarrow\hookrightarrow.
\end{equation}
The sets 
\begin{equation}\tag{(standard) evolutionary mappings}\label{eq:rnot10}
 L_{\textnormal{ev},\nu}(\s X,\s Y),\ \text{and }L_{\textnormal{sev}}(\s X,\s Y)
\end{equation}
 of evolutionary and standard evolutionary mappings are introduced in Definition \ref{d:em} and \ref{d:scd_sev}, respectively. The set 
 \begin{equation}\tag{closable evolutionary mappings}\label{eq:rnot11}
  C_{\textnormal{ev},\nu}(\s X)
 \end{equation}
of closable evolutionary mappings is defined in Definition \ref{d:cev}. 

Unless expressed explicitly otherwise, all vector spaces discussed in this manu\-script have the complex numbers as underlying scalar field. Particularly important for the representation theorem to be proved in Section \ref{s:Art} is that all scalar products are anti-linear in the \emph{first} and linear in the \emph{second} factor. The domain, kernel and range of a mapping $T$ are respectively denoted by 
\begin{equation}\tag{domain, kernel, range}\label{eq:rnot12}
\dom(T),\ \kar(T),\text{ and }\ran(T).
 \end{equation}
 The image of a set $D$ under $T$ is written as 
 \begin{equation}\tag{image}\label{eq:rnot13}
T[D].
 \end{equation}
  A linear operator $T$ with domain $\dom(T)$ mapping from a Hilbert space $\s X$ to a Hilbert space $\s Y$ is also written as 
  \begin{equation}\tag{operator}\label{eq:rnot14}
   T\colon \dom(T)\subseteq \s X\to \s Y.
  \end{equation}
 We will occasionally employ the custom of identifying $T$ with its graph 
 \begin{equation}\tag{graph}
 \label{eq:rnot15}T\subseteq \s X\times \s Y.
  \end{equation}
 $T$ is closed, if $T\subseteq \s X\times \s Y$ is closed and $T$ is closable, if $\overline{T}\subseteq \s X\times \s Y$ is right-unique in the sense that for all $(x,y_1),(x,y_2)\in \overline{T}$, we always have $y_1=y_2$. If $T\colon\dom(T)\subseteq \s X\to \s Y$ and $S\colon \dom(S)\subseteq \s Z\to \s W$ for some Hilbert spaces $\s Z$ and $\s W$, we write 
 \begin{equation}\tag{extensions}\label{eq:rnot16}
 T=S\text{ or }T\subseteq S,  
 \end{equation}
 if the graphs of $T$ and $S$ coincide or if the graph of $S$ is a superset of $T$. We write 
 \begin{equation}\tag{restriction}\label{eq:rnot18}
  T|_D
 \end{equation}
 for the restriction of $T$ to a set $D$ and say $T=S$ on $D$, if $T|_D=S|_D$. The sum $T+S$ and composition $ST$ of two operators $T$ and $S$ are defined on the respective natural domains, that is, $\dom(T)\cap\dom(S)$ and $\{ x\in \dom(T); Tx\in \dom(S)\}$. We will also use
\[
   S\cap T \colon \{ x\in \dom(S)\cap\dom(T); Sx=Tx\}\subseteq \s X\cap\s Z \to \s Y\cap \s W, x\mapsto Sx.
\]
The spectrum and the resolvent set of $T$ are respectively denoted by 
\begin{equation}\tag{spectrum and resolvent}\label{eq:rnot17}
  \sigma(T)\text{ and }\rho(T).
\end{equation}
The dual of a Hilbert space $\s X$ is denoted by 
\begin{equation}\tag{dual}\label{eq:rnot20}
 \s X^*. 
\end{equation}
 The index set of a net $(x_\iota)_\iota$ will generically be a directed set $I$; similarly, a sequence $(x_n)_n$ is thought of as a mapping $\N\ni n\mapsto x_n$. Limits are denoted by \begin{equation}\tag{limits}\label{eq:notlim}
 \lim_{n\to\infty} x_n,\, \lim_{\iota} x_\iota.
 \end{equation}
 We also write $x_\iota\stackrel{\iota}{\to}\lim_{\iota} x_\iota$ or $(x_\iota)_\iota\to \lim_{\iota} x_\iota$. We will leave it to the context to stress the particular topology. Unless explicitly expressed otherwise $d$ is a positive integer. 
If there is a risk of ambiguity, we will add a subscript to the scalar products $\langle\cdot,\cdot\rangle$ and norms $\|\cdot\|$ to identify the particular Hilbert or Banach space the computations are carried out. The closed unit ball of a Banach space $\s Z$ is written as 
\begin{equation}\tag{unit ball}\label{eq:not20}
B_{\s Z}.
\end{equation}
 The support of $f$ will be denoted by 
 \begin{equation}\tag{support}\label{eq:spt}
 \spt f. 
 \end{equation}

We will also use the standard notation to denote Sobolev spaces, Lebesgue spaces, spaces of continuous functions etc. We shall also use self-explanatory notation as in $\R_{>0}$, $\R_{\geq 0}$, $\mathbb{C}_{\Re>0}$, etc.

A definition will be ended by $\triangle$, a remark by \scalebox{1.5}[1.5]{$\diamond$}, the end-of-proof symbol is $\square$.
\chapter{Time-shift Invariant Operators}\label{ch:TSO}

In this chapter, we shall present the basic results needed for the remaining parts. In particular, we introduce the (time) derivative operator on weighted vector-valued $L^2$-type spaces. We provide an explicit spectral representation for the time derivative via the Fourier--Laplace transformation. Later on, we present a well-known representation result for time-shift invariant, causal operators. In the setting presented in this exposition, time-shift invariant, causal operators are equivalently described by bounded operator-valued, analytic functions of the time derivative. As such, time-shift invariant causal operators form an important class of so-called evolutionary mappings to be introduced in the next chapter.

\renewcommand{\baselinestretch}{0.65}\normalsize\mysection{The Time Derivative}{The Time Derivative}{time derivative $\cdot$ Young's inequality $\cdot$  Fourier transformation $\cdot$ Fourier--Laplace transformation $\cdot$ Theorem \ref{t:td_inv} $\cdot$ Theorem \ref{th:sp_td_inv}}\label{se:TD}

\renewcommand{\baselinestretch}{1}\normalsize
 We start out with the definition of the time derivative.  For this, we recall the $L^2$-variant of the Morgenstern norm known from the classical proof of the Picard--Lindelöf theorem (see \cite{Morgenstern1952}): Let $\nu\in \mathbb{R}$. We define
\[
   L_\nu^2(\mathbb{R})\coloneqq \{ f\in L_{\textrm{loc}}^2(\mathbb{R}); t\mapsto e^{-\nu t}f(t)\in L^2(\mathbb{R})\}.
\] 
With $\lambda_\nu$ being the Lebesgue measure on $\mathbb{R}$ with Radon--Nikodym derivative $t\mapsto e^{-2\nu t}$, we immediately verify the equality
\[
   L_\nu^2(\mathbb{R})=L^2(\lambda_\nu).
\]
Hence, $L_\nu^2(\mathbb{R})$ is a Hilbert space. Henceforth, we will encounter Hilbert space valued $L_\nu^2$-functions. That is, for a Hilbert space $\s X$, we let
\begin{align*}
   L_\nu^2(\s X)&\coloneqq L_\nu^2(\mathbb{R};\s X)
   \\     & \coloneqq L^2(\lambda_\nu;\s X)
   \\     & \coloneqq \{ f\colon \mathbb{R}\to \s X; f\text{ measurable, }t\mapsto \|f(t)\|_{\s X}\in L_\nu^2(\mathbb{R})\}.
\end{align*}
Note that $L_\nu^2(\mathbb{R};\s X)$ and $L_0^2(\mathbb{R};{\s X})=L^2(\mathbb{R};{\s X})$ are unitarily equivalent. Indeed, the mapping of multiplication by $t\mapsto e^{-\nu t}$,  denoted by $\textnormal{m}_\nu$, is the desired unitary operator from the weighted to the unweighted space: Take $f\in L_\nu^2(\mathbb{R};{\s X})$ and compute
\[
    \| \textnormal{m}_\nu f\|_{L^2}^2=\int_{\mathbb{R}} \|e^{-\nu t}f(t)\|_{\s X}^2 d t = \int_{\mathbb{R}} \|f(t)\|_{\s X}^2 e^{-2\nu t} dt = \| f\|_{L_\nu^2}^2.
\]
Realizing that $\textnormal{m}_{\nu}$ is a bijection on compactly supported measurable functions, we infer that the range of $\textnormal{m}_\nu$ is dense in $L^2(\mathbb{R};{\s X})$. Thus, $\textnormal{m}_\nu$ is unitary.

We recall the Sobolev space $H^1(\mathbb{R};{\s X})$ of Hilbert space ${\s X}$-valued $L^2$-functions with distributional derivative representable as $L^2$-function. The exponentially weighted version of this Sobolev space can be defined in two ways: either as the unitary image of $H^1$ or as the space of $L_\nu^2$-functions with distributional derivative lying in $L_\nu^2$:

\begin{Proposition}\label{p:sobsp} Let $\nu\in \mathbb{R}$, ${\s X}$ Hilbert space. Then 
\[
    H_\nu^1(\mathbb{R};{\s X})\coloneqq \{f\in L_\nu^2(\mathbb{R};{\s X}); f'\in L_\nu^2(\mathbb{R};{\s X})\} = \textnormal{m}_{-\nu} [H^1(\mathbb{R};{\s X})].
\]
\end{Proposition}
\begin{Proof}
  Let $f\in H_\nu^1(\mathbb{R};{\s X})$. Then $\textnormal{m}_\nu f \in L^2(\mathbb{R};{\s X})$. Moreover, for any $\phi\in C_c^\infty(\mathbb{R})$ (the space of smooth functions with compact support), we compute
 \begin{align*}
      \int_{\mathbb{R}}\phi'\textnormal{m}_\nu f &=\int_{\mathbb{R}} \phi'(t) e^{-\nu t}f(t) \dd t 
      \\ & = \int_{\mathbb{R}} ((\textnormal{m}_\nu\phi)'(t)+\nu e^{-\nu t}\phi(t))f(t) \dd t
      \\ & = \int_{\mathbb{R}} (\textnormal{m}_\nu\phi)'(t)f(t)+\nu e^{-\nu t}\phi(t)f(t) \dd t
      \\ & = -\int_{\mathbb{R}} (\textnormal{m}_\nu\phi)(t)f'(t)-\nu e^{-\nu t}\phi(t)f(t) \dd t.      
 \end{align*}
Hence, $\textnormal{m}_\nu f \in H^1(\mathbb{R};{\s X})$ and $(\textnormal{m}_\nu f)'=\textnormal{m}_\nu f'-\nu \textnormal{m}_\nu f$. Thus, we obtain the inclusion $H_\nu^1(\mathbb{R};{\s X})\subseteq \textnormal{m}_{-\nu}[H^1(\mathbb{R};{\s X})]$. On the other hand, let $f\in \textnormal{m}_{-\nu}[H^1(\mathbb{R};{\s X})]$ and let $g\in H^1(\mathbb{R};\mathcal{X})$ be such that $f=\textnormal{m}_{-\nu}g$. Then, for every $\phi\in C_c^\infty(\mathbb{R})$, we get
\begin{align*}
    \int_{\mathbb{R}} \phi'f& = \int_{\mathbb{R}} \phi' \textnormal{m}_{-\nu}g 
    \\ & =\int_{\mathbb{R}} (\textnormal{m}_{-\nu}\phi)' g-\nu \textnormal{m}_{-\nu}\phi g
    \\ & =\int_{\mathbb{R}} -\textnormal{m}_{-\nu}\phi g'-\nu \phi \textnormal{m}_{-\nu}g
    \\ & =\int_{\mathbb{R}} -\phi \textnormal{m}_{-\nu}g'-\nu \phi \textnormal{m}_{-\nu}g.
\end{align*}
Consequently, $f'=\textnormal{m}_{-\nu}g'+\nu \textnormal{m}_{-\nu}g\in L_\nu^2(\mathbb{R};{\s X})$, which yields the assertion.
\end{Proof}

With Proposition \ref{p:sobsp} at hand, we can define the time derivative.

\begin{Definition} Let $\nu\in \mathbb{R}$, ${\s X}$ Hilbert space. We define the graph
\[
    \partial_{t,\nu} \coloneqq \{ (f,g)\in L_\nu^2(\mathbb{R};{\s X})\times L_\nu^2(\mathbb{R};{\s X}); g=f' \}\subseteq L_\nu^2(\mathbb{R};\mathcal{X})\times L_\nu^2(\mathbb{R};\mathcal{X}).
\]
\end{Definition}

By linearity of the distributional derivative, $\partial_{t,\nu}$ is a linear relation. Moreover, by the uniqueness of the weak derivative, $\partial_{t,\nu}$ is actually a linear operator. Next, if for $f\in H_\nu^1(\mathbb{R};{\s X})$ we have $f'=0$ then a short argument shows that $f$ is constant almost everywhere. The only $\lambda_\nu$-integrable constant, however, is $0$. Thus, $\partial_{t,\nu}$ is one-to-one. But even more is true:  One of the key facts used in the following is that $\partial_{t,\nu}$ is \emph{boundedly} invertible as an operator in $L_\nu^2(\mathbb{R};\mathcal{X})$. Before stating and proving such a result, we recall the notion of convolution of two functions. 

\begin{Definition}
  Let $h\colon \mathbb{R}\to \mathbb{C}$, $f\colon \mathbb{R}\to \s X$, where $\s X$ is a Hilbert space. Let $t\in \mathbb{R}$ and assume that 
$s\mapsto h(t-s)f(s) \in L^1(\mathbb{R};\s X)$. Then define
\[
   h* f (t) \coloneqq \int_{\mathbb{R}} h(t-s)f(s)\dd s.
\]
\end{Definition}

When it comes to convolutions an inequality due to Young is of prime importance. We state the inequality in our present context.

\begin{Theorem}[Young's inequality]\label{t:YI} Let $\nu\in \mathbb{R}$, $\s X$ Hilbert space, $h\in L^1(\lambda_\nu)$, and let $f\in L^2(\lambda_\nu;\s X)$. Then,
for almost every $t\in \mathbb{R}$, $s\mapsto h(t-s)f(s)\in L^1(\R;\s X)$ and 
\[
   \|h*f\|_{L^2(\lambda_\nu;\s X)}\leq \|h\|_{L^1(\lambda_\nu)}\|f\|_{L^2(\lambda_\nu;\s X)}.
\] 
\end{Theorem}
\begin{Proof}
  We compute
  \begin{align*}
   & \|h*f\|_{L^2(\lambda_\nu;\s X)}^2 
   \\ &  = \int_\R \|h*f(t)\|^2e^{-2\nu t} \dd t
   \\ & \leq \int_\R \Big( \int_\R \|h(t-s)f(s)\| \dd s\Big)^2e^{-2\nu t}\dd t
   \\ & = \int_\R \Big( \int_\R |h(t-s)|e^{-\nu (t-s)}\|f(s)\|e^{-\nu s} \dd s\Big)^2\dd t
   \\ & \leq \int_\R \Big( \int_\R |h(t-s)|e^{-\nu (t-s)}\dd s\Big)\Big(\int_\R|h(t-s)|e^{-\nu (t-s)}\|f(s)\|^2e^{-2\nu s} \dd s\Big)\dd t
   \\ & = \|h\|_{L^1(\lambda_\nu)}^2 \int_\R \|f(s)\|^2e^{-2\nu s} \dd s
   \\ & = \|h\|^2_{L^1(\lambda_\nu)}\|f\|^2_{L^2(\lambda_\nu;\s X)}.
  \end{align*}
\end{Proof}

\begin{Remark}\label{r:conv_op}
 We shall put Theorem \ref{t:YI} into an operator-theoretic perspective. The inequality asserted is the same as saying that the \emph{operator of convolution with $h$} is a continuous mapping in $L_\nu^2(\R;\s X)$ with operator norm bounded by $\|h\|_{L^1(\lambda_\nu)}$. Later on, we employ this fact, when we recall the convolution theorem in the context of the Fourier transformation.
\end{Remark}

Next, we provide the result of bounded invertibility of the time derivative, which may be dated back to \cite{Picard1989}.

\begin{Theorem}\label{t:td_inv} Let $\nu\in \mathbb{R}\setminus\{0\}$, ${\s X}$ a Hilbert space. Then $\partial_{t,\nu}$ is boundedly invertible in $L_\nu^2(\mathbb{R};\mathcal{X})$ and, for all $f\in L_\nu^2(\mathbb{R};{\s X})$, $t\in \mathbb{R}$,  the formula
\begin{equation}\label{eq:td_inv}
      \partial_{t,\nu}^{-1} f (t) =\begin{cases}
                               \int_{-\infty}^t f(s)\dd s,& \nu>0,
                               \\  -\int_t^\infty f(s)\dd s,& \nu<0,
                            \end{cases}
\end{equation}
as well as the inequality $\|\partial_{t,\nu}^{-1}\|\leq 1/|\nu|$ hold true.
\end{Theorem}
\begin{Proof}
    To prove that the integrals in \eqref{eq:td_inv} exist and that the right-hand side of \eqref{eq:td_inv} defines an element in $L_\nu^2(\mathbb{R};\s X)$, we use Theorem \ref{t:YI}. In fact, for $\nu>0$, it is sufficient to observe that $h\colon t\mapsto \1_{[0,\infty)}(t)\in L^1(\lambda_\nu)$ satisfies $\|h\|_{L^1(\lambda_\nu)}=1/\nu$ as well as
\begin{align*}
    h*f(t) &=\int_{\R} h(t-s)f(s)\dd s
    \\     &=\int_{\R} \1_{[0,\infty)}(t-s)f(s)\dd s  
    \\     &=\int_{-\infty}^t f(s)\dd s\quad (t\in \mathbb{R}, f\in L_\nu^2(\s X)).
\end{align*}
 For $\nu<0$ use $h=-\1_{(-\infty,0]}$ instead. Once proved that $h*=\partial_{t,\nu}^{-1}$, the estimate for the operator norm follows from Theorem \ref{t:YI}. Next, by Fubini's Theorem the distributional derivative of $g\coloneqq \big(t\mapsto \int_{-\infty}^t f(s)\dd s\big)$ equals $f$, so $\partial_{t,\nu}\circ(h*)=\id_{L_\nu^2}$. Another application of Fubini's Theorem yields that for any $f\in \dom(\partial_{t,\nu})$ we have 
 \[h*\partial_{t,\nu}f = h*f'=(h*f)'=f.\]
 Therefore, $(h*)\circ \partial_{t,\nu}=\id_{\dom(\partial_{t,\nu})}$.
\end{Proof}

The main observation needed for the proof of Theorem \ref{t:td_inv} is
\begin{equation}\label{eq:td_conv}
 \partial_{t,\nu}^{-1}f = h*f,\text{ for }h=\begin{cases}
                                             \1_{[0,\infty)},&\nu>0,
                                             \\ -\1_{(-\infty,0]}, & \nu<0.
                                            \end{cases}
\end{equation}
  Exploiting equation \eqref{eq:td_conv} a bit further, we will move on proving an explicit spectral theorem for $\partial_{t,\nu}^{-1}$, that is, we will show that $\partial_{t,\nu}^{-1}$ is unitarily equivalent to a multiplication operator. The unitary operator yielding the spectral representation involves the Fourier transformation, which will be defined next.

\begin{Definition} Let $\s X$ Hilbert space. We define the \emph{Fourier transformation} \[\s F\colon L^2(\R;\s X)\to L^2(\R;\s X)\] as the continuous extension of the operator given by
\[
    \mathcal{F}\phi(\xi)\coloneqq \frac{1}{\sqrt{2\pi}}\int_\mathbb{R} e^{-i\xi t}\phi(t)\dd t,\quad \xi\in \mathbb{R},
\]
for all integrable measurable functions $\phi\colon\mathbb{R}\to {\s X}$. 
\end{Definition}

\begin{Remark}\label{r:PL}
(a) We recall that Plancherel's theorem (for Hilbert space valued functions) states that $\s F$ is not only continuous but unitary. 

(b) For $f\in L^1(\R)$ and $g\in L^2(\R;\s X)$ the so-called convolution theorem holds, that is, $\s F(f*g)=\sqrt{2\pi}(\s F f)\cdot(\s F g)$. In fact, the equality is an application of Fubini's theorem for $f\in L^1(\R)$ and $g\in L^1(\mathbb{R};\s X)\cap L^2(\mathbb{R};\s X)$. In particular, for $f\in L^1(\mathbb{R})$ and all $g\in L^2(\mathbb{R};\s X)$ it follows that
\begin{align*}
   \sqrt{2\pi}\|(\mathcal{F}f)\cdot(\mathcal{F}g)\|_{L^2(\mathbb{R};\s X)} &\leq \sqrt{2\pi}\|\mathcal{F}f\|_\infty \|\mathcal{F}g\|_{L^2(\mathbb{R};\s X)}
   \\&\leq\|f\|_{L^1(\mathbb{R})} \|g\|_{L^2(\mathbb{R};\s X)}
\end{align*}
where $\|\mathcal{F}f\|_\infty\coloneqq \sup_{\xi\in \mathbb{R}}|(\mathcal{F}f)(\xi)|$. We read off that the mapping
\[
   L^2(\mathbb{R};\s X)\ni g\mapsto \sqrt{2\pi}(\mathcal{F}f)\cdot (\mathcal{F}g) \in L^2(\mathbb{R};\s X)
\]
is continuous. Moreover, the latter mapping coincides with the continuous map
\[
   L^2(\mathbb{R};\s X)\ni g  \mapsto \mathcal{F}(f*g) \in L^2(\R;\s X)
\]
on the dense set $L^1(\mathbb{R};\s X)\cap L^2(\mathbb{R};\s X)\subseteq L^2(\mathbb{R};\s X)$; thus, $\s F(f*g)=\sqrt{2\pi}(\s F f)\cdot(\s F g)$ for all $f\in L^1(\mathbb{R})$ and $g\in L^2(\mathbb{R};\s X)$.
\end{Remark}

Being mainly interested in the case of exponentially weighted $L^2$-spaces, we introduce the \emph{Fourier--Laplace transformation} $\s L_\nu$ as an operator from $L_\nu^2(\R;\s X)$ to $L^2(\R;\s X)$, $\s X$ Hilbert space, $\nu\in \R$, as follows
\begin{equation}\label{e:FLT}
   \s L_\nu \coloneqq \s F \textnormal{m}_\nu.
\end{equation}
As a composition of unitary operators, $\s L_\nu$ is unitary itself. In the proof of Theorem \ref{t:td_inv}, we realized that $\partial_{t,\nu}^{-1}$ can be written as a convolution. From Remark \ref{r:PL}(b) we get that convolutions are multiplication operators on the Fourier transformed side. Putting these two facts together, we get the desired spectral representation for $\partial_{t,\nu}$ as multiplication operator with a certain function. For stating the result, we introduce ${T}_h$ the multiplication operator of multiplying by some measurable function $h\colon \R\to\mathbb{C}$ in $L^2(\mathbb{R};\s X)$, $\s X$ a Hilbert space.

\begin{Corollary}\label{co:spe_thm} Let ${\s X}$ Hilbert space, $\nu\in\mathbb{R}\setminus\{0\}$. Then the equality
 \[
     \partial_{t,\nu}^{-1}=\mathcal{L}_\nu^*{T}_{\hat h_\nu}\mathcal{L}_\nu
 \]
holds, where ${T}_{\hat h_\nu}$ is the operator of multiplication by the function $\hat h_\nu\colon \mathbb{R}\ni \xi\mapsto 1/(i\xi+\nu)$.
\end{Corollary}
\begin{Proof}
 Recall from equation \eqref{eq:td_conv}, the equality $\partial_{t,\nu}^{-1}g=h*g$ for $g\in L_\nu^2(\R;\s X)$ with $h=\1_{[0,\infty)}$ for $\nu>0$ or $h=-\1_{(-\infty,0]}$ for $\nu<0$. In either case we have $\sqrt{2\pi}\s F \textnormal{m}_\nu h = \hat h_\nu$.  Hence, with $\textnormal{m}_{-\nu}\s F^*=(\s F\textnormal{m}_{\nu})^*= \s L_\nu^*$ and using the convolution theorem (Remark \ref{r:PL}(b)) we get for $g\in L_\nu^2(\R;\s X)$:
\begin{align*}
   \partial_{t,\nu}^{-1}g & = h*g
   \\                     & = \textnormal{m}_{-\nu}\textnormal{m}_{\nu} (h* g)
   \\                     & = \textnormal{m}_{-\nu} (\textnormal{m}_{\nu}h* \textnormal{m}_{\nu}g)
   \\                     & = \textnormal{m}_{-\nu}\s F^* \s F (\textnormal{m}_{\nu}h* \textnormal{m}_{\nu}g)
   \\                     & = \s L_\nu^* \sqrt{2\pi} (\s F \textnormal{m}_{\nu}h)\cdot (\s F \textnormal{m}_{\nu}g)
   \\                     & = \s L_\nu^* {T}_{\hat h_\nu} \s L_{\nu}g.
\end{align*}
\end{Proof}

\begin{Remark}\label{r:spe_thm_dtnu}
  Note that Corollary \ref{co:spe_thm} also contains the spectral representation for $\partial_{t,\nu}$. Indeed, the operator $\s L_\nu$ being unitary it suffices to observe ${T}_{\hat h_\nu}^{-1}={T}_{1/\hat h_\nu}$, where ${T}_{1/\hat h_\nu}$ is the operator of multiplying by $\xi\mapsto i\xi+\nu$. Hence, we get
  \[
     \partial_{t,\nu}=\mathcal{L}_\nu^*{T}_{1/\hat h_\nu}\mathcal{L}_\nu.
  \]
  A consequence of the latter formula is
  \begin{equation}\label{eq:reder}
     \Re\langle \partial_{t,\nu}\phi,\phi\rangle = \nu \langle \phi,\phi\rangle \quad (\phi\in \dom(\partial_{t,\nu})).
  \end{equation}
  Indeed, take $\phi\in \dom(\partial_{t,\nu})$. Then
  \begin{align*}
     \Re\langle \partial_{t,\nu}\phi,\phi\rangle & =\Re\langle {T}_{1/\hat h_\nu}\mathcal{L}_\nu\phi,\mathcal{L}_\nu\phi\rangle 
     \\ & =\Re\langle (\xi\mapsto i\xi+\nu)\mathcal{L}_\nu\phi,\mathcal{L}_\nu\phi\rangle 
     \\ & = \nu \Re\langle\mathcal{L}_\nu \phi,\mathcal{L}_\nu\phi\rangle 
     \\ & = \nu \Re\langle\phi,\phi\rangle.
  \end{align*}
\end{Remark}

With the help of Corollary \ref{co:spe_thm} we can define a functional calculus for the operator $\partial_{t,\nu}^{-1}$. Before properly defining the functional calculus, we analyze the spectrum of $\partial_{t,\nu}^{-1}$ first. We restrict ourselves to the case $\nu>0$. Moreover, note that it suffices to state the scalar case only, as the respective operator defined on ${\s X}$-valued $L_\nu^2$-functions has the same spectrum as in the scalar-valued case as long as ${\s X}$ is at least one-dimensional. We denote by $B(a,\delta)$ the open ball in $\mathbb{C}$ centered at $a\in \mathbb{C}$ with radius $\delta>0$. The boundary of an open set $\Omega$ in some underlying topological space is denoted by $\partial\Omega$. The spectrum of $\partial_{t,\nu}^{-1}$, $\nu>0$, is as follows.

\begin{Theorem}\label{th:sp_td_inv} Let $\nu\in\mathbb{R}_{>0}$. Then $\sigma\left(\partial_{t,\nu}^{-1}\right)=\partial B(r,r)$ with $r\coloneqq1/(2\nu)$. 
\end{Theorem}

\begin{Lemma}\label{l:ran_h_nu} Let $\nu\in \mathbb{R}_{>0}$. Let $\hat h_\nu(\xi)=1/(i\xi+\nu)$, $\xi\in \mathbb{R}$. Then $\ran(\hat h_\nu)=\partial B(r,r)\setminus \{0\}$ with $r\coloneqq 1/(2\nu)$. 
\end{Lemma}
\begin{Proof}
 For proving $\ran(\hat h_\nu)\subseteq\partial B(r,r)\setminus\{0\}$, we take $\xi\in \mathbb{R}$ and compute
  \begin{align*}
      \left| \frac{1}{i\xi+\nu}-r\right|^2&=\left| \frac{-i\xi+\nu}{\xi^2+\nu^2}-\frac{1}{2\nu}\right|^2
      \\ &= \left| \frac{\nu}{\xi^2+\nu^2}-\frac{1}{2\nu}-i\frac{\xi}{\xi^2+\nu^2}\right|^2
      \\ &=  \left(\frac{\nu}{\xi^2+\nu^2}-\frac{1}{2\nu}\right)^2+\left(\frac{\xi}{\xi^2+\nu^2}\right)^2
      \\ &=  \left(\frac{\nu}{\xi^2+\nu^2}\right)^2+\left(\frac{1}{2\nu}\right)^2-\frac{1}{\nu}\frac{\nu}{\xi^2+\nu^2}+\left(\frac{\xi}{\xi^2+\nu^2}\right)^2
      \\ &=  \left(\frac{1}{2\nu}\right)^2=r^2.
  \end{align*}
Next, we show $\partial B(r,r)\setminus\{0\}\subseteq \ran(\hat h_\nu)$. For this, let $\theta\in (-\pi,\pi)$ and compute
\[
   \Re\Big(\frac{1}{re^{i\theta}+r}\Big)=\frac{1}{r}\Re \Big(\frac{e^{-i\theta}+1}{|e^{i\theta}+1|}\Big)=\frac{1}{r} \Big(\frac{\cos(\theta)+1}{2+2\cos(\theta)}\Big)=\nu.
\]
\end{Proof}

\begin{Proof}[of Theorem \ref{th:sp_td_inv}]
  By Corollary \ref{co:spe_thm}, the operator $\partial_{t,\nu}^{-1}$ is unitarily equivalent to the operator of multiplying by $\hat h_\nu\colon\xi\mapsto 1/(i\xi+\nu)$. The spectrum of a multiplication operator of multiplying with a continuous function coincides with the closure of the range of this function. Hence, it suffices to show that $\ran(\hat h_\nu)=\partial B(r,r)\setminus\{0\}$. The latter equality, however, has been settled in Lemma \ref{l:ran_h_nu}.
\end{Proof}

\renewcommand{\baselinestretch}{0.65}\normalsize\mysection{A Representation Theorem}{A Representation Theorem}{operator-valued functions of time derivative $\cdot$ causal continuous operators $\cdot$ translation-invariance $\cdot$ time-shift invariance $\cdot$ Paley--Wiener Theorem $\cdot$ Hardy--Lebesgue space $\cdot$ Corollary \ref{co:SW_weighted}}\label{s:Art}

\renewcommand{\baselinestretch}{1}\normalsize
In order to put the applications to be discussed later on into a broader perspective, we introduce \emph{operator-valued} functions of $\partial_{t,\nu}^{-1}$. For this, in principle, it would be sufficient to consider functions with domain $\partial B(r,r)$, $r=1/(2\nu)$. We will, however, also address causality in our solution theory, which necessitates the consideration of a somewhat different class of functions. We shall provide a proper definition of the notion of causality next. Afterwards, with the help of a well-established representation theorem, we motivate this specific class of (operator-valued) functions of $\partial_{t,\nu}^{-1}$. We highlight an important property of this class (Remark \ref{r:SW_OT}), which will serve as the basis for the definition of evolutionary mappings.

\begin{Definition}\label{d:cau} Let $\s X$, $\s Y$ Hilbert spaces, $\nu\in \mathbb{R}$, $\mathcal{M}\colon L_\nu^2(\mathbb{R};\s X)\to L_\nu^2(\mathbb{R};\s Y)$ bounded, linear.
\begin{enumerate}[label=(\alph*)]
 \item $\mathcal{M}$ is called \emph{(forward) causal}, if for all $t\in \mathbb{R}$ we have $P_t\mathcal{M}P_t=\mathcal{M}P_t$, where $P_t$ is the operator ${T}_{\1_{[t,\infty)}}$ of multiplication with $\1_{[t,\infty)}$;
 \item $\mathcal{M}$ is called \emph{translation-invariant} (or \emph{time-shift invariant}), if for all $h\in \mathbb{R}$ the equality $\tau_h \mathcal{M} =\mathcal{M}\tau_h$ holds true, where $\tau_hf(t)\coloneqq f(t+h)$, $t\in \mathbb{R}$.
\end{enumerate}  
\end{Definition}

\begin{Remark}\label{r:tinv_caus} Let $\s X$, $\s Y$ be Hilbert spaces, $\nu\in \mathbb{R}$, $\mathcal{M}\in L(L_\nu^2(\s X),L_\nu^2(\s Y))$. We define $L_\nu^2(0,\infty;\s X)$ as the space of $L_\nu^2(\mathbb{R};\s X)$-functions supported on $[0,\infty)$ only.

Assume that $\mathcal{M}$ is translation-invariant. Then $\mathcal{M}$ is causal iff $\mathcal{M}$ maps $L_\nu^2(0,\infty;\s X)$ into $L_\nu^2(0,\infty;\s Y)$. Indeed, the necessity being trivial, we show the sufficiency next. Let $t\in \mathbb{R}$, and let $f\in L_\nu^2(\mathbb{R};\s X)$. Then \[\tau_tP_tf=P_0\tau_t f\in L_\nu^2(0,\infty;\s X).\] Thus, $\mathcal{M}\tau_tP_tf\in L_\nu^2(0,\infty;\s Y)$ and, hence, 
  \begin{align*}
    \tau_t \mathcal{M}P_tf  & = \mathcal{M}\tau_tP_tf 
      \\                      & = P_0\mathcal{M}\tau_tP_tf
      \\                      & = P_0\tau_t \mathcal{M}P_tf
      \\                      & = \tau_tP_t\mathcal{M}P_tf.
  \end{align*}
So, $\mathcal{M}P_tf=P_t\mathcal{M}P_tf$ yielding the assertion.
\end{Remark}

Next, we give a well-known representation theorem for translation-invariant, causal mappings. We adopt the strategy given in \cite{Weiss1991} for the proof. 

\begin{Theorem}\label{th:SW} Let $\s X$, $\s Y$ Hilbert spaces, $\mathcal{M}\in L(L^2(\s X),L^2(\s Y))$. Assume that $\mathcal{M}$ is time-shift invariant and that $\mathcal{M}$ maps $L^2(0,\infty;\s X)$ into $L^2(0,\infty;\s Y)$. Then there is a unique, bounded and analytic function $M\colon \mathbb{C}_{\Re>0}\to L(\s X,\s Y)$ with the following property: For any $u\in L^2(0,\infty;\s X)$ we have
\begin{equation}\label{eq:SW_ptw}
          (\mathcal{L}_{\nu}\mathcal{M}u)(\xi)= M\left(i\xi+\nu\right)(\mathcal{L}_{\nu} u)(\xi),\quad \xi\in \mathbb{R},\nu>0.
\end{equation} 
\end{Theorem}

\begin{Remark}\label{r:SW_OT} In the situation of Theorem \ref{th:SW} define for $\nu>0$ the function $M_\nu \colon \mathbb{R}\to L(\s X, \s Y)$ by $M_\nu(\xi)\coloneqq M(i\xi+\nu)$, $\xi\in \mathbb{R}$. Then, using our convention for multiplication operators, we realize that equation \eqref{eq:SW_ptw} can be written as
\begin{equation}\label{eq:SW_OT}
   \mathcal{L}_\nu \mathcal{M} u =  {T}_{M_\nu}\mathcal{L}_\nu u\quad (\nu>0, u\in L^2(0,\infty;\s X)).
\end{equation}
Equivalently, equation \eqref{eq:SW_OT} may be written as
\[
   \mathcal{M} u =  \mathcal{L}_\nu^*{T}_{M_\nu}\mathcal{L}_\nu u\quad (\nu>0, u\in L^2(0,\infty;\s X)).
\]
 Next, as the left-hand side is translation-invariant, so is the right-hand side, which yields equality for all $u\in L^2_c(\mathbb{R};\s X)$, that is, for compactly supported $L^2$-functions. Furthermore, by the boundedness of $M$, the expression $\mathcal{L}_\nu^*{T}_{M_\nu}\mathcal{L}_\nu$ defines a bounded linear operator from $L_\nu^2(\s X)$ to $L_\nu^2(\s Y)$. Thus, $\mathcal{M}$ defined on $L^2_c(\R;\s X)$ admits a unique continuous extension as an operator on $L_\nu^2(\mathbb{R};\s X)$ and as such, it coincides with $\mathcal{L}_\nu^*{T}_{M_\nu}\mathcal{L}_\nu$. Since for any $u\in L_\nu^2(\s X)\cap L_\mu^2(\s X)$ with $\mu,\nu\in \mathbb{R}$ we can choose a sequence $(u_n)_{n\in \N}$ in $L_\textnormal{c}^2(\mathbb{R};\s X)$ converging in both the spaces $L_\nu^2(\mathbb{R};\s X)$ and $L_\mu^2(\mathbb{R};\s X)$ to $u$ (e.g.~take $\phi=\1_{[-1,1]}$ and set $u_n\coloneqq \phi(\cdot/n)u$, $n\in \N$), the respective extensions of $\mathcal{M}$ to $L_\nu^2(\mathbb{R};\s X)$ and $L_\mu^2(\mathbb{R};\s X)$ coincide on $L_\nu^2(\mathbb{R};\s X)\cap L_\mu^2(\mathbb{R};\s X)$, $\nu,\mu>0$. Write $\s M^\nu$ for the closure of $\s M$ as an operator in $L_\nu^2$. Thus, equation \eqref{eq:SW_OT} finally implies
\begin{equation}\label{eq:SW_OTs}
   \mathcal{M}^\nu = \mathcal{L}_\nu^* {T}_{M_\nu}\mathcal{L}_\nu\quad(\nu>0).
\end{equation}
We will discuss the property for operators of being continuously extendable to the spaces $L_\nu^2$ for all sufficiently large $\nu>0$ in the context of evolutionary mappings later on. 
\end{Remark}

The main corollary of Theorem \ref{th:SW} is the following.

\begin{Corollary}\label{co:SW_weighted} Let $\s X$, $\s Y$ Hilbert spaces, $\nu\in \mathbb{R}$, $\mathcal{M}\in L(L_{\nu}^2(\s X);L_\nu^2(\s Y))$ translation-invariant and causal. Then $\mathcal{M}$ extends to an operator in $L(L_{\mu}^2(\s X),L_{\mu}^2(\s X))$ for all $\mu>\nu$ and there is a unique, bounded, analytic $M\colon \mathbb{C}_{\Re>\nu}\to L(\s X,\s Y)$ satisfying
\begin{equation}\label{eq:SW_w1}
   (\mathcal{L}_{\mu}\mathcal{M}u)(\xi)=M\left(i\xi+\mu\right)\mathcal{L}_{\mu}u (\xi),\quad \xi\in \mathbb{R},\mu>\nu,
\end{equation}
 for all $u\in L_\mu^2(\mathbb{R};\s X)$. Moreover, 
 \[
    \|\mathcal{M}\|_{L(L_{\mu}^2(\s X),L_{\mu}^2(\s Y))}\leq \sup_{z\in \mathbb{C}_{\Re\geq \mu}}\|M(z)\|_{L(\s X,\s Y)}.
 \]
\end{Corollary}
\begin{Proof}
  First, we observe that $\mathcal{N}\coloneqq \textnormal{m}_{\nu}\mathcal{M}\textnormal{m}_{-\nu}\in L(L^2(\s X),L^2(\s Y))$ is causal. Moreover, ${\mathcal{N}}$ is translation-invariant.
  Thus, by Theorem \ref{th:SW}, there exists a mapping ${N}\colon \mathbb{C}_{\Re>0}\to L(\s X, \s Y)$ bounded and analytic with the property that for all $u\in L^2(0,\infty;\s X)$ we have
  \begin{equation}\label{eq:co1}
          (\mathcal{L}_{\eta}{\mathcal{N}}u)(\xi)= {N}\left(i\xi+\eta\right)(\mathcal{L}_{\eta} u)(\xi),\quad \xi\in \mathbb{R},\eta>0.
  \end{equation}
  The latter equation thus implies that for all $\xi\in \mathbb{R}$, $\eta>0$ and $u\in L_{\nu}^2(0,\infty;\s X)$ (use $\mathcal{L}_{\eta}\textnormal{m}_{\nu}=\mathcal{L}_{\eta+\nu}$)
  \[
        (\mathcal{L}_{\eta+\nu}\mathcal{M} u)(\xi)= {N}\left(i\xi+\eta\right)(\mathcal{L}_{\eta+\nu}  u)(\xi)
  \]
  Setting ${M}(i\xi+\mu)\coloneqq {N}(i\xi+\eta)$ for all $\mu=\eta+\nu$, $\eta>0$, we get the desired representation.
  
  For a proof of the uniqueness statement, let $M_1\colon \mathbb{C}_{\Re>\nu}\to L(\s X)$ be analytic, bounded and such that \eqref{eq:SW_w1} holds with $M$ replaced by $M_1$. Then, the mapping ${N}_1(i\xi+\eta)\coloneqq M_1(i\xi+\eta+\nu)$ for all $\eta>0$ satisfies \eqref{eq:co1} with $N$ replaced by $N_1$. By the uniqueness statement of Theorem \ref{th:SW} it follows that $N_1=N$ and, hence, $M_1=M$.
  
  The norm estimate follows from the unitarity of the Fourier--Laplace transformation.
\end{Proof}

Before we prove Theorem \ref{th:SW}, we want to state a converse to Corollary \ref{co:SW_weighted}. This requires the (vector-valued) Paley--Wiener theorem, see for instance \cite[Corollary 1.2.6]{Waurick2011} for a self-contained proof.

\begin{Theorem}[{{\cite[Chapter 19]{Futh1} and \cite[Corollary 2.7]{PicPhy}}}]\label{t:PW} Let $\s X$ be a Hilbert space. Define the \emph{Hardy--Lebesgue space} 
\[
   \s H^2(\s X)\coloneqq \{ f\colon \mathbb{C}_{\Re>0}\to \mathbb{C}; f \text{ analytic}, \|f\|_{\s H^2}\coloneqq \sup_{\eta>0}\|f(i(\cdot)+\eta)\|_{L^2(\s X)}<\infty\}
\]
endowed with the obvious norm. Then $\s H^2(\s X)$ is a Hilbert space and the mapping
\[
   L^2(0,\infty;\s X)\ni f\mapsto \left(i\xi+\eta\mapsto (\s L_{\eta}f)(\xi)\right)\in \s H^2(\s X)
\]
is unitary with inverse
\[
   \s H^2(\s X)\ni f\mapsto \lim_{\eta\to 0+} \s L_\eta^* f(i(\cdot)+\eta) \in L^2(0,\infty; \s X).
\]
\end{Theorem}

A converse to Corollary \ref{co:SW_weighted} reads as follows.

\begin{Theorem}\label{t:SW_conv} Let $\s X$, $\s Y$ Hilbert spaces, $\nu\in \mathbb{R}$. Let $G\colon \mathbb{C}_{\Re>\nu}\to L(\s X, \s Y)$ be bounded and analytic.
 Then, for any $\mu>\nu$ the operator given by
\begin{equation}\label{eq:SW_conv}
     \mathcal{L}_{\mu}^* {T}_{G_\mu} \mathcal{L}_{\mu}\text{ with }G_\mu(\xi)\coloneqq G(i\xi+\mu), \xi\in \mathbb{R},
\end{equation}
defines a causal, translation-invariant, bounded, linear operator on $L^2_\mu(\mathbb{R};\s X)$.
\end{Theorem}
\begin{Proof}Addressing translation-invariance first, we realize that for any $h\in \mathbb{R}$, we have for $f\in L^2_\textnormal{c}(\mathbb{R};\s X)$, $\mu>\nu$
\begin{align*}
    \mathcal{L}_{\mu}\tau_h f(\xi)& =  \frac{1}{\sqrt{2\pi}}\int_\R e^{-i\xi t-\mu t}f(t+h)\dd t 
    \\ & =\frac{1}{\sqrt{2\pi}}\int_\R e^{-i\xi (t-h)-\mu (t-h)}f(t)\dd t
    \\ & = e^{(i\xi+\mu) h}\mathcal{L}_{\mu}f(\xi), \quad(h\in \mathbb{R},\xi\in \mathbb{R}),
\end{align*}
and so $\tau_h = \mathcal{L}_{\mu}^* {T}_{t_\mu} \mathcal{L}_\mu$ with $t_\mu(\xi)=e^{(i\xi+\mu) h}$, $\xi\in \mathbb{R}, h\in \mathbb{R}$. Thus,
\[
   \mathcal{L}_\mu^* {T}_{G_\mu}\mathcal{L}_\mu \tau_h = \mathcal{L}_\mu^* {T}_{G_\mu}{T}_{t_\mu}\mathcal{L}_\mu =\mathcal{L}_\mu^* {T}_{t_\mu}{T}_{G_\mu}\mathcal{L}_\mu=\tau_h\mathcal{L}_\mu^* {T}_{G_\mu}\mathcal{L}_\mu.
\]
Next, we prove causality. By Remark \ref{r:tinv_caus}, it suffices to show that $\mathcal{L}_{\mu}^* {T}_{G_\mu} \mathcal{L}_{\mu}$ leaves functions supported on $[0,\infty)$ invariant. For this, take $f\in L_\mu^2(0,\infty;\s X)$. Then, by definition, $\textnormal{m}_\mu f\in L^2(0,\infty;\s X)$. Hence, using the Paley--Wiener Theorem \ref{t:PW}, we infer that \[
   \mathbb{C}_{\Re>0} \ni i\xi+\eta \mapsto (\mathcal{F}(\textnormal{m}_{\mu} f))(\xi-i\eta)\coloneqq (\mathcal{L}_{\eta}(\textnormal{m}_{\mu}f))(\xi)
\]
belongs to the Hardy--Lebesgue space $\s H^2(\s X)$. The boundedness and analyticity of $G$ implies that ${T}_{G_\mu}$ maps $\s H^2(\s X)$ into $\s H^2(\s Y)$ in the sense that if $g \in \s H^2(\s X)$ then 
\[\left({T}_{G_\mu} g \colon \mathbb{C}_{\Re>0} \to \s Y, (i\xi+\eta)\mapsto G(i\xi+\mu+\eta)g(i\xi+\eta)\right)\in \s H^2(\s Y).\] Again referring to Theorem \ref{t:PW}, we get $\s F^* ({T}_{G_\mu}(\mathcal{F}(\textnormal{m}_{\mu} f)))\in L^2(0,\infty;\s Y)$, or, equivalently,
\[
  \mathcal{L}_\mu^* {T}_{G_\mu} \mathcal{L}_{\mu} f \in L_\mu^2(0,\infty;\s Y),
\]
which yields the desired invariance property for the operator under consideration. The boundedness of $\mathcal{L}_{\mu}^* {T}_{G_\mu} \mathcal{L}_{\mu}$ is obvious from the boundedness of $G$.
\end{Proof}

\begin{Remark}\label{r:bdd_der} Let $\nu>0$. By Lemma \ref{l:ran_h_nu}, the sets $\mathbb{C}_{\Re>\nu}$ and $B(r,r)$ with $r=1/(2\nu)$ are biholomorphically mapped to one another by $z\mapsto 1/z$. Hence, the conclusion in Corollary \ref{co:SW_weighted} might also be written as that (apart from continuous extendability of $\mathcal{M}$) there exists a unique bounded analytic $M\colon B(r,r)\to L(\s X)$ satisfying
\begin{equation}\label{eq:SW_w}
   (\mathcal{L}_{\mu}\mathcal{M}^\mu u)(\xi)=M\left(1/(i\xi+\mu)\right)\mathcal{L}_{\mu}u (\xi),\quad \xi\in \mathbb{R},\mu>1/(2r),
\end{equation}
 for all $u\in L_\mu^2(\mathbb{R};\s X)$, or, with Remark \ref{r:SW_OT},
 \begin{equation}\label{eq:re_SW}
       \mathcal{M}^\mu u = \mathcal{L}_{\mu}^*M(\hat{h}_{\mu})\mathcal{L}_{\mu}u\quad (\mu>\nu, u\in L_{\mu}^2(\mathbb{R};\s X)),
 \end{equation}
 where $\hat{h}_{\mu}(\xi)=1/(i\xi+\mu)$ and $z\mapsto M(\hat{h}_{\mu})(z)\coloneqq M(\hat{h}_{\mu}(z))$, $\mu>\nu$.
\end{Remark}

Next, we prove Theorem \ref{th:SW} along the lines of \cite{Weiss1991}. The first preparatory result is the following lemma. For $s\in \mathbb{C}$ we put 
\[
   e_{-s} \colon \mathbb{R}\to \mathbb{C}, t\mapsto\begin{cases}
                                                      0,& t<0,\\
                                                      e^{-s t}, & t\geq 0.
                                                   \end{cases}
\]
\begin{Lemma}[{{\cite[Remark 1.1]{Weiss1991}}}]\label{le:prep} Let $s\in \mathbb{C}_{\Re>0}$, $z\in L^2(0,\infty)$. Assume that for all $h\geq 0$, we get $\tau_h z= e^{-s h}z$ on $[0,\infty)$. Then there is a unique $a\in \mathbb{C}$ such that $z=ae_{-s}$. 
\end{Lemma}
\begin{Proof}
 At first note that $z\in L^1(0,\infty)$. Indeed, we compute
\begin{align*}
    \int_0^\infty |z(t)|\dd t & = \sum_{n=0}^\infty \int_n^{n+1} |z(t)|\dd t  = \sum_{n=0}^\infty \int_0^1 |z(n+t)|\dd t 
    \\ & = \sum_{n=0}^\infty \int_0^1 |e^{-sn}z(t)|\dd t  = \sum_{n=0}^\infty \int_0^1 e^{(-\Re s)n}|z(t)|\dd t 
    \\ & \leq \sum_{n=0}^\infty \left(\int_0^1 e^{(-\Re s)2n}\dd t\right)^{1/2}\|z\|_{L^2}
    \\ & = \sum_{n=0}^\infty \left(e^{(-\Re s)}\right)^n\|z\|_{L^2} =\|z\|_{L^2} \frac{1}{1-e^{-\Re s}}<\infty.
\end{align*}
With the help of the dominated convergence theorem, we realize that \[
Z\colon t\mapsto \int_{t}^\infty z(\tau)\dd \tau                                                                     
                                                                    \]
 is continuous. Moreover, for $t\geq 0$, we infer
\[
   Z(t)=\int_t^\infty z(\tau)\dd \tau = \int_0^\infty z(\tau+t)\dd \tau = e^{-s t}Z(0).
\]
We read off that $Z$ is indefinitely differentiable. But, taking $\phi\in C_c^\infty(0,\infty)$, we also realize that
\[
   \int_{0}^\infty \phi' Z = \int_0^\infty \phi'(t)\int_t^\infty z(\tau)\dd \tau \dd t= \int_0^\infty \int_0^\tau \phi'(t) \dd t\ z(\tau) \dd \tau = \int_{0}^\infty \phi z 
\]
which yields $s e^{-s t}Z(0)=-Z'(t)=z(t)$ for almost every $t\in \mathbb{R}_{\geq0}$. 
\end{Proof}
The author is indebted to Sascha Trostorff for providing the latter short proof.

\begin{Lemma}\label{le:prep_cont} Let $\s X$, $\s Y$ Hilbert spaces, $\Phi\colon \s X\to L^2(0,\infty;\s Y)$ bounded and (anti-) linear, $s\in \mathbb{C}_{\Re>0}$. Assume that for every $w\in \s X$ there exists a unique $v_w\in \s Y$ such that $\Phi(w)=e_{-s}v_w$. Then $\s X\ni w\mapsto v_w\in \s Y$ is  (anti-) linear and continuous. 
\end{Lemma}
\begin{Proof}
  (Anti-) linearity of $w\mapsto v_{w}$ follows from the respective property of $\Phi$. Next, let $w\in \s X$. We estimate
\begin{align*}
   \|v_w\|^2 &= 2\Re s \frac{\|v_w\|^2}{2\Re s} = 2\Re s \int_0^\infty \|v_w e^{-s t}\|^2\dd t
   \\        &= 2\Re s \int_0^\infty \|\Phi(w)(t)\|^2\dd t = 2\Re s \|\Phi(w)\|^2 \leq 2\Re s \|\Phi\|^2\|w\|^2,
\end{align*}
which yields the assertion.
\end{Proof}

\begin{Lemma}[{{\cite[Lemma 2.1]{Weiss1991}}}]\label{le:prep_vv} Let $\s X$ Hilbert space, $s\in \mathbb{C}_{\Re>0}$ and $z\in L^2(0,\infty;\s X)$. Assume that for all $h\geq 0$, we get $\tau_{h}z=e^{-s h}z$ on $[0,\infty)$. Then there exists a unique $v\in \s X$ such that $z=e_{-s}v$. 
\end{Lemma}
\begin{Proof}
  Define $\Phi\colon \s X\times L^2(0,\infty;\s X)\to L^2(0,\infty)$ by $\Phi(w,y)(t)\coloneqq \langle w,y(t)\rangle_{\s X}$ for a.e.~$t\in \mathbb{R}$, $w\in \s X$, $y\in L^2(0,\infty;\s X)$. The Cauchy-Schwarz inequality implies the boundedness of the sesquilinear mapping $\Phi$ with bound $1$. Furthermore, note that for $h\in \mathbb{R}_{\geq0}$, we have
\[
    \tau_h\Phi(w,z)=\Phi(w,\tau_h z)=\Phi(w,e^{-s h}z)=e^{-s h}\Phi(w,z)\quad (w\in \s X, z\in L^2(0,\infty;\s X)).
\]
Hence, by Lemma \ref{le:prep}, for any $w\in \s X$ there exists $a_w\in \mathbb{C}$ such that $\Phi(w,z)=a_{w}e_{-s}$. By Lemma \ref{le:prep_cont} applied to $\Phi(\cdot,z)$ in the position of $\Phi$ and $a_{(\cdot)}$ in place of $v_{(\cdot)}$, we infer that $v\colon \s X\ni w\mapsto a_w$ defines a bounded, anti-linear functional. We identify $v\in \s X$ with its Riesz image. It now remains to show $z=e_{-s}v$. For this, we observe for all $w\in \s X$ the equality
\begin{equation}\label{eq:prep_vv}
   \Phi(w,z)(t)=a_we_{-s}(t)=\langle  w, v\rangle e_{-s}(t) = \langle w,(ve_{-s})(t) \rangle = \Phi(w,(ve_{-s}))(t)
\end{equation}
for a.e.~$t\in \mathbb{R}$. Next, we observe for $\alpha\in L^2(0,\infty), w\in \s X, y\in L^2(0,\infty;\s X)$
\[
    \langle \Phi(w,y),\alpha\rangle_{L^2(0,\infty)} =\int_0^\infty \langle w,y(t)\rangle_{\s X}^* \alpha(t)\dd t
     = \int_0^\infty \langle y(t),w\rangle_{\s X} \alpha(t)\dd t
     = \langle y,\alpha w\rangle_{L^2(0,\infty;\s X)}.
\]
In particular, putting $\alpha=\1_E$ for some $E\subseteq [0,\infty)$ bounded, measurable,  and $y=z$, we get for $w\in \s X$ with the help of equation \eqref{eq:prep_vv}
\[
    \langle  z,\1_E w\rangle = \langle \Phi(w,z),\1_E\rangle=\langle \Phi(w,ve_{-s}),\1_E\rangle=\langle ve_{-s},\1_Ew\rangle.
\]
As the set $\{\1_Ew;E\subseteq [0,\infty)\text{ bounded, measurable}, w\in \s X\}$ is separating for the space $L^2(0,\infty;\s X)$, the assertion follows.
\end{Proof}

We are now in the position to prove Theorem \ref{th:SW}.

\begin{Proof}[of Theorem \ref{th:SW}]
We note that as an operator in $L^2(\mathbb{R};\s X)$ we have that $\tau_{-h}^*=\tau_h$ for all $h\in \mathbb{R}$. Thus, as $\mathcal{M}$ is translation-invariant, so is $\mathcal{M}^*$. In particular, for $w\in \s Y$ and $h\in \mathbb{R}_{\geq 0}$ denoting by $P_t$ the multiplication by $\1_{[t,\infty)}$, we get
\begin{align*}
   \tau_h P_0\mathcal{M}^*e_{-s}w & = P_{-h}\tau_h\mathcal{M}^* e_{-s}w 
   \\& = P_{-h}\mathcal{M}^*\tau_h e_{-s}w 
   \\&= P_{-h}\mathcal{M}^*e^{-s h}e_{-s}w 
   \\&=e^{-s h}P_{-h}\mathcal{M}^*e_{-s}w \quad (s\in \mathbb{C}_{\Re>0}).
\end{align*}
Thus,
\[ 
   \tau_h P_0\mathcal{M}^*e_{-s}w = e^{-s h}P_0\mathcal{M}^*e_{-s}w\text{ on }[0,\infty).
\] 
Hence, by Lemma \ref{le:prep_vv} applied to $z=P_0\mathcal{M}^*e_{-s}w$, for $s\in \mathbb{C}_{\Re>0}$ and for all $w\in \s Y$ there exists a unique $v_w\in \s X$ with $P_0\mathcal{M}^*e_{-s}w=e_{-s}v_w$. Applying Lemma \ref{le:prep_cont} to $\Phi(w)=P_0\mathcal{M}^*e_{-s}w$, we get that $G(s)\colon w\mapsto v_w$ is continuous for all $s\in \mathbb{C}_{\Re>0}$. In particular, note that we get for all $w\in \s Y$, $s\in \mathbb{C}_{\Re>0}$,
\begin{align*}
   \|G(s)w\|_{\s X} & =\sqrt{2\Re(s)}|\|e_{-s}G(s)w\|_{L^2(0,\infty;\s X)}
   \\ & =\sqrt{2\Re(s)}\|\mathcal{M}^*e_{-s}w\|_{L^2(0,\infty;\s X)}
   \\& \leq \sqrt{2\Re(s)}\|\mathcal{M}^*\|_{L(L^2(\s Y),L^2(\s X))}\|e_{-s}w\|_{L^2(0,\infty;\s Y)}
   \\&=\|\mathcal{M}\|_{L(L^2(\s X),L^2(\s Y))}\|w\|_{\s Y},
\end{align*}
which implies the boundedness of $\mathbb{C}_{\Re>0}\ni s\mapsto G(s)\in L(\s Y,\s X)$.

Next, let $u\in L^2(0,\infty;\s X)$. We observe for $s=i\xi+\nu\in \mathbb{C}_{\Re>0}$ and $w\in \s Y$ by the first part of the proof that
\begin{align}
   \langle w,\mathcal{L}_{\nu} \mathcal{M}u(\xi)\rangle_{\s Y} &= \langle e_{i\xi-\nu}w,\mathcal{M}u\rangle_{L^2(\s Y)} \notag
   \\ & = \langle \mathcal{M}^*e_{-s^*}w,u\rangle_{L^2(\s X)} \notag
   \\ & = \langle \mathcal{M}^*e_{-s^*}w,P_0 u\rangle_{L^2(\s X)} \notag
   \\ & = \langle P_0\mathcal{M}^*e_{-s^*}w, u\rangle_{L^2(\s X)} \notag
   \\ & = \langle e_{-s^*}G(s^*)w, u\rangle_{L^2(\s X)} \notag
   \\ & = \langle e_{-s^*}w,(t\mapsto G(s^*)^*u(t))\rangle_{L^2(\s Y)} \notag
   \\ & = \langle w,(\mathcal{L}_\nu G(s^*)^*u)(\xi)\rangle_{\s Y}=\langle w,G(s^*)^*(\mathcal{L}_\nu u)(\xi)\rangle_{\s Y} \label{eq:SW_sc}
\end{align}
Setting $M(i\xi+\nu)\coloneqq G((i\xi+\nu)^*)^*$, we infer the equality asserted in the theorem to be proved. For analyticity, we realize by putting $u=e_{-1}v$ for some $v\in \s X$ into \eqref{eq:SW_sc} and using $\mathcal{L}_{\nu}e_{-1}(\xi)=\frac{1}{\sqrt{2\pi}}\frac{1}{i\xi+\nu+1}=\frac{1}{\sqrt{2\pi}}\frac{1}{s+1}$,
 \begin{equation}\label{eq:SW_anal}
    \langle w, G(s^*)^*v\rangle_{\s Y}={\sqrt{2\pi}(s+1)}\langle e_{-s^*}w,\mathcal{M}u\rangle_{L^2(\s Y)}.
 \end{equation}
Since $u$ is supported on $[0,\infty)$, thus, so is $\mathcal{M}u$ by hypothesis. Hence, the mapping $\mathbb{C}_{\Re>0}\ni s \mapsto \langle e_{-s^*}w,\mathcal{M}u\rangle_{L^2(\s Y)}$ is analytic, by Lebesgue's dominated convergence theorem. Thus, the right-hand side of \eqref{eq:SW_anal} is analytic in $s$. As this argument applies to all $w\in \s Y$, $v\in \s X$, the mapping $M$ is analytic on the $L(\s X,\s Y)$-norming set $\{ N\mapsto \langle w,Nv\rangle_{\s Y}; w\in \s Y,v\in \s X\} \subseteq L(\s X,\s Y)'$. This, together with the boundedness of $M$, where the boundedness we have proved already, is sufficient for analyticity of $M$, see \cite[Proposition A.3]{Arendt2011}. This finishes the proof.
\end{Proof}

\section{Comments}

As already mentioned the idea of discussing the time derivative in weighted spaces dates back to (at least as early as) the work of Morgenstern (\cite{Morgenstern1952}), where an exponentially weighted norm on the space of continuous functions was used to deduce existence and uniqueness of solutions of ordinary differential equations. Since then many people studied differential equations in weighted spaces. However, the core idea of discussing the derivative operator with this particular weight has been less prominent. In the late 1980's Rainer Picard (\cite{Picard1989}) studied integral transforms and their relation to explicit spectral theorems for differential operators. In these studies it then turned out that the Fourier--Laplace transformation (see \eqref{e:FLT}) is the unitary transformation realizing the spectral representation as a multiplication operator for the derivative on the weighted $L^2$-type spaces $L_\nu^2(\R;\s X)$, $\nu\in\R$. 

We will demonstrate that the freedom in the choice of the parameter $\nu$ together with the estimate $\|\partial_{t,\nu}^{-1}\|\leq1/|\nu|$ yields an easily accessible way of discussing ordinary differential equations in a $L^2$-type setting. We sketch a solution theory for possibly nonlinear ordinary differential equations as follows. Let $F\colon \mathbb{C}^d\to \mathbb{C}^d$ be Lipschitz continuous, with the property $F(0)=0$. Then for all $\nu\in \R$, it is easy to see that the Nemitskii-operator $N_F$ induced by $F$, that is,
\[
   N_F \colon L_\nu^2(\mathbb{R};\mathbb{C}^d)\to L_\nu^2(\mathbb{R};\mathbb{C}^d), \phi\mapsto (t\mapsto F(\phi(t)))
\]
is Lipschitz continuous as well. Next, let $f\in L^2(0,\infty;\mathbb{C}^d)\subseteq \bigcap_{\nu>0}L_\nu^2(\mathbb{R};\mathbb{C}^d)$, we want to find $u\colon \R\to\mathbb{C}^d$, which is locally weakly differentiable such that
\begin{equation}\label{eq:ode0}
    u'(t) = F(u(t))+f(t)\quad (t\in \mathbb{R}).
\end{equation}
The latter equation interpreted in $L_\nu^2(\R;\mathbb{C}^d)$ is the same as
\[
   \partial_{t,\nu} u = N_F(u)+f.
\]
Thus, multiplying by $\partial_{t,\nu}^{-1}$, we get
\[
      u = \partial_{t,\nu}^{-1}N_F(u)+\partial_{t,\nu}^{-1}f\eqqcolon \Phi_\nu(u),
\]
which is a fixed point problem for $u\in L_\nu^2(\mathbb{R};\mathbb{C}^d)$. Once we show that there exists $\nu>0$ such that $\Phi_\nu$ is a strict contraction, the existence of a weakly differentiable function $u$ for \eqref{eq:ode0} is warranted. But, if $|F|_{\textnormal{Lip}}$ is a Lipschitz constant for $F$, then $|F|_{\textnormal{Lip}}$ is also a Lipschitz constant for $N_F$. Hence, for $u,v\in L_\nu^2(\R;\mathbb{C}^d)$, we obtain
\begin{multline*}
  \|\Phi_\nu (u)-\Phi_\nu (v)\|_{L_\nu^2}= \|\partial_{t,\nu}^{-1}N_F(u)-\partial_{t,\nu}^{-1}N_F(v)\|_{L_\nu^2}
  \\ \leq \frac{1}{\nu}\|N_F(u)-N_F(v)\|_{L_\nu^2}\leq \frac{|F|_{\textnormal{Lip}}}{\nu}\|u-v\|_{L_\nu^2}.
\end{multline*}
Therefore, choosing $\nu$ large enough yields that $\Phi_\nu$ is a strict contraction and \eqref{eq:ode0} can be solved with the help of the contraction mapping theorem. The application presented is a first of many others. In fact, it is possible to extend these ideas in such a way that it leads to a unified solution theory for delay differential equations, that include discrete and continuous delays as well as neutral equations, see \cite{Kalauch}. 

This line of ideas has further applications. Indeed, starting out from the observation that $\partial_{t,\nu}^{-1}$ is the convolution with the Heaviside function (see Theorem \ref{t:td_inv}), it is possible to extend this unified way of looking at problems of the ordinary delay differential type to a Banach space setting, see \cite{Waurick2014OAM_DelayLp}. 

Another class tractable with this approach is the class of so-called (ordinary) integro-differential equations, that is, equations where the right-hand side of \eqref{eq:ode0} is replaced by some integral expression involving $u$. A prominent example are convolutions with respect to the time variable. These convolutions may be written as a multiplication by some function in Fourier space by the convolution theorem, see Remark \ref{r:PL}. Moreover, they can be represented as certain functions of $\partial_{t,\nu}^{-1}$ as in \eqref{eq:re_SW} in combination with Corollary \ref{co:spe_thm}, that is, in the form
\begin{equation}\label{eq:fnoftime}
   M(\partial_{t,\nu}^{-1})\coloneqq \s L_\nu^*M(\hat{h}_\nu)\s L_\nu
\end{equation}
for some appropriate $M$. We refer to \cite{Trostorff2015,Waurick2014MMAS_Frac,Waurick2014SIAM_HomFrac} for a thorough treatment of this type of functions of $\partial_{t,\nu}^{-1}$ and applications to partial differential equations. A specific integral operator is the fractional time derivative (\cite{Waurick2014MMAS_Frac,Waurick2014SIAM_HomFrac}), $\partial_{t,\nu}^\alpha$, which also falls into the class of operators that admit a representation as in \eqref{eq:fnoftime} for some bounded, analytic function $M$. Indeed, for $\alpha\in (-\infty,0]$ the operator $\partial_{t,\nu}^\alpha$ is translation-invariant and causal by Theorem \ref{t:SW_conv}. A detailed look at $\partial_{t,\nu}^\alpha$ shows that, in fact, if applied to functions being supported on $(0,\infty)$ only, the resulting operator is the Riemann-Liouville fractional derivative, see also \cite[p 3143]{Waurick2014MMAS_Frac} or \cite{Tsangi2014}. It should be noted that this fractional derivative has also been discussed with regards to numerical analysis, see \cite{Rubisch2015}.

The representation theorem proved in this chapter, Theorem \ref{th:SW}, has major applications in control and system theory. It is used in the representation theory for shift-invariant causal systems and is related to the description of the input-output relation of control systems by means of their so-called transfer function. We refer to the references in \cite{Weiss1991} for a more detailed account of that.

\cleardoublestandardpage

\chapter{Evolutionary Mappings and Causality}\label{ch:EMC}

In this chapter we provide the notion of evolutionary mappings. With regards to applications discussed later on, it is of interest whether the operators considered in certain $L_\nu^2$-type spaces are actually ``independent'' of $\nu$: Given an operator $\s S$, which is densely defined and continuous in both the spaces $L_\nu^2(\mathbb{R})$ and $L_\mu^2(\mathbb{R})$ for $\mu\neq \nu$, it is a priori unclear, whether the respective continuous extensions, denoted by $\s S^\nu$ and $\s S^\mu$ coincide on the intersection of their respective domains, that is, on $L_\nu^2(\mathbb{R})\cap L_\mu^2(\mathbb{R})$. The most prominent examples of such $\s S$ are solution operators to certain (abstract) partial differential equations. Defining evolutionary mappings in Section \ref{s:em}, we will obtain the desired independence result upon assuming one additional property for $\s S$. The additional property is (forward) causality, which is defined in Section \ref{s:cem} for evolutionary mappings.

\renewcommand{\baselinestretch}{0.65}\normalsize\mysection{Evolutionary Mappings}{Evolutionary Mappings}{standard causal domains $\s D(\s X)$, $\s D_\nu(\s X)$ $\cdot$ standard evolutionary mappings}\label{s:em}

\renewcommand{\baselinestretch}{1}\normalsize
In Remark \ref{r:SW_OT}, we pointed out an important property of causal, translation-invariant mappings, namely the property of being extendable to continuous operators on $L_{\mu}^2$ for all $\mu>\nu$. Though lacking the property of translation-invariance, multiplication operators are still causal and enjoy the same property of being extendable to the $L_\nu^2$-scale for all $\nu\in \mathbb{R}$:

\begin{Example}\label{ex:mult_op} Let $\s X$, $\s Y$ Hilbert spaces, $\nu\in \mathbb{R}$, and let $M\colon \mathbb{R}\to L(\s X, \s Y)$ be strongly measurable and bounded. Then 
\[
   \mathcal{M} \colon L_\nu^2(\mathbb{R};\s X) \to L_\nu^2(\mathbb{R};\s Y),f\mapsto (t\mapsto M(t)f(t))
\]
defines a bounded linear operator satisfying 
\[
 \|\mathcal{M}\|\leq \esssup_{t\in \mathbb{R}}\|M(t)\|.
\]
Assume, in addition, that $M$ is strongly differentiable, with bounded derivative, that is, for $u\in \s X$ the mapping $t\mapsto M(t)u$ is differentiable and $t\mapsto (M(\cdot)u)'(t)$ is bounded for every $u$. As a pointwise limit of measurable mappings, the mapping $(t\mapsto (M(\cdot)u)'(t))$ is measurable itself. Hence, $\s X \ni u\mapsto (t\mapsto (M(\cdot)u)'(t))\in L^\infty(\mathbb{R};\s Y)$ is well-defined. Moreover, by the mean-value inequality, this mapping is closed and continuous by the closed graph theorem. Thus, the mapping
\[
    \mathcal{M}'\colon L_\nu^2(\mathbb{R};\s X)\to L_\nu^2(\mathbb{R};\s Y), v\mapsto (t\mapsto (M(\cdot)v(t))'(t))
\]
is well-defined and continuous. Moreover, for $f\in \dom(\partial_{t,\nu})$ we get $\mathcal{M}f\in \dom(\partial_{t,\nu})$ and $\partial_{t,\nu}\mathcal{M}f=\mathcal{M}'f+\mathcal{M}\partial_{t,\nu}f$.
\end{Example}

Note that the operator norm of the multiplication operator introduced in Example \ref{ex:mult_op} is \emph{independent} of the chosen $\nu\in \mathbb{R}$. We recall that causal, translation-invariant mappings have an operator norm being decreasing in $\nu$. In fact, this is the upshot of Corollary \ref{co:SW_weighted}. Given an operator acting on the whole scale $(L_\nu^2(\mathbb{R}))_{\nu\in\mathbb{R}}$, the operator norm need not be decreasing in general, as the following example shows.

\begin{Example}\label{ex:tt} Let $h\in \mathbb{R}$ and for $\nu\in \mathbb{R}$ let $\tau_{h,\nu}\colon L_\nu^2(\mathbb{R})\to L_\nu^2(\mathbb{R}),f\mapsto f(\cdot+h)$. Then, $\|\tau_{h,\nu}\|=e^{h\nu}$. This implies that $\|\tau_{h,\nu}\|\to \infty$ as long as $\nu h>0$ and $|\nu|\to\infty$. Hence, in particular, $\|\tau_{h,\nu}\|\to\infty$ if $h>0$ and $\nu\to\infty$. Note that the operator $\tau_{h,\nu}$ is causal if and only if $h<0$. 
\end{Example}

We raise the idea of extendability to an operator acting continuously on $L_\nu^2(\mathbb{R};\s X)$ for all $\nu$ large enough with operator norm ``fairly independent'' of $\nu$ to the main definition of this section:

\begin{Definition}[evolutionary mappings]\label{d:em}
Let $\s X,\s Y$ Hilbert spaces, $\nu\in \mathbb{R}$. We call a linear
mapping 
\begin{equation}
\s S\colon \dom(\s S)\subseteq\bigcap_{\mu\geq \nu}L_{\mu}^{2}(\mathbb{R};\s X)\to\bigcap_{\mu\geq\nu}L_{\mu}^{2}(\mathbb{R};\s Y)\label{eq:evolutionary_map}
\end{equation}
\emph{evolutionary (at $\nu$)}, if $\dom(\s S)$ is dense in $L_{\mu}^{2}(\mathbb{R};\s X)$ for all $\mu\geq\nu$, $\s S$ extends to a bounded linear operator $\s S^\mu$ from $L_{\mu}^{2}(\mathbb{R};\s X)$
to $L_{\mu}^{2}(\mathbb{R};\s Y)$ for all $\mu\geq\nu$ and
is such that
\[
\limsup_{\mu\to\infty}\left\Vert\s S^\mu\right\Vert _{L(L_{\mu}^{2}(\mathbb{R};\s X),L_{\mu}^{2}(\mathbb{R};\s Y))}<\infty.
\]
 The continuous extension of $\s S$ to some $L_{\mu}^{2}$ will be denoted by $\s S^\mu$, and, if there is no risk of confusion, we will re-use the notation $\s S$.
We set
\begin{align*}
L_{\textnormal{ev},\nu}(\s X,\s Y)&\coloneqq\{\s S;\s S\text{ is as in (\ref{eq:evolutionary_map}) and is evolutionary at }\nu\};
\\L_{\textnormal{ev},\nu}(\s X)&\coloneqq L_{\textnormal{ev},\nu}(\s X,\s X).
\end{align*}
\end{Definition}
Note that $L_{\textnormal{ev},\nu}(\s X,\s Y)\subseteq L_{\textnormal{ev},\mu}(\s X,\s Y)$
for all $\mu\geq\nu$. 

\begin{Example}\label{ex:evo_map} Let $\s X$, $\s Y$ Hilbert spaces.

(a) Let $\nu\in \mathbb{R}$, $\mathcal{M}\in L(L_{\nu}^2(\s X),L_{\nu}^2(\s Y))$ be translation-invariant and causal. Then $\mathcal{M}$ with $\dom(\s M)=L_{\textnormal{c}}^2(\mathbb{R};\s X)$ is evolutionary at $\nu$. Indeed, this is Corollary \ref{co:SW_weighted}.

(b) Let $M\colon \mathbb{R}\to L(\s X)$ be strongly measurable and bounded. Then the induced multiplication operator $\mathcal{M}$ as introduced in Example \ref{ex:mult_op} is evolutionary if endowed with the domain $L_{\textnormal{c}}^2(\mathbb{R};\s X)$.

(c) Let $h\in \mathbb{R}$. Then the time-shift $\tau_h$ with the domain $L_{\textnormal{c}}^2(\mathbb{R};\s X)$ is evolutionary if and only if $h\leq 0$, see Example \ref{ex:tt}. 
\end{Example}

Evolutionary mappings will play the central role throughout this exposition. In applications, constitutive relations or material laws can be realized as evolutionary mappings. Moreover, we will show that for given (ordinary/partial) differential equations modeling physical processes the respective solution operators will be evolutionary mappings itself. So, it is of interest to study sum, product and inverses of evolutionary mappings and to address the question, whether the resulting operators are evolutionary again. However, we want to stress two subtleties in this context.

 (a) If we are given an evolutionary mapping $\s S$, its adjoint $\s S^*$ is, in general, not evolutionary again: Indeed, acting on different Hilbert space, it is a priori unclear what operator an adjoint of $\s S$ would be. So, in fact, the adjoint of $\s S$ computed in $L_\nu^2(\mathbb{R};\s X)$ and $L_\mu^2(\mathbb{R};\s X)$ for both $\mu,\nu$ large enough, $\mu\neq \nu$, might differ from one another. An example is the time-shift again, see Example \ref{ex:tt}: $\tau_{h,\nu}^*= \tau_{-h,\nu}e^{2h\nu}$, $\nu,h\in \mathbb{R}$.

 (b) The sum of two evolutionary mappings need not be densely defined any more. Indeed, take the time-shift $\tau_h$ for $h<0$. Consider $\tau^{(1)}_h$ and $\tau^{(2)}_h$ as the operator acting as $\tau_h$ but with 
 \begin{equation}\label{eq:domt1}
   \dom(\tau^{(1)}_h)=L_{\textnormal{c}}^2(\mathbb{R})
   \end{equation}
 and
 \begin{equation}\label{eq:domt2}
    \dom(\tau^{(2)}_h)=\lin \{t\mapsto t^n e^{-(t+\delta)^2/2};\delta\in \mathbb{R}, n\in \mathbb{N}_0\}. 
 \end{equation}
 Then, by the density of (linear combinations of) Hermite functions in $L^2(\mathbb{R})$, it is easy to see that $\dom(\tau^{(2)}_h)$ is dense in $L_{\nu}^2(\mathbb{R})$ for all $\nu\in \mathbb{R}$. Also note that only the zero element in $\dom(\tau^{(2)}_h)$ is compactly supported: For this let $g\in \dom(\tau^{(2)}_h)$ have compact support. Then there are polynomials $p_1,\ldots,p_n\colon \mathbb{R}\to\mathbb{C}$ and real numbers $-\infty<\delta_1<\cdots<\delta_n<\infty$ with 
 \[
   g(t) = \sum_{j=1}^n p_j(t) e^{-(t+\delta_j)^2/2}=\sum_{j=1}^n p_j(t) e^{-t^2/2-\delta_j t -\delta_j^2/2} \quad(t\in \mathbb{R}).
 \]
  For $t\in \mathbb{R}$ large enough, setting $q_j\coloneqq e^{-\delta_j^2/2} p_j$, $j\in \{1,\ldots,n\}$, we get 
  \[
    e^{t^2/2+\delta_1 t}g(t) = 0 = \sum_{j=1}^n q_j(t) e^{-(\delta_j-\delta_1) t}= q_1(t)+ \sum_{j=2}^nq_j(t) e^{-(\delta_j-\delta_1) t}.
  \]
  Hence, as $t\to\infty$ we have $\sum_{j=2}^nq_j(t) e^{-(\delta_j-\delta_1) t}\to 0$ and, thus, $q_1(t)=0$. Continuing in this manner, we infer $q_2=\cdots=q_n=0$. Thus, $g=0$. But, all functions in $L_{\textnormal{c}}^2(\mathbb{R})$ are compactly supported, and so $\dom(\tau_h^{(1)}+\tau_h^{(2)})=\{0\}$. Hence,
  \[
     \tau_h^{(1)}+\tau_h^{(2)}=0\subset 2\tau_h.
  \]
In order to circumvent the last problem, we will seek a possibility to endow an evolutionary mapping with a standard domain. It turns out that this can be done, if we assume causality for the evolutionary mapping under consideration. However, note that in Definition \ref{d:cau}, we have defined causality for closed continuous mappings only. The definition of an adapted version of causality for closable operators is postponed to Section \ref{s:ccm}. But, a first link of evolutionarity and causality can be given right away.

\begin{Remark}\label{evo=caus} Let $\s X$, $\s Y$ Hilbert spaces, $\nu\in \mathbb{R}$, $\s S\in  L_{\text{ev},\nu}(\s X,\s Y)$ and assume that for all $t\in \R$ the set $\dom(\s SQ_t)\cap \dom(\s S)$ is dense in $L_\mu^2(\R;\s X_0)$ for $\mu\geq\nu$\footnote{$\dom(\tau_h^{(1)})$ from \eqref{eq:domt1} meets and $\dom(\tau_h^{(2)})$ from \eqref{eq:domt2} violates this assumption.}, where $Q_t$ is the operator $T_{\1_{(-\infty,t]}}$ of  multiplication  by $\1_{(-\infty,t]}$. Then $\s S^\mu$ is causal for all $\mu\geq\nu$: Since $L_{\textnormal{ev},\nu}\subseteq L_{\textnormal{ev},\mu}$, it suffices to prove that $\s S^\nu$ is causal. For this, let $f\in L_\nu^2(\R;\s X)$, $t\in\R$ and assume that $Q_tf=0$. We choose $(f_n)_{n}$ in $\dom(\s SQ_t)\cap \dom(\s S)$ such that $f_n\to f$ in $L_\nu^2(\R;\s X)$ as $n\to\infty$. For all $n\in\N$ we have $Q_tf_n=f_n-P_tf_n \in \dom(\s S)$, $P_t\coloneqq 1-Q_t=T_{\1_{(t,\infty)}}$, and 
\[
   P_tf_n \to P_tf= f\quad (n\to\infty)
\]
 in $L_\nu^2(\R;\s X)$. In particular, the latter implies that $g_n\coloneqq P_tf_n$ approximates $f$ in $L_{\mu}^2(\R;\s X)$ for all $\mu\geq \nu$. Now, we follow the idea of \cite[Proof of Theorem 4.5]{Kalauch}. For this let $\phi\in L_{\textnormal{c}}^2(\mathbb{R};\s Y)$ with support bounded above by $t$. For $\mu\geq\nu$ we get that 
\begin{align*}
   \abs{\langle \s S^\nu f , \phi \rangle_{L^2(\R;\s Y)}} & = \lim_{n\to\infty}\abs{\langle \s S^\nu g_n,\phi\rangle } \\
                                            & = \lim_{n\to\infty}\abs{\langle\s S g_n,\phi\rangle } \\
                                            & \leq \lim_{n\to\infty}\|\s S g_n\|_{L_{\mu}^2}\|{\phi}\|_{L_{-\mu}^2} \\
                                            & \leq \Abs{\s S}_{L(L_{\mu}^2)} \|{f}\|_{L_\mu^2}\|{\phi}\|_{L^2}e^{\mu t} \\
					    & = \Abs{\s S}_{L(L_{\mu}^2)} \|{f(\cdot +t)}\|_{L_\mu^2}\|{\phi}\|_{L^2}  
\end{align*}
Letting $\mu\to\infty$, we deduce that $\langle \s S^\nu f , \phi \rangle_{L^2(\s Y)}=0$. Hence, $\s S^\nu f=0$ on $(-\infty,t]$, that is, $Q_t \s S^\nu f=0$ or $\s S^\nu P_tf=P_t\s S^\nu P_tf$ since $P_tf=f$. The latter yields causality, see Definition \ref{d:cau}.
\end{Remark}

The key observation of the latter remark is that certain conditions on the domain of evolutionary mappings result in the causality of the closure of these mappings. A prototype of such a domain is given next.

\begin{Definition}\label{d:scd_sev} Let $\s X$, $\s Y$ Hilbert spaces, $\nu\in \mathbb{R}$. Then the set
\[
   \mathcal{D}_{\nu}(\s X) \coloneqq \bigcap_{\mu\geq\nu} L_\mu^2(\mathbb{R};\s X)
\]
is called \emph{standard causal domain (at $\nu$)}; we set $\s D(\s X)\coloneqq \bigcap_{\nu\in \mathbb{R}} \s D_\nu(\s X)$. We call a map $\s S\in L_{\textnormal{ev},\nu}(\s X,\s Y)$ \emph{standard evolutionary (at $\nu$)}, if $\dom(\s S)$ is the standard causal domain (at $\nu$). We define the set of all standard evolutionary mappings
\begin{align*}
  L_{\textnormal{sev},\nu}(\s X,\s Y)&\coloneqq \{ \s S; \s S \text{ standard evolutionary at }\nu\}, \\
  L_{\textnormal{sev}}(\s X,\s Y)&\coloneqq \bigcup_{\nu\in \mathbb{R}}L_{\textnormal{sev},\nu}(\s X,\s Y),\\
  L_{\textnormal{sev},\nu}(\s X)&\coloneqq L_{\textnormal{sev},\nu}(\s X,\s X),\quad L_{\textnormal{sev}}(\s X)\coloneqq L_{\textnormal{sev}}(\s X,\s X).
\end{align*}
\end{Definition}

Standard evolutionary mappings are closed under vector space operations and composition:

\begin{Proposition}\label{p:sev_alg} Let $\s X$, $\s Y$, $\s Z$ Hilbert spaces, $\nu\in \mathbb{R}$, $\s S,\s T\in L_{\textnormal{sev},\nu}(\s X,\s Y)$, $\s U\in L_{\textnormal{sev},\nu}(\s Y,\s Z)$, $\alpha\in \mathbb{C}$. Then
\begin{enumerate}[label=(\alph*)]
 \item\label{alg1} $\s S+\alpha \s T\in L_{\textnormal{sev},\nu}(\s X,\s Y)$,
 \item\label{alg2} $\s U\s S \in L_{\textnormal{sev},\nu}(\s X,\s Z)$.
\end{enumerate}
\end{Proposition}
\begin{Proof}
 The statement in \ref{alg1} is easy. For \ref{alg2}, we observe that if $f\in \mathcal{D}_{\nu}(\s X)=\dom(\s S)$, then, by the evolutionarity of $\s S$, $\s S f\in \bigcap_{\mu\geq \nu}L_\nu^2(\s Y)=\mathcal{D}_{\nu}(\s Y)=\dom(\s U)$. The remaining norm estimate follows from the submultiplicativity of the operator norm.
\end{Proof}

We will elaborate more on evolutionary mappings once we discussed causality for closable mappings in the next section.

\renewcommand{\baselinestretch}{0.65}\normalsize\mysection{Causality for Closable Mappings}{Causality for Closable Mappings}{resolution space $\cdot$ causal $\cdot$ characterization of densely defined operators with causal closure $\cdot$ strongly causal $\cdot$ Theorem \ref{t:csc}}\label{s:ccm}

\renewcommand{\baselinestretch}{1}\normalsize
In order to motivate the upcoming notion of causality, we recall that for a Hilbert space $\s X$, and some $\nu\in \mathbb{R}$, we say a mapping $\s S \in L(L_\nu^2(\mathbb{R};\s X))$ is \emph{causal} (see Definition \ref{d:cau}), if 
\begin{equation} \label{eq:cau}
  P_t \s S P_t = \s S P_t \quad(t\in \mathbb{R}),
\end{equation}
 where $P_t$ is the operator of  multiplication by $\1_{(t,\infty)}$. If now $\s S$ is defined on a proper domain in $L_\nu^2(\mathbb{R};\s X)$ the way of defining causality just mentioned has the drawback that $\dom(\s SP_t)$ may only consist of the $0$ function, see \eqref{eq:domt2} for a possible domain. Hence, \emph{every} continuous mapping endowed with such a domain would be causal, if \eqref{eq:cau} characterized causality also for closable mappings. Anticipating the latter, we seek a different notion of causality, which for closed, continuous maps yields the same. For this, we discuss \eqref{eq:cau} in more detail: Introducing $Q_t\coloneqq 1-P_t$, we get that, equivalently to \eqref{eq:cau}, for all $t\in \mathbb{R}$  the equality
 \begin{equation}\label{eq:cau2}
  Q_t\s S = Q_t \s S Q_t
 \end{equation}
 or, for all $f\in L_\nu^2(\mathbb{R};\s X)$, the implication
 \begin{equation}\label{eq:cau3}
    Q_tf = 0\Longrightarrow Q_t \s S f =0
 \end{equation}
 hold true. 
 
 Yet another reformulation of \eqref{eq:cau3} or \eqref{eq:cau2} is that for all $\phi\in B_{L_\nu^2(\s X)}$, $B_{L_\nu^2(\s X)}$ the unit ball of $L_\nu^2(\s X)$, there exists $C\geq 0$ such that for all $f\in L_\nu^2(\mathbb{R};\s X)$ we have
 \begin{equation}\label{eq:cau4}
  |\langle Q_t \s S f,\phi\rangle|\leq C\|Q_t f\|.
 \end{equation}
 Indeed, the latter estimate follows from \eqref{eq:cau2} and implies \eqref{eq:cau3}, which, in turn, implies \eqref{eq:cau2}. The continuity estimate in \eqref{eq:cau4}, however, is the starting point for defining causality for closable mappings. 
 
 Beforehand, we introduce the concept of a resolution space. 

\begin{Definition}[{\cite{Saeks1970}}] Let $\s X$ be a Hilbert space, let $(Q_t)_{t\in \R}$ in $L(\s X)$ be a \emph{resolution of the identity}, that is, for all $t\in \R$ the operator $Q_t$ is an orthogonal projection, $\ran(Q_t)\subseteq \ran(Q_s)$ if and only if $t\leq s$ and $Q_t$ converges in the strong operator topology to $0$ and $1$ if $t\to -\infty$ and $t\to\infty$, respectively. The pair $(\s X,(Q_t)_t)$ is called \emph{resolution space}. 
\end{Definition}

In what follows we provide the notion of causality for closable mappings. We stick to the linear case here. For a possible way to define the respective concept for non-linear mappings as well, we refer to \cite{Waurick2015IM_Caus}.

\begin{Definition}\label{d:stc} Let $(\s X,(Q_t)_t)$, $(\s Y,(R_t)_t)$ resolution spaces, $D\subseteq \s Y$, $S\colon \dom(S)\subseteq \s X\to \s Y$ linear. We say that $S$ is \emph{causal} on $D$, if for all $r>0$, $t\in \R$, $\phi\in D$ there exists $C\geq 0$ such that for all $f\in B_S(0,r)$ we have
\[
   |\langle R_t S f,\phi\rangle|\leq C\| Q_t f\|,
\]where $B_S(0,r)\coloneqq \{ f\in \dom(S); \|f\|^2+\|Sf\|^2< r^2\}$.
If $D=\s Y$, then we say that $S$ is \emph{causal}.
\end{Definition}

\begin{Remark}\label{r:cau_cont} (a) Another way of expressing causality in Definition \ref{d:stc} is that the mapping
\begin{align*}
               \left( B_{S}(0,r), \abs{ Q_t\left(\cdot-\cdot\right)}\right) &\to \left(\s Y, \abs{\langle R_t\left(\cdot-\cdot\right),\phi\rangle}\right) \\
   f\mapsto Sf,
\end{align*} 
is Lipschitz continuous for all $\phi\in D$, $r>0$.

(b) If in Definition \ref{d:stc} the operator $S$ is also continuous, then $B_S(0,r)$ may be replaced by $B_{\s X}(0,r)\cap \dom(S)$, the open $r$-ball in $\s X$ intersected with the domain of $S$. Indeed, by continuity of $S$, we have $B_S(0,r)\subseteq B_{\s X}(0,r)\cap \dom(S)\subseteq B_S(0,Cr)$ for some $C>0$.
\end{Remark}

The notion of causality for closable operators coincides with the one for closed operators:

\begin{Theorem}[{{\cite[Theorem 1.7]{Waurick2015IM_Caus}}}]\label{thm:causal_cont_causal} Let $(\s X,(Q_t)_t)$, $(\s Y,(R_t)_t)$ resolution spaces. Let the operator $S\colon \dom(S)\subseteq \s X\to \s Y$ be linear and closable. Then the following statements are equivalent:
\begin{enumerate}[label=(\roman*)]
 \item\label{caus_cont_caus1} for all $f\in \dom(\overline{S})$, $t\in \mathbb{R}$, we have that $Q_tf=0$ implies $R_t\overline{S}f=0$;
 \item\label{caus_cont_caus2} $S$ is causal;
 \item\label{caus_cont_caus3} there exists $D\subseteq \s Y$ dense such that $S$ is causal on $D$.
\end{enumerate}
\end{Theorem}

Theorem \ref{thm:causal_cont_causal} also admits a generalization to the Banach space case. In this exposition, however, it is sufficient to consider  the Hilbert space case only. The somewhat more involved version of Theorem \ref{thm:causal_cont_causal} (including an adapted version of causality) for the general Banach space case can be found in \cite{Waurick2015IM_Caus}.

For the proof of Theorem \ref{thm:causal_cont_causal}, we need some prerequisites. Note that the next lemma has already been proven in the first few lines of this section for the case of continuous $S$ with $\dom(S)=\s X$:

\begin{Lemma}[{{\cite[Lemma 1.9]{Waurick2015IM_Caus}}}]\label{le:equiv_of_caus_for_weakly-closed} Let $(\s X,(Q_t)_t)$ and $(\s Y,(R_t)_t)$ be resolution spaces. Let $S\colon \dom(S)\subseteq \s X\to \s Y$ linear and closed.
Then the following assertions are equivalent:
\begin{enumerate}[label=(\roman*)]
 \item\label{equiv_cau_clo1} for all $f\in \dom(S)$, $t\in \mathbb{R}$, we have that $Q_tf=0$ implies $R_tSf=0$;
 \item\label{equiv_cau_clo2} $S$ is causal;
 \item\label{equiv_cau_clo3} there exists $D\subseteq \s Y$ dense such that $S$ is causal on $D$.
\end{enumerate} 
\end{Lemma}
\begin{Proof}
 The implication ``\ref{equiv_cau_clo2}$\Rightarrow$\ref{equiv_cau_clo3}'' is trivial. Assume that \ref{equiv_cau_clo3} holds. Let $f\in \dom(S)$ with $Q_tf=0$ for some $t\in \mathbb{R}$. By hypothesis, for all $t\in \mathbb{R}$ and $\phi\in D$, we find $C\geq 0$ such that
 \[
    |\langle R_tS g,\phi\rangle|\leq C\|Q_t g\| \quad (g\in B_{S}(0,\|f\|+1))
 \]  So, $\langle R_t S f, \phi\rangle=0$. As $\phi\in D$ is arbitrary, and $D$ is dense, we get $R_t S f=0$. Hence, \ref{equiv_cau_clo1} follows.

 For the sufficiency of \ref{equiv_cau_clo1} for \ref{equiv_cau_clo2}, we show that $S$ violates the condition stated in \ref{equiv_cau_clo1} provided $S$ is not causal. For this, let $r>0$, $t\in \R$, $\phi\in \s X$ and $\eps>0$ such that for all $n\in\N$ there is $f_n\in B_S(0,r)$ with 
\[
   \|Q_tf_n\|\leq\frac{1}{n}\text{ and }\abs{\langle R_tSf_n,\phi\rangle}\geq \eps.
\]
By boundedness of $(f_n)_n$ and $(Sf_n)_n$, there exists a subsequence $(n_k)_k$ of $(n)_n$, such that $(f_{n_k})_k$, and $(Sf_{n_k})_k$ weakly converge. By linearity and closedness of $S$, $S$ is weakly closed. Hence, we deduce that $f\coloneqq \textnormal{w-}\lim_{k\to\infty}f_{n_k}\in \dom(S)$ and $\textnormal{w-}\lim_{k\to\infty}Sf_{n_k}=Sf$. By (weak) continuity of $Q_t$ we get 
\[
  \|{Q_tf}\|\leq \liminf_{k\to\infty} \|Q_tf_{n_k}\|=0.
\]
 Now, from 
\[
 \abs{\langle R_tSf,\phi\rangle}=\lim_{k\to\infty}\abs{\langle R_tSf_{n_k},\phi\rangle}\geq \eps
\]
we read off that $S$ does not satisfy \ref{equiv_cau_clo1}.
\end{Proof}

\begin{Lemma}[{{\cite[Lemma 1.10]{Waurick2015IM_Caus}}}]\label{le:strongly_causal_is_equiv_strong_closure} Let $(\s X,(Q_t)_t), (\s Y,(R_t)_t)$ be resolution spaces, and $S\colon \dom(S)\subseteq \s X\to \s X$ closable, $D\subseteq \s X$. Then the following statements are equivalent.
\begin{enumerate}[label=(\roman*)]
 \item\label{strw1} $S$ is causal on $D$;
 \item\label{strw2} $\overline{S}$ is causal on $D$.
\end{enumerate} 
\end{Lemma}
\begin{Proof}
Let $r>0$, $t\in \R$. Then $B_S(0,r)$ is dense in $B_{\overline{S}}(0,r)$ with respect to $\abs{Q_t(\cdot-\cdot)}$. Indeed,  for $\eps>0$, $f\in B_{\overline{S}}(0,r)$ there exists $g\in B_S(0,r)$ such that 
\[
 \abs{f-g}+\abs{\overline{S}f-\overline{S}g}\leq \eps. 
\]
In particular, we have $\abs{Q_t(f-g)}\leq \Abs{Q_t}\eps$. Assuming the validity of \ref{strw1}, we see that
\[
   (B_{\overline{S}}(0,r),\abs{Q_t(\cdot-\cdot)})\to (\s X,\abs{\langle R_t (\cdot-\cdot),\phi\rangle}), f\mapsto \overline{S}f
\]
is Lipschitz continuous on the dense subset $B_{{S}}(0,r)$ for all $\phi\in D$. This implies \ref{strw2}, see also Remark \ref{r:cau_cont}. The converse is trivial.
\end{Proof}

\begin{Proof}[of Theorem \ref{thm:causal_cont_causal}]
 By Lemma \ref{le:equiv_of_caus_for_weakly-closed}, condition \ref{caus_cont_caus1} is equivalent to causality of $\overline{S}$ and to causality of $\overline{S}$ on some dense set, the latter two properties are, in turn, equivalent to causality of $S$ (condition \ref{caus_cont_caus2}) and causality of $S$ on some dense set $D$ (condition \ref{caus_cont_caus3}), respectively, by Lemma \ref{le:strongly_causal_is_equiv_strong_closure}.
\end{Proof}

For applications discussed later on, we give an instant of an example for closable causal mappings.

\begin{Proposition}\label{p:inverse_strongly_causal}
  Let $(\s X, (Q_t)_t)$ a resolution space, $S\colon \dom(S)\subseteq \s X\to \s X$ linear. Assume for every $t\in \mathbb{R}$ the inequality
\[
    |\langle Q_t Sf,f\rangle| \geq  c\langle Q_t f,f\rangle
\]
for some $c>0$ and all $f\in \dom(S)$. Then $S^{-1}$ defines a linear mapping and is causal. Furthermore, $S^{-1}$ satisfies the inequality
\begin{equation}\label{eq:isc}
  |\langle Q_t S^{-1} g ,\phi \rangle|\leq \frac{\|\phi\|}{c} \| Q_t g\| \quad(t\in \mathbb{R}, \phi\in \s X, g\in \ran(S)).
\end{equation}
\end{Proposition}
\begin{Proof} First of all, we verify that $S$ is one-to-one. For this, we have to verify that $\kar(S)=\{0\}$. Let $f\in \kar(S)$. Then for all $t\in \mathbb{R}$ we have
\[
   0 =    |\langle Q_t Sf,f\rangle| \geq  c\langle Q_t f,f\rangle.
\]
Letting $t\to\infty$, we infer $\|f\|^2=0$. We are left with showing \eqref{eq:isc}, which is also sufficient for causality of $S^{-1}$. For this, let $t\in \mathbb{R}$, $\phi\in \s X$, $g\in \ran(S)$. We compute, using the hypothesis,
\begin{align*}	
   |\langle Q_t g, Q_t S^{-1}g\rangle| & =|\langle Q_t^2 g, S^{-1}g\rangle|
   \\  & =  |\langle Q_t g, S^{-1}g\rangle|
   \\  &\geq c \langle Q_t S^{-1} g, S^{-1} g\rangle 
   \\  &= c \langle Q_t S^{-1} g,Q_t S^{-1} g \rangle=c\|Q_t S^{-1} g \|^2.
\end{align*}
Hence, 
\[
    |\langle Q_t S^{-1} g ,\phi \rangle|\leq  \|Q_t S^{-1} g \|\|\phi\|\leq \frac{\|\phi\|}{c} \| Q_t g\|.
\]
\end{Proof}

In view of the applications to follow, the inequality that has been shown for $S^{-1}$ is often satisfied. That is why, we introduce a concept slightly stronger than causality.

\begin{Definition}\label{d:sca} Let $(\s X, (Q_t)_t)$, $(\s Y, (R_t)_t)$ be resolution spaces, $S\colon \dom(S)\subseteq \s X\to \s Y$ linear. Then we call
$S$ \emph{strongly causal}, if for all $t\in \mathbb{R}$ there exists $C\geq 0$ such that for all $f\in \dom(S)$
\[
  \|R_t S f\|\leq C\| Q_t f\|.
\] 
\end{Definition}

The next theorem is the reason, why the notion just introduced is so important for applications. Namely, for continuous mappings, causality and strong causality are the same.

\begin{Theorem}\label{t:csc} Let $(\s X, (Q_t)_t)$, $(\s Y, (R_t)_t)$ resolution spaces, $S\colon \dom(S)\subseteq \s X\to \s Y$ densely defined, linear, continuous. Then the following conditions are equivalent:
\begin{enumerate}[label=(\roman*)]
 \item\label{csc1} for all $f\in \s X$, $t\in \mathbb{R}$ we have $Q_t f =0$ implies $R_t \overline{S} f=0$;
 \item\label{csc2} $S$ is causal;
 \item\label{csc2.5} $S$ is causal on $D$ for some $D\subseteq \s Y$ dense; 
 \item\label{csc3} $S$ is strongly causal;
 \item\label{csc4} for all $t\in \mathbb{R}$, we have $R_t \overline{S} Q_t = R_t \overline{S}$;
 \item\label{csc5} for all $t\in \mathbb{R}$, we have $(1-R_t)\overline{S}(1-Q_t)=\overline{S}(1-Q_t)$.
\end{enumerate} 
\end{Theorem}
\begin{Proof}
 Using that $\dom(\overline{S})=\s X$ as $S$ is continuous and densely defined, the equivalence of \ref{csc1}, \ref{csc2}, and \ref{csc2.5} has been established in Theorem \ref{thm:causal_cont_causal}. The equivalence of \ref{csc4} and \ref{csc5} is an easy computation. Moreover, linearity of $\overline{S}$ implies that \ref{csc1} is necessary for \ref{csc4}. But, if \ref{csc1} holds, then from $f=Q_tf+(1-Q_t)f$ and using that $Q_t(1-Q_t)=0$, we infer 
 \[
    R_t \overline{S}f =R_t \overline{S}(Q_t f + (1-Q_t)f)= R_t\overline{S}Q_t f + R_t\overline{S}(1-Q_t) f =  R_t\overline{S}Q_t f \quad(f\in \s X),
 \]
 which is \ref{csc4}.
 
 Next, clearly, \ref{csc3} is sufficient for \ref{csc2} and, assuming \ref{csc4}, we get for all $f\in \s X$
 \[
    \|R_t \overline{S} f\|= \|R_t \overline{S} Q_t f\|\leq \|R_t \overline{S}\| \| Q_t f\| \leq \| \overline{S}\| \| Q_t f\|,
 \]
 which implies \ref{csc3}.
\end{Proof}

\renewcommand{\baselinestretch}{0.65}\normalsize\mysection{Causality for Evolutionary Mappings}{Causality for Evolutionary Mappings}{causal evolutionary mappings are standard evolutionary $\cdot$ closures of standard evolutionary mappings are independent of $\nu$ $\cdot$ standard evolutionary mappings form a vector space $\cdot$ closable evolutionary $\cdot$ criterion for inverses being standard evolutionary $\cdot$ bounded sets of evolutionary mappings $\cdot$ Theorem \ref{t:inv_ce}}\label{s:cem}

\renewcommand{\baselinestretch}{1}\normalsize
In this section, we combine the results from Sections \ref{s:em} and \ref{s:ccm}. An important result of this section is Theorem \ref{t:evsc_sev} together with Proposition \ref{prop:indep_of_nu_for_bounded_mappings}, that is, roughly speaking,
\begin{itemize}
 \item causal evolutionary mappings are essentially the same as standard evolutionary mappings, (Theorem \ref{t:evsc_sev})
 \item the closure of a causal, evolutionary mapping is widely \emph{independent} of the exponential weight (Proposition \ref{prop:indep_of_nu_for_bounded_mappings}).
\end{itemize}
  To begin with we define causality for evolutionary mappings.

\begin{Definition}\label{d:sc} Let $\s X$, $\s Y$ Hilbert spaces, $Q_t\coloneqq T_{\1_{(-\infty,t)}}$, $\nu\in \mathbb{R}$. Then both the spaces $(L_\nu^2(\mathbb{R};\s X),(Q_t)_t)$ and $(L_\nu^2(\mathbb{R};\s Y),(Q_t)_t)$ are resolution spaces. Let $\s S\in L_{\textnormal{ev},\nu}(\s X,\s Y)$. We call $\s S$ \emph{causal}, if $\s S$ is causal considered as a mapping from $L_\nu^2(\mathbb{R};\s X)$ to $L_\nu^2(\mathbb{R};\s Y)$.
\end{Definition}

\begin{Proposition}\label{p:1cau} Let $\s X$, $\s Y$ be Hilbert spaces, $\nu\in \mathbb{R}$, $\s S\in L_{\textnormal{ev},\nu}(\s X,\s Y)$ causal. Then, for all $\mu\geq \nu$, $\s S\in L_{\textnormal{ev},\mu}(\s X,\s Y)$ is causal.   
\end{Proposition}
\begin{Proof}
  As $\s S$ is evolutionary, we are in the position to apply Theorem \ref{t:csc} to $\s S$ considered as a mapping from $L_\nu^2(\s X)$ to $L_\nu^2(\s Y)$. Thus, for all $t\in \mathbb{R}$, $\phi\in L_\textnormal{c}^2(\mathbb{R};\s Y)$ there exists $C\geq 0$ such that 
  \begin{equation}\label{eq:1cau}
     |\langle Q_t\s S f,\phi \rangle_{L_\nu^2}|\leq C \|Q_t f\|_{L_\nu^2} \quad (f\in \dom(\s S))
  \end{equation}
  But, as $\s S$ is evolutionary, $\dom(\s S)\subseteq L_\nu^2\cap L_\mu^2$. Hence, for $f\in \dom(\s S)$,
  \[
     \|Q_t f\|_{L_\nu^2}^2  =\int_{-\infty}^t \|f(\tau)\|^2e^{-2\nu \tau}d\tau 
     =\int_{-\infty}^t \|f(\tau)\|^2e^{-2\mu \tau}e^{2(\mu-\nu) \tau}d\tau 
     \leq \|Q_t f\|_{L_\mu^2}^2 e^{2(\mu-\nu)t}.
  \]
  Thus, using that multiplication by the exponential function is a bijection on $L_\textnormal{c}^2(\mathbb{R};\s Y)$, we infer from \eqref{eq:1cau} for all $t\in \mathbb{R}$, $\phi\in L_\textnormal{c}^2(\mathbb{R};\s Y)$, there exists $C\geq 0$ such that 
  \[
     |\langle Q_t \s S f,\phi\rangle_{L_\mu^2}|\leq Ce^{(\mu-\nu)t} \|Q_t f\|_{L_\mu^2} \quad (f\in \dom(\s S)),
  \]
  which implies that $\s S$ as a mapping in $L_\mu^2$ is causal on $L_\textnormal{c}^2(\mathbb{R};\s Y)$. Hence, $\s S$ is causal in $L_\mu^2$, by Theorem \ref{thm:causal_cont_causal}.
\end{Proof}

Recall that the condition of being standard evolutionary, that is, being evolutionary with the standard causal domain $\mathcal{D}_\nu(\s X)=\bigcap_{\mu\geq\nu}L_\mu^2(\s X)$ as underlying domain of definition, results in causality:

\begin{Proposition}\label{p:sev_sc} Let $\s X$, $\s Y$ Hilbert spaces, $\s S\in L_{\textnormal{sev}}(\s X,\s Y)$. Then $\s S\in L_{\textnormal{ev},\nu}(\s X,\s Y)$ is strongly causal for all $\nu\in \mathbb{R}$ large enough.
 \end{Proposition}
 \begin{Proof}
  By Remark \ref{evo=caus}, $\s S^\nu$ (the closure of $\s S$ in $L_\nu^2$) is causal for all $\nu$ large enough. Hence, the assertion follows from Proposition \ref{p:1cau} and Theorem \ref{t:csc}.
 \end{Proof}

 Next, we will seek to prove the following converse of the latter proposition.
 
 \begin{Theorem}\label{t:evsc_sev} Let $\s X$, $\s Y$ Hilbert spaces, $\nu\in \mathbb{R}$, $\s S\in L_{\textnormal{ev},\nu}(\s X,\s Y)$ causal. Then 
 \[
    \s T\coloneqq \{ (f,g);f\in \mathcal{D}_\nu(\s X) \land \exists \mu\geq  \nu: g=\s S^\mu f \}\in L_{\textnormal{sev},\nu}(\s X,\s Y).
 \]
 Moreover, for all $\mu\geq \nu$: $\s T^\mu=\s S^{\mu}$.
 \end{Theorem}
 \begin{Definition}\label{d:cr} Let $\s X$, $\s Y$ Hilbert spaces, $\nu\in \mathbb{R}$, $\s S\in L_{\textnormal{ev},\nu}(\s X,\s Y)$ causal. We call $\s T$ as defined in Theorem \ref{t:evsc_sev} the \emph{standard realization} of $\s S$.
\end{Definition}
We remark here that Theorem \ref{t:evsc_sev} implicitly asserts that $\s T$ defines a right-unique relation, that is, a mapping. Theorem \ref{t:evsc_sev} also serves as a justification for treating standard evolutionary mappings later on, only. Furthermore, in view of Theorem \ref{t:evsc_sev}, we shall even employ the custom to consider causal evolutionary mappings and standard evolutionary mappings as synonymous. The first step for the proof of Theorem \ref{t:evsc_sev} is to show that $\s T$ is right-unique.

\begin{Lemma}\label{le:pre_indep} Let $\mu, \nu\in \mathbb{R}$, $\mu\geq\nu$, $\s X, \s Y$ Hilbert spaces. Let both $\s S_\nu \in L(L_{\nu}^2(\s X),L_{\nu}^2(\s Y))$, $\s S_\mu \in L(L_{\mu}^2(\s X),L_{\mu}^2(\s Y))$ be causal. Assume there is $\s D\subseteq L_{\nu}^2(\s X)\cap L_{\mu}^2(\s X)$ dense in $L_{\mu}^2(\R;\s X)$ with the property that $\s S_\nu|_{\s D}=\s S_\mu|_{\s D}$. Then $\s S_\mu$ and $\s S_\nu$ coincide on the set $L_{\nu}^2(\s X)\cap L_{\mu}^2(\s X)$.
\end{Lemma}
\begin{Proof}
   Let $f\in L_{\nu}^2(\s X)\cap L_{\mu}^2(\s X)$, $t \in \R$. Let $Q_t$ be the operator of multiplication by $\1_{(-\infty,t]}$. By hypothesis, there exists $(f_n)_n$ in $\s D$ such that $f_n\to f$ in $L_\mu^2(\s X)$ as $n\to\infty$. Hence, $Q_t f_n \to Q_t f$ in $L_\nu^2(\s X)\cap L_\mu^2(\s X)$ as $n\to\infty$. Hence, for $n\in\N$,
   \[
       Q_t S_\mu Q_t f_n = Q_t S_\mu f_n = Q_t S_\nu f_n = Q_t S_\nu Q_t f_n.
   \]
   As both the left-hand and the right-hand side of the latter equality converge in $L^2_{\textnormal{loc}}(\s X)$, their respective limits coincide. So,
   \[
      Q_t S_\mu f = Q_t S_\mu Q_t f = Q_t S_\nu Q_t f = Q_t S_\nu f.
   \]
 Since the rationale presented applies to all $t\in \mathbb{R}$, the claim is proved.
\end{Proof}

\begin{Proposition}\label{prop:indep_of_nu_for_bounded_mappings} Let $\nu\in \mathbb{R}$, $\s X, \s Y$ Hilbert spaces. Let $\s S\in L_{\textnormal{ev},\nu}(\s X,\s Y)$ causal. Then $\s S^\nu$ and $\s S^\mu$, the closures of $\s S$ in $L_{\nu}^2$ and $L_{\mu}^2$, respectively, coincide on $L_{\nu}^2(\s X)\cap L_{\mu}^2(\s X)$, $\mu\geq \nu$. 
\end{Proposition}
\begin{Proof}
  As $\s S$ is evolutionary, $\dom(\s S)\subseteq L_{\nu}^2\cap L_\mu^2$ is dense in both $L_\mu^2$ and $L_\nu^2$. The respective continuous extensions $\s S^\nu$ and $\s S^\mu$ are causal, by Proposition \ref{p:1cau} and Theorem \ref{t:csc}. $\s S^\nu$ and $\s S^\mu$ coincide on $\dom(\s S)$ and, thus, on the intersection of the respective domains, by Lemma \ref{le:pre_indep}. 
\end{Proof}

\begin{Proof}[of Theorem \ref{t:evsc_sev}]
 Proposition \ref{prop:indep_of_nu_for_bounded_mappings} implies that $\s T$ is well-defined. Note that this also settles evolutionarity of $\s T$. In particular, we get that $\s T$ is standard evolutionary at $\nu$. The last assertion, $\s T^\mu=\s S^\mu$, $\mu\geq\nu$, can be seen as follows. By definition, for all $\mu\geq \nu$ we have
$ \s T^\mu|_{\s D_\nu}=\s S^\mu|_{\s D_\nu}$. So, $\s T^\mu=\s S^\mu$ as $\s D_\nu(\s X)$ is dense in $L_\mu^2(\s X)$.
\end{Proof}

Basically, Proposition \ref{prop:indep_of_nu_for_bounded_mappings} asserts that the closures of evolutionary (and causal) mappings do not depend on the particular realization in some $L_\nu^2$, that is, on the exponential weight parametrized by $\nu$. Theorem \ref{t:evsc_sev} contains the prototype of domains causal evolutionary mappings may be endowed with. 

Another consequence of Theorem \ref{t:evsc_sev} is the inclusion $L_{\textnormal{sev},\nu}\subseteq L_{\textnormal{sev},\mu}$, $\mu\geq \nu$, in the following sense:

\begin{Corollary}\label{c:sev_dir} Let $\s X$, $\s Y$ Hilbert spaces, $\nu,\mu\in \mathbb{R}$, $\mu\geq\nu$, $\s S\in L_{\textnormal{sev},\nu}(\s X,\s Y)$. Then $\s S\subseteq \s T$, where $\s T\in L_{\textnormal{sev},\mu}(\s X,\s Y)$ is the standard realization of $\s S$ considered as evolutionary at $\mu$. 
\end{Corollary}

In the sense of Corollary \ref{c:sev_dir}, any standard evolutionary mapping at $\nu$ may be considered as standard evolutionary at $\mu\geq \nu$, we shall do so in the following. This custom enables us to ease the formulations of several statements. For instance, a reformulation of Proposition \ref{p:sev_alg} reads as follows.

\begin{Proposition}\label{p:sev_alg2} Let $\s X$, $\s Y$, $\s Z$ Hilbert spaces, $\s S,\s T\in L_{\textnormal{sev}}(\s X,\s Y)$, $\s U\in L_{\textnormal{sev}}(\s Y,\s Z)$, $\alpha\in \mathbb{C}$. Then $\s S+\alpha \s T\in L_{\textnormal{sev}}(\s X,\s Y)$ and $\s U\s S\in L_{\textnormal{sev}}(\s X,\s Z)$. 
\end{Proposition}

In order to develop a solution theory for certain differential equations, apart from the Hadamard requirements of unique existence of solutions that depend continuously on the data, we ask for causality of the solution operator. Moreover, the solution operator should be widely independent of the exponential weight, which results in the requirement of evolutionarity for the solution operator. For (abstract) ordinary differential equations with potentially infinite-dimensional state space $\s X$, the solution operator can be computed explicitly as a composition (of inverses) of causal evolutionary mappings, see Chapter \ref{ch:ST}. 

The next theorem should be viewed in the context of Proposition \ref{p:sev_alg2}. Indeed, with additional regards to Theorem \ref{t:evsc_sev}, causal evolutionary mappings are closed under composition and addition. The next theorem complements these statements, by giving a criterion for an inverse being causal and evolutionary. Before giving the precise statement, we need to introduce a slightly more general concept than that of evolutionary mappings.

\begin{Definition}\label{d:cev} Let $\s X$, $\s Y$ Hilbert spaces, $\nu\in \mathbb{R}$. We call \[\s B\colon \dom(\s B)\subseteq \bigcap_{\mu\geq \nu} L_\mu^2(\mathbb{R};\s X)\to \bigcap_{\mu\geq \nu} L_\mu^2(\mathbb{R};\s Y)\] \emph{closable evolutionary (at $\nu$)}, if $\s B$ is linear, $\dom(\s B)$ dense in $L_\mu^2(\mathbb{R};\s X)$ and that for all $\mu\geq \nu$
\[
   \dom(\s B) \subseteq L_\mu^2(\mathbb{R};\s X)\to L_\mu^2(\mathbb{R};\s Y),\ f\mapsto \s Bf
\]
is closable. Again, we denote by $\s B^\mu$ the respective closure of $\s B$, $\mu\geq \nu$. We set
\[
    C_{\textnormal{ev},\nu}(\s X,\s Y)\coloneqq \{ \s B;\s B \text{ closable evolutionary at }\nu\};\ C_{\textnormal{ev},\nu}(\s X)\coloneqq C_{\textnormal{ev},\nu}(\s X,\s X).
\]
\end{Definition}

\begin{Example}\label{ex:cev} Let $\s X$, $\s Y$ Hilbert spaces, $\nu\in\mathbb{R}$.

(a) Recall $\s D_\nu(\s X)=\bigcap_{\mu\geq \nu} L_\mu^2(\s X)$. Then 
\[
   \check{\partial}_{t,\nu}\coloneqq \bigcap_{\mu\geq\nu} \partial_{t,\nu} =\{ (f,g); f,g\in \s D_\nu(\s X), g=f'\} 
\]
is closable evolutionary at $\nu>0$. Note that $\partial_t\coloneqq\check{\partial}_{t}\coloneqq \bigcap_{\mu>0} \partial_{t,\mu}$ is also closable evolutionary at $\nu>0$.

(b) Let $A\colon \dom(A)\subseteq \s X\to\s Y$ be a densely defined,  closable, linear operator. Then the (abstract multiplication) operator $\s A$ given by $\s A f\coloneqq (t\mapsto Af(t))$ for functions $f$ with bounded support and assuming values in $\dom(A)$ is closable as an operator in $L_\mu^2$, $\mu\in \R$. Denoting by $\s A^\mu$ the respective closure, we see that
\[
   \check{\s A} \coloneqq \bigcap_{\mu\in \mathbb{R}} \s A^\mu
\]
is closable evolutionary at $\nu$. We will also just write $\s A$ for $\check{\s A}$.

(c) If $\s B\in C_{\textnormal{ev},\nu}(\s X,\s Y)$, then
\[
   \check{\s B}_\nu \coloneqq \bigcap_{\mu\geq \nu}\s B^\mu
\]
is closable evolutionary at $\nu$ as well. 
\end{Example}

\begin{Theorem}\label{t:inv_ce} Let $\s X$ Hilbert space, $\nu\in \mathbb{R}$, $\s B\in C_{\textnormal{ev},\nu}(\s X)$. Assume there exists $c>0$ such that for all $\mu\geq \nu$, $t\in \mathbb{R}$,
 \[
      \Re\langle Q_t\s Bf,f\rangle_{L_\mu^2(\s X)}\geq c\langle Q_tf,f\rangle_{L_\mu^2(\s X)} \quad (f\in \dom(\s B)),
 \]
where $Q_t$ is multiplication by $\1_{(-\infty,t)}$. Assume, in addition, that $\ran(\s B)$ is dense in $L_\mu^2(\mathbb{R})$ for all $\mu\geq\nu$.

Then $\s S\coloneqq \s B^{-1}$ is evolutionary at $\nu$ and causal. Moreover, for all $\mu_1,\mu_2\geq \nu$, we have that $\s S^{\mu_1}=\s S^{\mu_2}$ on $L_{\mu_1}^2(\s X)\cap L_{\mu_2}^2(\s X)$.
\end{Theorem}
\begin{Proof} Existence and causality of $\s B^{-1}$ as a mapping in $L_\mu^2(\mathbb{R};\s X)$ follow from Proposition \ref{p:inverse_strongly_causal}. In fact, for all $\mu\geq \nu$ the inequality
\[
  |\langle Q_t \s B^{-1}g,\phi\rangle|\leq \frac{\|\phi\|}{c}\|Q_t g\|\quad  (t\in \mathbb{R}, \phi\in L_\mu^2(\mathbb{R};\s X), g\in \dom(\s B^{-1})) 
\]
holds true. Hence, by letting $t\to \infty$ and computing the supremum over $\phi$ with norm $1$, we arrive at
\[
   \limsup_{\mu\to\infty}\|\s B^{-1}\|_{L_\mu^2(\s X)\to L_\mu^2(\s X)}\leq 1/c.
\]
This, together with the density of $\dom(\s S)=\ran(\s B)$ in $L_\mu^2(\mathbb{R};\s X)$ establishes evolutionarity and causality of $\s S$. 
Let $\mu_1\geq\mu_2\geq\nu$. Then, by Proposition \ref{p:1cau}, $\s S\in L_{\textnormal{ev},\mu_2}(\s X)$ is causal. Hence, $S^{\mu_1}$ coincides with $S^{\mu_2}$ on $L_{\mu_1}^2(\s X)\cap L_{\mu_2}^2(\s X)$, by Proposition \ref{prop:indep_of_nu_for_bounded_mappings}.
  \end{Proof} 

The condition on the density of the range of $\s B$ can be dropped, if $\s B$ is assumed to be continuous:

\begin{Corollary}\label{c:inv_ce} Let $\s X$ Hilbert space, $\nu\in \mathbb{R}$, $\s B\in L_{\textnormal{ev},\nu}(\s X)$. Assume there exists $c>0$ such that for all $\mu\geq \nu$, $t\in \mathbb{R}$,
 \[
      \Re\langle Q_t\s Bf,f\rangle_{L_\mu^2(\s X)}\geq c\langle Q_tf,f\rangle_{L_\mu^2(\s X)} \quad (f\in \dom(\s B)),
 \]
where $Q_t$ is multiplication by $\1_{(-\infty,t)}$.

Then $\s S\coloneqq \s B^{-1}$ is evolutionary at $\nu$ and causal. Moreover, for all $\mu_1,\mu_2\geq \nu$, we have that $\s S^{\mu_1}=\s S^{\mu_2}$ on $L_{\mu_1}^2(\mathbb{R};\s X)\cap L_{\mu_2}^2(\mathbb{R};\s X)$.
\end{Corollary}

For the proof of Corollary \ref{c:inv_ce}, it is sufficient to observe that the inequality assumed implies (as $t\to\infty$):
\[
      \Re\langle \s Bf,f\rangle_{L_\mu^2(\s X)}\geq c\langle f,f\rangle_{L_\mu^2(\s X)} \quad (f\in \dom(\s B)),
 \]
Thus, the needed density result for the application of Theorem \ref{t:inv_ce} follows from the following observation:
\begin{Proposition}\label{p:posd} Let $\s X$ be a Hilbert space, $B$ a densely defined, linear operator in $\s X$ satisfying
\begin{equation}\label{eq:p_posd}
    \Re\langle B \phi,\phi\rangle\geq c\langle \phi,\phi\rangle
\end{equation}
for all $\phi\in \dom( B)$ and some $c>0$.
Then $B$ is closable. 

If, in addition, 
\begin{equation}\label{eq:p_posd2}
|\langle  B^* \psi,\psi\rangle|\geq c\langle \psi,\psi\rangle 
\end{equation}
for all $\psi\in \dom( B^*)$, then $0\in \rho(\overline{B})$, $\|\overline{ B}^{-1}\|\leq 1/c$, and $\ran(B)$ is dense in $\s X$.

If $B$ is continuous, then \eqref{eq:p_posd} implies \eqref{eq:p_posd2}.
\end{Proposition}
\begin{Proof}
 We address closability first. The Banach space version of the closability result can be found in \cite[Theorem 4.2.5]{Beyer2007}. Assume $B$ not to be closable. Then we find $(\phi_n)_n$ in $\dom(B)$ converging to $0\in \s X$ with the property that $(B\phi_n)_n$ converges to some non-zero $\psi\in X$. Without restriction, $\|\psi\|=1$. By the density of $\dom(B)$, there exists $\zeta\in \dom(B)$ with $\|\psi-\zeta\|<1/2$. Hence, $\|\zeta\|>1/2$. Next, for all $\beta>0$ and $n\in\mathbb{N}$, we obtain due to \eqref{eq:p_posd}
 \[
         \big\|\zeta-\frac{1}{\beta}\phi_n\big\|\leq \big\|\zeta-\frac{1}{\beta}\phi_n+ \beta B\big(\zeta-\frac{1}{\beta}\phi_n\big)\big\|
         = \big\|\zeta-\frac{1}{\beta}\phi_n+ \beta B\zeta-B\phi_n\big\|.
 \]
 Letting $n\to\infty$ and afterwards $\beta\to0$, we obtain
 \[
    \frac{1}{2}<\|\zeta\|\leq \|\zeta-\psi\|<\frac{1}{2},
 \]
 a contradiction, yielding closability.
 
 Next, the inequality \eqref{eq:p_posd} for $B$ implies the same for $\overline{B}$ (for all $\phi\in \dom(\overline{B})$). Hence, $\overline{B}$ is one-to-one and the, thus, existing inverse has an operator-norm bounded by $1/c$. Moreover, the same inequality implies the closedness of the range of $\overline{B}$. Hence, for showing $0\in \rho(\overline{B})$ we are left with showing that $\overline{B}$ is, in fact, onto. For this, we recall the orthogonal decomposition $\s X=\overline{\ran(B)}\oplus \kar(B^*)$. So, the inequality assumed for $B^*$ implies that $B^*$ is one-to-one, and, hence, $\kar(B^*)$ is trivial, implying
 \[
   \s X\supseteq \ran(\overline{B})=\overline{\ran(\overline{B})}\supseteq \overline{\ran(B)}=\s X.
 \] The latter settles both the density of $\ran(B)$ in $\s X$ as well as $0\in \rho(\overline{B})$. 
 \end{Proof}
The author is indebted to Sascha Trostorff for spotting a flaw and stating a simpler argument compared to an earlier version of the latter proposition. Sebastian Mildner eventually found the source \cite{Beyer2007}, which settled an issue concerning the closability statement.

In order to apply the following concept right away in the beginning of Chapter \ref{ch:ST}, we conclude with some convergence aspects of evolutionary mappings. Later on, we will deal with these issues in more detail. We introduce the notion of boundedness and convergence of standard evolutionary mappings.

\begin{Definition}\label{d:sev_bdd}
 Let $\s X$, $\s Y$ be Hilbert spaces. A subset $\mathfrak{S}\subseteq L_{\textnormal{sev},\nu}(\s X,\s Y)$
is called \emph{bounded} if  
\[
   \sup_{\mu\geq\nu}\sup_{\s S\in\mathfrak{S}}\|\s S\|_{L_\mu^2}<\infty. 
\]
We say $\mathfrak{S}\subseteq L_{\textnormal{sev}}(\s X,\s Y)$ is \emph{bounded}, if there exists $\nu\in \mathbb{R}$ such that $\mathfrak{S}\subseteq L_{\textnormal{sev},\nu}(\s X,\s Y)$ is bounded. A family $(\s S_{\iota})_{\iota\in I}$ in $L_{\textnormal{sev}}(\s X,\s Y)$
is called \emph{bounded} if $\{\s S_{\iota};\iota\in I\}$ is bounded. 

Let $(\s S_{n})_{n\in \N}$ in $L_{\textnormal{sev}}(\s X,\s Y)$, $\s T\in L_{\textnormal{sev}}(\s X,\s Y)$.
We say that $(\s S_{n})_{n\in \N}$ is \emph{convergent} to $\s T$ (or $(\s S_n)_n$ \emph{converges to $\s T$}), $\s S_n{\to} \s T\ (n\to\infty)$, if there exists $\nu\in \mathbb{R}$ such that
\[
   \s S_n^\mu  \to \s T^\mu \text{ in }L(L_\mu^2(\s X),L_\mu^2(\s Y))\text{ as }n\to\infty,
\]
for all $\mu\geq \nu$.
\end{Definition}

\section{Comments}

The notion of evolutionary mappings is inspired by the term `evolutionary' introduced in \cite{Picard}: It roots in determining the `time-like' directions in a given partial differential expression with constant coefficients. For sake of presentation, we think of a polynomial $p$ in $d$ variables with formally inserted the partial derivatives $\partial_1, \ldots,\partial_d$ leading to
\[
   p(\partial_1,\ldots,\partial_d) = \sum_{\alpha\in \N_0^{d}} c_\alpha (\partial_1,\ldots,\partial_d)^\alpha,
\]
where we employed multiindex notation and assume that all but finitely many $c_\alpha\in \mathbb{C}$ are $0$. Next, for given $f$ consider the problem of finding $u$ such that
\begin{equation}\label{eq:pdecon}
   p(\partial_1,\ldots,\partial_d) u =f.
\end{equation}
We might try setting up a solution theory for \eqref{eq:pdecon}. Similar to the ordinary differential equations case in the previous comments section, we seek a solution $u$ in an exponentially weighted space for every variable, that is, in a tensor product space $\bigotimes_{j=1}^d L_{w_j}^2(\mathbb{R})$ for $w=(w_1,\ldots,w_d)\in\mathbb{R}^d$. So, applying the Fourier--Laplace transformation in each variable and using Remark \ref{r:spe_thm_dtnu}, we get that \eqref{eq:pdecon} reads
\begin{equation}\label{eq:pdecon2}
  p(i\xi_1+w_1,\ldots,i\xi_d+w_d) \hat{u}(\xi_1,\ldots,\xi_d) =\hat{f}(\xi_1,\ldots,\xi_d)\quad(\xi_1,\ldots,\xi_d\in \mathbb{R})
\end{equation}
 for appropriate $\hat{u}$ and $\hat{f}$. Therefore, solving for $u$ in \eqref{eq:pdecon} leads to inverting $p(i\xi_1+w_1,\ldots,i\xi_d+w_d)$ for all $\xi_1,\ldots,\xi_d\in \mathbb{R}$. To make this procedure well-defined, we ask $p(i\xi_1+w_1,\ldots,i\xi_d+w_d)\neq 0$ for one $w=(w_1,\ldots,w_d)\in\mathbb{R}^d$ and all $(\xi_1,\ldots,\xi_d)\in \mathbb{R}^d$. Following \cite[Definition 3.1.14]{Picard}, we let $w\in \mathbb{R}^d$, $|w|=1$, and call $p(\partial_1,\ldots,\partial_d)$ \emph{evolutionary in direction $w$}, if there exists $\nu\in \mathbb{R}$ such that for all $\mu\geq\nu$ the polynomial $\mathbb{R}^d\ni\xi \mapsto p(i\xi+\mu w)\in \mathbb{C}$ has no zeros. When discussing so-called `canonical forms' of differential expressions, the direction of evolutionarity is singled out as the direction of time, see \cite[Section 3.1.7]{Picard}. We shall also refer to \cite[Section 3.1.6]{Picard}, where evolutionarity is discussed in view of the classical classification of partial differential equations into elliptic, parabolic and hyperbolic. 
 
 In the framework presented in this exposition, the direction of time is already given and modeled by the direction of the real line in the first variable of the space $L_\nu^2(\mathbb{R};\s X)$. Thus, similar to \cite{Picard}, the remaining variables are thought of being contained in $\s X$, the Hilbert space describing `spatial coordinates'. Followed by the introduction of evolutionarity in \cite{Picard}, causality of evolutionary partial differential expressions has been discussed as well. The definition of causality is similar to the one in Definition \ref{d:cau} but formulated for mappings from the space of distributions to the space of distributions (see \cite[Definition 3.1.47]{Picard}), see also Remark \ref{evo=caus}. The theorem of Paley--Wiener (Theorem \ref{t:PW}) has proven to be useful in showing causality already in \cite[pp 134]{Picard}, or \cite{PicPhy}. It rests on the usage of the Fourier--Laplace transformation. For non-autonomous problems this strategy, however, may not be applicable any more. Hence, we developed a framework which enables us to discuss and prove causality without employing the Fourier--Laplace transformation for the general setting of evolutionary mappings discussed here.
 
 The question of whether a closable mappings admits a causal closure has been addressed in the time translation-invariant case in the community of control and systems theory, see \cite{Jac2000}. The method of choice for answering this question is the (Fourier--)Laplace transformation or the $z$-transformation for discrete-time settings. In \cite{Waurick2015IM_Caus}, we gave a possible characterization of operators admitting a causal closure without asking for time-shift invariance. We also developed a Banach space analogue for this characterization. 

 The independence of the exponential weight for certain solution operators of certain partial differential equations has been addressed in \cite[Theorem 6.1.4]{Picard}, \cite[Lemma 3.6]{Trostorff2013a} for the time-shift invariant case. Hence, naturally, the arguments in \cite{Picard,Trostorff2013a} employ the Fourier--Laplace transformation. In \cite[Theorem 4.6]{Kalauch} the same question is discussed for possibly non-linear ordinary delay differential equations. In \cite[Section 4]{Waurick2014MMAS_Non} the independence of exponential weight has been shown for a specific class of linear non-autonomous partial differential equations. The line of ideas in \cite{Waurick2014MMAS_Non} together with \cite{Waurick2015IM_Caus} have then eventually lead to the treatment developed here. 
 
 There is also a huge theory of causal differential equations with a focus on ordinary differential equations in a Banach space setting developed in \cite{Lakshmikantham2009} and the references given there. We shall also refer to the references given in \cite{Waurick2015IM_Caus} for the treatment of causal mappings in other settings.
\cleardoublestandardpage
\chapter{Solution Theory for Evolutionary Equations}\label{ch:ST}

The aim of this chapter is to prove two well-posedness statements for evolutionary equations. In very abstract terms, we will consider operator equations of the type
\begin{equation}\label{eq:ade}
   \s B u = f
\end{equation}
for $\s B$ being defined in some $L_\nu^2(\R;\s X)$, $\s X$ Hilbert space. So, we address conditions on the continuous invertibility of $\s B$.

Assuming conditions on the structure of $\s B$, we will consider both ordinary and partial differential equations. More precisely, by assuming that $\s B$ is a sum of products of certain operators, we will provide conditions on the constituents of this composure yielding a solution theory for \eqref{eq:ade}. Here, in a solution theory, we gather the three Hadamard requirements, that is, existence, uniqueness of solutions as well as continuous dependence on the data. Furthermore, we want the solution operator $\s S=\s B^{-1}$ once existent to be causal. Moreover, being realized in certain weighted $L_\nu^2$-spaces we want $\s S$ to be fairly independent of $\nu$, that is, if $\s S$ existed in $L_\nu^2$ and $L_\mu^2$ then $\s S$ should be a well-defined mapping on $L_\nu^2\cup L_\mu^2$. The latter fact is properly restated as that the realizations $\s S^\nu$ and $\s S^\mu$ of $\s S$ on $L_\nu^2$ and $L_\mu^2$, respectively, should coincide on $L_\nu^2\cap L_\mu^2$. Hence, asking for a solution theory of \eqref{eq:ade}, amounts to the question of when $\s B^{-1}=\s S$ is evolutionary at some $\nu$ and causal or, equivalently (cf.~Proposition \ref{p:sev_sc} and Theorem \ref{t:evsc_sev}), is \emph{standard evolutionary}. 

In consequence, both the treatments of ordinary and partial differential (evolutionary) equations discussed are similar to one another: In a preparatory step, we will show continuous invertibility of a certain operator in some $L_\nu^2(\s X)$-space. This settles the three Hadamard requirements of existence and uniqueness as well as the continuous dependence on the data. The concluding step will be to apply the results derived in the preparatory step to (standard) evolutionary mappings and to show that the corresponding solution operator is standard evolutionary itself.

\renewcommand{\baselinestretch}{0.65}\normalsize\mysection{Ordinary Differential Equations}{Ordinary Differential Equations}{solution theory for ordinary differential equations $\cdot$ Theorem \ref{t:wp_ode_evo}}\label{s:wpode}

\renewcommand{\baselinestretch}{1}\normalsize
In order to further illustrate the notion of evolutionary mappings and some of the main ideas of the solution theory to be developed for non-auto\-no\-mous evolutionary equations, we stick to a specific class of evolutionary equations first. The class to be discussed in this section is the one of abstract ordinary differential equations with possibly infinite-dimensional state space.

\begin{Theorem}\label{t:wp_ode} Let $\s X$, $\s Y$ Hilbert spaces, $\nu, c_0,c_1>0$. Let $\mathcal{M}\in L(L_\nu^2(\mathbb{R};\s X))$ and let $\mathcal{N}=(\mathcal{N}_{ij})_{i,j\in \{0,1\}}\in L(L_\nu^2(\s X\times \s Y))$. Assume
\begin{equation}\label{eq:wp_ode} 
   \Re \langle \mathcal{M} \phi,\phi\rangle_{L_\nu^2(\s X)}\geq c_0 \langle \phi,\phi\rangle_{L_\nu^2(\s X)},\quad \Re \langle \mathcal{N}_{11}\psi,\psi\rangle_{L_\nu^2(\s Y)}\geq c_1\langle\psi,\psi\rangle_{L_\nu^2(\s Y)}
\end{equation}
for all $(\phi,\psi)\in L_{\nu}^2(\s X\times \s Y)$.
Further assume the validity of the estimate
\begin{equation}\label{eq:wp_ode2}
  c_1 \|\mathcal{N}_{00}\|+\|\mathcal{N}_{01}\|\|\mathcal{N}_{01}\|<\nu c_0 c_1.
\end{equation}
Then the operator 
\[
   \mathcal{B}\coloneqq \partial_{t,\nu} \begin{pmatrix}
                                     \mathcal{M} & 0 \\ 0 & 0 
                                   \end{pmatrix} + \begin{pmatrix}
                                     \mathcal{N}_{00} & \mathcal{N}_{01} \\ \mathcal{N}_{10} & \mathcal{N}_{11}
                                   \end{pmatrix}
\]
defined in $L_\nu^2(\s X\times \s Y)$ with domain $\dom\left(\partial_{t,\nu}\left(\begin{smallmatrix}
                                     \mathcal{M} & 0 \\ 0 & 0 
                                   \end{smallmatrix}\right)\right)$
is continuously invertible. Moreover, the estimate 
\begin{align*}
   & \Big\|\big(\overline{\mathcal{B}^{-1}\begin{pmatrix}
                      \partial_{t,\nu} & 0 \\ 0 & 1 
                     \end{pmatrix}}-\begin{pmatrix}
                      \mathcal{M}^{-1} & 0 \\ -\mathcal{N}_{11}^{-1}\mathcal{N}_{10}\mathcal{M}^{-1}  & \mathcal{N}_{11}^{-1} 
                     \end{pmatrix} \big)
\Big\|
\\ &\leq \frac{c_1\|\mathcal{N}_{01}\|+\|\mathcal{N}_{01}\|\|\mathcal{N}_{10}\|}{c_0c_1^2\nu} + \frac{\theta}{1-\theta}\Big(\frac{1}{c_0}+\frac{\|\mathcal{N}_{01}\|}{c_0c_1\nu}  +\frac{\|\mathcal{N}_{10}\|}{c_0c_1}+\frac{\|\mathcal{N}_{01}\|\|\mathcal{N}_{10}\|}{c_0c_1^2\nu}\Big),
\end{align*}
with $\theta=(c_1 \|\mathcal{N}_{00}\|+\|\mathcal{N}_{01}\|\|\mathcal{N}_{01}\|)/(\nu c_0 c_1)$ holds true. 
\end{Theorem}

\begin{Remark}\label{r:prespr} The strategy of the proof of Theorem \ref{t:wp_ode} will be based on an explicit computation of $\mathcal{B}^{-1}$. In fact, we will show that
\begin{align*}
   & \mathcal{B}^{-1}\begin{pmatrix}
      \partial_{t,\nu}& 0 \\
      0& 1
   \end{pmatrix} 
   \\ & = \begin{pmatrix}
      \mathcal{M}^{-1} & 0 \\
      -\mathcal{N}_{11}^{-1}\mathcal{N}_{10}\mathcal{M}^{-1}& \mathcal{N}_{11}^{-1}
   \end{pmatrix} + \begin{pmatrix}
      0 & -\mathcal{M}^{-1}\partial_{t,\nu}^{-1}\mathcal{N}_{01}\mathcal{N}_{11}^{-1}  \\
      0 & \mathcal{N}_{11}^{-1}\mathcal{N}_{10}\mathcal{M}^{-1}\partial_{t,\nu}^{-1}\mathcal{N}_{01}\mathcal{N}_{11}^{-1}
   \end{pmatrix} \\
   & \quad\quad\quad\quad\;+\sum_{k=1}^\infty \begin{pmatrix}
      {\s T}^k\mathcal{M}^{-1} & -{\s T}^k\mathcal{M}^{-1}\partial_{t,\nu}^{-1}\mathcal{N}_{01}\mathcal{N}_{11}^{-1} \\ -\mathcal{N}_{11}^{-1}\mathcal{N}_{10}{\s T}^k\mathcal{M}^{-1} & 
   \mathcal{N}_{11}^{-1}\mathcal{N}_{10}{\s T}^k\mathcal{M}^{-1}\partial_{t,\nu}^{-1}\mathcal{N}_{01}\mathcal{N}_{11}^{-1}
   \end{pmatrix},
\end{align*}
where we set $\s T=-(\partial_{t,\nu}\mathcal{M})^{-1}\s R$ as well as  $\s R=\mathcal{N}_{00}-\mathcal{N}_{01}\mathcal{N}_{11}^{-1}\mathcal{N}_{10}$ and the series being convergent in operator norm.
\end{Remark}

In the next statement, we will treat the case $\s Y=\{0\}$ in a slightly more general setting.

\begin{Lemma}\label{le:wp_ode0} Let $\s X$ Hilbert space. Let $D$ be a densely defined, closed linear operator in $\s X$ mapping onto $\s X$, $M,N\in L(\s X)$, $c_d,c_m>0$. Assume for all $\phi\in \dom(D)$ 
\begin{equation}\label{eq:wp_ode0} 
\Re \langle D\phi,\phi\rangle\geq c_d\langle \phi,\phi\rangle,\quad \Re \langle M\phi,\phi\rangle \geq c_m\langle \phi,\phi\rangle.
\end{equation}
  If, in addition, the estimate
\begin{equation}\label{eq:wp_est}
   \|N\|< c_d c_m
\end{equation}
holds, then $B\coloneqq DM-N$ is continuously invertible in $\s X$,
\[
   B^{-1} = \sum_{k=0}^\infty \big((DM)^{-1}N\big)^k(DM)^{-1},\text{ and } \| B^{-1}\|\leq \frac{1}{c_d c_m-\|N\|}.
\]
\end{Lemma}
\begin{Proof} By hypothesis both $D$ and $M$ are continuously invertible. Moreover, the estimates $\|D^{-1}\|\leq 1/c_d$, $\|M^{-1}\|\leq 1/c_m$ hold, which is a consequence of \eqref{eq:wp_ode0}, see Proposition \ref{p:posd}. From \eqref{eq:wp_est} it follows that $\|(DM)^{-1}N\|\leq c_d^{-1}c_m^{-1}\|N\|\eqqcolon \theta<1$. Hence, with the help of the Neumann series we get continuous invertibility of $(1-(DM)^{-1}N)$. As a composition of continuously invertible operators, we infer continuous invertibility of $B=(DM)(1-(DM)^{-1}N)$. Moreover, we compute
\[
   B^{-1} = (1-(DM)^{-1}N)^{-1}(DM)^{-1} =\sum_{k=0}^\infty \big((DM)^{-1}N\big)^k(DM)^{-1}.
\]
In order to prove the estimate asserted in the lemma, we observe
\begin{align*}
 \| B^{-1}\|&=\big\|\sum_{k=0}^\infty \big((DM)^{-1}N\big)^k(DM)^{-1}\big\|
 \\ &\leq \sum_{k=0}^\infty \big\|\big((DM)^{-1}N\big)^k(DM)^{-1}\big\| 
 \\ &\leq \frac{1}{c_dc_m}\sum_{k=0}^\infty \theta^k 
 \\ &= \frac{1}{c_dc_m}\frac{1}{1-\theta} =\frac{1}{c_d c_m-\|N\|}.
\end{align*}
\end{Proof}
\begin{Proof}[of Theorem \ref{t:wp_ode}] To begin with, note that by Lemma \ref{le:wp_ode0}, applied to $D=\partial_{t,\nu}$, $M=\mathcal{M}$, $N=\s R=\mathcal{N}_{00}-\mathcal{N}_{01}\mathcal{N}_{11}^{-1}\mathcal{N}_{10}$ the operator
\[
   \tilde{\mathcal{B}} \coloneqq \partial_{t,\nu}\mathcal{M}+\mathcal{N}_{00}-\mathcal{N}_{01}\mathcal{N}_{11}^{-1}\mathcal{N}_{10}
\]
is continuously invertible in $L_\nu^2(\s X)$ since for $c_m=c_0$, $c_d=\nu$ (cf.~\eqref{eq:reder}) and 
\begin{equation}\label{eq:esti_N}
  \|\s R\|= \|\mathcal{N}_{00}-\mathcal{N}_{01}\mathcal{N}_{11}^{-1}\mathcal{N}_{10} \|\leq \|\mathcal{N}_{00}\|+ \|\mathcal{N}_{01}\|c_1^{-1}\|\mathcal{N}_{10}\|
\end{equation}
the needed estimate in Lemma \ref{le:wp_ode0} is warranted by hypothesis. 
One immediately verifies the computation
\[
   \mathcal{B}=\begin{pmatrix}
      \partial_{t,\nu} \mathcal{M}+\mathcal{N}_{00} & \mathcal{N}_{01} \\
      \mathcal{N}_{10} & \mathcal{N}_{11}
   \end{pmatrix} =
   \begin{pmatrix}
      1& \mathcal{N}_{01}\mathcal{N}_{11}^{-1} \\
      0& 1
   \end{pmatrix}\begin{pmatrix}
      \tilde{\mathcal{B}} & 0 \\
      0& \mathcal{N}_{11}
   \end{pmatrix}\begin{pmatrix}
      1&0 \\
      \mathcal{N}_{11}^{-1}\mathcal{N}_{10} & 1
   \end{pmatrix}.
\]
Hence, $\mathcal{B}$ is a composition of continuously invertible operators and, thus, $\mathcal{B}$ is continuously invertible. Next, we compute $\mathcal{B}^{-1}$. For this, we employ again Lemma \ref{le:wp_ode0} to get
\begin{align*}
   \tilde{\mathcal{B}}^{-1} & = \sum_{k=0}^\infty \big(-(\partial_{t,\nu}\mathcal{M})^{-1}\s R\big)^k(\partial_{t,\nu}\mathcal{M})^{-1} =\sum_{k=0}^\infty \big(-(\partial_{t,\nu}\mathcal{M})^{-1}\s R\big)^k\mathcal{M}^{-1}\partial_{t,\nu}^{-1}.
\end{align*}
Hence, we obtain 
\begin{align*}
  &\mathcal{B}^{-1}\\ & = \begin{pmatrix}
      1&0 \\
      -\mathcal{N}_{11}^{-1}\mathcal{N}_{10} & 1
   \end{pmatrix}\begin{pmatrix}
      \tilde{\mathcal{B}}^{-1} & 0 \\
      0& \mathcal{N}_{11}^{-1}
   \end{pmatrix}\begin{pmatrix}
      1& -\mathcal{N}_{01}\mathcal{N}_{11}^{-1} \\
      0& 1
   \end{pmatrix} 
   \\  & = \sum_{k=0}^\infty \begin{pmatrix}
      1&0 \\
      -\mathcal{N}_{11}^{-1}\mathcal{N}_{10} & 1
   \end{pmatrix}\begin{pmatrix}
      \big(-(\partial_{t,\nu}\mathcal{M})^{-1}\s R\big)^k\mathcal{M}^{-1}\partial_{t,\nu}^{-1} & 0 \\
      0& 0^k\mathcal{N}_{11}^{-1}
   \end{pmatrix}\begin{pmatrix}
      1& -\mathcal{N}_{01}\mathcal{N}_{11}^{-1} \\
      0& 1
   \end{pmatrix}
   \\ & = \sum_{k=0}^\infty \begin{pmatrix}
      1&0 \\
      -\mathcal{N}_{11}^{-1}\mathcal{N}_{10} & 1
   \end{pmatrix}\begin{pmatrix}
      \big(-(\partial_{t,\nu}\mathcal{M})^{-1}\s R\big)^k\mathcal{M}^{-1} & 0 \\
      0& 0^k\mathcal{N}_{11}^{-1}
   \end{pmatrix}\begin{pmatrix}
      \partial_{t,\nu}^{-1}& -\partial_{t,\nu}^{-1}\mathcal{N}_{01}\mathcal{N}_{11}^{-1} \\
      0& 1
   \end{pmatrix}
\\ & = \sum_{k=0}^\infty \begin{pmatrix}
      1&0 \\
      -\mathcal{N}_{11}^{-1}\mathcal{N}_{10} & 1
   \end{pmatrix}\begin{pmatrix}
      \big(-(\partial_{t,\nu}\mathcal{M})^{-1}\s R\big)^k\mathcal{M}^{-1} & 0 \\
      0& 0^k\mathcal{N}_{11}^{-1}
   \end{pmatrix}\begin{pmatrix}
      1& -\partial_{t,\nu}^{-1}\mathcal{N}_{01}\mathcal{N}_{11}^{-1} \\
      0& 1
   \end{pmatrix}
   \\ & \quad \times \begin{pmatrix}
      \partial_{t,\nu}^{-1}& 0 \\
      0& 1
   \end{pmatrix}.
\end{align*}
In the series expression for $\mathcal{B}^{-1}\begin{pmatrix}
      \partial_{t,\nu}& 0 \\
      0& 1
   \end{pmatrix}$ just derived, the summand for $k=0$ reads as
\begin{align}\notag
   & \begin{pmatrix}
      1&0 \\
      -\mathcal{N}_{11}^{-1}\mathcal{N}_{10} & 1
   \end{pmatrix}\begin{pmatrix}
      \mathcal{M}^{-1} & 0 \\
      0& \mathcal{N}_{11}^{-1}
   \end{pmatrix}\begin{pmatrix}
      1& -\partial_{t,\nu}^{-1}\mathcal{N}_{01}\mathcal{N}_{11}^{-1} \\
      0& 1
   \end{pmatrix}
   \\ & = \begin{pmatrix}
      \mathcal{M}^{-1} & 0 \\
      -\mathcal{N}_{11}^{-1}\mathcal{N}_{10}\mathcal{M}^{-1}& \mathcal{N}_{11}^{-1}
   \end{pmatrix}+\begin{pmatrix}
      0 & -\mathcal{M}^{-1}\partial_{t,\nu}^{-1}\mathcal{N}_{01}\mathcal{N}_{11}^{-1}  \\
      0 & \mathcal{N}_{11}^{-1}\mathcal{N}_{10}\mathcal{M}^{-1}\partial_{t,\nu}^{-1}\mathcal{N}_{01}\mathcal{N}_{11}^{-1}
   \end{pmatrix}.\label{eq:k0}
\end{align}
The norm of the second summand in \eqref{eq:k0} is bounded above by 
\[
    \frac{1}{c_0c_1\nu}\|\mathcal{N}_{01}\|+\frac{1}{c_0c_1^2\nu}\|\mathcal{N}_{01}\|\|\mathcal{N}_{10}\|.
\]
Next, for $k\geq 1$, we compute using $\s T=-(\partial_{t,\nu}\s M)^{-1}\s R$
\begin{align}
  &\begin{pmatrix}
      1&0 \\
      -\mathcal{N}_{11}^{-1}\mathcal{N}_{10} & 1
   \end{pmatrix}\begin{pmatrix}
      {\s T}^k\mathcal{M}^{-1} & 0 \\
      0& 0
   \end{pmatrix}\begin{pmatrix}
      1& -\partial_{t,\nu}^{-1}\mathcal{N}_{01}\mathcal{N}_{11}^{-1} \\
      0& 1
   \end{pmatrix} \notag
   \\ & =  \begin{pmatrix}
      1&0 \\
      -\mathcal{N}_{11}^{-1}\mathcal{N}_{10} & 1
   \end{pmatrix}\begin{pmatrix}
      {\s T}^k\mathcal{M}^{-1} & -{\s T}^k\mathcal{M}^{-1}\partial_{t,\nu}^{-1}\mathcal{N}_{01}\mathcal{N}_{11}^{-1} \\ 0 & 0
   \end{pmatrix} \notag
   \\ & =\begin{pmatrix}
      {\s T}^k\mathcal{M}^{-1} & -{\s T}^k\mathcal{M}^{-1}\partial_{t,\nu}^{-1}\mathcal{N}_{01}\mathcal{N}_{11}^{-1} \\ -\mathcal{N}_{11}^{-1}\mathcal{N}_{10}{\s T}^k\mathcal{M}^{-1} & 
   \mathcal{N}_{11}^{-1}\mathcal{N}_{10}{\s T}^k\mathcal{M}^{-1}\partial_{t,\nu}^{-1}\mathcal{N}_{01}\mathcal{N}_{11}^{-1}
   \end{pmatrix}.\label{eq:kthop}
\end{align}
Thus, with the help of estimate \eqref{eq:esti_N}, we obtain
\[
   \|\s T\|=\|(\partial_{t,\nu}\s M)^{-1}\s R\|\leq \frac{1}{\nu c_0}\left(\|\mathcal{N}_{00}\|+ \|\mathcal{N}_{01}\|c_1^{-1}\|\mathcal{N}_{10}\|\right)\eqqcolon\theta.
\]
Therefore, for any $k\in \mathbb{N}_{\geq1}$ an estimate for the operator norm of the matrix given in \eqref{eq:kthop} reads
\begin{multline*}
\left\|\begin{pmatrix}
      {\s T}^k\mathcal{M}^{-1} & -{\s T}^k\mathcal{M}^{-1}\partial_{t,\nu}^{-1}\mathcal{N}_{01}\mathcal{N}_{11}^{-1} \\ -\mathcal{N}_{11}^{-1}\mathcal{N}_{10}{\s T}^k\mathcal{M}^{-1} & 
   \mathcal{N}_{11}^{-1}\mathcal{N}_{10}{\s T}^k\mathcal{M}^{-1}\partial_{t,\nu}^{-1}\mathcal{N}_{01}\mathcal{N}_{11}^{-1}
   \end{pmatrix}\right\| 
   \\ \leq    \theta^k\Big(\frac{1}{c_0}+\frac{\|\mathcal{N}_{01}\|}{c_0c_1\nu}  +\frac{\|\mathcal{N}_{10}\|}{c_0c_1}+\frac{\|\mathcal{N}_{01}\|\|\mathcal{N}_{10}\|}{c_0c_1^2\nu}\Big).
\end{multline*}
So, for the expression
\begin{align*}
   &\big(\mathcal{B}^{-1}\begin{pmatrix}
      \partial_{t,\nu}& 0 \\
      0& 1
   \end{pmatrix}  - \begin{pmatrix}
      \mathcal{M}^{-1} & 0 \\
      -\mathcal{N}_{11}^{-1}\mathcal{N}_{10}\mathcal{M}^{-1}& \mathcal{N}_{11}^{-1}
   \end{pmatrix}\big)
   \\& \subseteq \begin{pmatrix}
      0 & -\mathcal{M}^{-1}\partial_{t,\nu}^{-1}\mathcal{N}_{01}\mathcal{N}_{11}^{-1}  \\
      0 & \mathcal{N}_{11}^{-1}\mathcal{N}_{10}\mathcal{M}^{-1}\partial_{t,\nu}^{-1}\mathcal{N}_{01}\mathcal{N}_{11}^{-1}
   \end{pmatrix} \\
   & +\sum_{k=1}^\infty \begin{pmatrix}
      {\s T}^k\mathcal{M}^{-1} & -{\s T}^k\mathcal{M}^{-1}\partial_{t,\nu}^{-1}\mathcal{N}_{01}\mathcal{N}_{11}^{-1} \\ -\mathcal{N}_{11}^{-1}\mathcal{N}_{10}{\s T}^k\mathcal{M}^{-1} & 
   \mathcal{N}_{11}^{-1}\mathcal{N}_{10}{\s T}^k\mathcal{M}^{-1}\partial_{t,\nu}^{-1}\mathcal{N}_{01}\mathcal{N}_{11}^{-1}
   \end{pmatrix},
\end{align*}
we find
\begin{align*}
   &\big\|\mathcal{B}^{-1}\begin{pmatrix}
      \partial_{t,\nu}& 0 \\
      0& 1
   \end{pmatrix}  - \begin{pmatrix}
      \mathcal{M}^{-1} & 0 \\
      -\mathcal{N}_{11}^{-1}\mathcal{N}_{10}\mathcal{M}^{-1}& \mathcal{N}_{11}^{-1}
   \end{pmatrix}\big\|
   \\& \leq \frac{1}{c_0c_1\nu}\|\mathcal{N}_{01}\|+\frac{1}{c_0c_1^2\nu}\|\mathcal{N}_{01}\|\|\mathcal{N}_{10}\| 
   \\ &\quad +\sum_{k=1}^\infty \theta^k\Big(\frac{1}{c_0}+\frac{\|\mathcal{N}_{01}\|}{c_0c_1\nu}  +\frac{\|\mathcal{N}_{10}\|}{c_0c_1}+\frac{\|\mathcal{N}_{01}\|\|\mathcal{N}_{10}\|}{c_0c_1^2\nu}\Big)
   \\ & = \frac{c_1\|\mathcal{N}_{01}\|+\|\mathcal{N}_{01}\|\|\mathcal{N}_{10}\|}{c_0c_1^2\nu} \\ & \quad + \frac{\theta}{1-\theta}\Big(\frac{1}{c_0}+\frac{\|\mathcal{N}_{01}\|}{c_0c_1\nu}  +\frac{\|\mathcal{N}_{10}\|}{c_0c_1}+\frac{\|\mathcal{N}_{01}\|\|\mathcal{N}_{10}\|}{c_0c_1^2\nu}\Big).
\end{align*}
\end{Proof}

We conclude this section with the solution theory for linear abstract ordinary differential equations in the context of evolutionary mappings:

\begin{Theorem}\label{t:wp_ode_evo} Let $\s X,\s Y$ Hilbert spaces, $\mathcal{M}\in L_{\textnormal{sev}}(\s X)$, and $\s N=(\mathcal{N}_{ij})_{i,j\in\{0,1\}}\in L_{\textnormal{sev}}(\s X\times \s Y)$. Assume
that there exists $c>0$ such that for all $(\phi,\psi)\in \dom(\s M)\times \dom(\s N_{11})$ the positive definiteness conditions
\begin{equation}\label{eq:wp_ode_evo1} 
   \Re \langle Q_t\mathcal{M} \phi,\phi\rangle_{L_\nu^2(\s X)}\geq c \langle Q_t\phi,\phi\rangle_{L_\nu^2(\s X)},\quad \Re \langle Q_t\mathcal{N}_{11}\psi,\psi\rangle_{L_\nu^2(\s Y)}\geq c\langle Q_t\psi,\psi\rangle_{L_\nu^2(\s Y)}
\end{equation}
hold for all $t\in \mathbb{R}$ and eventually all $\nu$ large enough, where $Q_t$ is multiplication by $\1_{(-\infty,t)}$.
Then the operator 
\[
   \mathcal{B}\coloneqq \partial_{t,\nu} \begin{pmatrix}
                                     \mathcal{M} & 0 \\ 0 & 0 
                                   \end{pmatrix} + \begin{pmatrix}
                                     \mathcal{N}_{00} & \mathcal{N}_{01} \\ \mathcal{N}_{10} & \mathcal{N}_{11}
                                   \end{pmatrix}
\]
is one-to-one in $L_\nu^2(\s X\times \s Y)$ for all $\nu$ satisfying \eqref{eq:wp_ode_evo1} and 
\[
  c\|\s N_{00}^\nu\|+\|\s N_{10}^\nu\|\|\s N_{01}^\nu\| < \nu c^2.
\]
Moreover, $\mathcal{B}^{-1} \in L_{\textnormal{sev}}(\s X\times \s Y)$ and we have
\begin{equation}\label{eq:wp_ode_evo2}
     \limsup_{\nu\to\infty} \Big\|\big(\overline{\mathcal{B}^{-1}\begin{pmatrix}
                      \partial_{t,\nu} & 0 \\ 0 & 1 
                     \end{pmatrix}}-\begin{pmatrix}
                      \mathcal{M}^{-1} & 0 \\ -\mathcal{N}_{11}^{-1}\mathcal{N}_{10}\mathcal{M}^{-1}  & \mathcal{N}_{11}^{-1} 
                     \end{pmatrix} \big)
\Big\|_{L(L_\nu^2(\s X\times \s Y))}=0.
\end{equation}
\end{Theorem}
\begin{Proof}
 By Corollary \ref{c:inv_ce}, both $\s M^{-1}$ and $\s N^{-1}_{11}$ are evolutionary and causal. Thus, by Theorem \ref{t:evsc_sev}, (the standard realization of) $\s M^{-1}$ and $\s N^{-1}_{11}$ are standard evolutionary. In order to apply Theorem \ref{t:wp_ode} in the present context, we need to warrant inequality \eqref{eq:wp_ode2}, that is, we need to show that
 \begin{equation}\label{eq:ghj}
    c\|\s N_{00}^\nu\|+\|\s N_{10}^\nu\|\|\s N_{01}^\nu\| < \nu c^2
 \end{equation}
 for eventually all $\nu$ large enough. But, note that the left-hand side of the latter inequality remains bounded as $\nu\to\infty$ since $\s N$ is evolutionary. The right-hand side, however, blows up as $\nu\to\infty$. Hence, for eventually all $\nu$ large enough the inequality corresponding to \eqref{eq:wp_ode2}, that is, \eqref{eq:ghj}, is satisfied. Next, by Remark \ref{r:prespr}, $\s B^{-1}$ can be represented as an (infinite) sum of compositions of standard evolutionary mappings. Observing that the partial sums are standard evolutionary by Proposition \ref{p:sev_alg2} and bounded, we infer together with the convergence in operator norm (see Remark \ref{r:prespr}) for every $\nu$ large enough that the series converges to a standard evolutionary mapping. The last assertion of the lemma, that is, \eqref{eq:wp_ode_evo2}, follows from the estimate in Theorem \ref{t:wp_ode}.
\end{Proof}

\begin{Remark}\label{r:prest} Appealing to Theorem \ref{t:wp_ode}, we find the following estimate to hold true
\begin{align*}
   & \Big\|\big(\overline{\mathcal{B}^{-1}\begin{pmatrix}
                      \partial_{t,\nu} & 0 \\ 0 & 1 
                     \end{pmatrix}}-\begin{pmatrix}
                      \mathcal{M}^{-1} & 0 \\ -\mathcal{N}_{11}^{-1}\mathcal{N}_{10}\mathcal{M}^{-1}  & \mathcal{N}_{11}^{-1} 
                     \end{pmatrix} \big)
\Big\|_{L(L_\nu^2(\s X\times \s Y))}
\\ &\leq \frac{c\|\mathcal{N}_{01}\|_{\geq\nu}+\|\mathcal{N}_{01}\|_{\geq\nu}\|\mathcal{N}_{10}\|_{\geq\nu}}{c^3\nu} 
\\ &\quad + \frac{\theta}{1-\theta}\Big(\frac{1}{c}+\frac{\|\mathcal{N}_{01}\|_{\geq\nu}}{c^2\nu}  +\frac{\|\mathcal{N}_{10}\|_{\geq\nu}}{c^2}+\frac{\|\mathcal{N}_{01}\|_{\geq\nu}\|\mathcal{N}_{10}\|_{\geq\nu}}{c^3\nu}\Big),
\end{align*}
with $\theta=(c \|\mathcal{N}_{00}\|_{\geq\nu}+\|\mathcal{N}_{01}\|_{\geq\nu}\|\mathcal{N}_{01}\|_{\geq\nu})/(\nu c^2)$ for all $\nu>0$ large enough, where we set $\|\s R\|_{\geq\nu}\coloneqq \sup_{\mu\geq\nu}\|\s R^\mu\|_{L(L_\mu^2)}$ for a suitable evolutionary mapping $\s R$.
\end{Remark}

\renewcommand{\baselinestretch}{0.65}\normalsize\mysection{Partial Differential Equations -- Preliminaries}{Partial Differential Equations -- Preliminaries}{Proposition \ref{p:com_resA} $\cdot$ Corollary \ref{c:tdcapAdense} $\cdot$ Proposition \ref{p:Atrinv} $\cdot$ Lemma \ref{le:commus_of_a}}\label{sec:wpr0}

\renewcommand{\baselinestretch}{1}\normalsize
In the abstract settings discussed here, the main difference of ordinary differential equations and partial differential equations is the occurrence of another unbounded linear operator apart from the time derivative. Indeed, in order to cope with many commonly known linear evolutionary equations from mathematical physics, one has to take into account spatial derivatives, as well. So, the general equation $\s Bu=f$ from \eqref{eq:ade} to be studied in the following admits the more precise form
\begin{equation}\label{eq:pde0}
   \left(\partial_{t,\nu} \s M + \s N + \s A\right) u = f.
\end{equation}
In this preliminary section, we will introduce some of the assumptions on the operators $\s M$, $\s N$ (the `material law') and $\s A$ (the `unbounded spatial operator') and some of its consequences. Moreover, we will have the occasion to provide some results of a more general nature to be used later on.

The operator $\s A$ (see Hypothesis \ref{hyp:A}) is thought of containing the spatial derivatives. In manifold applications $\s A$ is an (unbounded) skew-selfadjoint operator (in the underlying spatial Hilbert space). For incorporating more involved evolutionary problems, we will, however, relax this condition. The main assumption is roughly rephrased by both the numerical range of $\s A$ and of its adjoint lying in a right half plane of the complex numbers. Further, anticipating the fact that in applications $\s A$ contains the spatial derivatives, we will assume a compatibility condition for $\s A$ with $\partial_{t,\nu}$.  

The assumptions on the coefficient $\s M$ of $\partial_{t,\nu}$ -- in comparison to the ODE-case -- have to be strengthened in the way that they should boundedly commute with time-differ\-enti\-ation (see Hypothesis \ref{hyp:mat_law}), which reflects the fact that if treating multiplication operators with operators depending explicitly on time these operators should be Lip\-schitz continuous, see also Example \ref{ex:mult_op}.

Next, we introduce the hypothesis on $\s A$:

\begin{Hypothesis}[on the unbounded spatial operator]\label{hyp:A} Let $\s X$ be a Hilbert space, $\nu>0$. Assume $\s A\colon \dom(\s A)\subseteq L_\nu^2(\R;\s X)\to L_\nu^2(\R;\s X)$ to be densely defined, closed, linear and such that $\partial_{t,\nu}^{-1}\s A \subseteq \s A \partial_{t,\nu}^{-1}$.
\end{Hypothesis}

\begin{Proposition}\label{p:com_resA} Assume Hypothesis \ref{hyp:A} to be satisfied. Then, for all $\epsilon>0$, the following inclusion holds
\[
    \left(1+\epsilon \partial_{t,\nu}\right)^{-1}\s A \subseteq \s A\left(1+\epsilon\partial_{t,\nu}\right)^{-1}
\] 
\end{Proposition}

For the proof of Proposition \ref{p:com_resA} some preparations are in order.

\begin{Remark}\label{r:eps_st} (a) A reason for merely being interested in $\left(1+\epsilon\partial_{t,\nu}\right)^{-1}$ rather than $\partial_{t,\nu}^{-1}$ as in Hypothesis \ref{hyp:A} is the following. The resolvents $\left(1+\epsilon\partial_{t,\nu}\right)^{-1}$ converge strongly to the identity operator as $\epsilon \to 0$. Indeed, since $\Re (1+\epsilon\partial_{t,\nu})\geq 1$ the operators $\left(1+\epsilon\partial_{t,\nu}\right)^{-1}$ are contractions for all $\epsilon>0$, see equation \eqref{eq:reder} together with Proposition \ref{p:posd} (recall $\nu>0$). Hence, it suffices to verify
\begin{equation}\label{eq:r_eps}
   \left(1+\epsilon\partial_{t,\nu}\right)^{-1}\phi \to \phi \quad (\epsilon\to 0)
\end{equation}
for $\phi$ belonging to the dense set $\dom(\partial_{t,\nu})\subseteq L_\nu^2(\s X)$. In order to deduce \eqref{eq:r_eps}, it is sufficient to observe for $\epsilon>0$ and $\phi\in \dom(\partial_{t,\nu})$
\begin{align}
   \left(1+\epsilon\partial_{t,\nu}\right)^{-1}\phi& =\left(1+\epsilon\partial_{t,\nu}\right)^{-1}(\phi+\epsilon\partial_{t,\nu}\phi)-\epsilon\left(1+\epsilon\partial_{t,\nu}\right)^{-1}\partial_{t,\nu}\phi \notag
   \\ &=\phi-\epsilon\left(1+\epsilon\partial_{t,\nu}\right)^{-1}\partial_{t,\nu}\phi,\label{eq:e_eps2}
\end{align}
where the last summand converges to $0$ as $\epsilon\to 0$.

(b) Similar to the convergence result in \eqref{eq:r_eps}, we observe the following
\begin{equation}\label{eq:r_eps2}
   \epsilon\partial_{t,\nu}\left(1+\epsilon\partial_{t,\nu}\right)^{-1} \phi \to 0\quad  (\epsilon\to 0)
\end{equation}
for all $\phi \in L_{\nu}^2(\R;\s X)$. In fact, from \eqref{eq:e_eps2} we read off that
\[
   \left(1+\epsilon\partial_{t,\nu}\right)^{-1}-1 =\epsilon\partial_{t,\nu}\left(1+\epsilon\partial_{t,\nu}\right)^{-1}.
\]
Hence, in view of (a), the left-hand side converges strongly to $0$ as $\eps\to0$, thus, so does the right-hand side.
\end{Remark}

\begin{Lemma}\label{l:resol} Let $\nu>0$ and $\epsilon>0$. Then for $r\coloneqq 1/(2\nu)$, we have the expansion
\[
   \left(1+\epsilon\partial_{t,\nu}\right)^{-1} = \epsilon^{-1}\partial_{t,\nu}^{-1} \frac{1}{1+(r/\epsilon)}\sum_{k=0}^\infty \left(\frac{(r-\partial_{t,\nu}^{-1})/\epsilon}{1+(r/\epsilon)}\right)^k.
\] 
\end{Lemma}
\begin{Proof}
   For all $z \in \partial B(r,r)$, we compute
  \begin{align*}
     \frac{1}{1+(z/\epsilon)}& =\frac{1}{1+(r/\epsilon)+((z-r)/\epsilon)}
     \\ &=\frac{1}{1+(r/\epsilon)}\frac{1}{1+\left(\frac{(z-r)/\epsilon}{1+(r/\epsilon)}\right)}
     \\ &=\frac{1}{1+(r/\epsilon)}\sum_{k=0}^\infty \left(\frac{(r-z)/\epsilon}{1+(r/\epsilon)}\right)^k,
  \end{align*}
  where the series converges uniformly in $z$. By the functional calculus induced by the unitary equivalence stated in Corollary \ref{co:spe_thm} for $\partial_{t,\nu}^{-1}$ (see also Theorem \ref{th:sp_td_inv}) we deduce that
  \[
    \left(1+\epsilon\partial_{t,\nu}\right)^{-1}  =\epsilon^{-1}\partial_{t,\nu}^{-1} \frac{1}{1+(r/\epsilon)}\sum_{k=0}^\infty \left(\frac{(r-\partial_{t,\nu}^{-1})/\epsilon}{1+(r/\epsilon)}\right)^k
  \]
  with convergence in operator norm, yielding the assertion.
\end{Proof}

\begin{Remark}\label{r:rtc} Note that $(1+\eps\partial_{t,\nu})^{-1}$ is translation-invariant and causal. In fact, this is a straightforward consequence of Theorem \ref{t:SW_conv}. For the proof of Theorem \ref{t:SW_conv}, we used the Paley--Wiener theorem, which we stated without proof. For having a self-contained proof of translation-invariance and causality of $(1+\eps\partial_{t,\nu})^{-1}$ in this exposition, we argue as follows. Indeed, the claim is a consequence of the explicit formula for $\partial_{t,\nu}^{-1}$ in Lemma \ref{eq:td_inv} (recall $\nu>0$) and the representation in Lemma \ref{l:resol}: The latter representation immediately yields translation-invariance (since $\partial_{t,\nu}^{-1}$ is translation-invariant) and that $(1+\eps\partial_{t,\nu})^{-1}$ leaves functions supported on $[0,\infty)$ invariant as well (since so does $\partial_{t,\nu}^{-1}$ by Lemma \ref{eq:td_inv}). Hence, $(1+\eps\partial_{t,\nu})^{-1}$ is causal, by Remark \ref{r:tinv_caus}. 
\end{Remark}

\begin{Remark}\label{r:stron_op_con} (a) Aiming for a proof of Proposition \ref{p:com_resA}, we will employ Lemma \ref{l:resol}. For this observe the following elementary fact. In a Hilbert space $\s X$, let $A$ be a closed linear operator and assume that there exists a sequence $(T_n)_n$ of bounded linear operators in $\s X$ being strongly convergent to some $T\in L(\s X)$. If, for all $n\in \mathbb{N}$, $T_n A\subseteq  AT_n$, then $T  A\subseteq   A T$. Indeed, let $\phi\in \dom( A)$. Then, by hypothesis, $T_n\phi\in \dom(  A)$, $n\in \mathbb{N}$, and $T_n \phi\to T\phi$ as well as $  A T_n \phi = T_n A\phi \to T  A\phi$ as $n\to\infty$. By the closedness of $  A$, we infer $T\phi\in \dom(  A)$ and $  A T\phi=T  A\phi$, which is the claim.

(b) Another observation being used in the following is in order. Assume $A$ to be only closable in $\s X$, $T, T'\in L(\s X)$. If $TA\subseteq \overline{A}T+ T'$, then $T\overline{A}\subseteq \overline{A}T+ T'$. Indeed, for $\phi\in \dom(\overline{A})$, let $(\phi_n)_n$ in $\dom(A)$ converge to $\phi$ in $\s X$ with $(A\phi_n)_n\to \overline{A}\phi$ in $\s X$. Then, for all $n\in \N$, we have 
\[
   \overline{A}T\phi_n=TA\phi_n-T'\phi_n \to  T\overline{A}\phi - T'\phi\quad (n\to\infty).
\]
Hence, by the closedness of $\overline{A}$ and continuity of $T$ and $T'$, we infer $T\phi\in \dom(\overline{A})$ and $\overline{A}T\phi=T\overline{A}\phi-T'\phi$. Thus, in other words, $T\overline{A}\subseteq \overline{A}T+T'$.
\end{Remark}

\begin{Proof}[of Proposition \ref{p:com_resA}] Appealing to Lemma \ref{l:resol} and Remark \ref{r:stron_op_con}, we observe that the assertion is proved, if we show for $n\in \mathbb{N}$ that
\begin{multline*}
   \epsilon^{-1}\partial_{t,\nu}^{-1} \frac{1}{1+(r/\epsilon)}\sum_{k=0}^n \left(\frac{(r-\partial_{t,\nu}^{-1})/\epsilon}{1+(r/\epsilon)}\right)^k \s A 
   \\ \subseteq \s A \epsilon^{-1}\partial_{t,\nu}^{-1} \frac{1}{1+(r/\epsilon)}\sum_{k=0}^n \left(\frac{(r-\partial_{t,\nu}^{-1})/\epsilon}{1+(r/\epsilon)}\right)^k.
\end{multline*}
The latter inclusion, however, is a straightforward consequence of 
\[
 (\lambda+\partial_{t,\nu}^{-1})\s A  =\lambda\s A+\partial_{t,\nu}^{-1}\s A 
 \subseteq \s A \lambda+\s A \partial_{t,\nu}^{-1}
  \subseteq \s A (\lambda+ \partial_{t,\nu}^{-1}),
\]
for all $\lambda\in \R$, where we used Hypothesis \ref{hyp:A}.
\end{Proof}

\begin{Corollary}\label{c:tdcapAdense} Assume Hypotheses \ref{hyp:A} to be satisfied. Then $\dom(\partial_{t,\nu})\cap \dom(\s A)$ is dense in $\dom(\s A)$, endowed with the graph norm of $\s A$, as well as dense in $L_\nu^2(\R;\s X)$. 
\end{Corollary}
\begin{Proof}
 Let $\eps>0$, $f\in \dom(\s A)$. By Proposition \ref{p:com_resA}, we have $(1+\eps\partial_{t,\nu})^{-1}\s A \subseteq\s A(1+\eps\partial_{t,\nu})^{-1}$. Thus, $(1+\eps\partial_{t,\nu})^{-1}f \in \dom(\partial_{t,\nu})\cap \dom(\s A)$. Since, by Remark \ref{r:eps_st}, \[(1+\eps\partial_{t,\nu})^{-1}f\to f\text{ and }\s A(1+\eps\partial_{t,\nu})^{-1}f=(1+\eps\partial_{t,\nu})^{-1}\s Af\to \s Af\text{ as }\eps\to 0,\] we conclude that $\dom(\partial_{t,\nu})\cap \dom(\s A)$ is dense in $\dom(\s A)$ endowed with the graph norm of $\s A$. The density of $\dom(\s A)$ in $L_\nu^2(\R;\s X)$ yields the second assertion.
\end{Proof}

With the techniques just employed, we can also show that $\s A$ is actually translation-invariant:

\begin{Proposition}\label{p:Atrinv} Assume Hypotheses \ref{hyp:A} to be satisfied. Then for all $h\in \mathbb{R}$, we have
\[
   \tau_h \s A\subseteq \s A\tau_h.
\] 
\end{Proposition}

To begin with, we represent the time translation as a function of the (inverse) time derivative.

\begin{Lemma}\label{l:tt_ftd} Let $\nu>0$, $r\coloneqq 1/(2\nu)$, $h\in \mathbb{R}$, $\s X$ Hilbert space. Then:

(a) For all $\epsilon>0$ the series
\[
  \tau_{h,\eps}\coloneqq \sum_{k=0}^\infty \frac{1}{k!}\left(\frac{h}{r+\epsilon}\sum_{\ell=0}^\infty\left(\frac{r-\partial_{t,\nu}^{-1}}{r+\epsilon}\right)^\ell\right)^k
\]
converges absolutely in $L(L_\nu^2(\mathbb{R};\s X))$.

(b) $\tau_{h,\eps}\to \tau_h$ as $\eps\to 0$ in the strong operator topology of $L(L_\nu^2(\mathbb{R};\s X))$.

(c) There exists $(p_n)_n$, a sequence of polynomials, such that $p_n(\partial_{t,\nu}^{-1})\to \tau_h$ as $n\to\infty$ in the strong operator topology of $L(L_\nu^2(\mathbb{R};\s X))$.
\end{Lemma}
\begin{Proof}
 (a) Let $\eps>0$. The convergence of the series can be seen by Fourier--Laplace transformation. Indeed, the series 
\[
  t_{h,\eps}(z)\coloneqq  \sum_{k=0}^\infty \frac{1}{k!}\left(\frac{h}{r+\epsilon}\sum_{\ell=0}^\infty\left(\frac{r-z}{r+\epsilon}\right)^\ell\right)^k
\]
converges uniformly for all $z\in \partial B(r,r)$.

(b) For $\eps>0$ using the formula for the exponential function and Neumann's series, we get for $z\in \partial B(r,r)$
\begin{align*}
   t_{h,\eps}(z)& = \exp\left( \frac{h}{r+\epsilon}\sum_{\ell=0}^\infty\left(\frac{r-z}{r+\epsilon}\right)^\ell\right) 
   \\           & = \exp\left( \frac{h}{r+\epsilon}\frac{1}{1- \frac{r-z}{r+\epsilon}}\right)
   \\           & = \exp\left( \frac{h}{r+\epsilon + z-r}\right) = \exp\left( \frac{1}{\epsilon + z}h\right).
\end{align*}
 From (dominated) pointwise convergence of $t_{h,\eps}$ to $t_{h}\colon z\mapsto e^{h/z}$ on $\partial B(r,r)$, we get, using Lebegue's dominated convergence theorem, that the multiplication operators associated with $t_{h,\eps}$ converge in the strong operator topology of $L(L^2(\mathbb{R};\s X))$ to the respective multiplication operator associated with $t_h$ as $\eps\to 0$. Via Fourier--Laplace transformation, (the multiplication operator associated with) $t_h$ is unitarily equivalent to $\tau_h$. Thus, the assertion follows.
 
 (c) This is a combination of the absolute convergence asserted in part (a) and the strong operator convergence of part (b).
\end{Proof}

\begin{Proof}[of Proposition \ref{p:Atrinv}] By Lemma \ref{l:tt_ftd}(c), the operator $\tau_h$ of time translation can be approximated by a sequence of polynomials $(p_n)_n$ applied to $\partial_{t,\nu}^{-1}$ with respect to the strong operator topology. An application of Remark \ref{r:stron_op_con} with $T_n=p_n(\partial_{t,\nu}^{-1})$ and $A=\mathcal{A}$ yields the assertion.
\end{Proof}

We conclude this preliminary section, with the hypotheses on $\s M$ and $\s N$ in \eqref{eq:pde0}, and a small consequence thereof:

\begin{Hypothesis}[on the material law]\label{hyp:mat_law} Let $\s X$ Hilbert space, $\s M, \s N \in L(L_\nu^2(\R;\s X))$. Assume that there exists $\s M'\in L(L_\nu^2(\R;\s X))$ such that 
\[
   \s M \partial_{t,\nu} \subseteq \partial_{t,\nu} \s M - \s M'.
\] 
\end{Hypothesis}

We denote the commutator of two operators $A,B$ by $[A,B]=AB-BA$ with its natural domain. Hence, Hypothesis \ref{hyp:mat_law} properly rephrased reads $[\partial_{t,\nu},\s M]\subseteq \s M'\in L_{\nu}^2(\R;\s X)$.

\begin{Lemma}\label{le:commus_of_a} Assume that $\s M, \s N$ satisfy Hypotheses \ref{hyp:mat_law}. Then the following statements hold:

(a) For every $\eps>0$, we have:
\[\overline{[(1+\eps\partial_{t,\nu})^{-1},\partial_{t,\nu} \s M]}=-\eps\partial_{t,\nu} (1+\eps\partial_{t,\nu})^{-1}\s M'(1+\eps\partial_{t,\nu})^{-1}\] 
Moreover,
\[
   \overline{[(1+\eps\partial_{t,\nu})^{-1},\partial_{t,\nu} \s M]}\stackrel{\tau_{\textnormal{s}}}{\to} 0 \quad (\eps\to 0)
\] in $L_\nu^2(\R;\s X)$, that is, the closure of the commutator converges to $0$ in the strong operator topology $\tau_{\textnormal{s}}$.

(b) For $\eps>0$, we have
\[
    [(1+\eps\partial_{t,\nu})^{-1},\s N]\stackrel{\tau_{\textnormal{s}}}{\to} 0 \quad(\eps \to 0).
\]
 \end{Lemma}
\begin{Proof}
(a) Let $u\in \dom(\partial_{t,\nu})$. We compute
\begin{align*}
  &[(1+\eps\partial_{t,\nu})^{-1},\partial_{t,\nu} \s M]u
  \\&= (1+\eps\partial_{t,\nu})^{-1}\partial_{t,\nu} \s M u - \partial_{t,\nu} \s M (1+\eps\partial_{t,\nu})^{-1}u\\
            & =(1+\eps\partial_{t,\nu})^{-1}\left(\partial_{t,\nu} \s M(1+\eps\partial_{t,\nu})-(1+\eps\partial_{t,\nu})\partial_{t,\nu} \s M \right)(1+\eps\partial_{t,\nu})^{-1}u \\
            & =-\eps\partial_{t,\nu}(1+\eps\partial_{t,\nu})^{-1} \s M'(1+\eps\partial_{t,\nu})^{-1}u. 
\end{align*}
By Remark \ref{r:eps_st} (a) and (b) together with the formula just derived, we get the desired convergence result in (a).

(b) The second statement follows from the continuity of $\s N$ and $(1+\eps\partial_{t,\nu})^{-1}\to 1$ as $\eps\to 0$ in the strong operator topology, by Remark \ref{r:eps_st} (a).
\end{Proof}

\renewcommand{\baselinestretch}{0.65}\normalsize\mysection{Partial Differential Equations -- Invertibility}{Partial Differential Equations -- Invertibility}{Computation of $[(1+\eps\partial_{t,\nu})^{-1},\s B]$ $\cdot$ for $u\in \dom(\overline{\s B})$ we have $(1+\eps\partial_{t,\nu})^{-1}u\in \dom(\s B)$ $\cdot$ the adjoint of $\s B$ $\cdot$ Theorem \ref{thm:genSolth}}

\renewcommand{\baselinestretch}{1}\normalsize
Next, we come to the announced result concerning the continuous invertibility of the operator sum $\s B = \partial_{t,\nu}\s M+\s N+\s A$, that is, the well-posedness of the abstract partial differential equation as in \eqref{eq:pde0}. In the next section, strengthening the positive definiteness requirement stated in the forthcoming hypothesis, we will address both causality and evolutionarity. 

The continuous invertibility result will be formulated in the following situation.
\begin{Hypothesis}\label{h:gS} Let $\s M,\s N,\s A$ be as in Hypotheses \ref{hyp:mat_law} and \ref{hyp:A}, that is, $\s M,\s N,\s M'\in L(L_\nu^2(\R;\s X))$, $\s A$ densely defined, closed in $L_\nu^2(\R;\s X)$ for some Hilbert space $\s X$ and $\nu>0$ with
\[
   \s M\partial_{t,\nu}\subseteq \partial_{t,\nu}\s M  - \s M'\text{ and } \partial_{t,\nu}^{-1}\s A\subseteq \s A\partial_{t,\nu}^{-1}.
\]
Furthermore, assume there exists $c>0$ such that the positivity conditions
\begin{equation}\label{ineq:pos}
 \Re\langle \left(\partial_{t,\nu}\s M+\s N +\s A\right)\phi,\phi\rangle \geq c \langle \phi,\phi\rangle
\end{equation}
and 
\begin{equation}\label{ineq:adjoint}
   \Re\langle \left(\left(\partial_{t,\nu}\s M+\s N\right)^* +\s A^*\right)\psi,\psi\rangle \geq c \langle \psi,\psi\rangle
\end{equation}
  hold for all $\phi\in \dom(\partial_{t,\nu}\s M)\cap \dom(\s A)$, $\psi\in \dom(\partial_{t,\nu})\cap \dom(\s A^*)$.
 We define $\s B\coloneqq {\partial_{t,\nu}\s M+\s N+\s A}$ with $\dom(\s B)=\dom(\partial_{t,\nu}\s M)\cap \dom(\s A)$.
\end{Hypothesis}

\begin{Theorem}\label{thm:genSolth} Assume Hypothesis \ref{h:gS}.
Then $\s B$ is closable, $\ran(\s B)$ is dense in $L_\nu^2(\R;\s X)$; $\overline{\s B}$ is continuously invertible in $L_\nu^2(\R;\s X)$ with $\Abs{\overline{\s B}^{-1}}\leq 1/c$.
\end{Theorem}

The idea for proving Theorem \ref{thm:genSolth} is to invoke Proposition \ref{p:posd}. For this, we want to understand $\dom(\overline{\s B})$ in a better way:

\begin{Lemma}\label{l:com_B}
Assume Hypothesis \ref{h:gS}. Then, $\s B$ is closable, and for $\eps>0$ and $u\in \dom(\overline{\s B})$ we have $(1+\eps\partial_{t,\nu})^{-1}u\in \dom(\overline{\s B})$ as well as the formula
\begin{align}
   &\left(1+\eps\partial_{t,\nu}\right)^{-1}\overline{\s B} u \notag
   \\ &= \overline{\s B}\left(1+\eps\partial_{t,\nu}\right)^{-1}u - \eps\partial_{t,\nu}(1+\eps\partial_{t,\nu})^{-1}\s M'\left(1+\eps\partial_{t,\nu}\right)^{-1}u  \label{eq:doad}
  \\ &\quad \notag+ \left[ \left(1+\eps\partial_{t,\nu}\right)^{-1},\s N\right] u.
\end{align}
\end{Lemma}
\begin{Proof}
 First of all note that the closability follows from Proposition \ref{p:posd} by inequality \eqref{ineq:pos} and Corollary \ref{c:tdcapAdense} in order that $\dom(\partial_{t,\nu})\cap \dom(\s A)\subseteq \dom(\s B)$ is dense in $L_\nu^2(\R;\s X)$. We recall Lemma \ref{le:commus_of_a} for the formulas of the commutators of $(1+\eps\partial_{t,\nu})^{-1}$ with the operators $\partial_{t,\nu}\s M$ and $\s N$, and Proposition \ref{p:com_resA} for the respective one with $\s A$.
 Then, for $u\in \dom(\partial_{t,\nu}\s M)\cap \dom(\s A)=\dom(\s B)$, we compute
 \begin{align*}
  & (1+\eps\partial_{t,\nu})^{-1}\s B u
  \\ & = (1+\eps\partial_{t,\nu})^{-1}\left(\partial_{t,\nu}\s M +\s N+\s A\right) u
  \\ & = (1+\eps\partial_{t,\nu})^{-1}\partial_{t,\nu}\s M u  +(1+\eps\partial_{t,\nu})^{-1}\s N u+(1+\eps\partial_{t,\nu})^{-1}\s A u
  \\ & = [(1+\eps\partial_{t,\nu})^{-1},\partial_{t,\nu}\s M] u  +[(1+\eps\partial_{t,\nu})^{-1},\s N] u+[(1+\eps\partial_{t,\nu})^{-1},\s A] u
  \\ &\quad + \partial_{t,\nu}\s M (1+\eps\partial_{t,\nu})^{-1} u + \s N(1+\eps\partial_{t,\nu})^{-1}u + \s A(1+\eps\partial_{t,\nu})^{-1}u
  \\ & =- \eps\partial_{t,\nu}(1+\eps\partial_{t,\nu})^{-1}\s M'\left(1+\eps\partial_{t,\nu}\right)^{-1}u  +[(1+\eps\partial_{t,\nu})^{-1},\s N] u+0 
  \\ &\quad + \s B (1+\eps\partial_{t,\nu})^{-1} u.
 \end{align*}
 We conclude with applying Remark \ref{r:stron_op_con} (b) to
 \begin{align*}
  A& = \s B,
  \\ T&= (1+\eps\partial_{t,\nu})^{-1}, \text{ and }
  \\ T'&=-\eps\partial_{t,\nu}(1+\eps\partial_{t,\nu})^{-1}\s M'\left(1+\eps\partial_{t,\nu}\right)^{-1}  +[(1+\eps\partial_{t,\nu})^{-1},\s N].
 \end{align*}
 Hence,  \eqref{eq:doad} follows, which, in particular, implies that $(1+\eps\partial_{t,\nu})^{-1}u\in \dom(\overline{\s B})$ for all $\eps>0$ and $u\in \dom(\overline{\s B})$.
\end{Proof}

\begin{Remark}\label{r:epinf} We note that in the same manner, one gets for $u\in \dom(\overline{\s B})$ also $\partial_{t,\nu}^{-1}u\in \dom(\overline{\s B})$ and 
 \begin{equation}\label{eq:doad2}
  \partial_{t,\nu}^{-1}\overline{\s B} u 
   = \overline{\s B}\partial_{t,\nu}^{-1}u - \s M'\partial_{t,\nu}^{-1}u 
  + \left[ \partial_{t,\nu}^{-1},\s N\right] u. 
 \end{equation}
\end{Remark}

\begin{Lemma}\label{l:appmation} Assume Hypothesis \ref{h:gS} to hold. Then for all $\eps>0$ and $u\in \dom(\overline{\s B})$, we have $(1+\eps\partial_{t,\nu})^{-1}u\in \dom(\partial_{t,\nu})\cap \dom(\s A)$.  For all $\eps>0$ there exists $C\geq 0$ such that 
\[
   \|\s A(1+\eps\partial_{t,\nu})^{-1}u\|\leq C\left(\|\overline{\s B}u\|+\|u\|\right)\quad (u\in \dom(\overline{\s B}))
\]
Furthermore, 
\[
 \overline{\s B}\left(1+\eps\partial_{t,\nu}\right)^{-1}u \to \overline{\s B}u\text{ in }L_\nu^2(\mathbb{R};\s X)
\]
  as $\eps\to 0+$ for all $u\in \dom(\overline{\s B})$. 
\end{Lemma}
\begin{Proof}
 Let $u\in \dom(\s B)=\dom(\partial_{t,\nu}\s M)\cap \dom(\s A)$, $\eps>0$, and $u_\eps\coloneqq (1+\eps\partial_{t,\nu})^{-1}u$. By Proposition \ref{p:com_resA}, $u_\eps\in \dom(\s A)$ and, since $u_\eps\in \dom(\partial_{t,\nu})$, we also have $u_\eps\in \dom(\partial_{t,\nu}\s M)$, by Hypothesis \ref{hyp:mat_law}. Then we compute with the help of Lemma \ref{l:com_B}
 \begin{align*}
    \s Au_\eps & = \s Au_\eps + \s Nu_\eps+\partial_{t,\nu}\s Mu_\eps-(\s Nu_\eps+\partial_{t,\nu}\s Mu_\eps)
    \\ & = \s B u_\eps -\s N u_\eps- \s M\partial_{t,\nu} u_\eps - \s M'u_\eps
    \\ & = (1+\eps\partial_{t,\nu})^{-1}\overline{\s B}u+ (\eps\partial_{t,\nu}(1+\eps\partial_{t,\nu})^{-1}\s M'u_\eps- \left[ \left(1+\eps\partial_{t,\nu}\right)^{-1},\s N\right] u)
    \\&\quad -\s N u_\eps- \s M\partial_{t,\nu} u_\eps - \s M'u_\eps.
 \end{align*}
  Hence, using $\|\eps\partial_{t,\nu}(1+\eps\partial_{t,\nu})^{-1}\|\leq 1$ (by Remark \ref{r:spe_thm_dtnu}), $\|(1+\eps\partial_{t,\nu})^{-1}\|\leq 1$ as well as the boundedness of $\s M$, $\s M'$ and $\s N$, we infer the desired estimate for $u\in \dom(\s B)$. Since, by definition, $\dom(\s B)$ is dense in $\dom(\overline{\s B})$ with respect to the graph norm of $\overline{B}$, we obtain the estimate also for $u\in \dom(\overline{\s B})$.
  
  For the proof of the convergence result, it suffices to recall equation \eqref{eq:doad} and Lemma \ref{le:commus_of_a} as well as that $(1+\eps\partial_{t,\nu})\stackrel{\tau_\textnormal{s}}{\to} 1$ as $\eps\to 0$, by Remark \ref{r:eps_st}.
\end{Proof}

\begin{Remark}\label{r:prestA} A more detailed look at the computations in the latter proof reveals that we have the more precise estimate
 \[
   \|\s Au_\eps \|\leq \|\overline{\s B}u\|+(2\|\s M'\|+3\|\s N\|+ \frac{1}{\eps}\|\s M\|)\|u\|.
 \]
\end{Remark}

We recall that we want to apply Proposition \ref{p:posd} for proving Theorem \ref{thm:genSolth}. For this, we need to compute the adjoint of $\s B$. We note that $\dom(\partial_{t,\nu})\cap \dom(\s A)$ is dense in $L_\nu^2(\mathbb{R};\s X)$ by Corollary \ref{c:tdcapAdense}. So, from $\dom(\partial_{t,\nu})\subseteq \dom(\partial_{t,\nu}\s M)$ by Hypothesis \ref{hyp:mat_law}, we get $\dom(\s B)=\dom(\partial_{t,\nu}\s M)\cap \dom(\s A)$ is dense in $L_\nu^2(\R;\s X)$, which yields that $\s B^*$ is a well-defined linear operator.

\begin{Lemma}\label{l:Bstar_res} Assume Hypothesis \ref{h:gS} to be satisfied. Let $f\in \dom(\s B^*)$. Then, for all $\eps>0$, $f_\eps\coloneqq (1+\eps\partial_{t,\nu}^*)^{-1}f\in \dom(\s B^*)\cap \dom(\s A^*)$,  
\begin{align}
\label{eq:Bstar1}
  \s B^*f_\eps &= \left(1+\eps\partial_{t,\nu}^*\right)^{-1}\s B^*f +[\left(1+\eps\partial_{t,\nu}\right)^{-1},\partial_{t,\nu} \s M]^*f+[(1+\eps\partial_{t,\nu})^{-1},\s N]^*f,
  \\ &=\left(\left(\partial_{t,\nu}\s M+\s N\right)^*+\s A^*\right)f_\eps. \label{eq:Bstar2}
\end{align} 
\end{Lemma}
\begin{Proof}
  With the help of Lemma \ref{l:com_B}, for $u\in \dom(\s B)$ we have
\begin{align*}
  &\langle \s B u,(1+\eps\partial_{t,\nu}^*)^{-1}f\rangle \\
  &= \langle \s B (1+\eps\partial_{t,\nu})^{-1}u,f\rangle + \langle [\left(1+\eps\partial_{t,\nu}\right)^{-1},\partial_{t,\nu} \s M]u,f\rangle  + \langle [(1+\eps\partial_{t,\nu})^{-1},\s N]u,f\rangle\\
  &=  \langle  u,(1+\eps\partial_{t,\nu}^*)^{-1}\s B^*f\rangle + \langle u,[\left(1+\eps\partial_{t,\nu}\right)^{-1},\partial_{t,\nu} \s M]^*f\rangle   + \langle u,[(1+\eps\partial_{t,\nu})^{-1},\s N]^*f\rangle,
\end{align*}
proving that $f_\eps \in \dom(\s B^*)$ as well as \eqref{eq:Bstar1}. For \eqref{eq:Bstar2}, we deduce with the help of the boundedness of both $\s M$ and $\s N$ that $\dom(\partial_{t,\nu}^*)\subseteq \dom((\partial_{t,\nu}\s M+\s N)^*)$. So, $f_\eps\in \dom((\partial_{t,\nu}\s M+\s N)^*)$.  Next, for $u\in \dom(\partial_{t,\nu})\cap \dom(\s A)$, we compute
\[
   \langle u,\s B^* f_\eps\rangle  =    \langle \left(\partial_{t,\nu}\s M + \s N + \s A\right) u,f_\eps\rangle = \langle \s A u,f_\eps\rangle +\langle u,\left(\partial_{t,\nu}\s M + \s N\right)^*f_\eps\rangle.
\]
Since $\dom(\partial_{t,\nu})\cap \dom(\s A)$ is dense in $\dom(\s A)$ with respect to the graph norm of $\s A$ (see Corollary \ref{c:tdcapAdense}), we infer \eqref{eq:Bstar2}.
\end{Proof}

\begin{Corollary}\label{c:Bstar} Assume Hypothesis \ref{h:gS} to be satisfied. Then, the operator 
\[
   \s C \colon \dom(\partial_{t,\nu}^*)\cap \dom(\s A^*) \subseteq L_\nu^2(\s X)\to L_\nu^2(\s X), u\mapsto ((\partial_{t,\nu}\s M+\s N)^*+\s A^*)u
\]
is closable, and 
\[
  \s B^*= \overline{\s C}.
\] 
\end{Corollary}
\begin{Proof}
 The closability follows from $\s C\subseteq \s B^*$. Next, let $f\in \dom(\s B^*)$. By Lemma \ref{l:Bstar_res}, $f_\eps\coloneqq (1+\eps\partial_{t,\nu}^*)^{-1}f\in \dom(\s C)$ and, by \eqref{eq:Bstar2}, $\s B^*f_\eps = \s Cf_\eps$, $\eps>0$. By Lemma \ref{le:commus_of_a}, we have
 \[
    \overline{[(1+\eps\partial_{t,\nu})^{-1},\partial_{t,\nu}\s M]}, [(1+\eps\partial_{t,\nu})^{-1},\s N]\to 0\quad (\eps\to0)
 \] in the strong operator topology. Hence, 
 \[
    [(1+\eps\partial_{t,\nu})^{-1},\partial_{t,\nu}\s M]^*, [(1+\eps\partial_{t,\nu})^{-1},\s N]^*\to 0\quad (\eps\to0)
 \]
 in the weak operator topology. Thus, by Lemma \ref{l:Bstar_res}, $(\s Cf_\eps)_\eps$ weakly converges to $\s B^*f$ as $\eps\to 0$. Moreover, since $f_\eps\to f$ weakly as $\eps\to 0$, we infer $f\in \dom(\overline{\s C})$ and $\overline{\s C}f=\s B^*f$, which yields the assertion.
\end{Proof}

We come to the proof of Theorem \ref{thm:genSolth}:

\begin{Proof}[of Theorem \ref{thm:genSolth}]
  We apply Proposition \ref{p:posd} to the operator $\s B$ and the space $L_\nu^2(\mathbb{R};\s X)$ as underlying Hilbert space. For this, we note that \eqref{eq:p_posd} is guaranteed by \eqref{ineq:pos}. Next, since $\s C$ from Corollary \ref{c:Bstar} satisfies the analogous positivity estimate \eqref{ineq:adjoint}, by the closability of $\s C$, the inequality is valid for $\s C$ replaced by $\overline{\s C}$. By Corollary \ref{c:Bstar}, however, $\overline{\s C}=\s B^*$, which yields \eqref{eq:p_posd2}. Thus, the assertion indeed follows from Proposition \ref{p:posd}.
\end{Proof}

\renewcommand{\baselinestretch}{0.65}\normalsize\mysection{Partial Differential Equations -- Causality and the Independence of $\nu$}{Partial Differential Equations -- Causality and the Independence of $\nu$}{a solution theory of partial differential equations $\cdot$ Hypothesis \ref{h:gSe} $\cdot$ Theorem \ref{t:gSe}}\label{sec:ind_of_nu}

\renewcommand{\baselinestretch}{1}\normalsize
This section is devoted to a proof of an adapted version of Theorem \ref{thm:genSolth} including causality. Moreover, we will prove the independence of the solution operator of the parameter $\nu$ in the solution theory. So, as in the case of ordinary differential equations, the aim is to show that the solution operator $\s S=\s B^{-1}$ associated with \eqref{eq:pde0} is standard evolutionary. Beforehand, we will state a sufficient condition warranting both the inequalities \eqref{ineq:adjoint} and \eqref{ineq:pos}.

\begin{Hypothesis}\label{h:wgS} Let $\s X$ be a Hilbert space, $\nu>0$, $\s A\colon \dom(\s A)\subseteq L_\nu^2(\s X)\to L_\nu^2(\s X)$ linear, densely defined, closable, $\s M,\s N,\s M'\in L(L_\nu^2(\s X))$, and let $\s D\subseteq \dom(\partial_{t,\nu})$ be a core for $\partial_{t,\nu}$.
We assume 
\begin{equation}\label{eq:pA}
   \Re\langle \phi,\s A\phi\rangle, \Re \langle\psi,\s A^*\psi\rangle\geq 0 \quad (\phi\in \dom(\s A),\psi\in \dom(\s A^*)),
\end{equation}
Further, assume
\[
    \s M\partial_{t,\nu}|_{\s D} \subseteq \partial_{t,\nu}\s M -\s M'\quad\text{ and }\partial_{t,\nu}^{-1}\s A\subseteq \overline{\s A}\partial_{t,\nu}^{-1}
\]
and 
\begin{equation}\label{eq:dmn}
   \Re\langle (\partial_{t,\nu}\s M +\s N)\phi,\phi\rangle \geq c\langle\phi,\phi\rangle
\end{equation}
for all $\phi\in \s D$ and some $c>0$.
\end{Hypothesis}

\begin{Proposition}\label{p:whyp} Hypothesis \ref{h:wgS} implies Hypothesis \ref{h:gS} (with $\overline{\s A}$ instead of $\s A$).
\end{Proposition}
\begin{Proof} It follows from Remark \ref{r:stron_op_con}(b) that the Hypotheses \ref{hyp:A} and \ref{hyp:mat_law} are valid. Hence, it remains to prove \eqref{ineq:adjoint} and \eqref{ineq:pos}. For this, we observe 
\[
  \Re\langle (\partial_{t,\nu}\s M +\s N)\phi,\phi\rangle=\Re\langle (\partial_{t,\nu}\s M +\s N)^*\phi,\phi\rangle 
\]
for all $\phi\in \dom(\partial_{t,\nu})=\dom(\partial_{t,\nu}^*)$. Thus, we are left with showing that \eqref{eq:dmn} carries over to all $\phi\in \dom(\partial_{t,\nu}\s M)$. This, however, follows from $\partial_{t,\nu}\s M = \overline{\partial_{t,\nu}\s M|_{D}}$, which we show in the next proposition.
\end{Proof}

\begin{Proposition}\label{p:cdddm} Assume Hypothesis \ref{hyp:mat_law}. Let $\s D\subseteq \dom(\partial_{t,\nu})$ be a core for $\partial_{t,\nu}$. Then 
\[
   \partial_{t,\nu}\s M = \overline{\partial_{t,\nu}\s M|_{\s D}}.
\] 
\end{Proposition}
\begin{Proof}
 Endowed with the respective graph norms, we observe that the canonical embedding $\dom(\partial_{t,\nu})\hookrightarrow \dom(\partial_{t,\nu}\s M)$ is continuous by Hypothesis \ref{hyp:mat_law}. Hence, it suffices to show that $\dom(\partial_{t,\nu})$ is a core for $\partial_{t,\nu}\s M$, that is, the mentioned embedding is dense.  So, take $u\in \dom(\partial_{t,\nu}\s M)$. Then, we compute for $\eps>0$
 \[
    \partial_{t,\nu}\s M(1+\eps\partial_{t,\nu})^{-1}u = [\partial_{t,\nu}\s M,(1+\eps\partial_{t,\nu})^{-1}]u+(1+\eps\partial_{t,\nu})^{-1}\partial_{t,\nu}\s Mu.
 \]
 Hence, letting $\eps\to 0$ and recalling Lemma \ref{le:commus_of_a}, we read off that 
 \[
    \dom(\partial_{t,\nu})\ni (1+\eps\partial_{t,\nu})^{-1}u \to u\text{ and }\partial_{t,\nu}\s M(1+\eps\partial_{t,\nu})^{-1}u\to \partial_{t,\nu}\s M u.  
 \] 
\end{Proof}

One of the reasons of having introduced Hypothesis \ref{h:wgS} is as follows. Theorem \ref{thm:genSolth} has a natural analogue in the context of evolutionary mappings with a variant of Hypothesis \ref{h:wgS} as the set of assumptions as we shall see next. This version of Hypothesis \ref{h:wgS} reads as follows.

\begin{Hypothesis}\label{h:gSe} Let $\s X$ Hilbert space, $\nu>0$, $\s M,\s M',\s N\in L_{\textnormal{sev},\nu}(\s X)$, $\s A\in C_{\textnormal{ev},\nu}(\s X)$, $c>0$. Assume for all $\mu\geq \nu$:
\[
   \s M\partial_{t,\mu}\subseteq \partial_{t,\mu}\s M^{\mu}-(\s M')^\mu,\quad \partial_{t,\mu}^{-1}\s A\subseteq \s A^\mu\partial_{t,\mu}^{-1},
\]
\[
  \Re\langle Q_0 \s A\phi,\phi\rangle_{L_\mu^2(\s X)}\geq 0, \Re\langle(\s A^\mu)^*\psi,\psi\rangle_{L_\mu^2(\s X)}\geq 0 \quad (\phi\in \dom(\s A),\psi\in \dom((\s A^\mu)^*)),
\]
as well as
\[
   \Re \langle Q_t (\partial_{t,\mu}\s M+\s N)\phi,\phi\rangle_{L_\mu^2(\s X)}\geq c\langle Q_t\phi,\phi\rangle\quad (\phi\in\s D),
\]
where $Q_t$ denotes multiplication by $\1_{(-\infty,t)}$ and $\s D\subseteq \bigcap_{\eta\geq\nu} \dom(\partial_{t,\eta})$ is a core for $\partial_{t,\mu}$.
Define $\s B\coloneqq \check{\partial}_{t,\nu}\s M+\s N+\check{\s A}$, where $\check{\s A}\coloneqq \bigcap_{\eta\geq\nu} \s A^\eta$ and $\check{\partial}_{t,\nu}\coloneqq \bigcap_{\eta\geq\nu}\partial_{t,\eta}$, see Example \ref{ex:cev}.
\end{Hypothesis}

\begin{Remark}\label{r:acheck} The reason of introducing the operator $\check{\s A}$ in the latter hypothesis is that being closable evolutionary only, the operator $\s A$ might be endowed with a domain that intersects only trivially with the domain of $\check{\partial}_{t,\nu}\s M$, see also the discussion after Example \ref{ex:evo_map}. 
\end{Remark}

\begin{Theorem}\label{t:gSe} Assume Hypothesis \ref{h:gSe}. Then $\s S\coloneqq \s B^{-1}$ is evolutionary at $\nu$ and causal, $\sup_{\mu\geq\nu}\|\s S^\mu\|_{L(L_\mu^2)}\leq 1/c$. In particular, the solution operator does not depend on the exponential weight, that is, for all $\mu_1\geq\mu_2\geq\nu$ the operators $\s S^{\mu_1}$ and $\s S^{\mu_2}$ coincide on the intersection of the respective domains.  
\end{Theorem}
 
With the results of Chapter \ref{ch:EMC} in mind, apart from the norm estimate, the assertion in Theorem \ref{t:gSe} may be expressed as $\s S\in L_{\textnormal{sev},\nu}(\s X)$. The proof of Theorem \ref{t:gSe} relies on the Theorems \ref{t:inv_ce} and \ref{thm:genSolth}. We need two preparatory results.

\begin{Proposition}\label{p:acau} Let $\s X$ Hilbert space, $\nu>0$, denote $Q_t$ as multiplication by $\1_{(-\infty,t)}$, $\s A\in C_{\textnormal{ev},\nu}(\s X)$. Assume for all $\mu\geq \nu$
\[
   \Re\langle Q_0 \s A\phi,\phi\rangle_{L_\mu^2(\s X)}\geq 0\quad (\phi\in \dom(\s A))	\quad\text{and}\quad
  \partial_{t,\mu}^{-1}\s A\subseteq \s A^\mu\partial_{t,\mu}^{-1}.
 \]
Then
\[
   \Re\langle Q_t \s A^\mu\phi,\phi\rangle_{L_\mu^2(\s X)}\geq 0\quad(\phi\in \dom(\s A^\mu), t\in \mathbb{R}).
\]
\end{Proposition}
\begin{Proof}
 By Proposition \ref{p:Atrinv} in combination with Remark \ref{r:stron_op_con}(b), we get for all $\mu\geq\nu$
 \[
    \tau_t\s A^\mu\subseteq \s A^\mu \tau_t \quad (t\in \mathbb{R}),
 \]
 where we recall $\tau_t f = f(\cdot+t)$.
 Hence, for all $t\in \mathbb{R}$, $\phi\in \dom(\s A^\mu)$, we have $\tau_t\phi\in \dom(\s A^\mu)$ and
 \begin{align*}
   0 &\leq \Re \langle Q_0 \s A^\mu \tau_t\phi,\tau_t\phi\rangle_{L_\mu^2} 
   = \Re \langle Q_0 \tau_t\s A^\mu \phi,\tau_t\phi\rangle_{L_\mu^2}
   \\ &= \Re \langle \tau_tQ_t \s A^\mu \phi,\tau_t\phi\rangle_{L_\mu^2}
   = \Re \langle Q_t \s A^\mu \phi,\tau_t^*\tau_t\phi\rangle_{L_\mu^2}
   \\ &=e^{2t\mu} \Re \langle Q_t \s A^\mu \phi,\tau_{-t}\tau_t\phi\rangle_{L_\mu^2}
   =e^{2t\mu} \Re \langle Q_t \s A^\mu \phi,\phi\rangle_{L_\mu^2},
 \end{align*}
where we used $\tau_{t}^*=\tau_{-t}e^{2t\mu}\in L(L_\mu^2(\s X))$.
\end{Proof}

\begin{Lemma}\label{l:Bdd} Assume Hypothesis \ref{h:gSe}. Then the operator $\s B$ is densely defined.  
\end{Lemma}
\begin{Proof}
  For all $\eta\geq \nu$, $\dom(\s A)$ is dense in $L_\eta^2(\mathbb{R};\s X)$ by hypothesis. Next, the operator $(1+\eps\partial_{t,\nu})^{-1}$ is translation-invariant and causal by Remark \ref{r:rtc}. So, we infer that $(1+\eps\partial_{t,\nu})^{-1}$ leaves $L_\eta^2(\mathbb{R};\s X)$ invariant for all $\eta\geq \nu$, by Remark \ref{r:SW_OT}. By Remark \ref{r:stron_op_con}(b), we have $\partial_{t,\eta}^{-1}\s A^\eta\subseteq \s A^\eta \partial_{t,\eta}^{-1}$. Thus, by Proposition \ref{p:com_resA}, we get for $\eta\geq\nu$ and $\eps>0$
  \[
     (1+\eps\partial_{t,\nu})^{-1}\s A \subseteq (1+\eps\partial_{t,\eta})^{-1}\s A^\eta \subseteq \s A^\eta(1+\eps\partial_{t,\eta})^{-1}.
  \]
Hence, 
 \[
    (1+\eps\partial_{t,\nu})^{-1}\s A \subseteq \bigcap_{\eta\geq \nu} (\s A^\eta(1+\eps\partial_{t,\eta})^{-1}).
 \]
 We read off that if $\phi\in \dom(\s A)$, then for $\eps>0$, we get \[(1+\eps\partial_{t,\nu})^{-1}\phi \in \bigcap_{\eta\geq\nu} \dom(\partial_{t,\eta})\cap \dom(\s A^\eta).\] Thus,
\[
   \lin \bigcup_{\eps>0} (1+\eps\partial_{t,\nu})^{-1}[\dom(\s A)] \subseteq \left(\bigcap_{\eta\geq \nu} \dom(\partial_{t,\eta})\cap \dom(\s A^\eta)\right) \subseteq \dom(\s B).
\]
Using Remark \ref{r:eps_st} and the density of $\dom(\s A)$, we realize that the left-hand side is dense in $L_\eta^2(\s X)$, $\eta\geq\nu$, hence, so is the right-hand side.
\end{Proof}

\begin{Proof}[of Theorem \ref{t:gSe}] We will apply Theorem \ref{t:inv_ce} to $\s B$. For this, we establish the positivity estimate required in Theorem \ref{t:inv_ce} first. Let $\mu\geq\nu$. By Proposition \ref{p:cdddm}, $\s D$ is a core for $\partial_{t,\mu}\s M^\mu$. Hence, for all $t\in \mathbb{R}$, we have
\[
  \Re \langle Q_t(\partial_{t,\mu}\s M^\mu +\s N^\mu)\phi,\phi\rangle\geq c\langle Q_t\phi,\phi\rangle\quad(\phi\in \dom(\partial_{t,\mu}\s M^\mu))
\]
Moreover, by Proposition \ref{p:acau}, we get
\[
   \Re\langle Q_t \s A^\mu\phi,\phi\rangle\geq 0.
\]
Thus, 
\begin{equation}\label{eq:pe}
  \Re \langle Q_t\s B \phi,\phi\rangle_{L_\mu^2} =\Re \langle Q_t(\partial_{t,\mu}\s M^\mu +\s N^\mu+\s A^\mu)\phi,\phi\rangle\geq c\langle Q_t\phi,\phi\rangle
\end{equation}
for all $\phi\in \dom(\s B)$, $t\in \R$. Therefore, the estimate required in Theorem \ref{t:inv_ce} is shown. 

Next, we show that $\s B\in C_{\textnormal{ev},\nu}(\s X)$. For this, we realize that $\s B$ is densely defined by Lemma \ref{l:Bdd}, so only closability is the issue here. But, if we let $t\to\infty$ in \eqref{eq:pe}, we get $\s B$ is closable by Proposition \ref{p:posd}.

To conclude, we are left with showing that $\ran(\s B)$ is dense in $L_\mu^2(\s X)$ for all $\mu\geq \nu$. This, however, follows from Proposition \ref{p:whyp} and Theorem \ref{thm:genSolth}. Indeed, for $\mu\geq\nu$, the estimates
\[
   \Re \langle (\partial_{t,\nu}\s M+\s N)\phi,\phi\rangle_{L_\mu^2(\s X)}\geq c\langle\phi,\phi\rangle_{L_\mu^2} \quad(\phi\in \s D)
\]
and 
\[ 
 \Re\langle \phi,\s A\phi\rangle_{L_\mu^2}\geq 0\quad (\phi\in \dom(\s A))
\]
follow either from Hypothesis \ref{h:gSe} or Proposition \ref{p:acau} by letting $t\to\infty$ in order that $Q_t\to 1$ strongly.
\end{Proof}

\section{Comments}

As it has been mentioned already, for ordinary differential equations, there is a wider class of problems, that may be studied in this $L^2$-type setting: Delay differential equations covering a class of functional differential equations, equations of neutral type, integro-differential equations or differential-algebraic equations. But note that, $\s M$ and $\s N$ in Theorem \ref{t:wp_ode_evo} are linear operators in space-time. Assuming time translation-invariance for $\s M$ and $\s N$ these are operators of convolution type. Hence, the set of equations treated in Theorem \ref{t:wp_ode_evo} may already be summarized by ``integro-differential-algebraic''. 

Theorem \ref{t:gSe} has its roots in \cite[Solution Theory]{PicPhy}. In \cite{PicPhy}, the problem of solving 
\[
   \s B u = (\partial_{t,\nu} M(\partial_{t,\nu}^{-1}) + \s A)u = f
\] for some bounded and analytic function $M$ of $\partial_{t,\nu}^{-1}$ with values in $L(\s X)$ and a skew-selfadjoint operator $A$ in $\s X$ has been addressed, where  $\s A$ is the lift of $A$ to $L_\nu^2(\R;\s X)$, as in Example \ref{ex:cev}(b). As noted in \cite[Section 3.1]{Waurick2014MMAS_Non}, Theorem \ref{t:gSe} covers the class discussed in \cite{PicPhy} by putting $\s M=M(\partial_{t,\nu}^{-1})$ and $\s N=0$. There are plenty of equations already covered by this class: A treatment of electro-seismic waves is included in \cite{McGhee2011} (one needs to involve fractional time derivatives in $M(\partial_{t,\nu}^{-1})$), a general class of fractional partial differential equations \cite{Waurick2014MMAS_Frac} (e.g.~fractional Fokker-Planck equations, or super and subdiffusion problems (\cite[Theorem 4.5, Remark 4.6]{Waurick2013AA_Hom})). For more examples, we refer to the list given in the introduction. We further remark here that the assumptions on $\s A$ in Theorem \ref{t:gSe} allow for spatial operators with certain differential equations as boundary conditions. A prominent example are boundary conditions of impedance type. We refer to \cite{Picard2012} and \cite{PSSM_GD}: In both these references the assumptions on $\s A$ being asked for in Hypothesis \ref{h:gSe} have been established for impedance type boundary conditions for the wave equation and for boundary conditions of Leontovitch type in the area of Maxwell's equation, respectively. A structural point of view treating possible boundary condition in a slightly more abstract setting can be found in \cite{Trostorff2014}.

We will treat some (standard) problems of mathematical physics in the Sections \ref{s:eca}, \ref{s:cdcnt} and \ref{s:homaw}.

The method of proof of the \cite[Solution Theory]{PicPhy} relied on the spectral representation of $\partial_{t,\nu}$. Later on, still using the explicit spectral theorem for $\partial_{t,\nu}$, in \cite{Picard2012}, the method has been generalized to include operators $\s A$ that commute with the inverse of the time derivative, as in Hypothesis \ref{hyp:A}. The latter comes in handy, when discussing the aforementioned problems with impedance type boundary conditions. For non-autonomous problems, techniques as the Fourier--Laplace transformation have limited applicability. Hence, the regularization technique presented in this exposition might be the method of choice. This technique was used in \cite{Waurick2013JEE_Non}. As shown in \cite{Waurick2014MMAS_Non}, this technique applies to a broader class. The strategy of proof of Theorem \ref{t:gSe} developed here has its roots in the methods from \cite{Waurick2013JEE_Non,Waurick2014MMAS_Non}. We note here that problems with changing type, that is, problems that are hyperbolic, parabolic and elliptic on different space-time regions may be addressed as well, see \cite[pp 765]{Waurick2013JEE_Non}, \cite[Remark 6.2]{Waurick2014IMAJMCI_ComprehContr}, or \cite{Waurick2016_StH}.

Another possible way to deduce a solution theory for problems of the type \eqref{eq:pde0} is the usage of the theory of maximal monotone relations. In this line of reasoning one treats \eqref{eq:pde0} as a sum of the two maximal monotone relations $\partial_{t,\nu}\s M$ and $\s A$. This strategy has been successfully applied to time translation-invariant $\s M$ and non-linear $\s A$ in \cite{Trostorff2011,Trostorff2012,Trostorff2013,Trostorff2015} eventually yielding a solution theory for partial differential \emph{inclusions}. For time dependent operators $\s M$, an adapted form can be found in \cite{Trostorff2014a}. For a possible extension to a Banach space setting, we refer to \cite{Wehowski2015}.

Going back to linear problems, one may address the minimal $\nu\in \mathbb{R}$ the solution operator $\s S$ in Theorem \ref{t:gSe} is evolutionary at. In fact, if this $\nu$ was negative, it is possible to address the question of exponential stability in this framework as well, see \cite{Trostorff2013a,Trostorff2015a}. 

Some remarks on the comparison to other strategies of finding solutions to this type of partial differential equations are in order. The overall strategy may be thought of as a particular instance of discussing sums of unbounded operators similar to the seminal paper \cite{Prato1975}. We refer to \cite[Section 5]{Prato1975}, where a hyperbolic type case is treated. The strategies developed in \cite{Prato1975} are thorough and deep and do also cover the Banach space case. Note that, however, restricting ourselves to a Hilbert space setting and employing the particular role of the time derivative yields a particularly accessible way of discussing continuous invertibility of evolutionary partial differential equations. Indeed, the method solely relies on emphasizing the special role of the time derivative and the well-known observation that strict positive definiteness eventually leads to continuous invertibility as demonstrated in Proposition \ref{p:posd}.

Putting $\s M=1$ and $\s N=0$ in \eqref{eq:pde0} with $A$ being quasi-m-accretive, we realize that the corresponding equation ($\s A$ is given as in Example \ref{ex:cev})
\[
  (\partial_{t,\nu} + \s A)u=f
\]
may well be treatable with $C_0$-semi-group theory. We refer to \cite{EngNag,Arendt2011,Pazy1983} for a thorough treatment of semi-groups with regards to the solution theory of differential equations. Being genuinely developed in a Banach space setting, $C_0$-semi-groups may give more insight on particular properties of the corresponding solution, that is, for instance, boundedness, $p$-integrability, positivity, or stability. We refer to \cite{EngNag,Arendt2011,Pazy1983} again for an account on that. Leading to a continuous-in-time-solution, semi-group theory may also be viewed as a regularity theory of evolutionary equations as in \eqref{eq:pde0}. 

The idea of introducing semi-groups or its second order analogue cosine families as a solution concept for partial differential equations has its roots in the finite-dimensional case: Given $A\in \mathbb{R}^{d\times d}$ the fundamental solution to $\partial_{t} u=-A u$ or $\partial_{t}^2 u =- Au$ may respectively be written as $t\mapsto e^{-tA}$  or $t\mapsto \cos(-tA)$. Keeping this idea in mind, non-autonomous equations are solved by finding a corresponding solution family associated to $\partial_t u=-A(t)u$, say. Solution families generalize the concept of the fundamental matrix for non-autonomous ordinary differential equations to the infinite-dimensional setting \cite{Kato1953,Pazy1983,Tanabe1979,Sohr1978}. In particular, when $A(t)$ is an unbounded operator for every $t$, one needs to cope with varying domains of $A(t)$. Again, we view the concept of evolution families as a certain regularity theory. As a particular example, we mention the non-autonomous heat equation, formally given by
\[
   \partial_t u(t,x) - \dive a(t,x)\grad u(t,x)= f(t,x).
\]
As it will be demonstrated in Section \ref{s:cdcnt}, we reformulate this equation into a problem of first order. This reformulation enables us to apply the solution theory given in Theorem \ref{t:gSe} without the need of coping with subtleties of possibly changing domains. We refer to \cite{Auscher2015} for a deep and thorough treatment of these kind of problems in a Banach space type setting.

Semi-groups, cosine families and evolution families basically provide solutions for initial value problems. Non-homogeneous problems are then solved using some sort of variation of constants type formulas. The solution concept developed here treats the non-homogeneous problem first. We sketch an adapted treatment of initial value problems as follows. We refer to \cite[Section 6.2.5]{Picard} for more details. A corresponding theory for initial value problems can be obtained by formally putting
\begin{equation}\label{eq:ivp0}
 \partial_{t,\nu} u + \s A u =\delta_0 u_0,
\end{equation}
where $\delta_0$ denotes the Dirac-$\delta$-distribution at $0$ and $u_0\in \dom(A)\subseteq \s X$ with $\s A$ being the lift of a quasi-m-accretive $A$ to $L_\nu^2(\s X)$, $\nu>0$. Recalling that $\partial_{t,\nu}\1_{[0,\infty)}=\delta_0$, we obtain from \eqref{eq:ivp0} the equation
\begin{equation}
\label{eq:ivp1}
 \partial_{t,\nu} (u-\1_{[0,\infty)}u_0) + \s A (u-\1_{[0,\infty)}u_0) = -\1_{[0,\infty)} Au_0.
\end{equation}
Hence, solving equation \eqref{eq:ivp1} for $v= u-\1_{[0,\infty)}u_0$, we obtain a solution to \eqref{eq:ivp0} by putting $u=v+\1_{[0,\infty)}u_0$ .  Indeed, causality of $(\partial_{t,\nu}+\s A)^{-1}$ yields that $v$ is supported on $[0,\infty)$ only. If, in addition, $v$ is weakly differentiable, then $v$ is continuous on $\R$, by Sobolev's embedding theorem (see \cite[Lemma 5.2]{Kalauch}). Hence, 
\[
   0= v(0) = \lim_{t\to 0-} v(t) =\lim_{t\to 0+} v(t) = \lim_{t\to 0+} (u(t)-\1_{[0,\infty)}(t)u_0)=\lim_{t\to 0+}u(t)- u_0.
\]
Thus, $u$ attains its initial value and satisfies \eqref{eq:ivp0} on $(0,\infty)$. For a more detailed treatment of initial value problems, we also refer to \cite[Theorem 5.4]{Kalauch} and \cite[Example 2.18]{Waurick2014OAM_DelayLp} for the ordinary differential equations case.

We summarize that semi-groups, cosine families and evolution families are the fundamental solutions or abstract Green's functions to certain (partial) differential equations. The solution theory developed here complements these treatments to partial differential equations as the existence of (a sufficiently regular) fundamental solution is a priori not needed. However, in the present approach -- roughly speaking -- $L^2$ right-hand sides are mapped to $L^2$-solutions only. The strategy of solving partial differential equations by means of semi-groups, cosine families or evolution families leads to more regular solutions. So, take a partial differential equation, where it is possible to apply both the approach discussed in the previous section as well as one of the three approaches semi-groups, cosine or evolution families. Then either of the latter three solution strategy may be viewed as a regularity theory for the approach advanced in the present exposition.

\cleardoublestandardpage
\chapter[Convergence of Evolutionary Mappings]{Convergence of Evolutionary Mappings and an Application to Ordinary Differential Equations}\label{Chapter:ODE}

In this chapter we will address convergence issues of evolutionary mappings. We have occasion to discuss the norm, the strong and the weak operator topology in the light of evolutionary mappings. As a first application of these concepts, we will consider abstract ordinary differential equations as discussed in Section \ref{s:wpode}.

\renewcommand{\baselinestretch}{0.65}\normalsize\mysection{Topologies on Evolutionary Mappings}{Topologies on Evolutionary Mappings}{norm topology $\cdot$ strong operator topology $\cdot$ weak operator topology $\cdot$ characterization of the topologies on bounded sets $\cdot$ compactness result for the weak operator topology $\cdot$ metrizability for bounded sets under the weak operator topology $\cdot$ continuity of multiplication and inversion under the norm and strong operator topologies}\label{s:nt}

\renewcommand{\baselinestretch}{1}\normalsize
We start out with the definition of the topologies we are interested in on standard evolutionary mappings at some $\nu$, see Definition \ref{d:scd_sev}. 

\begin{Definition}\label{d:top_sev} Let $\s X,\s Y$ be Hilbert spaces, $\nu\in \mathbb{R}$. The \emph{norm topology} $\tau_{\textnormal{n},\nu}$ on $L_{\textnormal{sev},\nu}(\s X,\s Y)$ is defined as the initial topology induced by
\[
   \s S\mapsto \s S^\mu \in L(L_\mu^2(\mathbb{R};\s X),L_\mu^2(\mathbb{R};\s Y))
\]
for all $\mu\geq \nu$. The \emph{strong operator topology} $\tau_{\textnormal{s},\nu}$ (\emph{weak operator topology} $\tau_{\textnormal{w},\nu}$) on the space $L_{\textnormal{sev},\nu}(\s X,\s Y)$ is defined as the initial topology induced by
\[
   \s S\mapsto \s S^\mu \in (L(L_\mu^2(\mathbb{R};\s X),L_\mu^2(\mathbb{R};\s Y)),\tau_{\textnormal{s}})\ \Big(\s S\mapsto \s S^\mu \in (L(L_\mu^2(\mathbb{R};\s X),L_\mu^2(\mathbb{R};\s Y)),\tau_{\textnormal{w}})\Big)
\]
for all $\mu\geq\nu$. Denote $L^{\textnormal{n}}_{\textnormal{sev},\nu}(\s X,\s Y)\coloneqq \big(L_{\textnormal{sev},\nu}(\s X,\s Y),\tau_{\textnormal{n},\nu}\big)$, $L^{\textnormal{s}}_{\textnormal{sev},\nu}(\s X,\s Y)\coloneqq \big(L_{\textnormal{sev},\nu}(\s X,\s Y),\tau_{\textnormal{s},\nu}\big)$ and $L^{\textnormal{w}}_{\textnormal{sev},\nu}(\s X,\s Y)\coloneqq \big(L_{\textnormal{sev},\nu}(\s X,\s Y),\tau_{\textnormal{w},\nu}\big)$ the corresponding topological spaces. 

Moreover, the \emph{norm}, \emph{strong operator}, and \emph{weak operator topology} on $L_{\textnormal{sev}}(\s X,\s Y)$ are defined as the final topologies induced by 
\begin{multline*}
   L^{\textnormal{n}}_{\textnormal{sev},\mu}(\s X,\s Y)\hookrightarrow L_{\textnormal{sev}}(\s X,\s Y), L^{\textnormal{s}}_{\textnormal{sev},\mu}(\s X,\s Y)\hookrightarrow L_{\textnormal{sev}}(\s X,\s Y)\\ \text{ and }L^{\textnormal{w}}_{\textnormal{sev},\mu}(\s X,\s Y)\hookrightarrow L_{\textnormal{sev}}(\s X,\s Y)
\end{multline*}
for all $\mu\in\mathbb{R}$, respectively. Denote the respective topological spaces by $L^{\textnormal{n}}_{\textnormal{sev}}(\s X,\s Y)$, $L^{\textnormal{s}}_{\textnormal{sev}}(\s X,\s Y)$ and $L^{\textnormal{w}}_{\textnormal{sev}}(\s X,\s Y)$. 
\end{Definition}

Note that as an immediate consequence of the latter definition, for $\mu\geq \nu$, we have the continuous (canonical) embeddings
\begin{align}
   L^{\textnormal{n}}_{\textnormal{sev},\nu}(\s X,\s Y) &\hookrightarrow L^{\textnormal{n}}_{\textnormal{sev},\mu}(\s X,\s Y), \notag
   \\ L^{\textnormal{s}}_{\textnormal{sev},\nu}(\s X,\s Y) &\hookrightarrow L^{\textnormal{s}}_{\textnormal{sev},\mu}(\s X,\s Y),\label{embmon}
   \\ L^{\textnormal{w}}_{\textnormal{sev},\nu}(\s X,\s Y) &\hookrightarrow L^{\textnormal{w}}_{\textnormal{sev},\mu}(\s X,\s Y)\notag
\end{align}
as well as the continuous embeddings
\begin{align*}
 L^{\textnormal{n}}_{\textnormal{sev},\nu}(\s X,\s Y)\hookrightarrow L^{\textnormal{s}}_{\textnormal{sev},\nu}(\s X,\s Y)\hookrightarrow L^{\textnormal{w}}_{\textnormal{sev},\nu}(\s X,\s Y).
\end{align*}
In particular, the embeddings 
\begin{equation}\label{eq:trivemb}
   L^{\textnormal{n}}_{\textnormal{sev}}(\s X,\s Y)\hookrightarrow L^{\textnormal{s}}_{\textnormal{sev}}(\s X,\s Y)\hookrightarrow L^{\textnormal{w}}_{\textnormal{sev}}(\s X,\s Y)
\end{equation}are continuous.

\begin{Remark}\label{r:pt} (a) We note that in Definition \ref{d:top_sev}, the final topology on $L_{\textnormal{sev}}^\textnormal{n}(\s X,\s Y)$, $L_{\textnormal{sev}}^\textnormal{s}(\s X,\s Y)$ or $L_{\textnormal{sev}}^\textnormal{w}(\s X,\s Y)$ is \emph{not} defined as the linear or locally convex final topology on $L_{\textnormal{sev}}(\s X,\s Y)$. More precisely, for $\textnormal{t}\in \{\textnormal{n},\textnormal{s},\textnormal{w}\}$ the topology $\tau^\textnormal{t}$ on $L_{\textnormal{sev}}^\textnormal{t}(\s X,\s Y)$ is given by
\[
   \tau^\textnormal{t}=\{ \mathfrak{U}\subseteq L_{\textnormal{sev}}(\s X,\s Y); f_\nu^{-1}(\mathfrak{U})\subseteq L_{\textnormal{sev},\nu}^\textnormal{t}(\s X,\s Y) \text{ open for all }\nu\in\R\},
\]
where $f_\nu \colon L_{\textnormal{sev},\nu}(\s X,\s Y)\hookrightarrow L_{\textnormal{sev}}(\s X,\s Y)$ is the canonical embedding, $\nu\in \R$. In particular, we do not view any of the spaces $L_{\textnormal{sev}}^\textnormal{n}(\s X,\s Y)$, $L_{\textnormal{sev}}^\textnormal{s}(\s X,\s Y)$ or $L_{\textnormal{sev}}^\textnormal{w}(\s X,\s Y)$ as topological vector spaces.

(b) Unless specified otherwise, for product spaces of the form $L_{\textnormal{sev}}^{\textnormal{t}}(\s X,\s Y)\times L_{\textnormal{sev}}^{\textnormal{t}}(\s X,\s Y)$ we will use the final topology induced by 
\[
   L_{\textnormal{sev},\mu}^\textnormal{t}(\s X,\s Y)\times L_{\textnormal{sev},\nu}^\textnormal{t}(\s X,\s Y) \hookrightarrow L_{\textnormal{sev}}(\s X,\s Y)\times L_{\textnormal{sev}}(\s X,\s Y)\quad (\mu,\nu\in \R)
\] instead of the product topology, $\textnormal{t} \in \{\textnormal{n}, \textnormal{s}, \textnormal{w}\}$ (see also Theorem \ref{t:cim} below). In particular, the mapping 
\[
   p\colon L_{\textnormal{sev}}^\textnormal{t}(\s X,\s Y) \times L_{\textnormal{sev}}^\textnormal{t}(\s X,\s Y) \to L_{\textnormal{sev}}^\textnormal{t}(\s X,\s Y), (\s S,\s T)\mapsto \s S+\s T
\]
is continuous for all $\textnormal{t}\in \{\textnormal{n}, \textnormal{s}, \textnormal{w}\}$.
\end{Remark}

\begin{Example}\label{ex:el} Let $\s X,\s Y$ be Hilbert spaces. Let $(S_\iota)_{\iota\in I}$ be a net in $L(\s X,\s Y)$ with the property $\sup_\iota\|S_\iota\|_{L(\s X,\s Y)}<\infty$. Let $\nu\in \mathbb{R}$ and 
\[
   \s S^{(\nu)}_\iota \colon L_\nu^2(\s X)\to L_\nu^2(\s Y), f\mapsto (t\mapsto S_\iota f(t))\quad(\iota\in I).
\]
Then, by standard density arguments, we find (one might also consult Proposition \ref{p:top} below):
\begin{enumerate}[label=(\alph*)]
 \item If $(S_\iota)_\iota$ converges in $(L(\s X,\s Y),\|\cdot\|_{L(\s X,\s Y)})$, then $(\s S^{(\nu)}_\iota)_\iota$ converges in $L_{\textnormal{sev},\nu}^\textnormal{n}(\s X,\s Y)$.
 \item If $(S_\iota)_\iota$ converges in the strong operator topology of $L(\s X,\s Y)$, then $(\s S^{(\nu)}_\iota)_\iota$ converges in $L_{\textnormal{sev},\nu}^\textnormal{s}(\s X,\s Y)$.
 \item If $(S_\iota)_\iota$ converges in the weak operator topology of $L(\s X,\s Y)$, then $(\s S^{(\nu)}_\iota)_\iota$ converges in $L_{\textnormal{sev},\nu}^\textnormal{w}(\s X,\s Y)$.
\end{enumerate} 
\end{Example}

Following the general philosophy of this exposition to formulate the results in $\nu$--independent type as much as possible the following small observation is in order.

\begin{Proposition}\label{p:top} Let $\s X$, $\s Y$ be Hilbert spaces, $D\subseteq \s X$ and $E\subseteq \s Y$ dense, $B\subseteq L(\s X,\s Y)$ bounded.
\begin{enumerate}[label=(\alph*)]
 \item\label{pn} Let $S\colon \dom(S)\subseteq \s X\to \s Y$ linear and assume $\dom(S)=D$. Then
 \[
     \sup_{\phi \in D\cap B_{\s X}}\|S\phi\|<\infty \iff \overline{S}\in L(\s X,\s Y).
 \]
 If $\overline{\s S}\in L(\s X,\s Y)$, then $\sup_{\phi \in D\cap B_{\s X}}\|S\phi\|=\|\overline{S}\|.$
 \item\label{ps} Let $\tau_{\textnormal{s}}$ be the strong operator topology on $L(\s X,\s Y)$. Define $\tau_{D}$ to be the initial topology on $L(\s X,\s Y)$ induced by the mappings
 \[
    T\mapsto T\phi\quad (\phi\in D).
 \]
 Then $(B,\tau_{\textnormal{s}})\hookrightarrow(B,\tau_D)$ is a homeomorphism.
 \item\label{pw} Let $\tau_{\textnormal{w}}$ be the weak operator topology on $L(\s X,\s Y)$. Define $\tau_{D,E}$ to be the initial topology on $L(\s X,\s Y)$ induced by the mappings
 \[
    T\mapsto \langle \psi,T\phi\rangle \quad (\psi\in E,\phi\in D).
 \]
 Then $(B,\tau_{\textnormal{w}})\hookrightarrow(B,\tau_{E,D})$ is a homeomorphism.
\end{enumerate} 
\end{Proposition}
\begin{Proof} For the proof of \ref{pn}, note that $\overline{S}\in L(\s X,\s Y)$ implies the finiteness of the supremum. The equality follows from the density of $D$. On the other hand, the supremum being finite together with the linearity of $S$ imply the Lipschitz continuity of $S$. Hence, $\overline{S}$ is Lipschitz continuous as well, yielding $\overline{S}\in L(\s X,\s Y)$.

To prove \ref{ps} and \ref{pw}, note that both the mappings $j_\textnormal{s}\colon (B,\tau_{\textnormal{s}})\to (B,\tau_D), x\mapsto x$ and $j_\textnormal{w}\colon(B,\tau_{\textnormal{w}})\to (B,\tau_{E,D}), x\mapsto x$ are continuous. Hence, it remains to prove continuity of $j_{\textnormal{s}}^{-1}$ and $j_{\textnormal{w}}^{-1}$. The arguments will conceptually be the same for both these inverses. We will only show continuity of $j_{\textnormal{w}}^{-1}$. 

Denote $\kappa \coloneqq \sup_{S\in B}\|S\|$, and let $(S_\iota)_\iota$ be a net in $(B,\tau_{E,D})$ convergent to some $T\in B$. Let $x\in \s X$, $y\in \s Y$, $\eps>0$. Choose $\phi\in D$, $\psi\in E$ such that $\|x-\phi\|+\|y-\psi\|\leq \eps$. There exists $\iota_0$ such that for all $\iota\geq\iota_0$ we have 
\begin{align*}
   |\langle y,(S_\iota-T)x\rangle|&\leq |\langle \psi,(S_\iota-T)\phi\rangle|+2\kappa (\|x\|\|y-\psi\|+\|y\|\|x-\phi\|)
   \\ &\leq \eps+2\kappa\eps(\|x\|+\|y\|).
\end{align*}
\end{Proof}

With the latter proposition at hand, we can formulate a description of the strong and weak operator topology on bounded subsets of $L_{\textnormal{sev}}$. We recall from Definition \ref{d:sev_bdd} that $\mathfrak{B}\subseteq L_{\textnormal{sev}}$ is bounded, if there is $\nu\in \mathbb{R}$ with $\mathfrak{B}\subseteq L_{\textnormal{sev},\nu}$ and $\sup_{\mu\geq\nu}\sup_{\s S\in \mathfrak{B}}\|\s S^\mu\|<\infty$. For a subset $\mathfrak{B}\subseteq L_{\textnormal{sev}}$, we denote the relative topology of $L_{\textnormal{sev}}^\textnormal{n}$, $L_{\textnormal{sev}}^\textnormal{s}$ and $L_{\textnormal{sev}}^\textnormal{w}$ by $\mathfrak{B}^\textnormal{n}$, $\mathfrak{B}^\textnormal{s}$ and $\mathfrak{B}^\textnormal{w}$, respectively.

\begin{Theorem}\label{t:ds} Let $\s X$, $\s Y$ Hilbert spaces, $\mathfrak{B}\subseteq L_{\textnormal{sev}}(\s X,\s Y)$ bounded. Recall the space $\mathcal{D}(\s X)\coloneqq \bigcap_{\nu\in \mathbb{R}}L_\nu^2(\mathbb{R};\s X)$. Then the following holds.

 A net $(\s S_\iota)_\iota$ in $\mathfrak{B}^\textnormal{s}$ is convergent to some $\s T\in \mathfrak{B}^\textnormal{s}$ if and only if there exists $\nu\in \mathbb{R}$ such that for all $\phi\in \s D(\s X)$ and $\mu\geq \nu$
 \begin{equation}\label{cscon}
     \s S_\iota \phi \stackrel{\iota}{\to} \s T\phi \text{ in }L_\mu^2(\R;\s Y).
 \end{equation}
\end{Theorem}
\begin{Proof}
 Let $(\s S_\iota)_{\iota\in I}$ a net in $\mathfrak{B}$, $\s T\in \mathfrak{B}$. Then, by definition, $\s S_\iota\stackrel{\iota}{\to}\s T$ in $\mathfrak{B}^\textnormal{s}$ if there exists $\nu\in \R$ such that $\{\s S_\iota;\iota\in I\}\subseteq L_{\textnormal{sev},\nu}(\s X,\s Y)$ and for all $\mu\geq \nu$ we have that $\s S^\mu_\iota\stackrel{\iota}{\to}\s T^\mu$ in $(L(L_\mu^2(\s X),L_\mu^2(\s Y)),\tau_{\textnormal{s}})$, that is, in the strong operator topology of $L(L_\mu^2(\s X),L_\mu^2(\s Y))$. But, the latter convergence implies that $\s S_\iota\stackrel{\iota}{\to}\s T$  in $(L(L_\mu^2(\s X),L_\mu^2(\s Y)),\tau_{\s D(\s X)})$, where we adopted the notation from Proposition \ref{p:top}\ref{ps} for the topology $\tau_{\s D(\s X)}$. Hence, $(\s S_\iota)_\iota$ converges to $\s T$ as in \eqref{cscon}.
 
 On the other hand, let $(\s S_\iota)_\iota$ converge to $\s T$ as in \eqref{cscon}. By the boundedness of $\mathfrak{B}\subseteq L_{\textnormal{sev}}(\s X,\s Y)$, there exists $\nu_0\in \mathbb{R}$ such that $\mathfrak{B}\subseteq L_{\textnormal{sev},\nu_0}(\s X,\s Y)$ with $\sup_{\mu\geq\nu_0,\s S\in \mathfrak{B}} \|\s S^\mu\|<\infty$. Moreover, by hypothesis, there exists $\nu_1\in \mathbb{R}$ such that for all $\mu\geq \nu_1$ we have that $(\s S_\iota)_\iota{\to}\s T$ in $(L(L_{\mu}^2(\s X),L_\mu^2(\s Y)),\tau_{\s D(\s X)})$. The latter is the same as $(\s S_\iota^\mu)_\iota {\to}\s T^\mu$ in the space $(L(L_{\mu}^2(\s X),L_\mu^2(\s Y)),\tau_{\s D(\s X)})$. Hence, for $\mu\geq \max\{\nu_0,\nu_1\}$ we get by Proposition \ref{p:top}\ref{ps} that $(\s S_\iota^\mu)_\iota$ converges to $\s T^\mu$ in the strong operator topology of $L(L_\mu^2(\s X),L_\mu^2(\s Y))$. Thus, as $L_{\textnormal{sev},\nu_0}^\textnormal{s}\hookrightarrow L_{\textnormal{sev},\max\{\nu_0,\nu_1\}}^\textnormal{s}$, the assertion follows. 
\end{Proof}

\begin{Theorem}\label{t:dw} Let $\s X$, $\s Y$ Hilbert spaces, $\mathfrak{B}\subseteq L_{\textnormal{sev}}(\s X,\s Y)$ bounded. Define the space $\mathcal{D}(\s X)\coloneqq \bigcap_{\nu\in \mathbb{R}}L_\nu^2(\mathbb{R};\s X)$ and correspondingly $\mathcal{D}(\s Y)$. Define the topology $\tau_{\mathcal{D}(\s Y),\mathcal{D}(\s X)}$ on $\mathfrak{B}$ induced by
 \begin{equation}\label{cwcon}
  \mathfrak{B}\ni \s S \mapsto \int_\R \langle \psi(t),\s S\phi(t)\rangle_{\s Y}\dd t \quad (\phi\in \mathcal{D}(\s X), \psi\in \mathcal{D}(\s Y)). 
 \end{equation}
 Then $\mathfrak{B}^{\textnormal{w}}\hookrightarrow (\mathfrak{B},\tau_{\s D(\s Y),\s D(\s X)})$ is a homeomorphism. 
\end{Theorem}

For the proof of the latter theorem, we need two elementary observations. The first one will be stated without proof.

\begin{Lemma}\label{l:intop} Let $\s X$ be a vector space, $\mathfrak{B}\subseteq \s X$, $\s D\subseteq \s X'$, $f\colon \s D\to \s D$ bijective. On $\mathfrak{B}$, the topology $\tau_{\s D}$ induced by
\[
   \mathfrak{B}\ni S\mapsto  x'(S)\quad(x'\in \s D)
\]
coincides with $\tau_{f[\s D]}$ induced by
\[
   \mathfrak{B}\ni S\mapsto  (f(x'))(S)\quad(x'\in \s D).
\] 
\end{Lemma}

\begin{Lemma}\label{l:topind} Let $\s X$, $\s Y$ Hilbert spaces, $\mathfrak{B}\subseteq L_{\textnormal{sev}}(\s X,\s Y)$, $\s D\subseteq \bigcap_{\nu\in \mathbb{R}}L_\nu^2(\s X)$, $\s E\subseteq \bigcap_{\nu\in \R} L_\nu^2(\s Y)$, $\mu\in \mathbb{R}$. Assume 
\[
   \{ t\mapsto e^{-2 \mu t}\psi(t); \psi\in \s E\}=\s E.
\]
Let $\tau_0$ be the initial topology on $\mathfrak{B}$ induced by
\begin{equation}\label{eq:tau0}
   \s S\mapsto \int_{\R} \langle \psi(t),\s S \phi(t)\rangle\dd t\quad (\psi\in \s E,\phi\in \s D)
\end{equation}
and $\tau_\mu$ be the initial topology on $\mathfrak{B}$ induced by
\begin{equation}\label{eq:taunu}
   \s S\mapsto \int_{\R} \langle \psi(t),\s S \phi(t)\rangle e^{-2\mu t}\dd t\quad (\psi\in \s E,\phi\in \s D)
\end{equation}
Then $\tau_0 = \tau_\mu$.
\end{Lemma}
\begin{Proof}
 Note that for any $\s S\in L_{\textnormal{sev}}(\s X,\s Y)$ the integrals in \eqref{eq:tau0} and \eqref{eq:taunu} exist: Indeed, for $\s S$ there exists $\eta\in \mathbb{R}$ such that $\s S^\eta\in L(L_{\eta}^2(\s X),L_{\eta}^2(\s Y))$. Thus, using \begin{multline*}
  \s D\subseteq \bigcap_{\nu\in \mathbb{R}}L_\nu^2(\s X)\subseteq L_\eta^2(\s X)\text{ and }\\e^{2 \eta (\cdot)}[\s E]\subseteq e^{2 \eta (\cdot)}[\bigcap_{\nu\in \mathbb{R}}L_\nu^2(\s Y)]=\bigcap_{\nu\in \mathbb{R}}L_\nu^2(\s Y)\subseteq L_\eta^2(\s Y),
 \end{multline*}we get
 \begin{align*}
   \int_\R |\langle \psi(t),\s S \phi(t)\rangle|\dd t& = \int_\R|\langle e^{2\eta t} \psi(t),\s S \phi(t)\rangle| e^{-2\eta t}\dd t
   \\ & \leq \int_\R\|e^{2\eta t} \psi(t)\|_{\s Y}\|\s S \phi(t)\|_{\s Y} e^{-2\eta t}\dd t
   \\ & \leq \|e^{2\eta (\cdot)} \psi(\cdot)\|_{L_\eta^2(\s Y)} \|\s S \phi\|_{L_\eta^2(\s Y)}.
 \end{align*}
 The same argument applies to \eqref{eq:taunu}. The assertion, $\tau_0= \tau_{\mu}$, follows from Lemma \ref{l:intop}.
\end{Proof}

\begin{Proof}[of Theorem \ref{t:dw}]
  By the boundedness of $\mathfrak{B}\subseteq L_{\textnormal{sev}}(\s X,\s Y)$ there exists $\nu\in \mathbb{R}$ with $\mathfrak{B}\subseteq L_{\textnormal{sev},\nu}(\s X,\s Y)$ and
  \[
     \sup_{\mu\geq \nu}\sup_{\s S\in \mathfrak{B}}\|\s S^\mu\|<\infty.
  \]  Next, note that $\s D(\s X)$ is dense in $L_\mu^2(\mathbb{R};\s X)$ for all $\mu\geq\nu$ and for all $\mu\in \R$, we have $e^{-2\mu(\cdot)}[\s D(\s Y)]=\s D(\s Y)$.
  Hence, by Lemma \ref{l:topind},
  \[
     (\mathfrak{B},\tau_{\s D(\s Y),\s D(\s X)})\hookrightarrow(\mathfrak{B},\tau_{\s D(\s Y),\s D(\s X),\mu})
  \]
  is a homeomorphism, where $\tau_{\s D(\s Y),\s D(\s X),\mu}$ is the initial topology induced by
 \[
    \s S\mapsto \langle \psi,\s S\phi\rangle_{L_{\mu}^2(\s Y)}\quad(\psi\in \s D(\s Y),\phi\in \s D(\s X)).
 \]
 Next, denoting $\mathfrak{B}_{\mu}\coloneqq \{ \s S^\mu; \s S\in \mathfrak{B}\}$, we realize that the spaces  $(\mathfrak{B},\tau_{\s D(\s Y),\s D(\s X),\mu})$ and $(\mathfrak{B}_\mu,\tau_{\s D(\s Y),\s D(\s X),\mu})$ are (obviously) homeomorphic via $\s S\mapsto \s S^\mu$. Hence, by Proposition \ref{p:top}\ref{pw}, we infer for all $\mu\geq \nu$ that
 \begin{equation}\label{eq:ho}
    (\mathfrak{B},\tau_{\s D(\s Y),\s D(\s X)})\ni \s S \mapsto \s S^\mu \in (\mathfrak{B}_\mu,\tau_{\textnormal{w}})\text{ is a homeomorphism,}
 \end{equation}
 where $\tau_{\textnormal{w}}$ is the weak operator topology on $L(L_\mu^2(\s X),L_\mu^2(\s Y))$. So, if a net $(\s S_\iota)_\iota$ converges to $\s T\in\mathfrak{B}^\textnormal{w}$, then there exists $\mu\geq \nu$ such that $(\s S^\mu_\iota)_\iota$ converges to $\s T^\mu$ in $(\mathfrak{B}_\mu,\tau_{\textnormal{w}})$. Hence, by \eqref{eq:ho}, $(\s S_\iota)_\iota$ converges to $\s T$ in $(\mathfrak{B},\tau_{\s D(\s Y),\s D(\s X)})$. On the other hand, if $(\s S_\iota)_\iota$ converges to some $\s T\in \mathfrak{B}$ in $(\mathfrak{B},\tau_{\s D(\s Y),\s D(\s X)})$, then, by \eqref{eq:ho} $(\s S_\iota^\mu)_{\iota}$ converges to $\s T^\mu$ in $(\mathfrak{B}_\mu,\tau_{\textnormal{w}})$ for all $\mu\geq \nu$. We infer that $(\s S_\iota)_\iota$ converges to $\s T$ in $\mathfrak{B}^\textnormal{w}$.
\end{Proof}

Another application of the almost trivial observation in Lemma \ref{l:topind} can be found in the proof of the next lemma, which will lead us to the proof of a compactness property for the weak operator topology of standard evolutionary mappings.

\begin{Lemma}\label{l:bhssev} Let $\s X, \s Y$ Hilbert spaces, $\nu\in \R$, $\mathfrak{B}\subseteq L_{\textnormal{sev}}(\s X,\s Y)$ with $\sup_{\mu\geq\nu,\s S\in \mathfrak{B}} \|\s S^\mu\|<\infty$. Then $\mathfrak{B}^\textnormal{w}\subseteq L_{\textnormal{sev}}^\textnormal{w}\hookrightarrow\mathfrak{B}\subseteq L_{\textnormal{sev},\nu}^\textnormal{w}(\s X,\s Y)$ is a homeomorphism.
\end{Lemma}
\begin{Proof}
 By definition $L_{\textnormal{sev},\nu}^\textnormal{w}(\s X,\s Y)\cap\mathfrak{B}$ is continuously embedded into $\mathfrak{B}^\textnormal{w}$. On the other hand, if $(\s S_\iota)_\iota$ converges to $\s T\in \mathfrak{B}^\textnormal{w}$, then there exists $\eta\in \mathbb{R}$ such that $(\s S_\iota^\eta)_\iota$ converges to $\s T^\eta$ in the weak operator topology of $L(L_\eta^2(\s X),L_{\eta}^2(\s Y))$. Hence, $(\s S_\iota)_\iota$ converges to $\s T$ in $(\mathfrak{B},\tau_\eta)$, where $\tau_\eta$ is given as in \eqref{eq:taunu} with $\eta$ instead of $\mu$ and $\s E=\s D(\s Y)$, $\s D=\s D(\s X)$. By Lemma \ref{l:topind}, $(\s S_\iota)_\iota$ converges to $\s T$ in $(\mathfrak{B},\tau_\mu)$ for all $\mu\in \mathbb{R}$. Thus, by Proposition \ref{p:top} and the boundedness of $\mathfrak{B}$ considered in $L(L_\mu^2(\s X),L_\mu^2(\s Y))$, we get that $(\s S_\iota^\mu)_\iota$ converges to $\s T^\mu$ in the weak operator topology of $L(L_\mu^2(\s X),L_\mu^2(\s Y))$ for all $\mu\geq\nu$. Thus, $(S_\iota)_\iota$ converges to $\s T$ in $L_{\textnormal{sev},\nu}^\textnormal{w}(\s X,\s Y)$.
\end{Proof}

\begin{Lemma}\label{l:ttp} Let $\s X,\s Y$ Hilbert spaces, $\nu\in \R$. Define the topological space
 \[\mathfrak{R}\coloneqq \prod_{\mu\geq\nu}\Big(\prod_{\substack{\phi\in L_\mu^2(\s X)\\ \psi\in L_\mu^2(\s Y)}} \overline{B_{\mathbb{C}}(0,\|\phi\|\|\psi\|)}\cap \{ \s S; \s S\colon L_\mu^2(\s X)\times L_\mu^2(\s Y)\to \mathbb{C}\text{ sesquilinear}\}\Big)\]
 endowed with the product topology.
 
 Then $\mathfrak{R}$ is compact.
\end{Lemma}
\begin{Proof}
We use Tikhonov's Theorem to deduce that
 \[
  \mathfrak{R}_{\mu}\coloneqq \prod_{\phi\in L_\mu^2(\s X),\psi\in L_\mu^2(\s Y)} \overline{B_{\mathbb{C}}(0,\|\phi\|\|\psi\|)}
 \]
 is compact. Moreover, it is easy to see that
 \[
   \mathfrak{R}_{\mu,\textnormal{ses}}\coloneqq \mathfrak{R}_{\mu} \cap \{ \s S; \s S\colon L_\mu^2(\s X)\times L_\mu^2(\s Y)\to \mathbb{C}\text{ sesquilinear}\}\subseteq\mathfrak{R}_{\mu}
  \]
  is closed. Hence, $\mathfrak{R}_{\mu,\textnormal{ses}}$ is compact, and, consequently, so is $\mathfrak{R}=\prod_{\mu\geq\nu} \mathfrak{R}_{\mu,\textnormal{ses}}$ by Tikhonov's Theorem again.
\end{Proof}
\begin{Theorem}\label{t:wco} Let $\s X$, $\s Y$ Hilbert spaces, $\mathfrak{B}\subseteq L_{\textnormal{sev}}(\s X,\s Y)$ bounded.
Then $\mathfrak{B}^\textnormal{w}\subseteq L_{\textnormal{sev}}^\textnormal{w}(\s X,\s Y)$ is relatively compact.
\end{Theorem}
\begin{Proof}
 Without loss of generality $\mathfrak{B}=\{\s S\in L_{\textnormal{sev},\nu}(\s X,\s Y); \sup_{\mu\geq\nu} \|\s S^\mu\|\leq 1\}$ for some $\nu\in \mathbb{R}$.
 Recall the compact topological space 
 \[\mathfrak{R}\coloneqq \prod_{\mu\geq\nu}\Big(\prod_{\substack{\phi\in L_\mu^2(\s X)\\ \psi\in L_\mu^2(\s Y)}} \overline{B_{\mathbb{C}}(0,\|\phi\|\|\psi\|)}\cap \{ \s S; \s S\colon L_\mu^2(\s X)\times L_\mu^2(\s Y)\to \mathbb{C}\text{ sesquilinear}\}\Big)\] from Lemma \ref{l:ttp}.
  
  In order that $\mathfrak{B}^{\textnormal{w}}$ is compact, we show that $\mathfrak{B}^\textnormal{w}$ can be identified with a closed subspace of $\mathfrak{R}$. Recalling the Riesz--Frechet representation theorem, we observe that, for Hilbert spaces $\s W$, $\s Z$, any contraction $T\in L(\s W,\s Z)$ is a sesquilinear mapping $s_T$ on $\s Z\times \s W$ with bound $1$ and vice versa. The isomorphism is induced by $s_{(\cdot)}\colon T\mapsto \left((\phi,\psi)\mapsto \langle \phi, T\psi\rangle\right)$. Hence, 
  \[
     j\colon \mathfrak{B}\ni \s T \mapsto (s_{\s T^\mu})_{\mu\geq\nu} \in \mathfrak{R}
  \]
  is a well-defined one-to-one mapping. By Lemma \ref{l:bhssev}, $\mathfrak{B}^\textnormal{w}$ carries the topology of $L_{\textnormal{sev},\nu}^\textnormal{w}$.
  Hence, by definition of the topology on $L_{\textnormal{sev},\nu}^\textnormal{w}$, the mapping $j\colon \mathfrak{B}^\textnormal{w}\to \mathfrak{R}$ is continuous. 
  
  For proving that $j$ is a homeomorphism and that $j[\mathfrak{B}^\textnormal{w}]\subseteq \mathfrak{R}$ is closed, we are left with showing that for any closed $\mathfrak{A}\subseteq \mathfrak{B}^\textnormal{w}$ the set $j[\mathfrak{A}]\subseteq \mathfrak{R}$ is closed as well. For this, let $\mathfrak{A}\subseteq \mathfrak{B}^\textnormal{w}$ closed, and $(\s T_\iota)_\iota$ be a net in $\mathfrak{A}$ such that $(j(\s T_\iota))_\iota$ converges in $\mathfrak{R}$ to some $(s_\mu)_{\mu\geq\nu}$. We have to show that $(\s T_\iota)_\iota$ converges in $\mathfrak{B}^\textnormal{w}$ to some $\s S\in \mathfrak{A}$ with $j(\s S)=\lim_{\iota} j(\s T_\iota)$.
  
  Employing the Riesz--Frechet theorem again, we infer -- by the definition of $\mathfrak{R}$ -- that there exists $\s S_\mu \in L(L_\mu^2(\s X),L_\mu^2(\s Y))$, $\|\s S_\mu\|\leq 1$, with $s_\mu = s_{\s S_\mu}$ for all $\mu\geq\nu$. Moreover, $\s T_\iota^\mu\stackrel{\iota}{\to} \s S_\mu$ in the weak operator topology. Hence, if we show that there exists $\s S\in \mathfrak{B}^\textnormal{w}$ with $\s S^\mu=\s S_\mu$ for all $\mu\geq \nu$, we infer $\lim_\iota \s T_\iota = \s S\in \mathfrak{A}$, by the closedness of $\mathfrak{A}$.
  
  In order that $\s S^\mu=\s S_\mu$ for all $\mu\geq \nu$ for some $\s S\in \mathfrak{B}$, it suffices to show that $\s S_{\mu_1}=\s S_{\mu_2}$ on $L_{\mu_1}^2(\s X)\cap L_{\mu_2}^2(\s X)$ for all $\mu_1,\mu_2\geq \nu$. So, let $f\in L_{\mu_1}^2(\s X)\cap L_{\mu_2}^2(\s X)$. Since both $\s S_{\mu_1}f,\s S_{\mu_2}f$ are measurable there exists a nullset $N_0$ such that for $J\coloneqq \mathbb{R}\setminus N_0$ we have $(\s S_{\mu_1}f)[J]\cup (\s S_{\mu_2}f)[J]\subseteq \s Y_0$ for some closed separable subspace $\s Y_0\subseteq \s Y$. Let $\s D_0\subseteq \s Y_0$ be countable and dense. Next, it is plain that $\s D\coloneqq \{h y; h\in C_c^\infty(\mathbb{R}), y\in \s D_0\}\subseteq \s D(\s Y)$. But, for $\psi=hy\in \s D$  we have
  \begin{multline*}
     \int_\R \langle\s S_{\mu_1} f(t),y\rangle_{\s Y}h(t) \dd t =\lim_\iota \int_\R \langle \s T_\iota^{\mu_1}f(t),\psi(t) \rangle_{\s Y} \dd t \\=\lim_\iota \int_\R \langle \s T_\iota^{\mu_2}f(t),\psi(t) \rangle_{\s Y} \dd t =\int_\R \langle\s S_{\mu_2} f(t),y\rangle_{\s Y}h(t) \dd t.
  \end{multline*}
  Hence, by the fundamental lemma of the calculus of variations and the countability of $\s D_0$, we infer that there exists a nullset $N_1\supseteq N_0$ such that for all $t\in \R\setminus N_1$ and $y\in \s D_0$, we have $\langle\s S_{\mu_1} f(t),y\rangle =\langle\s S_{\mu_2} f(t),y\rangle$. Thus, by the density of $\s D_0\subseteq \s Y_0$, we conclude $\s S_{\mu_1} f(t) =\s S_{\mu_2} f(t)$ for all $t\in \R\setminus N_1$.
 \end{Proof}
 
 \begin{Corollary}\label{c:wmet} Let $\s X$, $\s Y$ be separable Hilbert spaces, $\mathfrak{B}\subseteq L_{\textnormal{sev}}(\s X,\s Y)$ bounded. Then $\mathfrak{B}^\textnormal{w}$ is metrizable.  
 \end{Corollary}
 \begin{Proof} Without loss of generality $\mathfrak{B}=\{\s S\in L_{\textnormal{sev},\nu}(\s X,\s Y); \sup_{\mu\geq\nu} \|\s S^\mu\|\leq 1\}$ for some $\nu\in \mathbb{R}$.
 Since $\s X$ and $\s Y$ are separable, then so are $L_\nu^2(\s X)$ and $L_\nu^2(\s Y)$. In particular, the (standard causal) domains $\s D(\s X)=\bigcap_{\mu \in\R}L_{\mu}^2(\s X)$ and $\s D(\s Y)$ considered as respective (metric) subspaces of $L_\nu^2(\s X)$ and $L_\nu^2(\s Y)$ are separable, as well. So, take countable dense sets $\s D\subseteq \s D(\s X)$ and $\s E\subseteq \s D(\s Y)$. Then the mapping
 \[
    j \colon \mathfrak{B}^\textnormal{w}\ni \s S \mapsto \big((\psi,\phi)\mapsto \langle \psi,\s S\phi\rangle_{L_\nu^2(\s Y)}\big) \in \prod_{\psi\in \s E,\phi\in \s D} \overline{B_{\mathbb{C}}(0,\|\phi\|\|\psi\|)}
 \]
 is continuous and one-to-one. For injectivity of $j$ use that $\s E$ and $\s D$ are respectively dense in $L_\nu^2(\s Y)$ and $L_\nu^2(\s X)$. As $\mathfrak{B}^\textnormal{w}$ is compact, we infer that $j$ is a homeomorphism onto its image. Hence, as $\prod_{\psi\in \s E,\phi\in \s D} \overline{B_{\mathbb{C}}(0,\|\phi\|\|\psi\|)}$ is metrizable by the countability of $\s E\times \s D$, the metrizability of $\mathfrak{B}^\textnormal{w}$ follows.
\end{Proof}

 In the concluding parts of this section, we address continuity of multiplication and inversion in $L_{\textnormal{sev}}$. We recall that both the latter operations are \emph{not} continuous under the weak operator topology of $L(\s X,\s Y)$ for infinite-dimensional Hilbert spaces $\s X$ and $\s Y$. Thus, we cannot expect continuity of multiplication and inversion in $L_{\textnormal{sev}}^\textnormal{w}$, leading to more subtle statements, when the weak topology is involved. The results on continuity of multiplication read as follows, for which we recall from Definition \ref{d:sev_bdd} that a set  $\mathfrak{B}\subseteq L_{\textnormal{sev}}(\s X,\s Y)$ is called bounded, if there exists $\nu\in \mathbb{R}$ with $\mathfrak{B}\subseteq L_{\textnormal{sev},\nu}(\s X,\s Y)$ and that
\[
   \sup_{\mu\geq\nu}\sup_{\s S\in \mathfrak{B}} \|\s S^\mu\|<\infty.
\]
We say that $\mathfrak{D}\subseteq L_{\textnormal{sev}}(\s X,\s Y)\times L_{\textnormal{sev}}(\s Y,\s Z)$ is \emph{bounded}, if there exists bounded $\mathfrak{B}_1\subseteq L_{\textnormal{sev}}(\s X,\s Y)$, $\mathfrak{B}_2\subseteq L_{\textnormal{sev}}(\s Y,\s Z)$ with $\mathfrak{D}\subseteq \mathfrak{B}_1\times \mathfrak{B}_2$. We also recall Remark \ref{r:pt}.

\begin{Theorem}\label{t:cim} Let $\s X$, $\s Y$, $\s Z$ Hilbert spaces. Consider the multiplication $(\s S,\s T) \mapsto \s T\s S$ as a mapping in the following underlying topological spaces:
\begin{enumerate}[label=(\alph*)]
 \item\label{nn}$ L^\textnormal{n}_{\textnormal{sev}}(\s X,\s Y)\times L^\textnormal{n}_{\textnormal{sev}}(\s Y,\s Z)  \to L^\textnormal{n}_{\textnormal{sev}}(\s X,\s Z),$
 \item\label{ss}$ L^\textnormal{s}_{\textnormal{sev}}(\s X,\s Y)\times L^\textnormal{s}_{\textnormal{sev}}(\s Y,\s Z)  \to L^\textnormal{s}_{\textnormal{sev}}(\s X,\s Z),$
 \item\label{ww}$ L^\textnormal{w}_{\textnormal{sev}}(\s X,\s Y)\times L^\textnormal{w}_{\textnormal{sev}}(\s Y,\s Z)  \to L^\textnormal{w}_{\textnormal{sev}}(\s X,\s Z).$
\end{enumerate}
Then the multiplication in \ref{nn} is continuous, and, on bounded subsets, multiplication in \ref{ss} is continuous, whereas in \ref{ww} multiplication is only separately continuous, that is, for every $\s V\in L_{\textnormal{sev}}(\s X,\s Y)$ and $\s W\in L_{\textnormal{sev}}(\s Y,\s Z)$ the mappings $\s S\mapsto \s W\s S$ and $\s T\mapsto \s T\s V$ are continuous.
\end{Theorem}
\begin{Proof}
 Taking the continuous inclusions \eqref{embmon} into account, the results are straightforward consequences of the corresponding statements, when $L_{\textnormal{sev}}$ is replaced by $L(\s X,\s Y)$. As an example, we treat \ref{ss}: Let $((\s S_\iota,\s T_\iota))_\iota$ be a convergent net in the space $L^\textnormal{s}_{\textnormal{sev}}(\s X,\s Y)\times L^\textnormal{s}_{\textnormal{sev}}(\s Y,\s Z)$; denote the respective limit by $(\s U,\s W)$. By definition, there exists $\nu_0$ and $\nu_1$ such that $((\s S_\iota,\s T_\iota))_\iota$ converges in $L^\textnormal{s}_{\textnormal{sev},\nu_0}\times L^\textnormal{s}_{\textnormal{sev},\nu_1}\subseteq L^\textnormal{s}_{\textnormal{sev},\nu}\times L^\textnormal{s}_{\textnormal{sev},\nu}$, $\nu\coloneqq \max\{\nu_0,\nu_1\}$ (see \eqref{embmon}). Hence, in particular, by the submultiplicativity of the norm, we deduce that $\s T_\iota\s S_\iota \in L_{\textnormal{sev},\nu}$. And, by the corresponding statement of \ref{ss} when $L_{\textnormal{sev}}(\s X,\s Y)$ is replaced by $L(\s X,\s Y)$, the assertion follows.
\end{Proof}

\begin{Remark}\label{r:cim_fnu} Let $\nu\in \mathbb{R}$. We note that, if, in Theorem \ref{t:cim}, one replaces all $L_{\textnormal{sev}}(\s X,\s Y)$ by $L_{\textnormal{sev},\nu}(\s X,\s Y)$ the respective assertions hold true as well. In fact, it is even a more direct consequence of the analogous statements for the operator topologies on $L(\s W,\s Z)$ for Hilbert spaces $\s W$, $\s Z$. 
\end{Remark}

We conclude this section with addressing the continuity of computing the inverse of evolutionary mappings. Beforehand, for a Hilbert space $\s X$, $\kappa\geq 0$, we introduce the sets
\begin{align*}
   GL_{\textnormal{sev},\nu}(\s X) & \coloneqq \{ \s S\in L_{\textnormal{sev}}(\s X); \s S^{-1} \in L_{\textnormal{sev},\nu}(\s X)\},
   \\ GL_{\kappa,\textnormal{sev},\nu}(\s X) & \coloneqq \{ \s S\in GL_{\textnormal{sev},\nu}(\s X); \sup_{\mu\geq \nu}\|(\s S^{-1})^\mu\|\leq \kappa\}.
\end{align*}
We recall our convention to denote by $\mathfrak{B}^\textnormal{n}$ the topological subspace of $L_{\textnormal{sev}}^\textnormal{n}$ for a subset $\mathfrak{B}\subseteq L_{\textnormal{sev}}$ and similarly for $\mathfrak{B}^\textnormal{s}$.

\begin{Theorem}\label{t:cinv} Let $\s X$ Hilbert space, $\kappa\geq 0$. The mapping $\s S\mapsto \s S^{-1}$ is continuous as a mapping in the following underlying topological spaces
\begin{enumerate}[label=(\alph*)]
 \item\label{n} $GL_{\textnormal{sev},\nu}(\s X)^\textnormal{n}\to L_{\textnormal{sev}}^\textnormal{n}(\s X)$,
 \item\label{s} $GL_{\kappa,\textnormal{sev},\nu}(\s X)^\textnormal{s}\to L_{\textnormal{sev}}^\textnormal{s}(\s X)$.
\end{enumerate}
\end{Theorem}
\begin{Proof}
 For the proof of \ref{n}, we let $(\s S_\iota)_\iota$ be a convergent net in $GL_{\textnormal{sev},\nu}^\textnormal{n}$; denote $\s T\coloneqq \lim_\iota \s S_\iota$. Hence, there exists $\eta\geq\nu$ such that $\s S_\iota,\s T,\s S^{-1}_\iota, \s T^{-1}\in L_{\textnormal{sev},\eta}(\s X)$ for all $\iota$. Hence, for all $\iota$, we have 
 $\s S_\iota^{-1}- \s T^{-1}=\s S_\iota^{-1}(\s T-\s S_\iota)\s T^{-1} \in L_{\textnormal{sev},\eta}(\s X)$, by Proposition \ref{p:sev_alg}. Next, let $\mu\geq \eta$. Then $(\s S^{\mu}_\iota)^{-1}(\s T^\mu-\s S^{\mu}_\iota)(\s T^{\mu})^{-1}\to 0$ since  $\|\s T^\mu-\s S^{\mu}_\iota\|\to 0$ and $\|(\s S^{\mu}_\iota)^{-1}\|=\|(\s T^{\mu} + (\s S^{\mu}_\iota-\s T^\mu))^{-1}\|\leq 2\|(\s T^{\mu})^{-1}\|$, once $\|(\s T^{\mu})^{-1}(\s S^{\mu}_\iota-\s T^\mu)\|\leq 1/2$.
 
 In order to show the corresponding statement for \ref{s}, we similarly let $(\s S_\iota)_\iota$ be convergent in $GL_{\kappa,\textnormal{sev},\nu}^\textnormal{s}$; with $\s T$ the corresponding limit. As in \ref{n}, we find $\eta\geq \nu$, such that $\s S_\iota^{-1}(\s T-\s S_\iota)\s T^{-1} \in L_{\textnormal{sev},\eta}(\s X)$ for all $\iota$. Moreover, by definition of $GL_{\kappa,\textnormal{sev},\nu}(\s X)$, we have $\sup_{\mu\geq \eta, \iota}\|(\s S_\iota^{-1})^\mu\|\leq \kappa$. Hence, $(\s S_\iota^{-1})^\mu(\s T^\mu-\s S_\iota^\mu)(\s T^{-1})^\mu$ converges strongly to $0$ in $(L(L_\mu^2(\s X)),\tau_{\textnormal{s}})$, which yields the assertion.
\end{Proof}

\renewcommand{\baselinestretch}{0.65}\normalsize\mysection{Norm and Strong Topologies and Ordinary Differential Equations}{The Norm and the Strong Operator Topologies and Ordinary Differential Equations}{the set $\textnormal{SO}_{c,\nu}(\s X,\s Y)$ $\cdot$ continuity of the solution operator for ordinary differential equations in the coefficients $\cdot$ norm convergent series' and the norm, strong and weak operator topologies $\cdot$ Theorem \ref{t:cdoden} $\cdot$ Theorem \ref{t:cdodes}}\label{s:nsode}

\renewcommand{\baselinestretch}{1}\normalsize
In this section and the next section, we are aiming for continuity results of ordinary differential equations on the coefficients. More precisely, we will focus on equations of the form
\begin{equation}\label{eq:nsode}
   \Big(\partial_{t,\nu} \begin{pmatrix}
                                     \mathcal{M} & 0 \\ 0 & 0 
                                   \end{pmatrix} + \begin{pmatrix}
                                     \mathcal{N}_{00} & \mathcal{N}_{01} \\ \mathcal{N}_{10} & \mathcal{N}_{11}
                                   \end{pmatrix}\Big) U= F,
\end{equation}
where $\s M$, $\mathcal{N}_{00}$, $\mathcal{N}_{01}$, $\mathcal{N}_{10}$ and $\mathcal{N}_{11}$ are standard evolutionary mappings acting in suitable spaces. Assuming the well-posedness conditions as in Theorem \ref{t:wp_ode_evo}, we address the question of whether assigning a solution operator to \eqref{eq:nsode} is continuous. We recall the conditions that lead to a solution theory for \eqref{eq:nsode}: Let $\s X$, $\s Y$ be Hilbert spaces, $\nu,c>0$. Then $\s M\in L_{\textnormal{sev}}(\s X)$ needed to satisfy
\begin{equation}\label{eq:wpcondm}
   \Re\langle Q_t \s M \phi,\phi\rangle_{L_\mu^2(\s X)}\geq c\langle Q_t\phi,\phi\rangle_{L_\mu^2(\s X)}\quad (\phi\in \s D(\s X), t\in \mathbb{R},\mu\geq\nu).
\end{equation}
The condition for $\s N=(\s N_{ij})_{i,j\in\{0,1\}}\in L_{\textnormal{sev}}(\s X\times \s Y)$ read
\begin{equation}\label{eq:wpcondn}
   \Re\langle Q_t \s N_{11} \phi,\phi\rangle_{L_\mu^2(\s Y)}\geq c\langle Q_t\phi,\phi\rangle_{L_\mu^2(\s Y)}\quad (\phi\in \s D(\s Y), t\in \mathbb{R},\mu\geq\nu).
\end{equation}
We also recall that $Q_t$ is multiplication by $\1_{(-\infty,t)}$ and $\s D(\s X)=\bigcap_{\nu\in \R}L_{\nu}^2(\s X)$ and similarly for $\s D(\s Y)$.
For Hilbert spaces $\s X$ and $\s Y$, $c>0$, we define
\begin{align*}
   &\textnormal{SO}_{c,\nu}(\s X,\s Y)
   \\ &\coloneqq \{ \mathcal{M}\in L_{\textnormal{sev}}(\s X); \s M\text{ satisfies }\eqref{eq:wpcondm}\}\times \{ \mathcal{N}\in L_{\textnormal{sev}}(\s X\times \s Y); \s N\text{ satisfies }\eqref{eq:wpcondn}\}
   \\&\subseteq L_{\textnormal{sev}}(\s X)\times L_{\textnormal{sev}}(\s X\times\s Y)\subseteq L_{\textnormal{sev}}(\s X\times \s X\times\s Y).
\end{align*}
In the notation $\textnormal{SO}_{c,\nu}(\s X,\s Y)$ the letter `S' is a reminder of `solution theory', the `O' stands for `ordinary differential equations'. Being a subset of $L_{\textnormal{sev}}(\s X\times \s X\times\s Y)$, we may endow $\textnormal{SO}_{c,\nu}(\s X,\s Y)$ with the norm, the strong operator or the weak operator topology. It is easy to see that $\textnormal{SO}_{c,\nu}(\s X,\s Y)$ is a closed subset of each $L_{\textnormal{sev}}^\textnormal{n}(\s X\times \s X\times\s Y)$, $L_{\textnormal{sev}}^\textnormal{s}(\s X\times \s X\times\s Y)$ and $L_{\textnormal{sev}}^\textnormal{w}(\s X\times \s X\times\s Y)$. With the notation just introduced and recalling $\check{\partial}_{t}\coloneqq \bigcap_{\nu>0}\partial_{t,\nu}$, we may rephrase Theorem \ref{t:wp_ode_evo} next. For this, we adopt the custom that evolutionary and causal mappings are, in fact, standard evolutionary (and vice versa) in the sense of Theorem \ref{t:evsc_sev}.

\begin{Theorem}[Theorem \ref{t:wp_ode_evo}]\label{t:wp_ode_evoprm} Let $\s X$, $\s Y$ Hilbert spaces, $\nu,c>0$. Then the mapping
\begin{alignat}{1}\label{eq:ode_solm}
\begin{aligned}
     \textnormal{sol} \colon \textnormal{SO}_{c,\nu}(\s X,\s Y) &\to L_{\textnormal{sev}}(\s X\times \s Y)
     \\  (\s M,\s N) & \mapsto \Big(\check{\partial}_{t} \begin{pmatrix}
                                     \mathcal{M} & 0 \\ 0 & 0 
                                   \end{pmatrix} + \begin{pmatrix}
                                     \mathcal{N}_{00} & \mathcal{N}_{01} \\ \mathcal{N}_{10} & \mathcal{N}_{11}
                                   \end{pmatrix}\Big)^{-1}
\end{aligned}
\end{alignat}
is well-defined.
\end{Theorem}

So, the aim of this section is to establish the following two theorems.

\begin{Theorem}\label{t:cdoden} Let $\s X$, $\s Y$ Hilbert spaces, $\nu,c>0$. Then \textnormal{sol} given in \eqref{eq:ode_solm} is continuous on bounded sets as a mapping from $\textnormal{SO}_{\nu,c}(\s X,\s Y)^\textnormal{n}$ to $L_{\textnormal{sev}}^\textnormal{n}(\s X\times\s Y)$. 
\end{Theorem}

\begin{Theorem}\label{t:cdodes} Let $\s X$, $\s Y$ Hilbert spaces, $c>0$. Then \textnormal{sol} given in \eqref{eq:ode_solm} is continuous on bounded sets as a mapping from
$\textnormal{SO}_{\nu,c}(\s X,\s Y)^\textnormal{s}$ to $L_{\textnormal{sev}}^\textnormal{s}(\s X\times\s Y)$. 
\end{Theorem}

For the respective proofs of Theorem \ref{t:cdoden} and Theorem \ref{t:cdodes}, we need to analyze the solution operator. By Remark \ref{r:prespr}, for $(\s M,\s N)=(\s M,(\s N_{ij})_{ij\in \{0,1\}})\in \textnormal{SO}_{c,\nu}(\s X,\s Y)$ the formula
\begin{alignat}{1}\label{eq:solop}
\begin{aligned}
 &\textnormal{sol}(\s M,\s N)
\\
   & = \Big(\check{\partial}_{t} \begin{pmatrix}
                                     \mathcal{M} & 0 \\ 0 & 0 
                                   \end{pmatrix} + \begin{pmatrix}
                                     \mathcal{N}_{00} & \mathcal{N}_{01} \\ \mathcal{N}_{10} & \mathcal{N}_{11}
                                   \end{pmatrix}\Big)^{-1}
 \\ & = 
   \begin{pmatrix}
      \mathcal{M}^{-1} & 0 \\
      -\mathcal{N}_{11}^{-1}\mathcal{N}_{10}\mathcal{M}^{-1}& \mathcal{N}_{11}^{-1}
   \end{pmatrix}\begin{pmatrix}
      \check{\partial}_{t}^{-1}& 0 \\
      0& 1
   \end{pmatrix}  
   \\ &\quad\quad+ \begin{pmatrix}
      0 & -\mathcal{M}^{-1}\check{\partial}_{t}^{-1}\mathcal{N}_{01}\mathcal{N}_{11}^{-1}  \\
      0 & \mathcal{N}_{11}^{-1}\mathcal{N}_{10}\mathcal{M}^{-1}\check{\partial}_{t}^{-1}\mathcal{N}_{01}\mathcal{N}_{11}^{-1}
   \end{pmatrix}\begin{pmatrix}
      \check{\partial}_{t}^{-1}& 0 \\
      0& 1
   \end{pmatrix} \\
   & \quad\quad+\sum_{k=1}^\infty \begin{pmatrix}
      {\s T}^k\mathcal{M}^{-1} & -{\s T}^k\mathcal{M}^{-1}\check{\partial}_{t}^{-1}\mathcal{N}_{01}\mathcal{N}_{11}^{-1} \\ -\mathcal{N}_{11}^{-1}\mathcal{N}_{10}{\s T}^k\mathcal{M}^{-1} & 
   \mathcal{N}_{11}^{-1}\mathcal{N}_{10}{\s T}^k\mathcal{M}^{-1}\check{\partial}_{t}^{-1}\mathcal{N}_{01}\mathcal{N}_{11}^{-1}
   \end{pmatrix}\begin{pmatrix}
      \check{\partial}_{t}^{-1}& 0 \\
      0& 1
   \end{pmatrix},
\end{aligned}
\end{alignat}
with $\s T=-(\check{\partial}_{t}\mathcal{M})^{-1}\s R$ as well as  $\s R=\mathcal{N}_{00}-\mathcal{N}_{01}\mathcal{N}_{11}^{-1}\mathcal{N}_{10}$ holds true. Recall that the sum converges in $L_{\textnormal{sev}}^\textnormal{n}(\s X\times \s Y)$. Indeed, this follows from Theorem \ref{t:wp_ode_evo} \eqref{eq:wp_ode_evo2}. Hence, before coming to the proofs of the Theorems \ref{t:cdoden} and \ref{t:cdodes}, we need to establish a statement shedding light on infinite sums in view of the different topologies introduced. We need the following prerequisite of general nature.

\begin{Lemma}\label{l:dc} Let $\s Z$ be a Banach space, $(\alpha_k)_{k\in \mathbb{N}}$ in $(0,\infty)$ with $\sum_{k\in \mathbb{N}}\alpha_k<\infty$. Let $(z_{k,\iota})_{\iota\in I}$ be a convergent net in $\s Z$ with limit $w_k\in \s Z$ and assume that
\[
  \|z_{k,\iota}\|\leq \alpha_k\quad(k\in \mathbb{N},\iota\in I).
\]
Then $\sum_{k=1}^\infty w_k \in \s Z$ and
\[
   \lim_{\iota\in I} \sum_{k=1}^\infty z_{k,\iota} = \sum_{k=1}^\infty w_k.
\] 
\end{Lemma}
\begin{Proof}
 First of all note that $(1/\alpha_k)z_{k,\iota} \in B_{\s Z}$ implies that $(1/\alpha_k)w_k\in B_{\s Z}$ as the unit ball $B_{\s Z}=\{ z\in \s Z; \|z\|\leq 1\}$ is closed in $\s Z$. Thus, $\sum_{k=1}^\infty w_k \in \s Z$ since $\s Z$ is a Banach space. Let $\epsilon>0$. Then there exists $K\in \mathbb{N}$ such that $\sum_{k=K+1}^\infty \alpha_k\leq \epsilon$. We find $\iota_0\in I$ such that for all $k\in \{1,\ldots,K\}$ and $\iota\geq \iota_0$ we get
 \[
    \| z_{k,\iota} - w_k\| \leq \frac{\epsilon}{K}.
 \]
 Therefore, for all $\iota\geq\iota_0$
 \begin{align*}
    \big\|\sum_{k=1}^\infty z_{k,\iota}-\sum_{k=1}^\infty w_k \big\|
    & \leq \sum_{k=1}^\infty \|z_{k,\iota}- w_k \|
    \\&\leq \sum_{k=1}^K \|z_{k,\iota}-w_k \|+\sum_{k=K+1}^\infty \|z_{k,\iota}- w_k \|
    \\&\leq K\frac{\epsilon}{K}+\sum_{k=K+1}^\infty 2\alpha_k\leq 3\epsilon.
 \end{align*}
\end{Proof}

\begin{Remark} Even for the case $\s Z=\mathbb{C}$ the assertion of Lemma \ref{l:dc} does \emph{not} follow from Lebesgue's dominated convergence theorem, as the dominated convergence theorem asserts something on sequences instead of nets. 
\end{Remark}

\begin{Proposition}\label{p:sumdc} Let $\s X$, $\s Y$ Hilbert spaces, for $k\in \mathbb{N}$ let $(\s S_{k,\iota})_{\iota\in I}$ be a net in $L_{\textnormal{sev}}(\s X,\s Y)$, $(\s T_k)_{k\in\mathbb{N}}$ a sequence in $L_{\textnormal{sev}}(\s X,\s Y)$. Assume there exists $(\alpha_k)_{k\in\mathbb{N}}$ in $(0,\infty)$ with $\sum_{k\in \mathbb{N}}\alpha_k<\infty$ satisfying
\[
   \sup_{\mu\geq \nu}\sup_{\iota\in I} \|\s S_{k,\iota}^\mu\|\leq \alpha_k\quad (k\in \mathbb{N})
\]
for some $\nu\in \mathbb{R}$.
\begin{enumerate}[label=(\alph*)]
 \item\label{sumdcw} If for all $k\in \mathbb{N}$ we have $\s S_{k,\iota} \stackrel{\iota}{\to} \s T_k$ in $L_{\textnormal{sev},\nu}^\textnormal{w}(\s X,\s Y)$ then $\sum_{k=1}^\infty \s T_k$ converges in norm and 
 \[
    \sum_{k=1}^\infty \s S_{k,\iota} \stackrel{\iota}{\to} \sum_{k=1}^\infty \s T_k\text{ in }L_{\textnormal{sev}}^\textnormal{w}(\s X,\s Y).
 \]
 \item\label{sumdcs} If $\s S_{k,\iota} \stackrel{\iota}{\to} \s T_k$ in $L_{\textnormal{sev},\nu}^\textnormal{s}(\s X,\s Y)$ for all $k\in \mathbb{N}$ then 
 \[
    \sum_{k=1}^\infty \s S_{k,\iota} \stackrel{\iota}{\to} \sum_{k=1}^\infty \s T_k\text{ in }L_{\textnormal{sev}}^\textnormal{s}(\s X,\s Y).
 \]
 \item\label{sumdcn} If $\s S_{k,\iota} \stackrel{\iota}{\to} \s T_k$ in $L_{\textnormal{sev},\nu}^\textnormal{n}(\s X,\s Y)$ for all $k\in \mathbb{N}$ then 
 \[
    \sum_{k=1}^\infty \s S_{k,\iota} \stackrel{\iota}{\to} \sum_{k=1}^\infty \s T_k\text{ in }L_{\textnormal{sev}}^\textnormal{n}(\s X,\s Y).
 \]
\end{enumerate}
\end{Proposition}
\begin{Proof} For the proof of \ref{sumdcw}, we observe that for all $\mu\geq\nu$ and $k\in \mathbb{N}$, we have $\|\s S_{k,\iota}^\mu\|\leq \alpha_k$; hence $\|\s T_k^\mu\|\leq \alpha_k$ as the closed (norm) ball of $L(\s W,\s Z)$ for Hilbert spaces $\s W$, $\s Z$ is also closed with respect to the weak operator topology. Thus, for all $\mu\geq \nu$, the sum $\sum_{k=1}^\infty \|\s T_k^\mu\|$ is finite, yielding convergence of $\sum_{k=1}^\infty \s T_k$ in $L_{\textnormal{sev},\nu}^\textnormal{n}$.

For the proof of interchanging the limits, let $\mu\geq\nu$, $\phi\in L_\mu^2(\s X)$, $\psi\in L_\mu^2(\s Y)$ and apply Lemma \ref{l:dc} to $\s Z=\mathbb{C}$, $z_{k,\iota}=\langle \psi,\s S_{k,\iota}^\mu\phi\rangle_{L_\mu^2(\s Y)}$ and $w_k=\langle\psi,\s T_k^\mu\phi\rangle_{L_\mu^2(\s Y)}$ for $k\in \mathbb{N}$, $\iota\in I$.

In order to prove \ref{sumdcs}, let $\mu\geq\nu$, $\phi\in L_\mu^2(\s X)$ and apply Lemma \ref{l:dc} to $\s Z=L_\mu^2(\s Y)$, $z_{k,\iota}=\s S_{k,\iota}^\mu\phi$ and $w_k=\s T_k^\mu\phi$ for $k\in \mathbb{N}$, $\iota\in I$.

The proof of \ref{sumdcn} is an application of Lemma \ref{l:dc} to $\s Z=L(L_\mu^2(\s X),L_\mu^2(\s Y))$, $z_{k,\iota}=\s S_{k,\iota}^\mu$, $w_k=\s T_k^\mu$ for $k\in \mathbb{N}$, $\iota\in I$, $\mu\geq \nu$. 
\end{Proof}

We conclude this section with the proofs of Theorem \ref{t:cdoden} and \ref{t:cdodes}. The respective proofs follow similar lines, so we will provide these all in once.

\begin{Proof}[of Theorems \ref{t:cdoden} and \ref{t:cdodes}] Let $((\s M_\iota,\s N_\iota))_{\iota\in I}$ be a bounded and convergent net in $\textnormal{SO}_{c,\nu}^\textnormal{n}$ (or $\textnormal{SO}_{c,\nu}^\textnormal{s}$); let $(\s O,\s P)$ be the respective limit. Note that $\s M_{\iota}^{-1}\in GL_{1/c,\textnormal{sev},\nu}(\s X)$ and $\s N_{11,\iota}^{-1}\in GL_{1/c,\textnormal{sev},\nu}(\s Y)$, see Corollary \ref{c:inv_ce}. Hence, by Theorem \ref{t:cinv}, both the nets $(\s M_{\iota}^{-1})_\iota$ and $(\s N_{11,\iota}^{-1})_\iota$ converge in $L_{\textnormal{sev},\eta}^\textnormal{n}$ (or $L_{\textnormal{sev},\eta}^\textnormal{s}$) for some $\eta\geq\nu$ to $\s O^{-1}$ and $\s P_{11}^{-1}$, respectively. Hence, using Remark \ref{r:cim_fnu}, we get for $\s T_\iota=-(\check{\partial}_{t}\mathcal{M}_\iota)^{-1}\s R_\iota$ with  $\s R_\iota=\mathcal{N}_{\iota,00}-\mathcal{N}_{\iota,01}\mathcal{N}_{\iota,11}^{-1}\mathcal{N}_{\iota,10}$, $\iota\in I$, that there exists $\eta'\geq\eta\geq\nu$ with the property that for all $k\in \mathbb{N}$
\[
   (\s T_\iota^k)_\iota\text{ converges in } L_{\textnormal{sev},\eta'}^\textnormal{n} \text{(or } L_{\textnormal{sev},\eta'}^\textnormal{s}\text{) to} \left(-(\check{\partial}_{t}\s O)^{-1}(\mathcal{P}_{00}-\mathcal{P}_{01}\mathcal{P}_{11}^{-1}\mathcal{P}_{10})\right)^k.
\]
Hence, with the help of Remark \ref{r:cim_fnu} again and Proposition \ref{p:sumdc}, we infer the assertion taking the representation of the solution operator \eqref{eq:solop} into account. 
\end{Proof}

\renewcommand{\baselinestretch}{0.65}\normalsize\mysection{Weak Topology and Ordinary Differential Equations}{The Weak Operator Topology and Ordinary Differential Equations}{`failure' of continuity $\cdot$ the set $\textnormal{SO}_{c,\nu}(\s X)$ $\cdot$ estimates for limits in the weak operator topology $\cdot$ stability of positive definiteness under weak convergence of inverses $\cdot$ existence of convergent subsequences of solution operators $\cdot$ Theorem \ref{t:cdodew} $\cdot$ Theorem \ref{t:cdodewsol}}\label{s:wot}

\renewcommand{\baselinestretch}{1}\normalsize
Due to a lack of a version of Theorem \ref{t:cim} (\ref{nn} and \ref{ss}) and of Theorem \ref{t:cinv} for the weak operator topology, the result corresponding to the Theorems \ref{t:cdoden} and \ref{t:cdodes} is more involved. In consequence, the limiting equation is more involved. In applications, this is reflected in so-called `memory effects' occurring after `homogenization' of ordinary differential equations. For relations to homogenization problems and a thorough discussion for the occurrence of memory effects and the relationship to Young-measures, we refer the reader to the discussion in \cite{Waurick2014JAA_G} and to the comments at the end of this chapter.

We shall describe the question to be answered in the following in a bit more detail:

Let $((\s M_\iota,\s N_\iota))_\iota$ be a net in $\textnormal{SO}_{c,\nu}(\s X,\s Y)$. Are there conditions on (some type of convergence of) $((\s M_\iota,\s N_\iota))_\iota$ in terms of the weak operator topology such that $(\textnormal{sol}(\s M_\iota,\s N_{\iota}))_\iota$ converges in the weak operator topology to $\textnormal{sol}(\tilde{\s M},\tilde{\s N})$ for suitable $\tilde{\s M},\tilde{\s N}$? The first part of this section will be concerned with the first part of the question, the second part gives a description of $\tilde{\s M},\tilde{\s N}$, or, expressed differently, the second part provides the limit equation for a special case.

\begin{Theorem}\label{t:cdodew} Let $\s X$, $\s Y$ be Hilbert spaces, $c,\nu>0$, $((\s M_\iota,\s N_\iota))_\iota$ a bounded net in $\textnormal{SO}_{c,\nu}(\s X,\s Y)$. Assume that there exists $\eta\in \mathbb{R}$ such that
  \begin{align*}
   & (\mathcal{M}_\iota^{-1})_\iota,\quad\quad  (\mathcal{N}_{\iota,11}^{-1}\mathcal{N}_{\iota,10}\mathcal{M}_\iota^{-1})_\iota, \quad\quad (\mathcal{N}_{\iota,11}^{-1})_\iota, \quad\quad (\mathcal{M}_\iota^{-1}\check{\partial}_{t}^{-1}\mathcal{N}_{\iota,01}\mathcal{N}_{\iota,11}^{-1})_\iota, \\
  & (\mathcal{N}_{\iota,11}^{-1}\mathcal{N}_{\iota,10}\mathcal{M}_\iota^{-1}\check{\partial}_{t}^{-1}\mathcal{N}_{\iota,01}\mathcal{N}_{\iota,11}^{-1})_\iota, \quad\quad 
        ({\s T_\iota}^k\mathcal{M}_\iota^{-1})_\iota, \quad\quad ({\s T_\iota}^k\mathcal{M}_\iota^{-1}\check{\partial}_{t}^{-1}\mathcal{N}_{\iota,01}\mathcal{N}_{\iota,11}^{-1})_\iota, \\ & (\mathcal{N}_{\iota,11}^{-1}\mathcal{N}_{\iota,10}{\s T_\iota}^k\mathcal{M}_\iota^{-1})_\iota,
   \quad\quad (\mathcal{N}_{\iota,11}^{-1}\mathcal{N}_{\iota,10}{\s T_\iota}^k\mathcal{M}_\iota^{-1}\check{\partial}_{t}^{-1}\mathcal{N}_{\iota,01}\mathcal{N}_{\iota,11}^{-1})_\iota\quad\quad(k\in \mathbb{N}),
   \end{align*}
   converge in $L_{\textnormal{sev},\eta}^\textnormal{w}$, where $\s T_\iota=-(\check{\partial}_t\s M_\iota)^{-1}\s R_\iota$ and $\s R_\iota=\s N_{\iota,00}-\s N_{\iota,01}\s N_{\iota,11}^{-1}\s N_{\iota,10}^{-1}$\; for all $\iota$.
   
   Then $(\textnormal{sol}(\s M_\iota,\s N_\iota))_\iota$  converges in $L_{\textnormal{sev}}^\textnormal{w}(\s X\times \s Y)$ with $\textnormal{sol}$ being defined in \eqref{eq:ode_solm}.
\end{Theorem}
\begin{Proof}
 We use the representation of $\textnormal{sol}(\s M_\iota,\s N_\iota)$ as given in \eqref{eq:ode_solm}, see also \eqref{eq:solop}.
 Multiplication is separately continuous in the weak operator topology (Theorem \ref{t:cim}\ref{ww}). Next, interchanging the limits and summation is possible by Proposition \ref{p:sumdc}, which settles the assertion.
\end{Proof}

Next, we account for the computation of the limiting equation, that is, we seek to compute $\s M_\infty,\s N_\infty$ such that $\lim_\iota\textnormal{sol}(\s M_\iota,\s N_\iota)=\textnormal{sol}(\s M_\infty,\s N_\infty)$. This has been addressed for the case of (convergent) sequences already in \cite{Waurick2014JAA_G}. Note that without any further assumptions on $\s M_\iota$ and $\s N_\iota$ such as time translation-invariance, the representation of the limit equation given in \cite{Waurick2014JAA_G} does not fit into the representation $\textnormal{sol}(\s M_\infty,\s N_\infty)$ for some bounded evolutionary $(\s M_\infty,\s N_\infty)$. As we shall see later on, this is not true for the case $\s Y=\{0\}$. Hence, we focus on the case $\s Y=\{0\}$ in the following. For studying the respective solution operator, we introduce
\[
  \textnormal{SO}_{c,\nu}(\s X)
  \coloneqq \{ (\s M,\s N)\in L_{\textnormal{sev}}(\s X)^2; \s M \text{ satisfies }\eqref{eq:wpcondm}\}\subseteq L_{\textnormal{sev}}(\s X\times \s X).
\]
as well as
\begin{alignat}{1}\label{eq:solopd}
\begin{aligned}
   \textnormal{sol}\colon \textnormal{SO}_{c,\nu}(\s X)&\to L_{\textnormal{sev}}(\s X)
   \\ (\s M,\s N) & \mapsto \left(\check{\partial}_{t}\s M+\s N\right)^{-1}.
\end{aligned}
\end{alignat}
By equation \eqref{eq:solop}, we get for all $(\s M,\s N)\in \textnormal{SO}_{c,\nu}(\s X)$:
\begin{equation}\label{eq:solop2}
\textnormal{sol}(\s M,\s N) = \mathcal{M}^{-1}\check{\partial}_{t}^{-1}+\sum_{k=1}^\infty 
      (-\mathcal{M}^{-1}\check{\partial}_{t}^{-1}\mathcal{N})^k\mathcal{M}^{-1} \check{\partial}_{t}^{-1}.
\end{equation}
Hence, the corresponding version of the statement in Theorem \ref{t:cdodew} for $\s Y=\{0\}$ reads:
\begin{Corollary}\label{c:cdodew2} Let $\s X$ be a Hilbert space, $c,\nu>0$, $((\s M_\iota,\s N_\iota))_\iota$ a bounded net in $\textnormal{SO}_{c,\nu}(\s X)$. Assume that there exists $\eta\in \mathbb{R}$ such that
  \[
  (\mathcal{M}_\iota^{-1})_\iota\text{ and } 
        (\big(-(\check{\partial}_t\s M_\iota)^{-1}\s N_\iota\big)^k\mathcal{M}_\iota^{-1})_\iota\quad(k\in \mathbb{N})\]
   converge in $L_{\textnormal{sev},\eta}^\textnormal{w}$.
   
   Then $(\textnormal{sol}(\s M_\iota,\s N_\iota))_\iota$  converges in $L_{\textnormal{sev}}^\textnormal{w}(\s X)$ with $\textnormal{sol}$ being defined in \eqref{eq:solopd}. 
\end{Corollary}

For computing the limiting equation in the case just discussed in Corollary \ref{c:cdodew2}, we use the representation of the solution operator $\textnormal{sol}(\s M,\s N)$ in \eqref{eq:solop2}, the estimates given in Remark \ref{r:prest} and a Neumann series argument. First of all, we state a closedness result for the weak operator topology, leading to estimates, which will come in handy for applying the Neumann series argument.

\begin{Proposition}\label{p:bdd_cld_weak} Let $\s X,\s Y$ Hilbert spaces, $\nu\in \mathbb{R}$. For $r>0$ let
\[
  \mathfrak{C}_\nu(r)\coloneqq \{ \s S\in L_{\textnormal{sev},\nu}(\s X,\s Y); \sup_{\mu\geq\nu} \|\s S^\mu\|\leq r\}.
\] Then $\mathfrak {C}_\nu(r)\subseteq L_{\textnormal{sev}}^{\textnormal{w}}(\s X,\s Y)$ is closed.
\end{Proposition}
\begin{Proof} An application of Lemma \ref{l:bhssev} shows that $\mathfrak {C}_\nu(r)\subseteq L_{\textnormal{sev}}^{\textnormal{w}}(\s X,\s Y)$ is closed if and only if $\mathfrak {C}_\nu(r)\subseteq L_{\textnormal{sev},\nu}^{\textnormal{w}}(\s X,\s Y)$ is closed. The latter, however, is easily seen using that the unit ball of $L(L_\mu^2(\s X);L_\mu^2(\s Y))$ is closed under the weak operator topology for all $\mu\geq\nu$.
\end{Proof}

\begin{Remark}\label{r:up} In the sequel, Proposition \ref{p:bdd_cld_weak} may be applied as follows. Let $(\s S_\iota)_\iota$ be a convergent net in $L_{\textnormal{sev}}^{\textnormal{w}}(\s X,\s Y)$ with the property that for all $\iota$ we have $\s S_\iota\in \mathfrak{C}_\nu(r)$ for some $r,\nu>0$. Then $\lim_\iota \s S_\iota\in \mathfrak{C}_\nu(r)$. 
\end{Remark}

The next proposition gives a more precise estimate of the limit operator in relation to its net converging to it. For this, we define for a net of non-negative real numbers $(s_\iota)_\iota$ the number
\[
  \liminf_\iota s_\iota \coloneqq \inf\{ t\in [0,\infty]; t \text{ accumulation value of }(s_\iota)_\iota\}.
\]
\begin{Proposition}\label{p:liminfine} Let $\s X,\s Y$ Hilbert spaces, $(S_\iota)_\iota$ be a convergent net in $L(\s X,\s Y)$ with respect to the weak operator topology. Then
\[
   \|\lim_\iota  S_\iota \|\leq \liminf_\iota \|S_\iota\|.
\] 
\end{Proposition}
\begin{Proof}
  Put $T\coloneqq \lim_\iota S_\iota$. For $\phi\in B_{\s X},\psi\in B_{\s Y}$, that is, $\phi\in \s X,\psi\in \s Y$ with $\|\phi\|=\|\psi\|=1$,  we get for $\iota$
  \[
     |\langle \psi, S_\iota \phi\rangle|\leq \|S_\iota\|.
  \]
  Hence, as the only accumulation value of $(|\langle \psi, S_\iota \phi\rangle|)_\iota$ is $|\langle \psi, T \phi\rangle|$, we get 
  \[
   \|T\|=\sup_{\phi\in B_\s X,\psi\in B_\s Y}|\langle \psi, T \phi\rangle|\leq \liminf_\iota \|S_\iota\|.
  \]
\end{Proof}

 For having a solution theory for the limiting equation as well, we need to warrant the conditions imposed in Theorem \ref{t:wp_ode_evo} for $\s M$ and $\s N$ respectively interchanged by $\s M_\infty$ and $\s N_\infty$. As a part of this, we will need to study the inverse of the limit of $(\s M_\iota^{-1})_\iota$ in more detail.  In particular, we want the respective inverse to satisfy an estimate analogous to \eqref{eq:wpcondm} for some suitable $c>0$. This will be addressed next.
 
 \begin{Theorem}\label{t:invcau} Let $\s X$ be a Hilbert space, $c,\nu>0$, $(\s M_\iota)_{\iota\in I}$ a bounded net in $L_{\textnormal{sev},\nu}(\s X)$. Assume 
 \begin{equation}\label{eq:invcau1}
      \Re\langle Q_t \s M_\iota \phi,\phi\rangle_{L_\mu^2}\geq c\langle Q_t \phi,\phi\rangle_{L_\mu^2}\quad (\phi\in \s D(\s X), t\in \mathbb{R}, \mu\geq\nu, \iota\in I).
 \end{equation}
 If $(\s M_\iota^{-1})_\iota$ converges to some $\s O\in L_{\textnormal{sev}}^\textnormal{w}(\s X)$, then $\s O\in L_{\textnormal{sev},\nu}(\s X)$, $\sup_{\mu\geq\nu}\|\s O^\mu\|\leq 1/c$. Moreover, $\s M_\infty\coloneqq \s O^{-1} \in L_{\textnormal{sev},\nu}(\s X)$ with 
 \begin{equation}\label{eq:invcau0}
   \sup_{\mu\geq \nu} \|\s M_\infty^\mu\|\leq \frac{1}{c}\sup_{\mu\geq\nu,\iota\in I}\|\s M_\iota^\mu\|^2 
 \end{equation}
and
 \begin{equation}\label{eq:invcau2}
     \Re\langle Q_t \s M_\infty \phi,\phi\rangle_{L_\mu^2}\geq c\langle Q_t \phi,\phi\rangle_{L_\mu^2}\quad (\phi\in \s D(\s X), t\in \mathbb{R}, \mu\geq\nu),
 \end{equation}
 where $Q_t$ is multiplication by $\1_{(-\infty,t)}$.
 \end{Theorem}
 \begin{Proof} From \eqref{eq:invcau1}, we get for all $\iota\in I$ that $\s M_\iota^{-1}\in L_{\textnormal{sev},\nu}(\s X)$ with $\|(\s M_\iota^{-1})^\mu\|\leq 1/c$ for all $\mu\geq\nu$, by Corollary \ref{c:inv_ce} and Proposition \ref{p:posd}. In particular, by Theorem \ref{t:csc}, $Q_t(\s M_\iota^{-1})^\mu=Q_t(\s M_\iota^{-1})^\mu Q_t$, $t\in \mathbb{R}$, $\mu\geq\nu$. Further, $\s O\in \mathfrak{C}_\nu(1/c)$, by Proposition \ref{p:bdd_cld_weak}. Moreover, with $\psi=\s M_\iota \phi$, from 
 \[
  \Re\langle Q_t\s M_\iota \phi,\phi\rangle_{L_\mu^2}\geq c\langle Q_t\phi,\phi\rangle_{L_\mu^2}\quad (\phi\in \s D(\s X), t\in \mathbb{R}, \mu\geq\nu, \iota\in I),
 \]
 we read off
 \begin{align*}
  & \Re\langle Q_t \psi,\s M_\iota^{-1}\psi\rangle_{L_\mu^2} \geq  c\Re\langle Q_t \s M_\iota^{-1}\psi,\s M_\iota^{-1} \psi\rangle_{L_\mu^2}
  \\& = c\langle Q_t(\s M_\iota^{-1})^\mu Q_t\psi,Q_t(\s M_\iota^{-1})^\mu Q_t\psi\rangle_{L_\mu^2} \geq \frac{c}{\|\s M_\iota^\mu\|^2}\langle Q_t\psi,Q_t\psi\rangle_{L_\mu^2}
  \\&\geq \frac{c}{\sup_{\mu\geq\nu,\iota\in I}\|\s M_\iota^\mu\|^2}\langle Q_t\psi,\psi\rangle_{L_\mu^2}  \quad (\psi\in \s D(\s X), \mu\geq\nu, \iota\in I),
 \end{align*}
 where we used \[\|Q_t\psi\|=\|Q_t \s M_\iota^\mu(\s M_\iota^{-1})^\mu \psi\|=\|Q_t\s M_\iota^\mu Q_t (\s M_\iota^{-1})^\mu Q_t\psi\|\leq \|\s M_\iota^\mu\|\|Q_t(\s M_\iota^{-1})^\mu Q_t\psi\|.\] Hence, we obtain $\s M_\infty=\s O^{-1} \in L_{\textnormal{sev},\nu}(\s X)$ and \eqref{eq:invcau0}, by Corollary \ref{c:inv_ce} and Proposition \ref{p:posd}. Next, using \eqref{eq:invcau1} again, we get for all $t\in \R$, $\psi\in \s D(\s X)$, $\iota\in I$, $\mu\geq\nu$, 
 \[
    \Re\langle Q_t \psi, \s M_\iota^{-1}\psi\rangle_{L_\mu^2}\geq c\langle Q_t \s M_\iota^{-1}\psi,Q_t\s M_\iota^{-1}\psi\rangle_{L_\mu^2}.
 \]
 Hence, computing the limit in $\iota$ and using Proposition \ref{p:liminfine} for $\s Y=\mathbb{C}$ and $S_\iota=Q_t \s M_\iota^{-1}\psi$, we arrive at
 \[
    \Re\langle Q_t \psi, \s O\psi\rangle_{L_\mu^2}\geq c\liminf_\iota \|Q_t \s M_\iota^{-1}\psi\|_{L_\mu^2}^2\geq c\|Q_t\s O\psi\|_{L_\mu^2}^2.
  \]
  Substituting $\s O\psi =\phi$ we arrive at
  \[
    \Re\langle Q_t \s M_\infty \phi, \phi\rangle_{L_\mu^2}\geq c\|Q_t\phi\|_{L_\mu^2}^2\quad(\mu\geq\nu).
  \] 
 \end{Proof}
 We turn back to the derivation of the limit of $(\textnormal{sol}(\s M_\iota,\s N_\iota))_\iota$ in Corollary \ref{c:cdodew2}. So, assume the hypothesis of Corollary \ref{c:cdodew2} to be in effect and denote for all $k\in \mathbb{N}$:
\begin{alignat}{1}\label{eq:defwlimit}
\begin{aligned}
   \s O &\coloneqq \lim_\iota \mathcal{M}_\iota^{-1},
  \\  \s P_k & \coloneqq \lim_\iota\big(-(\check{\partial}_t\s M_\iota)^{-1}\s N_\iota\big)^k\mathcal{M}_\iota^{-1}
\end{aligned}
\end{alignat}
With these definitions at hand, we can describe the limit of $(\textnormal{sol}(\s M_\iota,\s N_\iota))_\iota$. Thus, the following theorem may be read as a sequel to Corollary \ref{c:cdodew2}.

\begin{Theorem}\label{t:cdodewsol} Assume all conditions in Corollary \ref{c:cdodew2} to be in effect. Then, with the help of the definitions \eqref{eq:defwlimit}, we get
\[
   \textnormal{sol}(\s M_\iota,\s N_\iota)\stackrel{\iota}{\to}\textnormal{sol}(\s M_\infty,0)
\]
in $L_{\textnormal{sev}}^\textnormal{w}(\s X)$, where
\begin{align*}
 \s M_\infty &= \s O^{-1}+\s O^{-1}\sum_{\ell=1}^\infty \s R^\ell,\\
 \s R & = -\sum_{k=1}^\infty \s P_k \s O^{-1}.
\end{align*}
\end{Theorem}
\begin{Proof} Without loss of generality we will assume $c<1$ in the following.
Using the representation in \eqref{eq:solop2} and recalling the argument in Corollary \ref{c:cdodew2} (or Theorem \ref{t:cdodew}), we get 
 \[
    \textnormal{sol}(\s M_\iota,\s N_\iota)\to  \s O\check{\partial}_{t}^{-1}+\sum_{k=1}^\infty 
      \s P_k \check{\partial}_{t}^{-1}\text{ in }L_{\textnormal{sev}}^\textnormal{w}(\s X).
 \]
 In the following we will show that $\|\s R\|_{L(L_\mu^2)}<1$ for eventually all $\mu$ large enough. Moreover, in order that $\textnormal{sol}$ being well-defined, we need to show that $\s M_\infty$ satisfies the positive definiteness estimate \eqref{eq:wpcondm} for some (possibly different) $c>0$. Having shown all these statements, it is then a straightforward computation to verify that
 \[
    \left(\s O\check{\partial}_{t}^{-1}+\sum_{k=1}^\infty 
      \s P_k \check{\partial}_{t}^{-1}\right)^{-1}=\check{\partial}_t\s M_\infty.
 \] 
 We address the norm estimate for $\s R$ first. From Remark \ref{r:prest}, we get for all sufficiently large $\eta\geq 0$, $\iota\in I$
 \[
    \|\textnormal{sol}(\s M_\iota,\s N_\iota)\check{\partial}_{t}-\s M_\iota^{-1}\|_{L(L_\eta^2(\s X))}\leq \frac{\theta}{c(1-\theta)},
 \]
 where $\theta=\sup_{\iota\in I}\|\s N_\iota\|_{\geq\eta}/(\nu c)$ and $\|\s N_\iota\|_{\geq\eta}\coloneqq \sup_{\mu\geq\eta}\|\s N_\iota^\mu\|$. Hence, for the limit in $\iota$, we obtain with Remark \ref{r:up}
 \[
   \Big\| \sum_{k=1}^\infty \s P_k \Big\|_{L(L_\eta^2)}\leq \frac{\theta}{c(1-\theta)}.
 \]
 So, there exists $\eta\geq\nu$ such that for every $\mu\geq\eta$, we obtain 
 \[
    \|\s R\|_{L(L_\mu^2)}=\Big\| \sum_{k=1}^\infty \s P_k \s O^{-1} \Big\|_{L(L_\mu^2)}\leq \frac{c}{2} (<1).
 \]
 We are left with showing the positive definiteness type estimate \eqref{eq:wpcondm} for $\s M_\infty$. For this, we compute
 \[
    \|\sum_{\ell=1}^\infty \s R^\ell\|_{L(L_\mu^2)}\leq \sum_{\ell=1}^\infty \big(\frac{c}{2}\big)^\ell=\frac{\frac{c}{2}}{1-\frac{c}{2}}=\frac{c}{2-c}<c \quad(\mu\geq\eta).
 \]
 Next, $\s T\coloneqq \sum_{\ell=1}^\infty \s R^\ell$ is a norm convergent limit of causal and evolutionary operators, and, thus, $\s T$ is causal and evolutionary itself. In particular, for $Q_t$ denoting multiplication by $\1_{(-\infty,t)}$ and using Theorem \ref{t:csc}, we get for every $\phi\in \s D(\s X)$ and $\mu\geq\nu$ 
 \[
    \|Q_t\s T\phi\|_{L_\mu^2}=\|Q_t\s TQ_t\phi\|_{L_\mu^2}\leq \|\s T\|_{L(L_\mu^2)}\|Q_t\phi\|_{L_\mu^2}\leq \frac{c}{2-c}\|Q_t\phi\|_{L_\mu^2}\quad(t\in \mathbb{R}).
 \]
 Next, by Theorem \ref{t:invcau} applied to $\s O^{-1}$ together with the estimate just derived, we obtain for all $t\in \mathbb{R}$, $\phi\in \s D(\s X)$ and $\mu\geq\eta$
 \begin{align*}
    \Re\langle Q_t \s M_\infty\phi,\phi\rangle_{L_\mu^2}& =\Re\langle Q_t (\s O^{-1}+\s T)\phi,\phi\rangle_{L_\mu^2}
    \\ &=\Re\langle Q_t \s O^{-1}\phi,\phi\rangle_{L_\mu^2} +\Re\langle Q_t\s T\phi,Q_t\phi\rangle_{L_\mu^2}
    \\ &\geq c\langle Q_t \phi,\phi\rangle_{L_\mu^2} - \frac{c}{2-c}\|Q_t\phi\|_{L_\mu^2}^2
    \\ & = \frac{(1-c)c}{2-c}\langle Q_t \phi,\phi\rangle_{L_\mu^2}.
 \end{align*}
\end{Proof}

\begin{Remark}\label{r:comcase} In the statement of Theorem \ref{t:cdodewsol} it is somewhat awkward that the solution operator converges to $\textnormal{sol}(\s M_\infty,\s N_\infty)$ with $\s N_\infty=0$. In fact, if $\partial_t^{-1}$ commuted with $\s M_\iota$ and $\s N_\iota$ as -- for instance -- in the case of time translation-invariant coefficients (see Section \ref{s:Art}), we have the more natural statement that 
\[
   \textnormal{sol}(\s M_\iota,\s N_\iota)\stackrel{\iota}{\to} \textnormal{sol}(\s M_\infty,\s N_\infty)
\]
with $\s M_\infty= \s O^{-1}$ and $\s N_\infty=\s O^{-1}\sum_{\ell=1}^\infty \tilde{\s R}^\ell$ with $\tilde{\s R} = -\sum_{k=1}^\infty \tilde{\s P}_k \s O^{-1}$, where
\begin{alignat}{1}\label{eq:defwlimit2}
\begin{aligned}
   \s O &\coloneqq \lim_\iota \mathcal{M}_\iota^{-1},
  \\  \tilde{\s P}_k & \coloneqq (\check{\partial}_t^{-1})^{k-1}\lim_\iota\big((-(\s M_\iota)^{-1}\s N_\iota\big)^k\mathcal{M}_\iota^{-1}.
\end{aligned}
\end{alignat}
 In fact, this representation is indeed more natural, as $(\s M_\infty,\s N_\infty)\in \textnormal{SO}_{c,\nu}(\s X)$, by Theorem \ref{t:invcau} (see \eqref{eq:invcau2} in particular), whereas the positive definiteness constant for $\s M_\infty$ given in Theorem \ref{t:cdodewsol} has to be adjusted (cf.~the concluding lines of the proof of Theorem \ref{t:cdodewsol}).
\end{Remark}
 
 If the Hilbert space $\s X$ is separable, a combination of the compactness and metrizability result for the weak operator topology, that is, Theorem \ref{t:wco} and Corollary \ref{c:wmet}, immediately yields the following statement for sequences of coefficients instead of nets.
 
 \begin{Theorem}\label{t:cdsep} Let $\s X$, $\s Y$ be separable Hilbert spaces, $c,\nu>0$, $((\s M_n,\s N_n))_n$ a bounded sequence in $\textnormal{SO}_{c,\nu}(\s X,\s Y)$. 
   
   Then there exists a strictly increasing sequence $(n_k)_k$ in $\mathbb{N}$ such that $(\textnormal{sol}(\s M_{n_k},\s N_{n_k}))_k$ converges in $L_{\textnormal{sev}}^\textnormal{w}(\s X\times \s Y)$ with $\textnormal{sol}$ being defined in \eqref{eq:ode_solm}.  
 \end{Theorem}
 \begin{Proof}
   It suffices to observe that, by relative sequential compactness of bounded sets in $L_{\textnormal{sev}}^\textnormal{w}$ (combine Theorem \ref{t:wco} and Corollary \ref{c:wmet}), we may choose a strictly increasing sequence $(n_k)_k$ in $\mathbb{N}$ such that
  \begin{align*}
   & (\mathcal{M}_{n_k}^{-1})_k,\quad  (\mathcal{N}_{{n_k},11}^{-1}\mathcal{N}_{{n_k},10}\mathcal{M}_{n_k}^{-1})_k, \quad (\mathcal{N}_{{n_k},11}^{-1})_k, \quad (\mathcal{M}_{n_k}^{-1}\check{\partial}_{t}^{-1}\mathcal{N}_{{n_k},01}\mathcal{N}_{{n_k},11}^{-1})_k, \\
  & (\mathcal{N}_{{n_k},11}^{-1}\mathcal{N}_{{n_k},10}\mathcal{M}_{n_k}^{-1}\check{\partial}_{t}^{-1}\mathcal{N}_{{n_k},01}\mathcal{N}_{{n_k},11}^{-1})_k, \quad
        (\s T_{n_k}^\ell\mathcal{M}_{n_k}^{-1})_k, \\ & (\s T_{n_k}^\ell\mathcal{M}_{n_k}^{-1}\check{\partial}_{t}^{-1}\mathcal{N}_{{n_k},01}\mathcal{N}_{{n_k},11}^{-1})_k, \quad (\mathcal{N}_{{n_k},11}^{-1}\mathcal{N}_{{n_k},10}\s T_{n_k}^\ell\mathcal{M}_{n_k}^{-1})_k, \\
   & (\mathcal{N}_{{n_k},11}^{-1}\mathcal{N}_{{n_k},10}\s T_{n_k}^\ell\mathcal{M}_{n_k}^{-1}\check{\partial}_{t}^{-1}\mathcal{N}_{{n_k},01}\mathcal{N}_{{n_k},11}^{-1})_k\quad(\ell\in \mathbb{N}),
   \end{align*}
   converge in $L_{\textnormal{sev},\eta}^\textnormal{w}$, $\s T_{n_k}=-(\check{\partial}_t\s M_{n_k})^{-1}\s R_{n_k}$ and $\s R_{n_k}=\s N_{{n_k},00}-\s N_{{n_k},01}\s N_{{n_k},11}^{-1}\s N_{{n_k},10}$\; for all $k\in \mathbb{N}$. Hence, the assertions follows from Theorem \ref{t:cdodew}.
 \end{Proof}

\renewcommand{\baselinestretch}{0.65}\normalsize\mysection{Drude--Born--Fedorov Model}{The Drude--Born--Fedorov Model}{Maxwell's equations $\cdot$ $\curl$ $\cdot$ admissible domain $\cdot$ convergence of multiplication operators and the strong operator topology $\cdot$ Theorem \ref{t:dbfns} $\cdot$ Theorem \ref{t:dbfw}}\label{s:DBF}

\renewcommand{\baselinestretch}{1}\normalsize
As a first major application of the theory developed so far, we will treat the Drude--Born--Fedorov model for electromagnetism. As it has been found in \cite{Picard2013a}, this formulation of Maxwell's equations may be written as an \emph{ordinary} differential equation in an infinite-dimensional Hilbert space. Hence, in the present chapter, it is of interest to apply the results of the preceding sections to this particular example. Before, however, applying the abstract theory of evolutionary equations, we need to frame the Drude--Born--Fedorov model into a proper functional analytic setting. We note that the results in  \cite{Picard2013a} on the well-posedness for the Drude--Born--Fedorov model apply to a more general situation than the one discussed in the present exposition.

Throughout this section, let $\Omega\subseteq \mathbb{R}^3$ be open. Formally, the equations may be written as follows.
\begin{alignat}{1}\label{eq:dbf0}
\begin{aligned}
   \partial_t (1+\beta\curl)\eps E - \curl H &= J
   \\ \partial_t (1+\beta\curl)\mu H + \curl E& = 0,
\end{aligned}
\end{alignat}
on $[0,\infty)\times \Omega$, where for simplicity, we assume homogeneous initial conditions. Some comments on the constituents of \eqref{eq:dbf0} are in order. The unknowns of \eqref{eq:dbf0} are the two components of the electromagnetic vector field $(E,H)\colon \mathbb{R}_{\geq 0}\times \Omega\to \R^3\times \R^3$. The mapping $J\colon \mathbb{R}_{\geq 0}\times \Omega\to \R^3$ models the source term, that is, external electric currents. The $3\times 3$-matrix valued functions $\eps$ and $\mu$ defined on $\mathbb{R}_{\geq0}\times \Omega$ are the dielectricity and the magnetic permeability of the underlying medium, $\beta$ is a non-zero real number. The expression $\curl$ is the differential operator acting on the spatial variables of $E$ and $H$ only, which is formally given by
\[
 \phi  \mapsto\left(\begin{array}{ccc}
0 & -\partial_{3} & \partial_{2}\\
\partial_{3} & 0 & -\partial_{1}\\
-\partial_{2} & \partial_{1} & 0
\end{array}\right)\phi,
\]
for any smooth vector field $\phi\colon \Omega\to \mathbb{R}^3$. We will use $\curl \phi$ also for $L^2$-vector fields $\phi$ in the distributional sense.

We will need the following assumptions on the ingredients of \eqref{eq:dbf0}.
\begin{Definition}\label{d:om} We say that $\Omega$ is an \emph{admissible domain}, if there exists 
\[
   D \subseteq \dom(\curl)\coloneqq \{ \phi\in L^2(\Omega)^3; \curl \phi \in L^2(\Omega)^3\}
\]
such that the operator
\[
   \curl_{\diamond} \colon D \subseteq L^2(\Omega)^3\to L^2(\Omega)^3, \phi\mapsto \curl \phi
\]
is self-adjoint and $-1/\beta \in \rho(\curl_\diamond)$. The operator $\curl_\diamond$ is called an \emph{admissible realization (of $\curl$)}.
\end{Definition}

We shall elaborate on the relationship of admissible domains to the Dru\-de--Born--Fe\-dorov model as follows.

\begin{Remark}\label{r:filpi} A particular realization of the $\curl$-operator can be found in \cite{Picard1998,Picard1998a}. The operator $\curl$ with this boundary condition is used for the description of the Drude--Born--Fedorov model, see \cite{Picard2013a}. It can be shown that this realization is selfadjoint provided $\Omega$ has finite Lebesgue measure, see \cite{Filonov2000}. Another consequence is that the spectrum of this particular realization of $\curl$ is countable. In particular, for bounded open sets $\Omega$ there are uncountably many $\beta\in \mathbb{R}$ such that $\Omega$ is an admissible domain. We shall also refer to \cite{Picard1998,Picard1998a} for a corresponding treatment of unbounded $\Omega$. 
\end{Remark}

We emphasize that the precise selfadjoint realization $\curl_\diamond$ of $\curl$ is not important for the analysis to follow as long as $-1/\beta\in \rho(\curl_\diamond)$. This observation together with Remark \ref{r:filpi} leads to a generalized treatment of the Drude--Born--Fedorov model.

\begin{Hypothesis}[on $\eps$ and $\mu$]\label{h:cof} We assume there exists $\hat{\eps},\hat{\mu}\in L^\infty(\R\times\Omega)^{3\times 3}$ with \[\hat{\eps}|_{[0,\infty)\times \Omega}=\eps, \quad \hat{\mu}|_{[0,\infty)\times \Omega}=\mu\] and $c>0$ such that 
\begin{equation}\label{eq:hcof}
   \Re \langle \hat{\eps}(t,x)\xi,\xi\rangle\geq c\langle\xi,\xi\rangle, \quad \Re \langle \hat{\mu}(t,x)\xi,\xi\rangle\geq c\langle\xi,\xi\rangle
\end{equation}
for almost every $(t,x)\in \mathbb{R}\times \Omega$.
\end{Hypothesis}

\begin{Remark}\label{r:inex} In the following, we will not distinguish between $\eps$, $\mu$ and its respective extensions $\hat{\eps}$, $\hat{\mu}$ to the whole real line. Thanks to causality, we will see that a solution to \eqref{eq:dbf0} will be \emph{independent} of the extension of the coefficients to the negative reals. 
\end{Remark}

In what follows, in order to ease readability considerably, we will identify $\eps$ with the corresponding multiplication operator on $L^2_\nu(\R;L^2(\Omega)^3)$ for all $\nu\in \mathbb{R}$. Likewise, we shall do so for $\mu$. Moreover, we identify $\curl_\diamond$ with its lifting to $L^2_\nu(\R;L^2(\Omega)^3)$ as an (abstract) multiplication operator with $L^2_\nu(\R;\dom(\curl_\diamond))$ as domain  for all $\nu\in \mathbb{R}$, see also Example \ref{ex:cev}.

The well-posedness theorem corresponding to \eqref{eq:dbf0} reads as follows.
\begin{Theorem}\label{t:wpdbf} Let $\Omega$ be an admissible domain and $\curl_\diamond$ an admissible realization. Assume that $\eps$ and $\mu$ satisfy Hypothesis \ref{h:cof}. Then
\begin{multline*}
   \s S_{\textnormal{DBF}}(\eps,\mu)\coloneqq\\ \Big(\partial_{t}\begin{pmatrix}
                 \eps & 0 \\ 0 & \mu 
               \end{pmatrix} + \begin{pmatrix} 0 & -\curl_\diamond \\ \curl_\diamond & 0 \end{pmatrix}(1+\beta\curl_\diamond)^{-1}\Big)^{-1}\in L_{\textnormal{sev}}(L^2(\Omega)^6)
\end{multline*}
\end{Theorem}
\begin{Proof}
  We want to apply Theorem \ref{t:wp_ode_evo}. More precisely, take \begin{multline*}\s X=L^2(\Omega)^6, \s Y=\{0\}, \s M=\begin{pmatrix}
                 \eps & 0 \\ 0 & \mu 
               \end{pmatrix}\text{ and }\\ \s N=\begin{pmatrix} 0 & -\curl_\diamond \\ \curl_\diamond & 0 \end{pmatrix}(1+\beta\curl_\diamond)^{-1}.\end{multline*} Hypothesis \ref{h:cof} guarantees inequality \eqref{eq:wpcondm}. Indeed, this follows from Hypothesis \ref{h:cof}, the fact that $Q_t$ commutes with $\s M$ and Example \ref{ex:mult_op}. Example \ref{ex:mult_op} also yields (standard) evolutionarity of $\s N$, since, as $\curl_\diamond$ is an admissible realization, the operator \[\curl_\diamond(1+\beta\curl_\diamond)^{-1}\] considered in $L^2(\Omega)^6$ is bounded.
\end{Proof}

\begin{Remark}\label{r:posreal}
  The solution theory obtained in Theorem \ref{t:wpdbf} says that for any $\tilde{J}\in L_{\nu}^2([0,\infty);L^2(\Omega)^3)$ (for eventually all large enough $\nu$), the equations \eqref{eq:dbf0} admit a unique solution $(E,H)\in L_{\nu}^2(\R;L^2(\Omega)^6)$ in the sense that
\begin{alignat}{1}\label{eq:dbf1}
\begin{aligned}
   \partial_t \eps E - \curl_\diamond(1+\beta\curl_\diamond)^{-1} H &= \tilde{J}
   \\ \partial_t \mu H + \curl_\diamond (1+\beta\curl_\diamond)^{-1}E& = 0,
\end{aligned}	
\end{alignat}  
holds, where $\curl_\diamond$ is an admissible realization of $\curl$. Moreover, note that, as $\tilde{J}$ is supported on $[0,\infty)$ only, so is the electromagnetic field $(E,H)$, by causality. In particular, equation \eqref{eq:dbf1} is trivial on $(-\infty,0)$, making the solution operator independent of the chosen extensions $\hat{\eps}$ and $\hat{\mu}$ in Hypothesis \ref{h:cof}.
\end{Remark}

Next, we will apply the results on the continuous dependence to the Drude--Born--Fedorov model. We will need a prerequisite, which is particularly useful for the result corresponding to the strong operator topology. 

\begin{Proposition}\label{p:stronmult} Let $p,d\in\mathbb{N}$, $(\eps_n)_n$ a sequence in $L^\infty(\R\times \Sigma)^{p\times p}$ for some measurable $\Sigma\subseteq \mathbb{R}^d$, $T\in L(L^2(\R\times \Sigma)^p)$. Assume that $T_{\eps_n}\to T$ as $n\to\infty$ in the strong operator topology of $L(L^2(\R\times \Sigma)^p)$, where we recall that $T_{\eps_n}$ denotes the associated multiplication operator of $\eps_n$. Then
\begin{enumerate}[label=(\alph*)]
 \item\label{sm1} $(T_{\eps_n})_n$ is bounded, there exists $\eps\in L^\infty(\R\times \Sigma)^{p\times p}$ such that $T=T_\eps$, and
 \item\label{sm2} $T_{\eps_n}\to T_\eps$ in $L_{\textnormal{sev}}^\textnormal{s}(L^2(\Sigma)^p)$.
\end{enumerate} 
\end{Proposition}
\begin{Proof}
We start out with \ref{sm1}. Being strongly convergent, the sequence $(T_{\eps_n})_n$ is bounded in $L(L^2(\R\times\Sigma)^p)$, by the uniform boundedness principle. Next, the convergence asserted implies convergence in the weak operator topology, which, in turn, for multiplication operators is easily seen to be equivalent to convergence in the weak* topology of $L^\infty$. But, by separability of $L^1(\R\times\Sigma)$, the unit ball of $L^\infty(\mathbb{R}\times \Sigma)^{p\times p}$ is sequentially compact under the weak* topology. Hence, there exists $\eps\in L^\infty(\R\times\Sigma)^{p\times p}$ being the limit of a weakly* convergent subsequence of $(\eps_n)_n$. Since (any subsequence of) $(T_{\eps_n})_n$ also converges in the weak operator topology, the strong operator 
 topology limit of $(T_{\eps_n})_n$ is induced by multiplication by $\eps$.
 
 For the proof of \ref{sm2}, we will use Theorem \ref{t:ds}. For this, observe that for all $\phi\in \s D(L^2(\Sigma)^p)=\bigcap_{\nu\in \mathbb{R}}L_\nu^2(L^2(\Sigma)^p)$, we have that
  \[
     \eps_n\phi\to {\eps}\phi\text{ as }n\to\infty \text{ in }L^2(\R;L^2(\Sigma)^p).
  \]
  Since, $\eps_n$ and ${\eps}$ commute with multiplication by functions of the type $t\mapsto e^{\xi t}$, $\xi\in \mathbb{R}$, which is a bijection on $\s D(L^2(\Sigma)^p)$, we get that for all $\nu\in\mathbb{R}$
  \[
     \eps_n\phi\to {\eps}\phi \text{ as }n\to\infty \text{ in }L^2_{\nu}(\R;L^2(\Sigma)^p).
  \]
  Thus, by Theorem \ref{t:ds} employing the boundedness of $(T_{\eps_n})_n$ again, we infer $T_{\eps_n}\to T_{\eps}$ in $L_{\textnormal{sev}}^\textnormal{s}(L^2(\Sigma)^p)$ as $n\to\infty$.
\end{Proof}

\begin{Theorem}\label{t:dbfns} Let $\Omega$ an admissible domain, $\curl_\diamond$ an admissible realization, $(\eps_n)_n, (\mu_n)_n$ in $L^\infty(\R\times \Omega)^{3\times 3}$. Assume that $\eps_n$ and $\mu_n$ satisfy Hypothesis \ref{h:cof} and that there exists $c>0$ such that for all $n\in \mathbb{N}$ the mappings $\eps_n$ and $\mu_n$ satisfy the inequalities \eqref{eq:hcof}. Let $\s S_{\textnormal{DBF}}$ be as in Theorem \ref{t:wpdbf}.
\begin{enumerate}[label=(\alph*)]
 \item\label{dbfn} If $(\eps_n)_n$, $(\mu_n)_n$ converge in $L^\infty(\mathbb{R}\times\Omega)^{3\times 3}$, then 
 \[
    \lim_{n\to\infty}\s S_{\textnormal{DBF}}(\eps_n,\mu_n) = \s S_{\textnormal{DBF}}(\lim_{n\to\infty}\eps_n,\lim_{n\to\infty}\mu_n)\in L_{\textnormal{sev}}^\textnormal{n}(L^2(\Omega)^6).
 \]
 \item\label{dbfs} If $(\eps_n)_n$, $(\mu_n)_n$ are such that the associated multiplication operators on $L(L^2(\mathbb{R}\times\Omega)^3)$ converge in the strong operator topology, then 
 \[
    \lim_{n\to\infty}\s S_{\textnormal{DBF}}(\eps_n,\mu_n) = \s S_{\textnormal{DBF}}(\lim_{n\to\infty}\eps_n,\lim_{n\to\infty}\mu_n)\in L_{\textnormal{sev}}^\textnormal{s}(L^2(\Omega)^6).
 \]
\end{enumerate}
\end{Theorem}
\begin{Proof}
 For the proof of \ref{dbfn}, we want to apply Theorem \ref{t:cdoden}. For this, we observe that convergence in $L^\infty(\mathbb{R}\times\Omega)^{3\times 3}$ implies convergence of the associated multiplication operators in  $L(L_\nu^2(\R;L(\Omega)^3))$ for all $\nu\in\mathbb{R}$. In particular, the sequences $(\eps_n)_n$ and $(\mu_n)_n$ are bounded and
 \begin{align*}
    \s M_n &\coloneqq \begin{pmatrix}
     \eps_n  & 0 \\ 0 & \mu_n
    \end{pmatrix}\to \begin{pmatrix}
     \lim_{n\to\infty}\eps_n  & 0 \\ 0 & \lim_{n\to\infty}\mu_n
    \end{pmatrix}\in L_{\textnormal{sev}}^{\textnormal{n}}(L^2(\Omega)^6);
    \\ \s N_n&\coloneqq \begin{pmatrix} 0 & -\curl_\diamond \\ \curl_\diamond & 0 \end{pmatrix}(1+\beta\curl_\diamond)^{-1}.
 \end{align*} Thus, the assertion indeed follows from Theorem \ref{t:cdoden}.
 
 In order to prove \ref{dbfs}, we deduce by Proposition \ref{p:stronmult} that (the multiplication operators induced by) $(\eps_n)_n$ and $(\mu_n)_n$ are bounded and converge in $L_{\textnormal{sev}}^\textnormal{s}(L^2(\Omega)^3)$. The assertion follows from Theorem \ref{t:cdodes}.
\end{Proof}

Our next aim is to derive a result corresponding to Theorem \ref{t:dbfns} for the weak operator topology. We are aiming at a result, which provides a complete description of the limiting equation. For this, we need to confine ourselves with a restricted class of multiplication operators.

\begin{Theorem}\label{t:dbfw} Let $\Omega$ an admissible domain, $\curl_\diamond$ an admissible realization. Let $\eps,\mu\colon\mathbb{R}\to \mathbb{C}$ bounded, measurable, $1$-periodic functions with $\Re\eps(t),\Re\mu(t)\geq c$ for some $c>0$. For $n\in \mathbb{N}$ define $\eps_n(t)\coloneqq \eps(nt)$, $\mu_n(t)\coloneqq \mu(nt)$, $t\in \mathbb{R}$. Then, with $\s S_{\textnormal{DBF}}$ as in Theorem \ref{t:wpdbf},
\[
   \lim_{n\to\infty} \s S_{\textnormal{DBF}}(\eps_n,\mu_n) = \s S_{\textnormal{DBF}}\Big(\big(\int_0^1\frac{1}{\eps}\big)^{-1},\big(\int_0^1\frac{1}{\mu}\big)^{-1}\Big)\in L_{\textnormal{sev}}^\textnormal{w}(L^2(\Omega)^6).
\] 
\end{Theorem}

The proof of Theorem \ref{t:dbfw} is an application of Theorem \ref{t:cdodewsol}. For this, we need to compute the limits in \eqref{eq:defwlimit}. First of all, however, we recall a well-known statement of more general nature, which is related to periodic mappings and will be stated without proof.

\begin{Theorem}[{{see e.g.~\cite[Theorem 2.6]{CioDon}}}]\label{t:perm} Let $d\in \mathbb{N}$, $\eps\colon \mathbb{R}^d\to \mathbb{C}$ bounded, measurable and $(0,1)^d$-periodic, that is, for all $z\in \mathbb{Z}^d$ and almost every $x\in \mathbb{R}^d$, we have $\eps(x+z)=\eps(x)$. Then $\eps_n\coloneqq \eps(n\cdot)\to\int_{(0,1)^d}\eps(x)\dd x$ in the weak* topology of $L^\infty(\R^d)$ as $n\to\infty$. 
\end{Theorem}

Next, we recall a result from \cite{Waurick2014JAA_G}.
\begin{Theorem}[{{\cite[Theorem 5.7]{Waurick2014JAA_G}}}]\label{t:timprod} Let $a_{1},\ldots,a_{k}\in L^{\infty}(\mathbb{R})$ be $1$-periodic, $k\in \N$. For $n\in\mathbb{N}$ let $T_{j,n}$ denote the multiplication operator induced by $t\mapsto a_j(nt)$, $j\in \{1,\ldots,k\}$. Then we have 
\[
 T_{1,n}\check{\partial}_{t}^{-1}T_{2,n}\check{\partial}_{t}^{-1}T_{3,n}\cdots\check{\partial}_{t}^{-1}T_{k,n} \to \left(\check{\partial}_{t}^{-1}\right)^{k-1}\prod_{j=1}^{k}\int_{0}^{1}a_{j}(y)\dd y
\]
in $L_{\textnormal{sev}}^{\textnormal{w}}(\s X)$ for any Hilbert space $\s X$.
\end{Theorem}
\begin{Proof}
Without loss of generality, we may assume $\s X=\mathbb{C}$. Let $\nu>0$. For $n\in\mathbb{N}$ let $T_n\coloneqq T_{1,n}\partial_{t,\nu}^{-1}T_{2,n}\partial_{t,\nu}^{-1}T_{3,n}\cdots\partial_{t,\nu}^{-1}T_{k,n}$. For $K,L\subseteq\mathbb{R}$ bounded, measurable, we compute with the help of Theorem \ref{t:td_inv}{
\begin{align*}
&\langle\1_{K},T_{n}\1_{L}\rangle_{L_\nu^2} 
\\& =\int_{K}a_{1}(nt_{1})\int_{-\infty}^{t_{1}}a_{2}(nt_{2})\int_{-\infty}^{t_{2}}\cdots\int_{-\infty}^{t_{k-1}}a_{k}(nt_{k})\1_{L}(t_{k})\dd t_{k}\cdots \dd t_{2}e^{-2\nu t_{1}}\dd t_{1}\\
 & =\int_{K}\int_{-\infty}^{t_{1}}\int_{-\infty}^{t_{2}}\cdots\int_{-\infty}^{t_{k-1}}\left(\prod_{j=1}^{k}a_{j}(nt_{j})\right)\1_{L}(t_{k})e^{-2\nu t_{1}}\dd t_{k}\cdots \dd t_{1}\\
 & =\underbrace{\int_{\mathbb{R}}\cdots\int_{\mathbb{R}}}_{k\text{-times}}\left(\prod_{j=1}^{k}a_{j}(nt_{j})\right)\1_{K}(t_{1})\left(\prod_{j=2}^{k}\1_{(0,\infty)}(t_{j-1}-t_{j})\right)\1_{L}(t_{k})e^{-2\nu t_{1}}\dd t_{k}\cdots \dd t_{1}.
\end{align*}}
Next, observe that \[(t_{1},\ldots,t_{k})\mapsto\1_{K}(t_{1})\left(\prod_{j=2}^{k}\1_{(0,\infty)}(t_{j-1}-t_{j})\right)\1_{L}(t_{k})e^{-2\nu t_{1}}\in L^{1}(\mathbb{R}^{k}).\]
Moreover, the mapping $(t_{1},\ldots,t_{k})\mapsto\prod_{j=1}^{k}a_{j}(t_{j})$
is $(0,1)^{k}$-periodic. Thus, by Theorem \ref{t:perm}, we conclude
that 
\[
\langle\1_{K},T_{n}\1_{L}\rangle_{L_\nu^2}\to\left\langle \1_{K},\left(\partial_{t,\nu}^{-1}\right)^{k-1}\prod_{j=1}^{k}\int_{0}^{1}a_{j}(y)\dd y\1_{L}\right\rangle_{L_\nu^2}
\]
 as $n\to\infty$ for all $K,L\subseteq\mathbb{R}$ bounded and measurable.	
A density argument yields
\[
\langle\psi,T_{n}\phi\rangle_{L_\nu^2}\to\left\langle \psi,\left(\partial_{t,\nu}^{-1}\right)^{k-1}\prod_{j=1}^{k}\int_{0}^{1}a_{j}(y)dy\phi\right\rangle_{L_\nu^2}
\]
for all $\phi,\psi\in \s D = \bigcap_{\nu\in\mathbb{R}}L_\mu^2(\R)$. Thus, the assertion follows from Theorem \ref{t:dw}.
\end{Proof}	

With Theorem \ref{t:timprod} at hand, we can conclude this chapter with a proof of Theorem \ref{t:dbfw}.

\begin{Proof}[of Theorem \ref{t:dbfw}] We apply Theorem \ref{t:cdodewsol}. In the course of doing so, with 
\[
   \s M_n = \begin{pmatrix}
              \eps_n & 0 \\ 0 & \mu_n
            \end{pmatrix},\quad \s N_n=\s N = \curl_\diamond(1+\beta\curl_\diamond)^{-1}\begin{pmatrix}
              0 & -1 \\ 1 & 0
            \end{pmatrix}\quad(n\in \mathbb{N}),
\]
we want to compute $\s O$ and $\s P_k$, $k\in \mathbb{N}$, as defined in \eqref{eq:defwlimit}.	By Theorem \ref{t:perm}, we obtain that
\[
   \s O =\lim_{n\to\infty} \s M_n^{-1}=\begin{pmatrix}
              \int_0^1\frac{1}{\eps} & 0 \\ 0 & \int_0^1\frac{1}{\mu}
            \end{pmatrix}.
\]
Next, observe that $\s C\coloneqq \curl_\diamond(1+\beta\curl_\diamond)^{-1}$ commutes with $\s M_n$, $\s O$ and $\check{\partial}_t^{-1}$ as the former only acts on the spatial variables and the latter only act on the temporal ones. Hence, for $k\in \mathbb{N}$,
\begin{align*}
  \s P_k &= \lim_{n\to\infty} (-\s M_n^{-1}\check{\partial}^{-1}_t\s N)^k\s M_n^{-1} = \lim_{n\to\infty} (-\s C)^k(\s M_n^{-1}\check{\partial}^{-1}_t\big(\begin{smallmatrix}
                                                                        0 & -1 \\ 1 & 0
                                                                       \end{smallmatrix}\big)
)^k\s M_n^{-1} 
\\ &=  (-\s C)^k(\s O \check{\partial}^{-1}_t\big(\begin{smallmatrix}
                                                                        0 & -1 \\ 1 & 0
                                                                       \end{smallmatrix}\big)
)^k\s O
=   (-\s O \check{\partial}^{-1}_t\s N)^k\s O,
\end{align*}
where we used Theorem \ref{t:timprod}. Thus, we may apply Theorem \ref{t:cdodewsol} to get 
\[
  \textnormal{sol}(\s M_n,\s N)\to \textnormal{sol}(\s M_\infty,0) \text{ in }L_{\textnormal{sev}}^\textnormal{w}(L^2(\Omega)^6)\text{ as }n\to\infty,
\]
with
\[ \s M_\infty = \s O^{-1}+\s O^{-1}\sum_{\ell=1}^\infty \s R^\ell,\quad\s R  = -\sum_{k=1}^\infty \s P_k \s O^{-1}.\]
So, employing $ \s P_k=(-\s O \check{\partial}^{-1}_t\s N)^k\s O$ for all $k\in \mathbb{N}$, we arrive at
\begin{align*}
  \s M_\infty & = \s O^{-1}\sum_{\ell=0}^\infty \s R^\ell
   = \s O^{-1}\sum_{\ell=0}^\infty \big(-\sum_{k=1}^\infty \s P_k \s O^{-1}\big)^\ell
  \\ & = \s O^{-1}\sum_{\ell=0}^\infty \big(-\sum_{k=1}^\infty (-\s O \check{\partial}^{-1}_t\s N)^k\s O \s O^{-1}\big)^\ell
   = \s O^{-1}\sum_{\ell=0}^\infty \big(1-\sum_{k=0}^\infty (-\s O \check{\partial}^{-1}_t\s N)^k\big)^\ell
  \\ & = \s O^{-1}\sum_{\ell=0}^\infty \big(1-(1+\s O\check{\partial}_t^{-1}\s N)^{-1}\big)^\ell
   = \s O^{-1}(1-\big(1-(1+\s O\check{\partial}_t^{-1}\s N)^{-1}\big))^{-1}
  \\ & = \s O^{-1}((1+\s O\check{\partial}_t^{-1}\s N)^{-1})^{-1}
   = \s O^{-1}+\check{\partial}_t^{-1}\s N.
\end{align*}
So, $\check{\partial}_t\s M_\infty= \check{\partial}_t\s O^{-1} + \s N$ yields the assertion. 
\end{Proof}

\section{Comments}

In \cite{Waurick2011}, we have introduced a topology on possible coefficients $\s M$ on equations of the form
\[
   (\partial_{t,\nu}\s M + \s A)u=f\in L_\nu^2(\mathbb{R};\s X),
\]
where $\s M$ was assumed to be translation-invariant and causal. Hence, $\s M=M(\partial_{t,\nu}^{-1})$ for an operator-valued, bounded, analytic function $M\colon B(r,r)\to L(\s X)$, $r>\frac{1}{2\nu}$, $\s X$ Hilbert space. We gather the coefficients treated in the set
\[
   \s H^\infty \coloneqq \s H^\infty(B(r,r);L(\s X))\coloneqq \{ M\colon B(r,r)\to L(\s X); M \textnormal{ bounded, analytic}\}.
\]
We endowed $\s H^\infty$ with the topology induced by
\[
   \s H^\infty \ni M\mapsto (z\mapsto \langle \psi,M(z)\phi\rangle_{\s X})\in \s H(B(r,r)),
\]
where $\s H(B(r,r))$ is the space of scalar valued analytic functions on $B(r,r)$ endowed with the compact open topology, that is, the topology of uniform convergence on compacts, $\phi,\psi\in \s X$. It turns out that \[B_{\s H^\infty}\coloneqq \{ M\in \s H^\infty; \sup_z\|M(z)\|\leq 1\}\] is compact under this topology, see \cite[Theorem 4.3]{Waurick2014SIAM_HomFrac}. Moreover, one can show that, if $(M_\iota)_\iota$  in $B_{\s H^\infty}$ converges to some $N\in B_{\s H^\infty}$, then, for all $\mu>1/2r$, we obtain $(M_\iota(\partial_{t,\mu}^{-1}))_\iota$ converges to $N(\partial_{t,\mu}^{-1})$ in the weak operator topology of $L(L_\mu^2(\s X))$, see \cite[Lemma 3.5]{Waurick2012MMAS_ODES}. Thus, 
\[
   B_{\s H^\infty} \ni M \mapsto M(\partial_{t,\nu}^{-1})\in L_{\textnormal{sev},\nu}^\textnormal{w}(\s X)
\]
is continuous. And, by compactness of $B_{\s H^\infty}$, the mapping just defined is even a homeomorphism on its image. Thus, the results derived in the previous chapter and the results to follow are proper generalizations of the results being initially restricted to $\s H^\infty$. 

The study in \cite{Waurick2011} has been developed for treating problems in homogenization theory. We address the general idea of homogenization theory by means of an example in Section \ref{s:homaw}. For now, we mention that, as a by-product of the functional analytic point of view developed, it is possible to explain memory effects occurring due to the process of homogenization: We consider for $\eps>0$ the solution $u_\eps\in L_\nu^2(\R;L^2(\R))$ of
\begin{equation}\label{eq:hom0}
   \partial_{t,\nu} u_\eps(t,x) + \sin((2\pi x)/\eps)u_\eps(t,x)= f(t,x)\quad (t,x\in \mathbb{R})
\end{equation}
for some $f\in C_c(\R\times\R)$. By the variation of constants formula we obtain
\[
   u_\eps(t,x) =  \int_{-\infty}^t   e^{-(t-s)\sin(2\pi x/ \eps)}f(s,x)\dd s\quad(t,x\in \mathbb{R}).
\]
Hence, multiplying the latter formula by $\phi\in L^2(\R)$ and integrating over $x\in \R$, by Theorem \ref{t:perm}, we infer
\[
   u_\eps \to u \coloneqq \int_{-\infty}^t J_0(i(t-s)) f(s,x)\dd s \text{ weakly in }L_\nu^2(\mathbb{R};L^2(\R)) \text{ as }\eps\to 0,
\]
where $J_0(z)=\sum_{k=0}^\infty \frac{\left(-\frac{1}{4}z^2\right)^k}{(k!)^2}$, $z\in \mathbb{C}\setminus(-\infty,0]$, is the $0$th order Bessel function of the first kind  (Note that $\int_0^1 e^{-(t-s)\sin(2\pi x)}\dd x = J_0(i(t-s))$).  Is it possible to find a differential equation, which is solved by $u$ and which has $f$ as a given source term? In fact, using the Theorems \ref{t:perm} and \ref{t:cdodewsol}, we obtain with a rather lengthy but straightforward computation
\[
   \partial_{t,\nu} u+\sum_{j=1}^\infty \partial_{t,\nu}\left(-\sum_{\ell=1}^\infty\frac{\Gamma(\frac{1}{2}+\ell)}{\sqrt{\pi}\Gamma(1+\ell)}\partial_{t,\nu}^{-2\ell}\right)^j u=f.
\]
We mention here that causal, translation-invariant coefficients for ordinary differential equations have been dealt with intensively in \cite{Waurick2014JAA_G,Waurick2012MMAS_ODES}; in \cite{Waurick2014JAA_G} we also treated causal evolutionary coefficients. However, we focused merely on sequences converging in the weak operator topology and did not choose the general perspective of discussing the continuity of the solution operator in the coefficients. The limiting equation is an equation of intergro-differential type. Hence, memory effects occur. A functional-analytic explanation is that computing the inverse of an operator is not a continuous process in the weak operator topology.

For an account of homogenization theory with regards to ordinary differential equations, we refer to \cite{Mascarenhas1984,TartarMemEff,Petrini1998,Antoni'c1993} for a non-exhaustive list. In the bulk of these studies, however, the description of the limiting equation uses the notion of Young measures, see also \cite[Remark 3.8]{Waurick2014JAA_G}. We briefly elaborate on Young measures as follows. Given a bounded sequence $(a_n)_n$ of $[\alpha,\beta]$-valued $L^\infty(\R)$-functions for some real $\alpha<\beta$. For an appropriately chosen subsequence $(a_{n_k})_k$ it is possible to describe the weak*-limit of $(g\circ a_{n_k})_k$ for any continuous real function $g\in C(\R)$ by means of a family of measures $(\nu_x)_{x\in \mathbb{R}}$ in the way that
\[
   g\circ a_{n_k} \to \big(x\mapsto \int_\alpha^\beta g (s)\dd \nu_x(s)\big)\text{ in }\sigma(L^\infty,L^1)\text{ as }k\to\infty.
\]
The derivation of the \emph{Young measure} $(\nu_x)_{x\in\R}$ (associated to $(a_{n_k})_k$) is not constructive and relies on a compactness theorem for the weak* topology for measures, see also \cite[Theorem 2]{Ball}. The limit equation of $(\partial_{t,\nu}-a_{n_k})u_k = f$ for some appropriate $f$ can then be described by
\[
   \partial_{t,\nu} u(t,x) + a^0(x)u-\int_0^tK(t-s,x)u(s,x)\dd s = f(t,x)
\]
with $a^0$ being the weak*-limit of $(a_{n_k})_k$ and 
\[
   K(t,x)= \int_0^\infty e^{-\mu t}\dd \nu_x(\mu),
\]
see \cite[p 930]{TartarNonlHom}.

For other treatments of homogenization theory for ordinary differential equations, we refer to the list of references in \cite{Waurick2014JAA_G}.

\cleardoublestandardpage
\chapter[Continuous Dependence on Coefficients in Partial Differential Equations]{The Continuous Dependence on the Coefficients in Partial Differential Equations}\label{ch:pde}

In this chapter we will treat partial differential equations with regards to varying coefficients. More precisely, in Chapter \ref{ch:ST} Section \ref{sec:wpr0}, see \eqref{eq:pde0}, we discussed equations of the form
\begin{equation}\label{eq:pde1}
   (\partial_{t,\nu}\s M+\s N+\s A)u=f
\end{equation}
for some bounded $\s M$, $\s N$ and a possibly unbounded $\s A$. So, in view of applications discussed later on, we ask for continuous dependence on $\s M$ and $\s N$ under the various topologies introduced in Section \ref{s:nt}. In the first section to follow we will study both the norm and the strong operator topology. The second section will be concerned with the weak operator topology. Similar to the case of ordinary differential equations, the result for the weak operator topology is somewhat more involved. This chapter is concluded with continuous dependence results for partial differential equations and an application to homogenization theory.

\renewcommand{\baselinestretch}{0.65}\normalsize\mysection{Norm and Strong Topologies and Partial Differential Equations}{The Norm and the Strong Operator Topologies and Partial Differential Equations}{the set $\textnormal{SP}_{c,\nu,r}^\textnormal{s}$ $\cdot$ continuity result for the strong operator topology $\cdot$ continuity estimate for the norm topology $\cdot$ Theorem \ref{t:cdpdes} $\cdot$ Theorem \ref{t:cdpden}}\label{s:nsotpde}

\renewcommand{\baselinestretch}{1}\normalsize
Similar to the case of ordinary differential equations in the previous chapter, we define the solution operator according to \eqref{eq:pde1}. For this, we recall the assumptions that lead to a solution theory for \eqref{eq:pde1}. We start out with the conditions on $\s A$.

\begin{Hypothesis}\label{h:ntpde} Let $\s X$ Hilbert space, $\nu>0$, $\s A\in C_{\textnormal{ev},\nu}(\s X)$.
Assume for all $\mu\geq \nu$:
\[
   \partial_{t,\mu}^{-1}\s A\subseteq \s A^\mu\partial_{t,\mu}^{-1},
\]
and
\begin{multline*}
  \Re\langle Q_0 \s A\phi,\phi\rangle_{L_\mu^2(\s X)}\geq 0,\; \Re\langle(\s A^\mu)^*\psi,\psi\rangle_{L_\mu^2(\s X)}\geq 0 \\ (\phi\in \dom(\s A),\psi\in \dom((\s A^\mu)^*)).
\end{multline*}
\end{Hypothesis}

Next, we define the solution operator. We note that we only show a continuous dependence result for the strong operator topology. The corresponding result for the norm topology is a mere variant of a continuity estimate.

\begin{Definition}\label{d:solpde} Assume Hypothesis \ref{h:ntpde}. Let $r,c>0$. We define the set
\begin{align*}
  \textnormal{SP}^\textnormal{s}_{c,\nu,r}(\s X) \coloneqq & \Big\{ (\s M,\s N)\in L_{\textnormal{sev},\nu}^\textnormal{s}(\s X)^2 ;
  \\& \quad\quad\quad\quad \text{there is }\s M'\in \mathfrak{C}_\nu(r),\s D\subseteq \dom(\partial_{t,\nu})\text{ such that for all }\mu\geq\nu 
  \\&\quad\quad\quad\quad  \s M\partial_{t,\mu}\subseteq \partial_{t,\mu}\s M^{\mu}-(\s M')^\mu,
  \\  &\quad\quad\quad\quad   \Re \langle Q_t (\partial_{t,\mu}\s M+\s N)\phi,\phi\rangle_{L_\mu^2(\s X)}\geq c\langle Q_t\phi,\phi\rangle\quad (t\in \R,\phi\in \s D),
  \\  & \quad\quad\quad\quad	 \s D\subseteq \dom(\partial_{t,\mu}) \text{ is a core for }\partial_{t,\mu} \Big\},
\end{align*}
where $Q_t$ denotes multiplication by $\1_{(-\infty,t)}$ and \[\mathfrak{C}_\nu(r)=\{\s S\in L_{\textnormal{sev},\nu}(\s X);\sup_{\mu\geq\nu}\|\s S^\mu\|\leq r\}.\]
With the help of Theorem \ref{t:gSe}, the following mapping is well-defined.
\begin{align*}
  \textnormal{sol}\colon \textnormal{SP}_{c,\nu,r}^\textnormal{s}(\s X) &\to L_{\textnormal{sev}}^\textnormal{s}(\s X)
  \\ (\s M,\s N)& \mapsto \left(\check{\partial}_{t,\nu}\s M+\s N+\check{\s A}\right)^{-1},
\end{align*}
where $\check{\s A}\coloneqq \bigcap_{\mu\geq\nu} \s A^\mu$ and $\check{\partial}_{t,\nu}\coloneqq \bigcap_{\mu\geq\nu}\partial_{t,\mu}$, see Example \ref{ex:cev}.
\end{Definition}

In this section, we aim for establishing the following two results on the continuous dependence on $\s M$ and $\s N$:

\begin{Theorem}\label{t:cdpdes} Assume Hypothesis \ref{h:ntpde}, $r,c>0$. Then the operator
\[
   \textnormal{sol}\colon \textnormal{SP}_{c,\nu,r}^\textnormal{s}(\s X)\to L_{\textnormal{sev}}^\textnormal{s}(\s X)
\]
is continuous, where $\textnormal{sol}$ is given in Definition \ref{d:solpde}. 
\end{Theorem}
For the norm topology, we have the following announced quantitative estimate.
\begin{Theorem}\label{t:cdpden} Assume Hypothesis \ref{h:ntpde} to be satisfied, $r,c>0$. Then, for $(\s M,\s N)$, $(\s O,\s P)\in \textnormal{SP}_{c,\nu,r}^\textnormal{s}(\s X)$, $\mu\geq \nu$, we have
 \begin{multline*}
    \| \big(\textnormal{sol}(\s M,\s N)\textnormal{sol}(\s O,\s P)^{-1}-1\big)\check{\partial}_{t,\nu}^{-1}\textnormal{sol}(\s O,\s P)\|_{L(L_{\mu}^2)}
    \\ \leq \frac{1}{c} \Big(\|(\s M-\s O)\textnormal{sol}(\s O,\s P)\|_{L(L_\mu^2)}
    \\+\|(\s M'-\s O')\check{\partial}_{t,\nu}^{-1}\textnormal{sol}(\s O,\s P)\|_{L(L_\mu^2)}+\|(\s N-\s P)\check{\partial}_{t,\nu}^{-1}\textnormal{sol}(\s O,\s P)\|_{L(L_\mu^2)}\Big).
 \end{multline*}
\end{Theorem}

\begin{Remark} Note that in the notation of $\textnormal{sol}(\s M,\s N)$, we did not keep explicit reference to $\s M'$. Once existent this operator is uniquely determined by the inclusion
\[
   \s M\partial_{t,\mu}\subseteq \partial_{t,\mu}\s M^{\mu}-(\s M')^\mu
\]
for $\mu\geq\nu$. 
\end{Remark}

Both the results Theorem \ref{t:cdpdes} and Theorem \ref{t:cdpden} have their roots in the following fundamental identity. We recall that we will suppress the superscript of evolutionary mappings indicating the spaces, where the closure is computed.

\begin{Proposition}\label{p:funid} Assume Hypothesis \ref{h:ntpde}, $r,c>0$ and $(\s M,\s N), (\s O,\s P)\in \textnormal{SP}_{c,\nu,r}^\textnormal{s}(\s X)$. Then we have
\begin{multline}\label{eq:funid}
   \big(\textnormal{sol}(\s M,\s N)-\textnormal{sol}(\s O,\s P)\big)\textnormal{sol}(\s O,\s P)^{-1}\check{\partial}_{t,\nu}^{-1}\textnormal{sol}(\s O,\s P)
   \\=\textnormal{sol}(\s M,\s N)\big((\s O-\s M)+(\s O'-\s M')\check{\partial}_{t,\nu}^{-1}+(\s P-\s N)\check{\partial}_{t,\nu}^{-1}\big)\textnormal{sol}(\s O,\s P).
\end{multline}
on $\s D(\s X)$, where $\textnormal{sol}$ is given in Definition \ref{d:solpde}.
\end{Proposition}
\begin{Proof}
  Let $\phi\in \s D(\s X)$ and define $\s B_1\coloneqq \textnormal{sol}(\s M,\s N)^{-1}$ and $\s B_2\coloneqq \textnormal{sol}(\s O,\s P)^{-1}$. We compute in the space $L_\mu^2$ for some $\mu\geq\nu$. Note that $\check{\partial}_{t,\nu}^{-1}\s B_j^{-1}\phi\in \dom(\check{\s A}^\nu)$ (see e.g.~Lemma \ref{l:appmation}). We compute with the help of Remark \ref{r:epinf}
  \begin{align*}
     & \big(\textnormal{sol}(\s M,\s N)-\textnormal{sol}(\s O,\s P)\big)\textnormal{sol}(\s O,\s P)^{-1}\check{\partial}_{t,\nu}^{-1}\textnormal{sol}(\s O,\s P)\phi
     \\ & = (\s B_1^{-1}\s B_2-1)\check{\partial}_{t,\nu}^{-1}\s B_2^{-1}\phi
     \\ & = (\s B_1^{-1}\s B_2\check{\partial}_{t,\nu}^{-1}\s B_2^{-1}-\check{\partial}_{t,\nu}^{-1}\s B_2^{-1})\phi
     \\ & = (\s B_1^{-1}\s B_2[\check{\partial}_{t,\nu}^{-1},\s B_2^{-1}]+\s B_1^{-1}\check{\partial}_{t,\nu}^{-1}-\check{\partial}_{t,\nu}^{-1}\s B_2^{-1})\phi
     \\ & = (\s B_1^{-1}[\s B_2,\check{\partial}_{t,\nu}^{-1}]\s B_2^{-1}+\s B_1^{-1}(\check{\partial}_{t,\nu}^{-1}\s B_2-\s B_1\check{\partial}_{t,\nu}^{-1})\s B_2^{-1})\phi
     \\ & = (\s B_1^{-1}[\s B_2,\check{\partial}_{t,\nu}^{-1}]\s B_2^{-1}+\s B_1^{-1}([\check{\partial}_{t,\nu}^{-1},\s B_2]+\s B_2\check{\partial}_{t,\nu}^{-1}-\s B_1\check{\partial}_{t,\nu}^{-1})\s B_2^{-1})\phi
     \\ & = \s B_1^{-1}\big(\s B_2\check{\partial}_{t,\nu}^{-1}-\s B_1\check{\partial}_{t,\nu}^{-1}\big)\s B_2^{-1}\phi
     \\ & = \s B_1^{-1}\big(\check{\partial}_{t,\nu}(\s O-\s M)\check{\partial}_{t,\nu}^{-1}+(\s P-\s N)\check{\partial}_{t,\nu}^{-1}\big)\s B_2^{-1}\phi
     \\ & = \s B_1^{-1}\big((\s O-\s M)+(\s O'-\s M')\check{\partial}_{t,\nu}^{-1}+(\s P-\s N)\check{\partial}_{t,\nu}^{-1}\big)\s B_2^{-1}\phi.
  \end{align*}
\end{Proof}

\begin{Proof}[of Theorem \ref{t:cdpden}] Using Proposition \ref{p:top} for computing the operator norm, the result is a direct consequence of equality \eqref{eq:funid} in Proposition \ref{p:funid}. 
\end{Proof}

The proof of Theorem \ref{t:cdpdes} needs yet another prerequisite. Indeed, equality \eqref{eq:funid} shows that, if $(\s M_\iota,\s N_\iota)_\iota$ converge in $\textnormal{SP}_{c,\nu,r}^\textnormal{s}(\s X)$ (and $(\s M_\iota')_\iota$) to some $(\s O,\s P)$ (and $\s O'$), then \[\big(\textnormal{sol}(\s M_\iota,\s N_\iota)-\textnormal{sol}(\s O,\s P)\big)\stackrel{\iota}{\to} 0\] strongly on the range of the operator \[\textnormal{sol}(\s O,\s P)^{-1}\check{\partial}_{t,\nu}^{-1}\textnormal{sol}(\s O,\s P),\] which is the same as the domain of $\textnormal{sol}(\s O,\s P)^{-1}\check{\partial}_{t,\nu}\textnormal{sol}(\s O,\s P)$. Hence, by the boundedness of $(\textnormal{sol}(\s M_\iota,\s N_\iota))_\iota$, in order to prove Theorem \ref{t:cdpdes}, one needs to show that \[\textnormal{sol}(\s O,\s P)^{-1}\check{\partial}_{t,\nu}\textnormal{sol}(\s O,\s P)\] is densely defined. When we come to the proof of Theorem \ref{t:cdpdes}, we will elaborate on the convergence of $(\s M_\iota')_\iota$.

\begin{Lemma}\label{l:dn} Assume Hypothesis \ref{h:ntpde}, $r,c>0$, $(\s M,\s N)\in \textnormal{SP}_{c,\nu,r}^\textnormal{s}(\s X)$. Then, for all $\mu\geq\nu$, the operator $(\textnormal{sol}(\s M,\s N)^{-1})^\mu{\partial}_{t,\mu}\textnormal{sol}(\s M,\s N)^\mu$ is densely defined in $L_{\mu}^2(\s X)$. 
\end{Lemma}
\begin{Proof} First of all note that for all $\eps>0$, we obtain 
\[
   \s S_\eps\coloneqq \left(1+\eps(\check{\partial}_{t,\nu}\s M+\s N+\check{\s A})\right)^{-1}\in L_{\textnormal{sev},\nu}(\s X).
\]
Indeed, $\eps\check{\s A}$ satisfies Hypothesis \ref{h:ntpde} and $(\eps\check{\partial}_{t,\nu}\s M,1+\eps\s N)\in \textnormal{SP}_{1,\nu,r}^\textnormal{s}(\s X)$.
Let $\mu\geq\nu$. We make all computations in the space $L_\mu^2(\s X)$ and consider the closures of all evolutionary operators involved in this space, without explicitly recording it in the notation. As in Remark \ref{r:eps_st}, it is readily seen that $\s S_\eps\to 1$ in the strong operator topology of $L(L_\mu^2(\s X))$. Moreover, for all $\phi\in L_\mu^2(\s X)$ we infer that $\s S_\eps\phi\in \dom(\textnormal{sol}(\s M,\s N)^{-1})$. For $\phi\in L_\mu^2(\s X)$ define for $\delta,\eps>0$ 
\[
   \phi_{\delta,\eps}\coloneqq \textnormal{sol}(\s M,\s N)^{-1}(1+\delta\partial_{t,\mu})^{-1}\s S_\eps \textnormal{sol}(\s M,\s N)\phi.
\]
Next, $\phi_{\delta,\eps}$ is well-defined, that is, $(1+\delta\partial_{t,\mu})^{-1}\s S_\eps$ maps into the domain of the operator $\textnormal{sol}(\s M,\s N)^{-1}$ by Lemma \ref{l:com_B}.
Next, $\phi_{\delta,\eps}\to \phi$ in $L_\mu^2(\s X)$ as $\delta,\eps\to 0$, by Lemma \ref{l:appmation} and Remark \ref{r:eps_st}. It is obvious, that $\phi_{\delta,\eps}\in \dom({\partial}_{t,\mu}\textnormal{sol}(\s M,\s N))$. Moreover, again by Lemma \ref{l:com_B}, we infer
\begin{align*}
 & {\partial}_{t,\mu}\textnormal{sol}(\s M,\s N)\phi_{\delta,\eps} 
\\ & =\partial_{t,\mu}(1+\delta\partial_{t,\mu})^{-1}\s S_\eps \textnormal{sol}(\s M,\s N)\phi
  \\ & = \frac{1}{\delta}(1+\delta\partial_{t,\mu}-1)(1+\delta\partial_{t,\mu})^{-1}\s S_\eps \textnormal{sol}(\s M,\s N)\phi
  \\ & = \frac{1}{\delta}\left(1-(1+\delta\partial_{t,\mu})^{-1}\right)\s S_\eps \textnormal{sol}(\s M,\s N)\phi
  \\ & = \frac{1}{\delta}\s S_\eps \textnormal{sol}(\s M,\s N)\phi-\frac{1}{\delta}(1+\delta\partial_{t,\mu})^{-1}\s S_\eps \textnormal{sol}(\s M,\s N)\phi\in \dom(\textnormal{sol}(\s M,\s N)^{-1}).
\end{align*}
\end{Proof}

We can now conclude this section with a proof of Theorem \ref{t:cdodes}.

\begin{Proof}[of Theorem \ref{t:cdodes}] Recalling the reasoning right before Lemma \ref{l:dn} and Lemma \ref{l:dn} itself, we are left with showing the following: Let $(\s M_{\iota})_\iota$ be a convergent net in $L_{\textnormal{sev}}^\textnormal{s}(\s X)$ with the property that there exists $r>0$ such that for all $\iota$ there is $\s M'_\iota\in \mathfrak{C}_\nu(r)$ with the property
\begin{equation}\label{eq:compr}
   \s M_\iota\partial_{t,\mu}\subseteq \partial_{t,\mu}\s M_\iota^\mu-(\s M_\iota')^\mu\quad(\mu\geq\nu).
\end{equation}
Then $\lim_\iota \s M'_\iota$ exists in $\mathfrak{C}_\nu(r)$ as a limit in $L_{\textnormal{sev}}^\textnormal{s}(\s X)$ and 
\begin{equation}\label{eq:compr1}
     \lim_\iota \s M_\iota\partial_{t,\mu} \subseteq \partial_{t,\mu}(\lim_\iota\s M_\iota)^\mu-\lim_\iota \s M'_\iota.
\end{equation}
First of all note that $(\s M'_\iota)_\iota$ lies $\mathfrak{C}_\nu(r)$. So, if $((\s M'_\iota)^\eta)_\iota$ converges in the strong operator topology for some $\eta\geq\nu$, then $((\s M'_\iota)^\eta)_\iota$ converges in the weak operator topology of $L(L_\eta^2(\s X))$. But, $\mathfrak{C}_\nu(r)^\textnormal{w}$ is closed in $L_{\textnormal{sev}}^\textnormal{w}(\s X)$ by Proposition \ref{p:bdd_cld_weak}, this implies that $(\s M'_\iota)_\iota$ converges in $\mathfrak{C}_\nu(r)^\textnormal{w}$ by Theorem \ref{t:dw}. Hence, any accumulation value of $(\s M'_\iota)_\iota$ under the strong operator topology lies in $\mathfrak{C}_\nu(r)$.

Next, from \eqref{eq:compr}, we read off that 
\[
     \partial_{t,\mu}^{-1}\s M_\iota\partial_{t,\mu}\subseteq \s M_\iota^\mu-\partial_{t,\mu}^{-1}(\s M_\iota')^\mu\quad(\mu\geq\nu,\iota\in I).
\]
Let $\eta\geq\nu$ be such that $(\s M_\iota^\mu)_\iota$ converges in the strong operator topology of $L(L_\mu^2(\s X))$ for all $\mu\geq \eta$. For all $\phi\in \dom(\partial_{t,\mu})$ we get
\begin{equation}\label{eq:compr2}
     \partial_{t,\mu}^{-1}\s M_\iota^\mu\partial_{t,\mu}\phi - \s M_\iota^\mu\phi=-\partial_{t,\mu}^{-1}(\s M_\iota')^\mu\phi\quad(\mu\geq\eta,\iota\in I).
\end{equation}
  The left hand side of equation \eqref{eq:compr2} converges in $L_\mu^2(\s X)$ since $(\s M_\iota^\mu)_\iota$ converges strongly. But, $(\partial_{t,\mu}^{-1}(\s M_\iota')^\mu)_\iota$ is a bounded net of linear operators from $L_\mu^2(\s X)$ to $\dom(\partial_{t,\mu})$, where the latter is endowed with the graph norm of $\partial_{t,\mu}$ and so $(\partial_{t,\mu}^{-1}(\s M_\iota')^\mu)_\iota$ converges in the strong operator topology of $L(L_\mu^2(\s X),\dom(\partial_{t,\mu}))$. Next, as $\partial_{t,\mu}^{-1}\colon L_\mu^2(\s X)\to \dom(\partial_{t,\mu})$ is a Banach space isomorphism, $((\s M_\iota')^\mu)_\iota$ converges in the strong operator topology of $L(L_\mu^2(\s X))$. Finally, from \eqref{eq:compr2}, by performing the limit in $\iota$, we get \eqref{eq:compr1}.
\end{Proof}

\renewcommand{\baselinestretch}{0.65}\normalsize\mysection{Weak Topology and Partial Differential Equations}{The Weak Operator Topology and Partial Differential Equations}{continuity result for the weak operator topology $\cdot$ the set $\textnormal{SP}_{c,\nu,r}^\textnormal{w}$ $\cdot$ Theorem of Aubin--Lions $\cdot$ weak-strong principle $\cdot$ Theorem \ref{t:cdpdew}}\label{s:wotpde}

\renewcommand{\baselinestretch}{1}\normalsize
Similar to our way of presenting the case of ordinary differential equations, we now seek a result analogous to the Theorems \ref{t:cdpdes} and \ref{t:cdpden} for the weak operator topology. As the rationale in the previous chapter shows, the weak operator topology is likely to be more involved due to the missing continuity statements in Theorem \ref{t:cim} \ref{nn} and \ref{ss} or a corresponding result of Theorem \ref{t:cinv} for the weak operator topology.

In this section, we will provide a criterion roughly saying the following:
\begin{alignat}{1}\label{eq:roughst}
\begin{aligned}
 &\text{If }(\s M_\iota)_\iota \text{ converges in the weak operator topology,}\\ &\text{then }((\partial_{t,\nu}\s M_\iota +\s A)^{-1})_\iota\text{ converges to }(\partial_{t,\nu}\lim_\iota\s M_\iota +\s A)^{-1}.
\end{aligned}
\end{alignat}
In view of Corollary \ref{c:cdodew2} or Theorem \ref{t:dbfw}, the limit to equal $(\partial_{t,\nu}\lim_\iota\s M_\iota +\s A)^{-1}$ might be somewhat unexpected: One might suspect that the limit should involve some sort of harmonic mean (as in Theorem \ref{t:dbfw}). Indeed, if we formally set 
\[
  \s A=\curl_\diamond (1+\beta\curl_\diamond)^{-1}\begin{pmatrix} 0 & -1 \\ 1 & 0 \end{pmatrix}\text{ and }\s M_n=\begin{pmatrix} \eps_n & 0 \\ 0 & \mu_n \end{pmatrix},
\]
Theorem \ref{t:dbfw} asserts 
\[
  \lim_{n\to\infty} (\partial_{t,\nu}\s M_n +\s A)^{-1} = (\partial_{t,\nu}(\lim_{n\to\infty}\s M_n^{-1})^{-1} +\s A)^{-1}.
\]
As $(\lim_{n\to\infty}\s M_n^{-1})^{-1}\neq \lim_{n\to\infty}\s M_n$ for convergence in the weak operator topology in general, a result of the type \eqref{eq:roughst} can only be true under additional assumptions on $\s A$. In fact, this is where the unboundedness of $\s A$ comes into play:
 
\begin{Hypothesis}\label{h:awpde} We assume the conditions in Hypothesis \ref{h:ntpde}, that is, $\s X$ Hilbert space, $\nu>0$, $\s A\in C_{\textnormal{ev},\nu}(\s X)$.
Assume for all $\mu\geq \nu$ and $\phi\in \dom(\s A),\psi\in \dom((\s A^\mu)^*)$,
\[
   \partial_{t,\mu}^{-1}\s A\subseteq \s A^\mu\partial_{t,\mu}^{-1},\;\Re\langle Q_0 \s A\phi,\phi\rangle_{L_\mu^2(\s X)}\geq 0,\; \Re\langle(\s A^\mu)^*\psi,\psi\rangle_{L_\mu^2(\s X)}\geq 0.
\]In addition, assume there exists a Hilbert space $\s Y$ compactly embedded into $\s X$ such that 
\[
   (\dom(\s A^\mu),\|\cdot\|_{\s A^\mu})\hookrightarrow L_\mu^2(\R; \s Y)\quad (\mu\geq\nu),
\]
where $\|\cdot\|_{\s A^\mu}$ is the graph norm of $\s A^\mu$.
\end{Hypothesis}

The solution operator to study in this section reads as follows.

\begin{Definition}\label{d:solpdew} Assume Hypothesis \ref{h:awpde}, $r,c>0$. Define the set
\begin{align*}
  &\textnormal{SP}^\textnormal{w}_{c,\nu,r}(\s X) \coloneqq \Big\{ \s M\in L_{\textnormal{sev},\nu}^\textnormal{w}(\s X) ;
  \\& \quad\quad\quad\quad\quad\quad\quad\quad \text{there is }\s M'\in \mathfrak{C}_\nu(r), \s D\subseteq \dom(\partial_{t,\nu})\text{ such that for all }\mu\geq\nu 
  \\&\quad\quad\quad\quad\quad\quad\quad\quad  \s M\partial_{t,\mu}\subseteq \partial_{t,\mu}\s M^{\mu}-(\s M')^\mu,
  \\  &\quad\quad\quad\quad\quad\quad\quad\quad   \Re \langle Q_t\partial_{t,\mu}\s M\phi,\phi\rangle_{L_\mu^2(\s X)}\geq c\langle Q_t\phi,\phi\rangle\quad (t\in \R,\phi\in \s D),
  \\  & \quad\quad\quad\quad\quad\quad\quad\quad \s D\subseteq \dom(\partial_{t,\mu}) \text{ is a core for }\partial_{t,\mu} \Big\},
\end{align*}
where $\mathfrak{C}_\nu(r)=\{\s S\in L_{\textnormal{sev},\nu}(\s X); \sup_{\mu\geq\nu}\|\s S^\mu\| \leq r\}$ and $Q_t$ denotes multiplication by $\1_{(-\infty,t)}$.
Due to Theorem \ref{t:gSe}, we may define
\begin{align*}
  \textnormal{sol}\colon \textnormal{SP}^\textnormal{w}_{c,\nu,r}(\s X) &\to L_{\textnormal{sev},\nu}^\textnormal{w}(\s X)
  \\ \s M&\mapsto \left(\check{\partial}_{t,\nu}\s M+\check{\s A}\right)^{-1}.
\end{align*}
where $\check{\s A}\coloneqq \bigcap_{\mu\geq\nu} \s A^\mu$ and $\check{\partial}_{t,\nu}\coloneqq \bigcap_{\mu\geq\nu}\partial_{t,\mu}$, see Example \ref{ex:cev}.
\end{Definition}

With all the relevant notions at hand, we may now state the main theorem of the current section:

\begin{Theorem}\label{t:cdpdew} Assume Hypothesis \ref{h:awpde}, $c,r>0$. Then 
\[
  \textnormal{sol}\colon \textnormal{SP}^\textnormal{w}_{c,\nu,r}(\s X)\subseteq L_{\textnormal{sev},\nu}^\textnormal{w}(\s X)  \to L_{\textnormal{sev},\nu}^\textnormal{w}(\s X) 
\]
 as introduced in Definition \ref{d:solpdew} is continuous on bounded sets.
\end{Theorem}

The proof of Theorem \ref{t:cdpdew} needs several prerequisites. We denote by $H^1(J;\s X)$ the Bochner--Sobolev space of Hilbert space $\s X$-valued weakly differentiable functions with $L^2(J;\s X)$-derivative, $J\subseteq \R$ open interval. First of all, we quote a well-known compactness theorem.

\begin{Theorem}[Theorem of Aubin--Lions, {\cite[p. 67, $2^\circ$]{Simon1987}}]\label{t:AuLi} Let $\s X, \s Y$ be Hilbert spaces, $J\subseteq \R$ bounded, open interval. Assume that $\s Y\hookrightarrow\hookrightarrow\s X$, that is, $\s Y$ is compactly embedded into $\s X$. Then 
\[
   H^1(J;\s X)\cap L^2(J;\s Y)\hookrightarrow\hookrightarrow L^2(J;\s X).
\] 
\end{Theorem}

Next, we will state the most important result for a proof of Theorem \ref{t:cdpdew}, which is essentially a consequence of the Aubin--Lions Theorem and causality. Beforehand, recall the notation $H^1_\mu(\R;\s X)$ for the domain of $\dom(\partial_{t,\mu})$ defined in $L_\mu^2(\R;\s X)$ endowed with the graph norm of $\partial_{t,\mu}$, $\s X$ Hilbert space. Moreover, for a bounded interval $J\subseteq \R$, we have that the canonical embedding $L^2(J;\s X)\hookrightarrow L_\mu^2(\R;\s X)$, where we extend any function in the left-hand side by zero to the whole real line, is continuous for all $\mu\in \mathbb{R}$.

\begin{Theorem}[weak-strong principle]\label{t:wsp} Let $\s X,\s Y$ Hilbert spaces with $\s Y\hookrightarrow\hookrightarrow \s X$. Let $\nu>0$ and $(v_\iota)_{\iota\in I}$ be a weakly convergent, bounded net in $L^2_\mu(\R;\s Y)\cap H_{\mu}^1(\R;\s X)$ for all $\mu\geq\nu$. Assume there exists $t_0\in \R$ with $Q_{t_0}v_\iota=0$ for all $\iota\in I$, where $Q_{t_0}$ denotes multiplication by $\1_{(-\infty,t_0)}$. If $(\s M_\iota)_{\iota\in I}$ is a bounded, convergent net in $L_{\textnormal{sev}}^\textnormal{w}(\s X)$, then $(\s M_\iota v_\iota)_\iota$ weakly converges in $L_\mu^2(\R;\s X)$ and
\[
   \lim_{\iota\in I} \s M_\iota^\mu v_\iota = \lim_{\iota\in I}\s M_\iota^\mu \lim_{\iota \in I}v_\iota \in L_\mu^2(\R;\s X)
\] for all sufficiently large $\mu$.
\end{Theorem}
\begin{Proof} Let $\eta\geq\nu$ be such that for all $\mu\geq\eta$ we have that $(v_\iota)_\iota$ is weakly convergent in $L^2_\mu(\R;\s Y)\cap H_{\mu}^1(\R;\s X)$ and that $(\s M_\iota)_{\iota}$ is bounded and convergent in $L_{\textnormal{sev},\eta}^\textnormal{w}(\s X)$.

Next, let $\phi \in L^2_{c}(\R;\s X)$ and define $t\coloneqq \sup \spt \phi$. Let $Q_t$ be multiplication by $\1_{(-\infty,t)}$ and denote $w\coloneqq \lim_{\iota}v_\iota\in \bigcap_{\mu\geq\nu}L^2_\mu(\R;\s Y)\cap H_{\mu}^1(\R;\s X)$ and $\s N\coloneqq \lim_{\iota}\s M_\iota$. Now, by Theorem \ref{t:AuLi}, we deduce that $(Q_t v_\iota)_\iota$ converges to $Q_t v$ in (norm in) $L^2((t_0,t);\s X)\subseteq L_\mu^2(\R;\s X)$ for all $\mu\geq\eta$. Next, for all $\iota$, we compute with the help of causality of $\s M_\iota$ and $\s N$ (Proposition \ref{p:sev_sc} and Theorem \ref{t:csc}) for all $\mu\geq\eta$:
\begin{align*}
  \langle \s M_\iota^\mu v_\iota,\phi\rangle_{L_\mu^2} & = \langle \s M_\iota^\mu v_\iota, Q_t\phi\rangle_{L_\mu^2}  = \langle Q_t \s M_\iota^\mu v_\iota, \phi\rangle_{L_\mu^2} \\ & = \langle Q_t \s M_\iota^\mu Q_t  v_\iota, \phi\rangle_{L_\mu^2} = \langle \s M_\iota^\mu Q_t  v_\iota, Q_t\phi\rangle_{L_\mu^2} \\
                                   & \to \langle \s N^\mu Q_t v, Q_t\phi\rangle_{L_\mu^2}  = \langle Q_t\s N^\mu Q_t v, \phi\rangle_{L_\mu^2}  = \langle \s N^\mu v, \phi\rangle_{L_\mu^2}.
\end{align*}By the boundedness of $(\s M_\iota^\mu v_\iota)_\iota$ and the density of $L_c^2(\R;\s X)$ in $L_\mu^2(\s X)$, we get the assertion.
\end{Proof}

Next, we prove an adapted version of an assertion in the proof of Theorem \ref{t:cdpdes}. Namely, the weak operator topology convergence of $(\s M_\iota)_\iota$ in $\textnormal{SP}_{c,\nu,r}^\textnormal{w}$ implies the same for $(\s M_\iota')_\iota$:

\begin{Lemma}\label{l:wlc} Let $\s X$ Hilbert space, $D$ densely defined, closed, linear operator in $\s X$ with $0\in\rho(D)$. Let $(M_\iota)_\iota$ be a convergent net in $L(\s X)$ under the weak operator topology. Assume that there exists $r>0$ such that for all $\iota\in I$ there is $M_\iota'\in L(\s X)$ with $\|M_\iota'\|\leq r$ and 
\[
   M_\iota D \subseteq D M_\iota - M_\iota'.
\]
Then $\lim_{\iota\in I} M_\iota'$ exists in the weak operator topology, $\|\lim_{\iota\in I} M_\iota'\|\leq r$ and 
\begin{equation}\label{eq:wlc}
    (\lim_{\iota\in I} M_\iota)D \subseteq D(\lim_{\iota\in I} M_\iota) - \lim_{\iota\in I}M_\iota'.
\end{equation}
 Moreover, for all $u\in \dom(D)$, $(D M_\iota u)_\iota$ is weakly convergent and
\[
   \lim_{\iota\in I} (D M_\iota u) = D (\lim_{\iota\in I}M_\iota) u.
\]
\end{Lemma}
\begin{Proof}
Define $N\coloneqq \lim_{\iota\in I} M_\iota$.
 For $\phi,\psi\in \s X$, $\iota\in I$ we compute
\begin{align*}
  \langle M_\iota'{D}^{-1}\phi,({D}^{-1})^*\psi \rangle &=\langle ({D}M_\iota-M_\iota {D}){D}^{-1}\phi,({D}^{-1})^*\psi\rangle \\
                                             &= \langle {D}^{-1}({D}M_\iota-M_\iota{ D}){D}^{-1}\phi,\psi\rangle \\
                                             &= \langle (M_\iota{D}^{-1}-{D}^{-1}M_\iota)\phi,\psi\rangle \\
                                             &\to \langle N{D}^{-1}\phi,\psi\rangle-\langle{D}^{-1}N \phi,\psi\rangle\\
                                             &= \langle (DN-ND){D}^{-1}\phi,({D}^{-1})^*\psi \rangle.  
\end{align*}
Thus, the boundedness of $(M_\iota')_\iota$ together with the density of both the domains $\dom({D})$ ($D$ is densely defined) and $\dom({D}^*)$ ($D$ is closed) implies that $(M_\iota')_\iota$ converges in the weak operator topology. The estimate for the operator norm of $\lim_\iota M_\iota'$ follows from Proposition \ref{p:liminfine}.  Moreover, for all $\phi\in \dom(D)$, we have 
\[
   \lim_{\iota\in I} M_\iota'\phi = (DN-ND)\phi=D(\lim_{\iota\in I}M_\iota)\phi-(\lim_{\iota\in I} M_\iota)D\phi.
\]
Hence, \eqref{eq:wlc} follows. The remaining assertion is a straightforward consequence of \eqref{eq:wlc} and the convergence of $(M_\iota')_\iota$.
\end{Proof}

\begin{Remark}\label{r:SPwcls} By Lemma \ref{l:wlc}, we infer the closedness of $\textnormal{SP}_{c,\nu,r}^\textnormal{w}(\s X)\subseteq L_{\textnormal{sev},\nu}^\textnormal{w}(\s X)$: Let $(\s M_\iota)_\iota$ be a net in $\textnormal{SP}_{c,\nu,r}^\textnormal{w}(\s X)$ convergent in $L_{\textnormal{sev},\nu}^\textnormal{w}(\s X)$. Then, for all $\mu\geq\nu$, $(\s M_\iota^\mu)_\iota$ converges in the weak operator topology of $L(L_\mu^2(\s X))$. Hence, applying Lemma \ref{l:wlc} to $D=\partial_{t,\mu}$, $M_\iota=\s M_\iota^\mu$, we infer that $(M_\iota')_\iota=((\s M_\iota')^\mu)_\iota$ converges in the weak operator topology of $L(L_\mu^2(\s X))$, $\|(\s M_\iota')^\mu\|\leq r$ as well as 
\[
   \lim_{\iota\in I}\s M_\iota \partial_{t,\mu}\subseteq \lim_{\iota\in I}\s M_\iota^\mu \partial_{t,\mu}\subseteq \partial_{t,\mu} \lim_{\iota\in I}\s M_\iota^\mu - \lim_{\iota\in I}(\s M_\iota')^\mu.
\]
Next, for $\iota\in I$, there exists $\s D\subseteq \dom(\partial_{t,\nu})$ a core for $\partial_{t,\mu}$ for all $\mu\geq\nu$ such that for all $t\in \mathbb{R}$ and $\phi\in \s D$ we have
\begin{equation}\label{eq:SPwcls}
 \Re\langle Q_t (\s M_\iota^\mu\partial_{t,\mu}+(\s M_\iota')^\mu)\phi,\phi\rangle_{L_\mu^2}=\Re\langle Q_t \partial_{t,\mu}\s M_\iota\phi,\phi\rangle_{L_\mu^2}\geq c\langle Q_t\phi,\phi\rangle_{L_\mu^2}.
\end{equation}
Thus, by continuous extension, \eqref{eq:SPwcls} holds for all $\phi\in \bigcap_{\eta\geq\nu}\dom(\partial_{t,\eta})$. Taking the limit in $\iota$ in \eqref{eq:SPwcls}, we infer 
\[
 \Re\langle Q_t ((\lim_{\iota}\s M_\iota^\mu)\partial_{t,\mu}+(\lim_{\iota}(\s M_\iota')^\mu))\phi,\phi\rangle_{L_\mu^2}\geq c\langle Q_t\phi,\phi\rangle_{L_\mu^2}.
\]
So, 
\[
  \Re\langle Q_t \partial_{t,\mu}(\lim_{\iota}\s M_\iota)\phi,\phi\rangle_{L_\mu^2}\geq c\langle Q_t\phi,\phi\rangle_{L_\mu^2}\quad(\mu\geq\nu)
\]
for all $\phi\in \bigcap_{\eta\geq\nu}\dom(\partial_{t,\eta})$, which is a core for $\partial_{t,\mu}$.
\end{Remark}

Before we come to the proof of Theorem \ref{t:cdpdew}, we briefly recall the `subnet argument' in topological spaces. For this also recall that $(x_{\phi(\kappa)})_{\kappa\in J}$ is a \emph{subnet} of a net $(x_\iota)_{\iota\in I}$, if $(J,\leq_J)$ is a directed set and $\phi\colon J\to I$ is \emph{cofinal in $I$}, that is, for all $\iota_0\in I$ there exists $\kappa_0\in J$ such that for all $\kappa\geq_J \kappa_0$ we have $\phi(\kappa)\geq_I \iota_0$.

\begin{Proposition}\label{p:sna} Let $(\Omega,\tau)$ be a Hausdorff topological space, $(x_\iota)_\iota$ a net in $\Omega$, $y\in \Omega$. Then $(x_\iota)_\iota$ converges to $y$ if and only if any subnet of $(x_\iota)_\iota$ contains a subnet, which converges to $y$. 
\end{Proposition}
\begin{Proof} The necessity is easy.

On the other hand, assume that $(x_\iota)_\iota$ does not converge to $y$. Then there exists a neighborhood $U\subseteq \Omega$ of $y$ such that for all $\iota$ we find $\phi(\iota)\geq \iota$ such that $x_{\phi(\iota)}\notin U$. The, thus, defined map $\phi\colon I\to I$ is cofinal in $I$. Indeed, take $\iota_0\in I$ then we get for all $\iota\geq \iota_0$ that $\phi(\iota)\geq\iota\geq\iota_0$. So, $(x_{\phi(\iota)})_{\iota\in I}$ is a subnet of $(x_\iota)_\iota$. Now, if $(x_{\psi(\kappa)})_{\kappa\in J}$ is a subnet of $(x_{\phi(\iota)})_{\iota\in I}$, we have $\{x_{\psi(\kappa)};\kappa\in J\}\cap U\subseteq \{x_{\phi(\iota)};\iota\in I\}\cap U=\emptyset$. Thus, $(x_{\psi(\kappa)})_{\kappa\in J}$  does not converge to $y$. 
\end{Proof}

\begin{Proof}[of Theorem \ref{t:cdpdew}] Let $(\s M_\iota)_{\iota\in I}$ be a bounded, convergent net in $\textnormal{SP}_{c,\nu,r}^\textnormal{w}(\s X)$, denote $\s N\coloneqq \lim_{\iota}\s M_\iota$. Lemma \ref{l:wlc} implies that $\s N'\coloneqq \lim_{\iota}\s M_\iota'$ converges in $L_{\textnormal{sev},\nu}^\textnormal{w}(\s X)$ and is such that
\[
   \s N\partial_{t,\mu}\subseteq \partial_{t,\mu}\s N^\mu -(\s N')^\mu\quad(\mu\geq\nu).
\]
By Theorem \ref{t:gSe}, we infer that $\sup_{\mu\geq\nu}\|\textnormal{sol}(\s M_\iota)^\mu\|\leq 1/c$. Let $\mu\geq\nu$. In order to reduce cluttered notation, we will drop the superindex $\mu$ in the notation of all evolutionary mappings for denoting the closure in $L_\mu^2(\s X)$; this also applies to $\textnormal{sol}(\s M_\iota)$(and its inverse), $\s A$ and $\s M_\iota'$ with $\s M_\iota \partial_{t,\mu}\subseteq \partial_{t,\mu}\s M_\iota - \s M_\iota'$, $\iota\in I$.

 Let $f\in L_c^2(\R;\s X)$. For $\iota\in I$ we set
\[
   u_\iota \coloneqq \textnormal{sol}(\s M_\iota) f.
\]
Our aim is to show that 
\begin{equation}\label{eq:s1}
  \lim_\iota u_\iota = \textnormal{sol}(\lim_\iota\s M_\iota) f.
\end{equation}
 Then, from $\|\textnormal{sol}(\s M_\iota)\|_{L(L_\mu^2)}\leq 1/c$ and the density of $L_c^2(\R;\s X)$ in $L_\mu^2(\R; \s X)$ it follows that the net $(\textnormal{sol}(\lim_\iota\s M_\iota))_\iota$ converges in the weak operator topology of $L(L_\mu^2(\s X))$.

Since $f$ is compactly supported, there exists $t_0\in \R$ such that $Q_{t_0} f=0$, where $Q_{t_0}$ is multiplication by $\1_{(-\infty,t_0)}$. Hence, by causality of $\textnormal{sol}(\s M_\iota)$, for all $\iota\in I$, we get  
\[
   Q_{t_0}u_\iota = Q_{t_0}\textnormal{sol}(\s M_\iota) f = Q_{t_0}\textnormal{sol}(\s M_\iota) Q_{t_0} f = 0.
\]
 Next, for $\eps>0,\iota\in I$ we set $u_{\eps,\iota} \coloneqq (1+\eps\partial_{t,\mu})^{-1}u_\iota$. Then $Q_{t_0}u_{\eps,\iota}=0$, by causality of $(1+\eps\partial_{t,\mu})^{-1}$ (see Remark \ref{r:rtc}). With the help of Remark \ref{r:prestA} together with Lemma \ref{l:appmation}, we get
 \begin{alignat}{1}\label{eq:abd}
 \begin{aligned}
   \|\s A  u_{\eps,\iota}\|_{L_\mu^2(\s X)} & \leq \|f\|_{L_\mu^2} + (2\|\s M_\iota'\|_{L(L_\mu^2)}+\frac{1}{\eps}\|\s M_\iota\|_{L(L_\mu^2)})\|u_{\iota}\|_{L_\mu^2}
   \\&\leq
   ((2r +\frac{1}{\eps}\|\s M_\iota\|_{L(L_\mu^2)})\frac{1}{c}+1)\|f\|_{L_\mu^2}.
 \end{aligned} 
 \end{alignat}
 Since $(u_\iota)_\iota$ is bounded in $L_\mu^2(\R;\s X)$, using the weak compactness of the unit ball of $L_\mu^2(\R;\s X)$, we find a weakly convergent subnet $(u_{\phi(\kappa)})_{\kappa\in J}$ and set $u_\kappa\coloneqq u_{\phi(\kappa)}$ and accordingly define $u_{\eps,\kappa}$, $\s M_\kappa$ and $\s M_\kappa'$, $\kappa\in J$, $\eps>0$. Let $v\coloneqq \lim_\kappa u_{\kappa}$. In the following, we will prove that 
 \begin{equation}\label{eq:s2}
   v=\lim_\kappa u_{\kappa}=\textnormal{sol}(\lim_\iota M_\iota)f.
  \end{equation}
 By Proposition \ref{p:sna}, we infer that then $\lim_\iota u_\iota = v$, which yields \eqref{eq:s1}, and, hence, the assertion of the Theorem.
 
 As the mapping $(1+\eps\partial_{t,\mu})^{-1}$ is continuous, we get
 \begin{equation}\label{eq:wH}
   \lim_\kappa u_{\eps,\kappa} = \lim_\kappa u_{\eps,\phi(\kappa)} = v_\eps\coloneqq (1+\eps\partial_{t,\mu})^{-1}v \text{ weakly in }H^1_\mu(\R;\s X)=\dom(\partial_{t,\mu}).
 \end{equation}
 Moreover, by estimate \eqref{eq:abd}, $(u_{\eps,\kappa})_\kappa$ is a bounded sequence in the graph space of $\s A$. As the latter is a Hilbert space as well, we may choose a weakly convergent subnet $(u_{\eps,\psi(o)})_o$ in the graph Hilbert space of $\s A$. But as this graph Hilbert space is continuously embedded into $L_\mu^2(\s X)$, we get that $\lim_o u_{\eps,\psi(o)}=\lim_\kappa u_{\eps,\kappa}$. So, 
 \begin{equation}\label{eq:wA}
  \lim_\kappa u_{\eps,\kappa} = v_\eps \text{ weakly in }(\dom(\s A);\|\cdot\|_{\s A}).
 \end{equation}
 The continuity of $\s A$ considered as a mapping from $(\dom(\s A);\|\cdot\|_{\s A})$ to $L_\mu^2(\s X)$ implies, given \eqref{eq:wA},
 \begin{equation}\label{eq:wAc}
  \lim_\kappa \s A u_{\eps,\kappa} = \s A v_\eps \text{ weakly in }L_\mu^2(\s X).
 \end{equation}
 Next, by \eqref{eq:wH}, \eqref{eq:wA} and since the embedding $(\dom(\s A);\|\cdot\|_{\s A})\hookrightarrow L_\mu^2(\R;\s Y)$ is continuous, we infer
 \begin{equation}\label{eq:wAL}
   \lim_\kappa u_{\eps,\kappa} = v_\eps \text{ weakly in }H^1_\mu(\R;\s X)\cap L_\mu^2(\R; \s Y).
 \end{equation}
 
 Setting $\s B_{\kappa}\coloneqq \textnormal{sol}(\s M_{\kappa})^{-1}$, we obtain from Lemma \ref{l:com_B} for all $\eps>0$ and $\kappa\in J$: 
 \begin{alignat}{1}\label{eq:kap}
 \begin{aligned}
    (1+\eps\partial_{t,\mu})^{-1}f&=(1+\eps\partial_{t,\mu})^{-1}\s B_\kappa u_{\kappa} 
    \\& = \s B_\kappa u_{\eps,\kappa} - \eps\partial_{t,\mu}(1+\eps\partial_{t,\mu})^{-1}\s M_{\kappa}' u_{\eps,\kappa} \\
     & = \partial_{t,\mu}\s M_{\kappa} u_{\eps,\kappa} + \s A u_{\eps,\kappa}  - \eps\partial_{t,\mu}(1+\eps\partial_{t,\mu})^{-1}\s M_{\kappa}' u_{\eps,\kappa}.
 \end{aligned}
 \end{alignat}
 We inspect the convergence of each summand of the right hand side in \eqref{eq:kap}. We apply Theorem \ref{t:wsp} to $(\s M_\kappa)_\kappa$ and $(u_{\eps,\kappa})_\kappa$ in place of $(\s M_\iota)_\iota$ and $(v_\iota)_\iota$. Since $Q_{t_0}u_{\eps,\kappa}=0$ as well as \eqref{eq:wAL}, Theorem \ref{t:wsp} is indeed applicable. Thus, we infer
 \[
    \lim_\kappa \s M_{\kappa} u_{\eps,\kappa} = \lim_\kappa \s M_\kappa \lim_\kappa u_{\eps,\kappa}=\lim_\iota \s M_\iota \lim_\kappa u_{\eps,\kappa} = \s N v_\eps.
 \]
 But, the net $(\partial_{t,\mu}\s M_{\kappa} u_{\eps,\kappa})_\kappa$ is bounded in $L_\mu^2$:
 \begin{align*}
    & \|\partial_{t,\mu}\s M_{\kappa} u_{\eps,\kappa}\|_{L_\mu^2(\s X)}
    \\ & \leq \|\s M_{\kappa} \partial_{t,\mu}\left(1+\eps\partial_{t,\mu}\right)^{-1} u_{\kappa}\|_{L_\mu^2(\s X)}+ \|\s M_{\kappa}' \left(1+\eps\partial_{t,\mu}\right)^{-1} u_{\kappa}\|_{L_\mu^2(\s X)}
    \\ & \leq \|\s M_{\kappa}\| \frac{1}{\eps} \|u_{\kappa}\|_{L_\mu^2(\s X)}+ \|\s M_{\kappa}'\| \|u_{\kappa}\|_{L_\mu^2(\s X)}
    \\ & \leq \sup_{\iota}\|\s M_{\iota}\| \frac{1}{c\eps} \|f\|_{L_\mu^2(\s X)}+ \frac{r}{c} \|f\|_{L_\mu^2(\s X)}.
 \end{align*}
 A compactness argument thus yields
 \begin{equation}\label{eq:f1}
    \lim_\kappa \partial_{t,\mu}\s M_{\kappa} u_{\eps,\kappa}=\partial_{t,\mu}\s N v_\eps\text{ weakly in }L_\mu^2(\R;\s X).
 \end{equation}
 Recall from the beginning of the proof that we may apply Theorem \ref{t:wsp} also to $(\s M_{\kappa}')_\kappa$. Hence, by the continuity of $\eps\partial_{t,\mu}(1+\eps\partial_{t,\mu})^{-1}$ in $L_\mu^2(\s X)$, we get
 \begin{equation}\label{eq:f2}
    \lim_\kappa \eps\partial_{t,\mu}(1+\eps\partial_{t,\mu})^{-1} \s M_\kappa'u_{\eps,\kappa}=\eps\partial_{t,\mu}(1+\eps\partial_{t,\mu})^{-1} \s N'v_{\eps} \text{ weakly in }L_\mu^2(\R;\s X).
 \end{equation}
 Thus, putting the convergence results \eqref{eq:f1}, \eqref{eq:wAc} and \eqref{eq:f2} together, we read off from equation \eqref{eq:kap} that for all $\eps>0$
 \begin{align*}
    &(1+\eps\partial_{t,\mu})^{-1}f\\ & =\lim_{\kappa}\left(\partial_{t,\mu}\s M_{\kappa} u_{\eps,\kappa} + \s A u_{\eps,\kappa}  - \eps\partial_{t,\mu}(1+\eps\partial_{t,\mu})^{-1}\s M_{\kappa}' u_{\eps,\kappa}\right)
    \\ & = \lim_{\kappa}(\partial_{t,\mu}\s M_{\kappa} u_{\eps,\kappa}) + \lim_\kappa (\s A u_{\eps,\kappa})  - \lim_\kappa (\eps\partial_{t,\mu}(1+\eps\partial_{t,\mu})^{-1}\s M_{\kappa}' u_{\eps,\kappa})
    \\ & = \partial_{t,\mu}\s N v_{\eps} + \s A v_{\eps}  -  \eps\partial_{t,\mu}(1+\eps\partial_{t,\mu})^{-1}\s N' v_{\eps}.
 \end{align*}
  Hence, with $\s B\coloneqq \textnormal{sol}(\s N)^{-1}$ we get for all $\eps>0$
  \[
    \s B v_\eps = (1+\eps\partial_{t,\mu})^{-1}f   +  \eps\partial_{t,\mu}(1+\eps\partial_{t,\mu})^{-1}\s N' v_{\eps}
  \]
  or, equivalently,
  \[
      v_\eps = \textnormal{sol}(\s N)((1+\eps\partial_{t,\mu})^{-1}f   +  \eps\partial_{t,\mu}(1+\eps\partial_{t,\mu})^{-1}\s N'(1+\eps\partial_{t,\mu})^{-1} v)
  \]
  Hence, by the continuity of $\textnormal{sol}(\s N)$ and from $(1+\eps\partial_{t,\mu})^{-1}f\to f$ as $\eps\to0$ (see Remark \ref{r:eps_st}) as well as $\eps\partial_{t,\mu}(1+\eps\partial_{t,\mu})^{-1}\s N'(1+\eps\partial_{t,\mu})^{-1} v\to 0$ as $\eps\to 0$ (see Lemma \ref{le:commus_of_a}) it follows
  \begin{align*}
      v &= \lim_{\eps\to0} (1+\eps\partial_{t,\mu})^{-1}v =\lim_{\eps\to0} v_\eps 
      \\ &= \textnormal{sol}(\s N)(\lim_{\eps\to0}((1+\eps\partial_{t,\mu})^{-1}f   +  \eps\partial_{t,\mu}(1+\eps\partial_{t,\mu})^{-1}\s N'(1+\eps\partial_{t,\mu})^{-1} v))
      \\ &=\textnormal{sol}(\s N) f =\textnormal{sol}(\lim_\iota \s M_\iota) f,
  \end{align*}
 which yields \eqref{eq:s2}, thus, \eqref{eq:s1} and, hence, the assertion.
\end{Proof}

\begin{Remark}\label{r:sharper} (a) We note that there is also a way of including the topology of $L_{\textnormal{sev}}^\textnormal{w}(\s X)$ as the underlying topology for $\textnormal{SP}_{c,\nu,r}^\textnormal{w}(\s X)$ and in the target space of $\textnormal{sol}$ in Theorem \ref{t:cdpdew} in the following sense. Let $(\s M_\iota)_{\iota\in I}$ be a bounded and convergent net in $L_{\textnormal{sev}}^\textnormal{w}(\s X)$ with the property that $\s M_\iota\in \textnormal{SP}_{c,\nu,r}^\textnormal{w}(\s X)$ for all $\iota\in I$. By convergence of $(\s M_\iota)_{\iota\in I}$, there exists $\eta\geq\nu$ such that $(\s M_\iota)_{\iota\in I}$ converges in $L_{\textnormal{sev},\eta}^\textnormal{w}(\s X)$. Observe that from $\textnormal{SP}_{c,\nu,r}^\textnormal{w}(\s X)\subseteq \textnormal{SP}_{c,\eta,r}^\textnormal{w}(\s X)$, we get that $\lim_\iota M_\iota\in \textnormal{SP}_{c,\eta,r}^\textnormal{w}(\s X)$, as the space $\textnormal{SP}_{c,\eta,r}^\textnormal{w}(\s X)$ is a closed subset of $L_{\textnormal{sev},\eta}^\textnormal{w}(\s X)$, by Remark \ref{r:SPwcls}. Hence, $\textnormal{sol}(\s M_\iota)$ converges to $\textnormal{sol}(\lim_\iota \s M_\iota)$ in $L_{\textnormal{sev}}^\textnormal{w}(\s X)$, by Theorem \ref{t:cdpdew} applied to $\textnormal{SP}_{c,\eta,r}^\textnormal{w}(\s X)$ instead of $\textnormal{SP}_{c,\nu,r}^\textnormal{w}(\s X)$.

(b) It should be noted that a slightly more detailed analysis than the one done in the proof of Theorem \ref{t:cdpdew} shows the following stronger statement. If $(f_\iota)_\iota$ in $L_\nu^2(\R;H)$ is a bounded, weakly convergent net with $\inf_\iota\inf\spt f_\iota>-\infty$, then we also have
\[
   \textnormal{sol}(\s M_\iota)^\mu f_\iota \to \textnormal{sol}(\lim_{\iota}\s M_\iota)^\mu (\lim_{\iota}f_\iota)
\]with weak convergence in $L_\mu^2(\s X)$ for all $\mu\geq\nu$, where $(\s M_\iota)_\iota$ is any bounded and convergent net in $\textnormal{SP}_{c,\nu,r}^\textnormal{w}(\s X)$. 
\end{Remark}

The concluding sections of this chapter are devoted to examples.

\renewcommand{\baselinestretch}{0.65}\normalsize\mysection{Eddy-Current Approximation in Electromagnetic Theory}{The Eddy-Current Approximation in Electromagnetic Theory}{Maxwell's equations $\cdot$ electric boundary condition $\cdot$ Theorem \ref{t:eca}}\label{s:eca}

\renewcommand{\baselinestretch}{1}\normalsize
In this section, we will consider Maxwell's equations in matter. For this, let throughout this section $\Omega\subseteq \mathbb{R}^3$ be an open set. The equations are formally given by
\begin{alignat}{1}\label{eq:max}
\begin{aligned}
   \partial_{t} \eps E + \sigma E - \curl H &= J
   \\ \partial_{t} \mu H  + \curl E &= 0
\end{aligned}
\end{alignat}
 on the space time cylinder $\R\times\Omega$ subject to homogeneous electric boundary conditions for $E$ of vanishing tangential components. As in the Drude--Born--Fedorov model discussed in Section \ref{s:DBF}, the unknowns are the two components of the electromagnetic field $(E,H)\colon \mathbb{R}\times \Omega\to \R^3\times \R^3$. The material's properties are gathered in the coefficients $\eps,\mu,\sigma \colon \R\times \Omega\to \mathbb{C}^{3\times 3}$, which respectively are the dielectricity, magnetic permeability and the electric conductivity of the underlying medium. The given right-hand side $J\colon \mathbb{R}\times\Omega \to \mathbb{R}^3$ is a source term modeling external currents. We think of the system \eqref{eq:max} of being given on the whole time line $\mathbb{R}$ bearing in mind that -- thanks to causality -- the consideration of the real half-line as time parameter space eventually merely results in a restriction on the support of $J$, see also Remarks \ref{r:inex} and \ref{r:posreal}.
 
 In the study of eddy-currents (see e.g.~\cite{Yousept2012}), the dielectricity $\eps$ is observed to be rather small compared to the other operators $\sigma$ and $\mu$ involved in \eqref{eq:max}. That is why -- for simplicity -- $\eps$ is often neglected to the effect that the resulting system, the so-called \emph{eddy-current approximation}, formally reads
\begin{alignat}{1}\label{eq:edd}
\begin{aligned}
    \sigma E - \curl H &= J
   \\ \partial_{t} \mu H  + \curl E &= 0
\end{aligned}
\end{alignat}
on $\R\times \Omega$ subject to the electric boundary condition. Substituting the equation for $\sigma E$ into the one of $\partial_t \mu H$, one obtains an equation of parabolic type.

As an application of our results in Section \ref{s:nsotpde}, we will study the ``distance'' of solutions from \eqref{eq:max} to \eqref{eq:edd}. Before, however, doing so we set up the functional analytic framework for both the equations \eqref{eq:max} and \eqref{eq:edd}. For this, we introduce the $L^2$-operator realization of the $\curl$-operator with homogeneous electric boundary condition:

\begin{Definition}\label{d:elc} Denote $C^1_c(\Omega)^3$ the set of continuously differentiable vector fields and define
\begin{align*}
   \curl_c \colon C^1_c(\Omega)^3\subseteq L^2(\Omega)^3 & \to L^2(\Omega)^3
   \\  \phi &\mapsto \left(\begin{array}{ccc}
0 & -\partial_{3} & \partial_{2}\\
\partial_{3} & 0 & -\partial_{1}\\
-\partial_{2} & \partial_{1} & 0
\end{array}\right)\phi.
\end{align*}
Note that $\curl_c\subseteq \curl_c^*= \curl$, where $\curl$ is the (maximal) $L^2(\Omega)$-realization of the (distributional) $\curl$ operator, that is, 
\[
   \curl\colon \{\phi\in L^2(\Omega)^3;\curl\phi\in L^2(\Omega)^3\}\subseteq L^2(\Omega)^3 \to L^2(\Omega)^3
\]
acting as the distributional $\curl$ operator, see also Definition \ref{d:om}. Define
\[
   \curl_0 \coloneqq \overline{\curl}_c.
\]
\end{Definition}

We note that $\curl_0^*=\curl$. As we do not assume any regularity of the boundary of $\Omega$, in general, there is no continuous (tangential) trace operator. So, the replacement of the homogeneous electric boundary condition for $E$ is that $E\in \dom(\curl_0)$. A first step towards a solution theory for \eqref{eq:edd} and \eqref{eq:max} is the following almost trivial observation:

\begin{Lemma}\label{l:maxopsk} The operator
\begin{multline*}
  A\coloneqq \begin{pmatrix}
     0 &-\curl \\
    \curl_0 &  0
   \end{pmatrix}\colon \dom(\curl_0)\times \dom(\curl)
   \\ \subseteq L^2(\Omega)^3\times L^2(\Omega)^3\to L^2(\Omega)^3\times L^2(\Omega)^3
\end{multline*}
is skew-selfadjoint, that is, $A^*=-A$.
\end{Lemma}
\begin{Proof}
  The result is an application of the observation that for densely defined, closed linear operators $B_1\colon \dom(B_1)\subseteq\s X\to \s Y$ and $B_2\colon \dom(B_2)\subseteq \s Y\to \s X$ for some Hilbert spaces $\s X$ and $\s Y$, we have
  \[
     \begin{pmatrix}
     0 & B_1 \\
    B_2 &  0
   \end{pmatrix}^* = \begin{pmatrix}
     0 & B_2^* \\
    B_1^* & 0 
   \end{pmatrix}.
  \]
\end{Proof}

As in the case of the Drude--Born--Fedorov model, we lift the operator $A$ defined in Lemma \ref{l:maxopsk} as an (abstract) multiplication operator to the space time setting discussed: We set
\begin{equation}\label{eq:dmaxlift}
  \s A^\nu\colon  L_\nu^2(\R;\dom(A)) \subseteq L_\nu^2(\R; L^2(\Omega)^6) \to L_\nu^2(\R; L^2(\Omega)^6), \Phi\mapsto (t\mapsto A\Phi(t))
\end{equation}
for all $\nu\in \mathbb{R}$. We record the following facts regarding $\s A^\nu$:
\begin{enumerate}[label=(l\arabic*)]
 \item\label{sksA} for all $\nu\in \mathbb{R}$ the operator $\s A^\nu$ is skew-selfadjoint;
 \item\label{comQ} for all $\nu\in \mathbb{R}$ and $\phi\in L_\nu^2(\R;\dom(A))$, we infer $Q_0\phi\in L_\nu^2(\R;\dom(A))=\dom(\s A^\nu)$, $Q_0$ multiplication by $\1_{(-\infty,0)}$, and
 \[
    \s A^\nu Q_0 \phi = Q_0 \s A^\nu \phi.
 \]
\end{enumerate}

Next, we study the assumptions on the operators of multiplying by $\eps$, $\mu$ and $\sigma$:

\begin{Hypothesis}\label{h:ems} Let $\eps,\mu,\sigma\colon \mathbb{R}\to L^\infty(\Omega)^{3\times 3}$ be bounded and measurable. Moreover, assume $\eps,\mu\in C^1_b(\mathbb{R};L^\infty(\Omega)_{\textnormal{pd}}^{3\times 3})$, that is, $\eps$ and $\mu$ are continuously differentiable with bounded derivatives and attaining values in the selfadjoint, positive definite matrices: 
\[
 \eps(t)(x)=\eps(t)(x)^*\geq 0,\; \mu(t)(x)=\mu(t)(x)^*\geq 0\text{ for all }t\in \mathbb{R}, \text{ a.e.~}x\in \Omega.
\]
 Next, assume that there is $\nu\in \R$, $c>0$, such that for all $\eta\geq\nu$ we have
\begin{equation}\label{eq:emse}
   \eta\langle \xi,\eps(t)(x)\xi\rangle + (1/2)\langle \xi,\eps'(t)(x)\xi\rangle\geq 0\quad(\xi\in \mathbb{C}^3)
\end{equation}
and
\begin{multline}\label{eq:emsms}
   \eta\langle \xi,\mu(t)(x)\xi\rangle + (1/2)\langle \xi,\mu'(t)(x)\xi\rangle\geq c\langle\xi,\xi\rangle,\\ \Re\langle\xi,\sigma(t)(x)\xi\rangle\geq c\langle\xi,\xi\rangle \quad(\xi\in \mathbb{C}^3)
\end{multline} 
for almost every $x\in \Omega$ and $t\in\R$.
\end{Hypothesis}

Before applying the continuous dependency results from Section \ref{s:nsotpde}, we show that the assumptions in Hypothesis \ref{h:ems} together with the operator $\s A^\nu$ introduced in \eqref{eq:dmaxlift} lead to a proper solution theory for both equations \eqref{eq:max} and \eqref{eq:edd}. For this, we provide the following two lemmas:

\begin{Lemma}\label{l:maxa} Let $\nu\in(0,\infty)$, $\s A^\nu$ given by \eqref{eq:dmaxlift}. Then the following assertions hold true.
 \begin{enumerate}[label=(\alph*)]
  \item\label{maxa} For all $\phi\in \dom(\s A^\nu)=\dom((\s A^\nu)^*)$, we have $\Re\langle Q_0\s A^\nu\phi,\phi\rangle=0$, where $Q_0$ is multiplication by $\1_{(-\infty,0)}$.
  \item\label{maxb} For all $\phi\in \dom(\s A^\nu)$, we have $\partial_{t,\nu}^{-1}\s A^\nu\phi=\s A^\nu\partial_{t,\nu}^{-1}\phi$.
 \end{enumerate}
\end{Lemma}
\begin{Proof} For \ref{maxa}, we observe that the skew-selfadjointness of $A$ (Lemma \ref{l:maxopsk}) implies the same for $\s A^\nu$ (see \ref{sksA}). Hence, $\dom(\s A^\nu)=\dom(-\s A^\nu)=\dom((\s A^\nu)^*)$. We recall that for a skew-selfadjoint operator, the respective real-part vanishes. This together with \ref{comQ} implies for all $\phi\in \dom(\s A^\nu)$
\[
  \Re \langle Q_0 \s A^\nu\phi,\phi\rangle =\Re \langle Q_0 \s A^\nu\phi,Q_0\phi\rangle=\Re \langle  \s A^\nu Q_0\phi,Q_0\phi\rangle=\langle  \Re\s A^\nu Q_0\phi,Q_0\phi\rangle=0.
\]
In order to prove part \ref{maxb}, we apply Theorem \ref{t:td_inv} to $\s X=A$, that is, the Hilbert space $A$ considered as a closed subspace of $L^2(\Omega)^6\times L^2(\Omega)^6$ (see Lemma \ref{l:maxopsk}). By Theorem \ref{t:td_inv}, we infer that $\partial_{t,\nu}^{-1}$ is a continuous linear operator from 
\[
  L_\nu^2(\R;A) = \{  t\mapsto (\phi(t),A\phi(t)); \phi\in L_\nu^2(\R;\dom(A))\} = \s A^\nu
\]
into itself. This, in fact, is the assertion.  
\end{Proof}

As in the section on the Drude--Born--Fedorov model discussed earlier, for the sake of readability, we identify $\eps$, $\mu$ and $\sigma$ with their respective multiplication operators in $L_\nu^2(\R;L^2(\Omega)^3))$.

\begin{Lemma}\label{l:ems} Let $\eps$ be as in Hypothesis \ref{h:ems} satisfying \eqref{eq:emse}, $\eta\geq\nu$. Then the following conditions are true.
\begin{enumerate}[label=(\alph*)]
 \item\label{ems0} For all $\phi\in \dom(\partial_{t,\eta})$ we have $\partial_{t,\eta}\eps\phi=\eps\partial_{t,\eta}\phi+\eps'\phi$, where $\eps'$ is the multiplication operator induced by $(t\mapsto \eps(t))'$.
 \item\label{ems1} For all $\phi\in C_c^1(\mathbb{R};L^2(\Omega)^3)$ we have
 \[
   \Re \langle Q_t\partial_{t,\eta}\eps\phi,\phi\rangle_{L_\eta^2} \geq 0,
 \]where $Q_t$ is multiplication by $\1_{(-\infty,t)}$.
\end{enumerate}
\end{Lemma}
\begin{Proof}
  Note that \ref{ems0} has already been proven in Example \ref{ex:mult_op}. Thus, we are left with proving \ref{ems1}: Taking $\phi\in \dom(\partial_{t,\eta})$, we compute using integration by parts
  \begin{align*}
    &\langle Q_t\partial_{t,\eta}\eps\phi,\phi\rangle_{L_\eta^2} 
    \\& = \int_{-\infty}^t \langle (\eps\phi)'(s),\phi(s)\rangle_{L^2(\Omega)^3} e^{-2\eta s}\dd s
    \\ & = \langle \eps(t)\phi(t),\phi(t)\rangle_{L^2(\Omega)^3} e^{-2\eta t} -\int_{-\infty}^t \langle \eps\phi(s),\phi'(s)\rangle_{L^2(\Omega)^3} e^{-2\eta s}\dd s 
    \\ &\quad +2\eta \int_{-\infty}^t \langle \eps\phi(t),\phi(t)\rangle_{L^2(\Omega)^3} e^{-2\eta t}\dd t
    \\ & \geq -\langle Q_t\eps\phi,\partial_{t,\eta}\phi\rangle_{L_\eta^2} +2\eta \langle Q_t\eps\phi,\phi\rangle_{L_\eta^2}
  \end{align*}
Hence,
 \begin{align*}
   & 2\Re \langle Q_t\partial_{t,\eta}\eps\phi,\phi\rangle_{L_\eta^2} 
   \\&= \langle Q_t\partial_{t,\eta}\eps\phi,\phi\rangle_{L_\eta^2}+ \langle \phi,Q_t\partial_{t,\eta}\eps\phi\rangle_{L_\eta^2}
   \\ & \geq -\langle Q_t\eps\phi,\partial_{t,\eta}\phi\rangle_{L_\eta^2} +2\eta \langle Q_t\eps\phi,\phi\rangle_{L_\eta^2} + \langle Q_t\phi,\partial_{t,\eta}\eps\phi\rangle_{L_\eta^2}
   \\ & = -\langle Q_t\eps\phi,\partial_{t,\eta}\phi\rangle_{L_\eta^2} +2\eta \langle Q_t\eps\phi,\phi\rangle_{L_\eta^2} 
   \\&\quad+ \langle Q_t\eps\phi,\partial_{t,\eta}\phi\rangle_{L_\eta^2}+\langle Q_t\phi,\eps'\phi\rangle_{L_\eta^2}
   \\ & = 2\Big(\eta \int_{-\infty}^t \langle\big(\eps(s)\phi(s)+(1/2)\eps'(s)\phi(s)\big),\phi(s)\rangle_{L^2(\Omega)^3} e^{-2\eta s}\dd s\Big)\geq 0.
 \end{align*}
 \end{Proof}
So, we come to a solution theory for both \eqref{eq:max} and \eqref{eq:edd}. Note that for the coefficients $\eps$, $\mu$ and $\sigma$ satisfying Hypothesis \ref{h:ems} both equations to study in this section are covered as $\eps=0$ satisfies inequality \eqref{eq:emse} as well as the regularity requirements asked for in Hypothesis \ref{h:ems}.
We stress that we will not record the parameter $\nu$ in the notation of the operators $\curl$ and $\curl_0$ of the respective liftings to $L_\nu^2(\R; L^2(\Omega)^6)$.

\begin{Theorem}\label{t:wpmax} Assume Hypothesis \ref{h:ems} to be satisfied. Then there exists $\nu\in\R$ such that 
\[
   \s S_{\textnormal{MAX}}(\eps,\mu,\sigma)\coloneqq \left( \partial_{t,\nu}\begin{pmatrix}
                                                                              \eps & 0 \\ 0 & \mu
                                                                            \end{pmatrix} + \begin{pmatrix}
                                                                              \sigma & 0 \\ 0 & 0
                                                                            \end{pmatrix}+ \begin{pmatrix}
                                                                              0 & -\curl \\ \curl_0 & 0
                                                                            \end{pmatrix}\right)^{-1}
\]
is standard evolutionary at $\nu$. More precisely, there exists $r>0$ such that \[\left(\begin{pmatrix}
                                                                              \eps & 0 \\ 0 & \mu
                                                                            \end{pmatrix},\begin{pmatrix}
                                                                              \sigma & 0 \\ 0 & 0
                                                                            \end{pmatrix}\right)\in \textnormal{SP}^\textnormal{s}_{c,\nu,r}(L^2(\Omega)^6),\] where the latter space is given in Definition \ref{d:solpde}.
\end{Theorem}
\begin{Proof}
  By Lemma \ref{l:maxa}, the operator $\s A^\nu$ given in \eqref{eq:dmaxlift} satisfies all assumptions on $\s A$ in Hypothesis \ref{h:gSe}. Moreover, Lemma \ref{l:ems} ensures that 
  \[
     \Re \left\langle Q_t\left(\partial_{t,\eta}\begin{pmatrix}
                                                                              \eps & 0 \\ 0 & \mu
                                                                            \end{pmatrix} + \begin{pmatrix}
                                                                              \sigma & 0 \\ 0 & 0
                                                                            \end{pmatrix}\right) \begin{pmatrix}\phi \\ \psi \end{pmatrix},\begin{pmatrix}\phi \\ \psi \end{pmatrix}\right\rangle\geq c \left\|Q_t\begin{pmatrix}\phi \\ \psi \end{pmatrix}\right\|^2_{L_\eta^2}
  \]for all $\eta\geq\nu$, $\phi,\psi\in C_c^1(\R;L^2(\Omega)^3), t\in\mathbb{R}$ (apply Lemma \ref{l:ems} to $\mu-c$ in place of $\eps$); $\Re Q_t\sigma \geq cQ_t$ is easy to see. The rest of the conditions needed for Theorem \ref{t:gSe}, that is, the remaining conditions in Hypothesis \ref{h:gSe}, have been established in Lemma \ref{l:ems} as well.
\end{Proof}

The continuous dependency result, now, justifies the eddy-current formulation of Max\-well's equation as a proper approximation of the equations given originally. We state the result not in its most general form, as we will leave $\mu$ and $\sigma$ being fixed. We focus on variations in the dielectricity only:

\begin{Theorem}\label{t:eca} Let $\mu,\sigma$ as in Hypothesis \ref{h:ems}, $(\eps_n)_n$ bounded in $C^1_b(\R;L^\infty(\Omega)^{3\times 3}_{\textnormal{pd}})$ such that for all $n\in \mathbb{N}$ the map $\eps_n$ satisfies inequality \eqref{eq:emse} for all $\eta\geq\nu$ for some $\nu\in\R$; assume that $\eps_n\to 0$ in $C_b(\R;L^\infty(\Omega)^{3\times 3}_{\textnormal{pd}})$. Let $\s S_{\textnormal{MAX}}$ be given as in Theorem \ref{t:wpmax}. Then
\begin{equation}\label{eq:eca}
    \s S_{\textnormal{MAX}}(\eps_n,\mu,\sigma)\to \s S_{\textnormal{MAX}}(0,\mu,\sigma) \text{ in }L_{\textnormal{sev}}^\textnormal{s}(\R;L^2(\Omega)^6)\text{ as }n\to\infty,
\end{equation}
and, for all $n\in\N$, we have for all $\eta\geq\nu$
\begin{multline*}
    \| \big(\s S_{\textnormal{MAX}}(\eps_n,\mu,\sigma)\s S_{\textnormal{MAX}}(0,\mu,\sigma)^{-1}-1\big)\check{\partial}_{t,\nu}^{-1}\s S_{\textnormal{MAX}}(0,\mu,\sigma)\|_{L(L_{\eta}^2)} \\ \leq \frac{1}{c^2} \Big(\|\eps_n\|_{\infty}+\frac{1}{\eta}\|\eps_n'\|_{\infty}\Big),
\end{multline*}
where $\|\eps_n\|_{\infty}\coloneqq \sup_{t\in\mathbb{R}} \|\eps_n(t)\|_{L^\infty(\Omega)^{3\times 3}}$.
\end{Theorem}
\begin{Proof} First of all observe that for all $\eta\in \mathbb{R}$ the operator norm of $T_f$ the multiplication operator associated to $f\in C_b(\R;L^\infty(\Omega)^{3\times 3})$ realized as operator in $L(L_{\eta}^2(\R;L^2(\Omega)^3))$ satisfies $\|T_f\|_{L(L_{\eta}^2(L^2(\Omega)^3))}\leq \sup_{t\in \mathbb{R}} \|f(t)\|_{L^\infty(\Omega)^{3\times 3}}.$
The latter inequality implies the continuity of the embedding 
\[
  C_b(\R;L^\infty(\Omega)^{3\times 3})\ni f\mapsto T_f \in L_{\textnormal{sev}}^\textnormal{n}(L^2(\Omega)^3). 
\]
 Hence, (the multiplication operators associated to) $\eps_n$ converge in $L_{\textnormal{sev}}^\textnormal{n}(L^2(\Omega)^3)$ to $0$ . Thus, $\eps_n\to 0$ in $L_{\textnormal{sev}}^\textnormal{s}(L^2(\Omega)^3)$ as $n\to\infty$, by \eqref{eq:trivemb}. Hence, 
\[
   \begin{pmatrix}\eps_n & 0 \\ 0 & \mu \end{pmatrix} \to \begin{pmatrix}
                                                                              0 & 0 \\ 0 & \mu
                                                                            \end{pmatrix}\in L_{\textnormal{sev}}^\textnormal{s}(L^2(\Omega)^6)
\]
as $n\to\infty$ and by Theorem \ref{t:wpmax}, we get
\[
   \left(\begin{pmatrix}\eps_n & 0 \\ 0 & \mu\end{pmatrix},\begin{pmatrix}
                                                                              \sigma & 0 \\ 0 & 0
                                                                            \end{pmatrix}\right) \to
                                                                            \left(\begin{pmatrix}
                                                                            0 & 0 \\ 0 & \mu
                                                                            \end{pmatrix},\begin{pmatrix}
                                                                              \sigma & 0 \\ 0 & 0
                                                                            \end{pmatrix}\right)\in \textnormal{SP}^\textnormal{s}_{c,\nu,r}(L^2(\Omega)^6)
\]
as $n\to\infty$ for $r\coloneqq \sup_{n\in\N}\|\eps_n'\|_\infty$. Therefore, we infer \eqref{eq:eca} from Theorem \ref{t:cdpdes}.

For the proof of the estimate, employing Theorem \ref{t:cdpden} and using that 
\[
   \textnormal{sol}\left(\begin{pmatrix}\eps_n & 0 \\ 0 & \mu\end{pmatrix},\begin{pmatrix}
                                                                              \sigma & 0 \\ 0 & 0
                                                                            \end{pmatrix}\right) = \s S_{\textnormal{MAX}}(\eps_n,\mu,\sigma),
\]
we compute for all $\eta\geq\nu$
\begin{align*}
  & \| \big(\s S_{\textnormal{MAX}}(\eps_n,\mu,\sigma)\s S_{\textnormal{MAX}}(0,\mu,\sigma)^{-1}-1\big)\check{\partial}_{t,\nu}^{-1}\s S_{\textnormal{MAX}}(0,\mu,\sigma)\|_{L(L_{\eta}^2)}
    \\ & \leq \frac{1}{c} \Big(\left\|\Big(\begin{pmatrix}\eps_n & 0 \\ 0 & \mu\end{pmatrix}-\begin{pmatrix}0 & 0 \\ 0 & \mu\end{pmatrix}\Big)\s S_{\textnormal{MAX}}(0,\mu,\sigma)\right\|_{L(L_\eta^2)}
    \\ & \quad +\left\|\Big(\begin{pmatrix}\eps_n' & 0 \\ 0 & \mu'\end{pmatrix}-\begin{pmatrix}0 & 0 \\ 0 & \mu'\end{pmatrix}\Big)\partial_{t,\nu}^{-1}\s S_{\textnormal{MAX}}(0,\mu,\sigma)\right\|_{L(L_\eta^2)}
    \\ & \quad+\left\|\Big(\begin{pmatrix}\sigma & 0 \\ 0 & 0\end{pmatrix}-\begin{pmatrix}
      \sigma & 0 \\ 0 & 0
       \end{pmatrix}\Big)\partial_{t,\nu}^{-1}\s S_{\textnormal{MAX}}(0,\mu,\sigma)\right\|_{L(L_\eta^2)}\Big)
       \\ & \leq \frac{1}{c} \Big(\left\|\Big(\begin{pmatrix}\eps_n & 0 \\ 0 & \mu\end{pmatrix}-\begin{pmatrix}0 & 0 \\ 0 & \mu\end{pmatrix}\Big)\right\|_{L(L_\eta^2)} \frac{1}{c}
     +\left\|\Big(\begin{pmatrix}\eps_n' & 0 \\ 0 & 0\end{pmatrix}-\begin{pmatrix}0 & 0 \\ 0 & 0\end{pmatrix}\Big)\right\|\frac{1}{\eta}\frac{1}{c}\Big)
    \\ &\leq \frac{1}{c^2} \Big(\left\|\eps_n\right\|_{L(L_\eta^2)} +\left\|\eps_n'\right\|_{L(L_\eta^2)}\frac{1}{\eta}\Big)
    \\ &\leq \frac{1}{c^2} \Big(\left\|\eps_n\right\|_{\infty} +\left\|\eps_n'\right\|_\infty\frac{1}{\eta}\Big)
\end{align*}
where we also used that $\|\partial_{t,\nu}^{-1}\|_{L(L_\eta^2)}\leq 1/\eta$ and $\|\s S_{\textnormal{MAX}}(0,\mu,\sigma)\|_{L(L_\eta^2)}\leq 1/c$.
\end{Proof}

\begin{Remark} We remark here that in order that the solution operators of the originally given Maxwell's equation strongly converge to its respective eddy-current approximation where the dielectricity is formally set to $0$, one only needs to ensure the dielectricity to be uniformly small in the norm in $C_b$ and the respective derivatives being bounded.  
\end{Remark}

\renewcommand{\baselinestretch}{0.65}\normalsize\mysection{Dependency on Conductivity in Non-autonomous Thermodynamics}{The Continuous Dependence on the Conductivity in Non-autonomous Thermodynamics}{non-autonomous heat equation with rough coefficients $\cdot$ $\dive$ $\cdot$ $\grad$ $\cdot$ homogeneous Dirichlet boundary conditions $\cdot$ Theorem \ref{t:cdheat}}\label{s:cdcnt}

\renewcommand{\baselinestretch}{1}\normalsize
In this section, we present an application of Theorem \ref{t:cdpdes}. We focus on the heat equation with time dependent, non-symmetric and rough coefficients on general domains $\Omega\subseteq \mathbb{R}^d$. For the rest of this section, let $\Omega\subseteq\mathbb{R}^d$ be open.

We mention that there exists a deep theory for this type of equations with equally general coefficients in an $L^p$-type setting on $\mathbb{R}^d$ as underlying spatial domain, see \cite{Auscher2015} and the profound list of references therein. Concerning the solution theory, we do not claim any originality here, however, we stress the comparatively easy way of deriving the present well-posedness result. The equations to be discussed in this section read as
\begin{alignat}{1}\label{eq:heat}
\begin{aligned}
  \partial_t \theta + \dive q & = f 
  \\ q & = - a \grad \theta
\end{aligned}
\end{alignat}
on $\mathbb{R}\times \Omega$. The vector analytic operators $\dive$ and $\grad$ are computed with respect to the spatial variables $x\in \Omega$ only. The map $a\colon \mathbb{R}\times\Omega \to \mathbb{C}^{d\times d}$ models the thermal conductivity, $\theta$ and $q$ are the heat and the heat flux, respectively. The right-hand side $f\colon\mathbb{R}\times\Omega\to \mathbb{R}$ is a given external heat source. Subject to homogeneous Dirichlet boundary conditions to be satisfied by $\theta$ on the boundary of $\Omega$, we try to solve equation \eqref{eq:heat} for $(\theta,q)$. Having done so, we address the question of continuity of the solution operator in the thermal conductivity subject to an appropriate topology. As in the Sections \ref{s:DBF} and \ref{s:eca}, we will build up a proper functional analytic framework for \eqref{eq:heat} first and provide a solution theory for this problem. After that, we will apply Theorem \ref{t:cdpdes} in order to address the continuity in the thermal conductivity.

For the functional analytic set up, we introduce the vector analytic operators $\dive$ and $\grad$:

\begin{Definition}\label{d:gd} We set 
\begin{align*}
    \grad_c \colon C_c^1(\Omega)\subseteq L^2(\Omega)&\to L^2(\Omega)^d
    \\ \phi\mapsto (\partial_i\phi)_{i\in\{1,\ldots,d\}}.
\end{align*}
Define 
\[
   \grad_0 \coloneqq \overline{\grad}_c
\]
and set $\dive\coloneqq -\grad_0^*$.
\end{Definition}

Note that $\dom(\grad_0)=H_0^1(\Omega)$. Parallel to the Lemmas \ref{l:maxopsk} and \ref{l:maxa} the corresponding statements for the block operator
\begin{alignat}{1}\label{eq:dheato}
\begin{aligned}
   A\coloneqq \Big(\begin{pmatrix}
     0 & \dive\\ \grad_0 & 0
   \end{pmatrix} \colon &\dom(\grad_0)\times\dom(\dive)\\ &\quad\quad\subseteq L^2(\Omega)\times L^2(\Omega)^d\to L^2(\Omega)\times L^2(\Omega)^d\Big)
\end{aligned}
\end{alignat}
is true. We gather these assertions in one lemma, which we will state without proof as the reasoning follows the lines of the one in the Lemmas \ref{l:maxopsk} and \ref{l:maxa} upon replacing $\curl_0$ by $\grad_0$ and $\curl$ by $\dive$.

\begin{Lemma}\label{l:heata} Let $\nu\in (0,\infty)$, $\s A^\nu$ be the lift of $A$ given in \eqref{eq:dheato} to $L_\nu^2(\R;L^2(\Omega)^{d+1})$ with $\dom(\s A^\nu)=L_\nu^2(\mathbb{R};\dom(A))$. Then the following assertions are true.
 \begin{enumerate}[label=(\alph*)]
  \item\label{heatz} The operator $\s A^\nu$ is skew-selfadjoint.
  \item\label{heata} For all $\phi\in \dom(\s A^\nu)$, we have $\Re\langle Q_0\s A^\nu\phi,\phi\rangle=0$, where $Q_0$ is multiplication by $\1_{(-\infty,0)}$.
  \item\label{heatb} For all $\phi\in \dom(\s A^\nu)$, we have $\partial_{t,\nu}^{-1}\s A^\nu\phi=\s A^\nu\partial_{t,\nu}^{-1}\phi$.
\end{enumerate} 
\end{Lemma}

The assumptions on the conductivity $a$ are gathered in the following hypothesis.

\begin{Hypothesis}\label{h:cond} Let $a\colon \mathbb{R}\times\Omega \to \mathbb{C}^{d\times d}$ be bounded and measurable, $c>0$. Assume that
\begin{equation}\label{eq:cond}
   \Re \langle \xi,a(t,x)\xi\rangle \geq c\langle\xi,\xi\rangle \quad(\xi\in \mathbb{C}^d)
\end{equation} 
for almost every $x\in \Omega$ and $t\in\R$.
\end{Hypothesis}

For concluding the solution theory for \eqref{eq:heat}, we need the following preparatory step.

\begin{Lemma}\label{l:ainv} Assume Hypothesis \ref{h:cond} to be satisfied. Then for almost every $(t,x)\in \mathbb{R}\times\Omega$, we have
\[
   \Re \langle \zeta,a(t,x)^{-1}\zeta\rangle\geq \frac{c}{\|a\|^2_\infty}\langle \zeta,\zeta\rangle\quad(\zeta\in \mathbb{C}^d),
\] 
where $\|a\|_\infty\coloneqq \esssup_{(t,x)\in \mathbb{R}\times \Omega} \|a(t,x)\|$
\end{Lemma}
\begin{Proof}
  Let $(t,x)\in \mathbb{R}\times\Omega$ be such that \eqref{eq:cond} is satisfied for $a(t,x)$. Then, $a(t,x)$ is an invertible matrix, and $\|\zeta\|=\|a(t,x)a(t,x)^{-1}\zeta\|\leq \|a(t,x)\| \|a(t,x)^{-1}\zeta\|$ yields for $\zeta\in \mathbb{C}^d$ with $\xi\coloneqq a(t,x)^{-1}\zeta$:
  \begin{multline*}
     \Re\langle \zeta,a(t,x)^{-1}\zeta\rangle =\Re\langle a(t,x)\xi,\xi\rangle \\ \geq c \Re\langle \xi,\xi\rangle\geq \frac{c}{\|a(t,x)\|^2}\|a(t,x)^{-1}\xi\|^2=\frac{c}{\|a(t,x)\|^2}\|\zeta\|^2,
  \end{multline*}
 which implies the assertion.
\end{Proof}

We are now in the position to provide a solution theory for \eqref{eq:heat}.

\begin{Theorem}\label{t:wpheat} Let Hypothesis \ref{h:cond} be satisfied. Then 
\[
   \s S_{\textnormal{HEAT}}(a)\coloneqq \left( \partial_{t,\nu}\begin{pmatrix} 1 & 0 \\ 0 & 0\end{pmatrix} + \begin{pmatrix} 0 & 0 \\ 0 & a^{-1} \end{pmatrix} + \begin{pmatrix} 0 & \dive \\ \grad_0 & 0\end{pmatrix}\right)^{-1}
\]
is standard evolutionary at $\nu>0$. Moreover, 
\[
   \left(\begin{pmatrix} 1 & 0 \\ 0 & 0\end{pmatrix},\begin{pmatrix} 0 & 0 \\ 0 & a^{-1} \end{pmatrix}\right)\in \textnormal{SP}_{c',c',0}^\textnormal{s}(L^2(\Omega)^{d+1}),
\]
where $c'\coloneqq c/(\|a\|_\infty^2)$ and $\textnormal{SP}_{c',c',0}^\textnormal{s}$ is given in Definition \ref{d:solpde}. 
\end{Theorem}
\begin{Proof}
 By Lemma \ref{l:heata}, the operator $\s A^\nu(=\left(\begin{smallmatrix} 0 & \dive \\ \grad_0 & 0\end{smallmatrix}\right))$ defined in Lemma \ref{l:heata} meets the requirements imposed on $\s A$ in Hypothesis \ref{h:gSe}. Moreover, integration by parts, shows that $\Re\langle Q_t \partial_{t,\nu} \phi,\phi\rangle_{L_\eta^2} \geq \nu \langle Q_t \phi,\phi\rangle_{L^2_\eta}$ for all $\phi\in  C_c^1(\mathbb{R};L^2(\Omega))$ and $\eta\geq\nu>0$, $t\in \mathbb{R}$. Further, from Lemma \ref{l:ainv} it follows that 
 \[
   \Re \langle Q_t a^{-1} \phi,\phi\rangle =  \Re \langle  a^{-1} Q_t\phi,Q_t\phi\rangle \geq \frac{c}{\|a\|^2_\infty}\langle Q_t \phi,\phi\rangle.
 \]
 Hence, 
 \begin{align*}
    & \Re \langle Q_t \left(\partial_{t,\nu}\begin{pmatrix} 1 & 0 \\ 0 & 0\end{pmatrix} + \begin{pmatrix} 0 & 0 \\ 0 & a^{-1} \end{pmatrix}\right)\begin{pmatrix} \theta \\ q \end{pmatrix} ,\begin{pmatrix} \theta \\ q \end{pmatrix} \rangle_{L_\eta^2} 
    \\&\geq \langle \begin{pmatrix} \nu & 0 \\ 0 & \frac{c}{\|a\|_\infty^2}\end{pmatrix}\begin{pmatrix} \theta \\ q \end{pmatrix}, Q_t \begin{pmatrix} \theta \\ q \end{pmatrix}\rangle_{L_\eta^2} 
    \\&\geq \min\{ \nu,c'\} \langle Q_t \begin{pmatrix} \theta \\ q \end{pmatrix},\begin{pmatrix} \theta \\ q \end{pmatrix}\rangle_{L_\eta^2}\quad(\eta\geq\nu).
 \end{align*}
 Thus, all requirements in Hypothesis \ref{h:gSe} are warranted. So, Theorem \ref{t:gSe} applies and $\s S_{\textnormal{HEAT}}$ is evolutionary at $\nu$ and causal, or, equivalently, standard evolutionary at $\nu$. The remaining assertion follows from the estimates just derived together with the fact that $\s M'=0$ for $\s M=\begin{pmatrix} 1 &0 \\ 0& 0 \end{pmatrix}$.
\end{Proof}

The continuous dependence result on the conductivity is presented next. We recall that we again identified $a\in L^\infty(\R\times \Omega)^{d\times d}$ with its associated multiplication operators on $L_\nu^2(\R;L^2(\Omega)^d)$ for all $\nu\in\R$.

\begin{Theorem}\label{t:cdheat} Let $(a_k)_k$ be a sequence in $L^\infty(\R\times \Omega)^{d\times d}$ and assume that $b\coloneqq \lim_{k\to\infty} a_k$ exists in the strong operator topology of $L(L^2(\mathbb{R}\times\Omega))$. Assume there exists $c>0$ such that $a_k$ satisfies Hypothesis \ref{h:cond} with this $c$ for all $k\in\mathbb{N}$. 
 
 Then
 \[
    \s S_{\textnormal{HEAT}}(a_k) \to \s S_{\textnormal{HEAT}}(b) \text{ in }L_{\textnormal{sev}}^\textnormal{s}(L^2(\Omega)^{d+1})\text{ as }k\to\infty,
 \]
 where $\s S_{\textnormal{HEAT}}$ is given in Theorem \ref{t:wpheat}.
\end{Theorem}
\begin{Proof}
  By Theorem \ref{t:wpheat}, we may apply Theorem \ref{t:cdpdes}. For this, we have to ensure that 
  \[
     \begin{pmatrix}  0 & 0 \\ 0 & a_k^{-1} \end{pmatrix} \to \begin{pmatrix}  0 & 0 \\ 0 & b^{-1} \end{pmatrix} \in L_{\textnormal{sev}}^\textnormal{s}(L^2(\Omega)^{d+1})\text{ as }k\to\infty.
  \]
   But, this is the same as saying that $(a_k^{-1})_k$ is convergent to $b^{-1}$ in $L_{\textnormal{sev}}^\textnormal{s}(L^2(\Omega)^d)$. By Proposition \ref{p:stronmult}, we infer boundedness and convergence of the sequence $(a_k)_k$ to $b$ in $L_{\textnormal{sev}}^\textnormal{s}(L^2(\Omega)^d)$. Note that also $\Re\langle Q_tb\phi,\phi\rangle_{L_\nu^2}\geq c\langle Q_t\phi,\phi\rangle_{L_\nu^2}$ for all $\phi\in \s D(L^2(\Omega)^d)=\bigcap_{\eta\in \mathbb{R}}L^2_\eta(L^2(\Omega)^d)$. Thus, from $\|a_k^{-1}\|_{L(L^2_\nu)},\|b^{-1}\|_{L(L^2_\nu)}\leq 1/c$ for all $\nu\in\mathbb{R}$, see Corollary \ref{c:inv_ce} and Proposition \ref{p:posd}, we infer from Theorem \ref{t:cinv} that $a_k^{-1}\to b^{-1}$ in $L_{\textnormal{sev}}^\textnormal{s}$, which together with Theorem \ref{t:cdpdes} eventually proves the assertion.
\end{Proof}

\begin{Remark} We may roughly rephrase the contents of Theorem \ref{t:cdheat} as follows. Under the assumptions given in the theorem, let  $f\in L^2(0,\infty;L^2(\Omega))\subseteq \bigcap_{\eta>0}L_\eta^2(L^2(\Omega))$ and for $k\in \mathbb{N}$ let $\theta_k\in L_\nu^2(\R;L^2(\Omega))$ be a solution of 
\[
   \partial_{t,\nu} \theta_k -\dive a_k \grad_0 \theta_k = f,
\]
then $\theta_k \to u$ in $L_\nu^2(\R;L^2(\Omega))$ as $k\to\infty$, where $\theta$ satisfies
\[
   \partial_{t,\nu} u -\dive b \grad_0 u = f.
\]
The only thing, which has not been proved in Theorem \ref{t:cdheat} is that $u$ satisfies the equation asserted in the strong sense, that is, $u\in \dom(\partial_{t,\nu})\cap \dom(\dive b\grad_0)$. Hence, the term `roughly'.
\end{Remark}

\renewcommand{\baselinestretch}{0.65}\normalsize\mysection{Homogenization of Acoustic Wave Propagation in Bounded Domains}{On the Homogenization of Acoustic Wave Propagation in Bounded Domains}{the set $M(\alpha,\beta)$ $\cdot$ $G$-convergence $\cdot$ solution theory for elliptic type equations $\cdot$ relationship to the weak operator topology $\cdot$ Theorem \ref{t:wpwcdw}}\label{s:homaw}

\renewcommand{\baselinestretch}{1}\normalsize
The last application of the results developed concerns the weak operator topology. The motivation of this kind of problems stems from homogenization theory. The idea is to consider heterogeneous materials that have highly oscillatory material coefficients. The aim of homogenization theory is to determine the `effective' properties of the material by looking at the behavior of the solutions of the respective equations when the frequency of oscillations becomes infinitely large. To be more precise, on a bounded domain $\Omega\subseteq \mathbb{R}^d$ consider the wave equation formally given as follows (see also \cite[Example 5.3]{CioDon})
\begin{equation}\label{eq:wave}
  \partial_t^2 u - \dive a \grad u  = f
\end{equation}
on $\mathbb{R}\times \Omega$, where $a\colon \Omega\to \mathbb{C}^{d\times d}_{\textnormal{pd}}$ (taking values in the symmetric, positive definite $d$ by $d$ matrices) describes the material properties, $f\colon \R\times \Omega\to \mathbb{C}$ is a given source term, $u\colon\mathbb{R}\times\Omega\to \mathbb{C}$ describes the unknown wave propagation subject to homogeneous Dirichlet boundary conditions. In the theory of homogenization, for $\eps>0$, one is interested in the problem 
\begin{equation}\label{eq:wave2}
  \partial_t^2 u_\eps - \dive a_\eps \grad u_\eps  = f
\end{equation}
often with $a_\eps(x)\coloneqq a(x/\eps)$ for $x\in \Omega$ and addresses the question, whether $(u_\eps)_\eps$ converges as $\epsilon\to 0$ and, if so, whether the respective limit solves an equation of 'similar' type as in \eqref{eq:wave}. There exists a vast literature on homogenization theory, we only refer to \cite{BenLiPap}, \cite{CioDon} and \cite{TarIntro} to mention a few. In order to tackle problems in homogenization theory Spagnolo \cite{Spagnolo1967,Spagnolo1968a} introduced the concept of $G$-convergence. We recall this concept here. For this, we use the vector analytic operators given in Definition \ref{d:gd}. The literature also gives an account on what conditions are needed for $a$ such that $a(\cdot/\eps)$ is `$G$-convergent' and provides formulas for the respective limits as well as quantitative estimates for the rate of convergence related to $\eps$.

\begin{Definition}[{see e.g.~{\cite[Definition 13.1]{CioDon}}}]\label{d:gc} Let $0<\alpha<\beta$, and let $(a_\eps)_{\eps>0}$ in 
\begin{multline*}
   M(\alpha,\beta)\coloneqq\\ \{ a\in L^\infty(\Omega)^{d\times d}; \|a\|_\infty\leq \beta, \langle a(x)\xi,\xi\rangle_{\mathbb{C}^d}\geq \alpha \|\xi\|^2\;(\xi\in \mathbb{C}^d, \text{a.e.~}x\in \Omega)\}.
\end{multline*}
Then $(a_\eps)_\eps$ is said to \emph{$G$-converge} to $b\in M(\alpha,\beta)$ as $\eps\to 0$, if for all $f\in H^{-1}(\Omega)=(H_0^1(\Omega))^*$ the solution $u_\eps\in \dom(\grad_0)$ of
\[
   -\dive a_\eps \grad_0 u_\eps = f
\]
is such that $(u_\eps)_\eps$ converges weakly in $H_0^1(\Omega)$ to $u_0$ satisfying
\[ 
    -\dive b \grad_0 u_0 = f.
\]
\end{Definition}

In order to put the latter definition into perspective of the results developed so far, we insert a short interlude on the solvability of elliptic problems, in particular to those mentioned in Definition \ref{d:gc}. Although a solution theory for equations discussed in Definition \ref{d:gc} is well-known, we like to point out a slightly more abstract point of view, which has proved useful for general (non-linear) elliptic type problems in \cite{Waurick2014MN_Ellip}. It rests on the following observations. Let $\s X,\s Y$ be Hilbert spaces and let $S\colon \s X\to \s Y$ be a continuous bijection. Then, by the closed graph theorem, $S^{-1}$ is continuous as well. Hence, $S$ is Banach space isomorphism from $\s X$ to $\s Y$. Thus, we change the scalar product in $\s Y$ to be $\s Y\times \s Y\ni (\phi,\psi)\mapsto \langle S^{-1}\phi,S^{-1}\psi\rangle_{\s X}$ resulting in a scalar product on $\s Y$, which is equivalent to the original one. We call $\s Y_S$ the Hilbert space endowed with this modified scalar product. It is easy to see that $S\colon \s X\to \s Y_S$ is unitary. Thus, the dual operator $S'\colon \s Y^*\to \s X^*$ is a Banach space isomorphism. By identifying $\s Y^*=\s Y$ via the unitary Riesz isomorphism $R\colon \s Y^*\to \s Y$, we infer that the (modified) dual
\[
   S^\diamond \colon \s Y \to \s X^*, \phi \mapsto S'R\phi
\] 
is a Banach space isomorphism again. 

We apply this rationale to (a modification) of the vector-analytic operators introduced in Definition \ref{d:gd}. For this, observe that, as $\Omega$ is bounded we have a Poincar\'e inequality, that is, we find $c>0$ such that
\[
   \|u\|_{L^2(\Omega)}\leq c\|\grad_0 u\|_{L^2(\Omega)^d} \quad(u\in H_0^1(\Omega)).
\]
The latter ensures two-fold, on the one hand $\grad_0$ is injective and on the other hand, by the closedness of $\grad_0$, the range of $\grad_0$ is a closed subspace of $L^2(\Omega)^d$. We denote $\pi\colon L^2(\Omega)^d\to \ran(\grad_0)$ the orthogonal projection onto $\ran(\grad_0)$. Applying the reasoning just developed to $\s X=H_0^1(\Omega)$ and $\s Y=\ran(\grad_0)$, we get that
\begin{equation}\label{eq:guni}
   \pi\grad_0 \colon H_0^1(\Omega) \to \s Y_{\pi\grad_0}, u\mapsto \pi\grad_0 u\text{ is unitary.}
\end{equation}
Moreover, using $H^{-1}(\Omega)=H_0^1(\Omega)^*$, we obtain,
\begin{equation}\label{eq:duni}
   -\dive\pi^*=(\pi\grad_0)^\diamond \colon \s Y \to H^{-1}(\Omega)\text{ is a Banach space isomorphism,}
\end{equation}
where $-\dive$ is the distributional divergence on $L^2(\Omega)^d$ with values in $H^{-1}(\Omega)$. Having made these observations, a solution theory for divergence form problems as discussed in Definition \ref{d:gc} is almost immediate:

\begin{Theorem}[{{\cite[Theorem 3.1.1]{Waurick2014MN_Ellip}}}]\label{t:wpet} Let $0<\alpha<\beta$ and let $M(\alpha,\beta)$ as in Definition \ref{d:gc}. For all $a\in M(\alpha,\beta)$ the mapping
\[
   -\dive a \grad_0 \colon H_0^1(\Omega)\to H^{-1}(\Omega)
\]
is continuously invertible with inverse given by
\[
    (-\dive a \grad_0)^{-1} = (\pi\grad_0)^{-1} (\pi a \pi^*)^{-1} (-\dive \pi^*)^{-1},
\]
where $\pi\grad_0$ and $-\dive \pi$ are given in \eqref{eq:guni} and \eqref{eq:duni}, respectively.
\end{Theorem}
\begin{Proof}
 First of all note that $\grad_0^\diamond = -\dive\colon L^2(\Omega)^d\to H^{-1}(\Omega)$. Moreover, recall that $\ran(\grad_0)^\bot=\kar(\dive)$. Hence, 
 \begin{multline*}
    \dive = \dive (\pi^*\pi+(1-\pi^*\pi))=\dive \pi^*\pi\text{ and }\\ \grad_0 = (\pi^*\pi+(1-\pi^*\pi))\grad_0 = \pi^*\pi \grad_0.
 \end{multline*}
 So,
 \[
    -\dive a \grad_0 = -\dive\pi^*\pi a\pi^*\pi \grad_0 = (-\dive\pi^*)(\pi a\pi^*)(\pi \grad_0).
 \]
 From, \eqref{eq:guni} and \eqref{eq:duni}, we read off that $\pi \grad_0$ and $-\dive\pi^*$ considered with the space $\ran(\grad_0)$ replacing $\s Y$ are Banach space isomorphisms. Thus, the assertion follows once we show that $\pi a \pi^*$ is continuously invertible in $\ran(\grad_0)$. But, using that $a$ is bounded together with the positive definiteness condition mentioned in the definition of $M(\alpha,\beta)$, we infer
 \begin{multline*}
    \langle \pi a \pi^* \phi,\phi\rangle_{\ran(\grad_0)} =\langle  a \pi^* \phi,\pi^*\phi\rangle_{L^2(\Omega)^n}\\ \geq \alpha\langle \pi^* \phi,\pi^*\phi\rangle_{L^2(\Omega)^n}=\alpha\langle\phi,\phi\rangle_{\ran(\grad_0)}\quad(\phi\in \ran(\grad_0)).
 \end{multline*}Hence, as $\pi a \pi^*$ is continuous in $\ran(\grad_0)$, the continuous invertibility, thus, follows from Proposition \ref{p:posd}.
\end{Proof}

Next, we relate Theorem \ref{t:wpet} to the concept of $G$-convergence:

\begin{Theorem}\label{t:gcwo} Let $0<\alpha<\beta$, $(a_\eps)_\eps$ in $M(\alpha,\beta)$ (see Definition \ref{d:gc}), $b\in M(\alpha,\beta)$. Then the following conditions are equivalent:
\begin{enumerate}[label=(\roman*)]
 \item\label{gc} $(a_\eps)_\eps$ $G$-converges to $b$;
 \item\label{w0} $(\dive a_\eps \grad_0)^{-1} \to (\dive b \grad_0)^{-1}$ in $L(H^{-1}(\Omega);H_0^1(\Omega))$ with respect to the weak operator topology;
 \item\label{wo} $((\pi a_\eps \pi^*)^{-1})_\eps$ converges in the weak operator topology of $L(\ran(\grad_0))$ to the mapping $(\pi b \pi^*)^{-1}$, where $\pi\colon L^2(\Omega)^d\to \ran(\grad_0)$ is the orthogonal projection onto $\ran(\grad_0)$.
\end{enumerate} 
\end{Theorem}
\begin{Proof}
  The conditions \ref{gc} and \ref{w0} are trivial reformulations of one another. Theorem \ref{t:wpet} together with \eqref{eq:guni} and \eqref{eq:duni} imply that \ref{w0} is true if and only if $((\pi a_\eps \pi^*)^{-1})_\eps$ converges to $(\pi b \pi^*)^{-1}$ in the weak operator topology of the space $L(\s Y_{\pi\grad_0})$. The latter, in turn, is equivalent to \ref{wo}.
\end{Proof}

Before being able to apply our result concerning the weak operator topology to the homogenization type problem in \eqref{eq:wave2}, we need to warrant the compactness condition in Hypothesis \ref{h:awpde}. For this, the following observation comes in handy. The proof of which stems from \cite{Waurick2013AA_Hom}.

\begin{Theorem}[{{\cite[Lemma 4.1]{Waurick2013AA_Hom}}}]\label{t:cc} Let $\s X,\s Y$ Hilbert spaces, $S:\dom(S)\subseteq \s X\to \s Y$ densely defined, closed  and assume that $(\dom(S),\|\cdot\|_S)\hookrightarrow \s X$ is compact. Then $(\dom(S^*)\cap \kar(S^*)^{\bot},\|\cdot\|_{S^*}) \hookrightarrow \s Y$ is compact. 
\end{Theorem}
\begin{Proof} 
 We use the polar decomposition \cite[p 334]{Kato1980} for densely defined, closed, linear operators. We have
 \begin{equation}\label{PoCo}
    S= U|S|,
 \end{equation}
 where $U: \overline{\ran(S^*)}\to \s X$ is a linear isometry from $\overline{\ran(S^*)}$ to $\overline{\ran(S)}$. Note that by equation \eqref{PoCo}, we see that 
 \[
    V: \overline{\ran(S^*)} \to \overline{\ran(S)}: x\mapsto Ux
 \]
 is a linear isometry with dense range. Thus, $V$ is unitary. Furthermore, we have $V^{-1}x=U^*x$ for all $x\in\kar(S^*)^\bot=\overline{\ran(S)}$. 
 
 Let $(x_n)_n$ be a bounded sequence in $(\dom(S^*)\cap \kar(S^*)^{\bot},\|\cdot\|_{S^*})$. Adjoining equation \eqref{PoCo} yields that $(U^*x_n)_n$ is a bounded sequence in $\dom(|S|)$. Since $\dom(|S|)=\dom(S)$, we may choose a convergent subsequence of $(U^*x_n)_n$, for which we use the same notation. Since $V$ is unitary, we have that $(VU^*x_n)_n$ also strongly converges and, thus, so does
 \[
    (x_n)_n =  (VV^{-1}x_n)_n = (VU^*x_n)_n.
 \]
\end{Proof}

In order to put Theorem \ref{t:cdpdew} into the perspective of homogenization theory, we reformulate equation \eqref{eq:wave} according to our general setting of evolutionary equations developed in this exposition. So, assuming that $a\in M(\alpha,\beta)$ for some $0<\alpha<\beta$ and recalling \eqref{eq:wave}
\begin{equation}\label{eq:wavem}
  \partial_t^2 u - \dive a \grad_0 u  = f,
\end{equation}
we set $v\coloneqq \partial_t u$ and $p \coloneqq \pi a \pi^*\pi \grad_0$ with $\pi\colon L^2(\Omega)^d\to \ran(\grad_0)$ being the orthogonal projection. Hence, we obtain
\begin{equation}\label{eq:waves}
  \left( \partial_t \begin{pmatrix}
             1 & 0 \\ 0 & (\pi a \pi^*)^{-1} 
         \end{pmatrix} - \begin{pmatrix}
             0 & \dive\pi^* \\ \pi\grad_0 & 0
         \end{pmatrix}\right)\begin{pmatrix} v \\ p \end{pmatrix} =\begin{pmatrix} f \\ 0 \end{pmatrix}.
\end{equation}
In the reformulated fashion, we can now apply Theorem \ref{t:gSe}, where we employ the custom of not keeping track of the particular $L_\nu^2$ the respective operators are realized in as well as of identifying elements $a\in M(\alpha,\beta)$ with their corresponding multiplication operators.

\begin{Theorem}\label{t:wpwave} Let $0<\alpha<\beta$, $a\in M(\alpha,\beta)$ (see Definition \ref{d:gc}); let $\pi\colon L^2(\Omega)^d\to \ran(\grad_0)$ be the orthogonal projection. Then
  \begin{equation}\label{eq:wpwave}
  \s S_{\textnormal{WAVE}}(a) \coloneqq \left( \partial_{t,\nu} \begin{pmatrix}
             1 & 0 \\ 0 & (\pi a \pi^*)^{-1} 
         \end{pmatrix} - \begin{pmatrix}
             0 & \dive\pi^* \\ \pi\grad_0 & 0
         \end{pmatrix}\right)^{-1}
 \end{equation}
 is standard evolutionary at $\nu>0$. Moreover, 
 \[
    \begin{pmatrix}
             1 & 0 \\ 0 & (\pi a \pi^*)^{-1} 
         \end{pmatrix} \in \textnormal{SP}_{c,\nu,0}^\textnormal{w}(L^2(\Omega)\times \ran(\grad_0)),
 \]
 for some $c>0$, where $\textnormal{SP}_{c,\nu,0}^\textnormal{w}$ is given in Definition \ref{d:solpdew}.
\end{Theorem}
\begin{Proof}
Employing the same rationale (see e.g.~Lemma \ref{l:heata}) of the previous two sections, we infer that (the lifting to $L_\nu^2$-functions of) $\begin{pmatrix} 0 & \dive\pi^* \\ \pi\grad_0 & 0 \end{pmatrix}$ commutes with $Q_t$, $\partial_{t,\nu}^{-1}$ ($\nu>0$) and is skew-selfadjoint as $\dive\pi^*=-(\pi\grad_0)^*$. Moreover, the strict positive definiteness of $\pi a\pi^*$ in $\ran(\grad_0)$ implies a similar one for $(\pi a\pi^*)^{-1}$ in $\ran(\grad_0)$. Hence, with $\s M=\left(\begin{smallmatrix}1 & 0 \\ 0 & (\pi a \pi^*)^{-1} 
         \end{smallmatrix}\right)$ the corresponding commutator $\s M'$ of $\s M$ with $\partial_{t,\nu}$ vanishes for all $\nu\in\R$, we infer with the help of integration by parts that there exists $c>0$ such that
         \[
            \Re \langle Q_t \partial_{t,\nu} \begin{pmatrix}
             1 & 0 \\ 0 & (\pi a \pi^*)^{-1} 
         \end{pmatrix} \begin{pmatrix} v \\ p \end{pmatrix},\begin{pmatrix} v \\ p \end{pmatrix}\rangle_{L_\eta^2} \geq c\langle Q_t \begin{pmatrix} v \\ p \end{pmatrix},\begin{pmatrix} v \\ p \end{pmatrix}\rangle_{L_\eta^2}
         \]
for all $\eta\geq\nu>0$, $v,p$ being compactly supported and continuously differentiable with values in $L^2(\Omega)$ and $\ran(\grad_0))$, respectively.
        Hence, Theorem \ref{t:gSe} applies and yields the assertion.
\end{Proof}

Finally, we come to the announced result on the homogenization of \eqref{eq:wave2}.

\begin{Theorem}\label{t:wpwcdw} Let $0<\alpha<\beta$, $(a_\eps)_\eps$ in $M(\alpha,\beta)$. Assume that $(a_\eps)_\eps$ $G$-converges to $b\in M(\alpha,\beta)$ (see Definition \ref{d:gc}). Then 
\[
   \s S_{\textnormal{WAVE}}(a_\eps) \to \s S_{\textnormal{WAVE}}(b)\text{ in }L_{\textnormal{sev}}^\textnormal{w}(L^2(\Omega)\times \ran(\grad_0))\text{ as }\eps\to 0,
\]
where $\s S_{\textnormal{WAVE}}$ is given in Theorem \ref{t:wpwave}.
\end{Theorem}
\begin{Proof}
  We want to apply Theorem \ref{t:cdpdew}. For this, we recall that if $(T_\eps)_\eps$ converges in the weak operator topology to $T$ of $L(\s X)$ for some Hilbert space $\s X$, then the corresponding liftings of $(T_\eps)_\eps$ to $L(L_\nu^2(\R;\s X))$ converge to the corresponding lifting of $T$ in the weak operator topology of $L(L_\nu^2(\R;\s X))$, see also Example \ref{ex:el}(c). Hence, by Theorem \ref{t:gcwo},
  \[
   \begin{pmatrix}
             1 & 0 \\ 0 & (\pi a_\eps \pi^*)^{-1} 
         \end{pmatrix} \to  \begin{pmatrix}
             1 & 0 \\ 0 & (\pi b \pi^*)^{-1} 
         \end{pmatrix}\text{ in }L_{\textnormal{sev},\nu}^\textnormal{w}\text{ as }\eps\to 0,
  \]
  for all $\nu\in \mathbb{R}$. The assertion, thus, follows from Theorem \ref{t:cdpdew} once we show the compactness condition in Hypothesis \ref{h:awpde} for 
  \begin{multline*}
   \s A^\nu \coloneqq \Big(\begin{pmatrix}
             0 & \dive\pi^* \\ \pi\grad_0 & 0
         \end{pmatrix} \colon L_\nu^2(\dom(\pi\grad_0)\times\dom(\dive\pi^*))
         \\ \subseteq L_\nu^2(L^2(\Omega)\times \ran(\grad_0)) \to L_\nu^2(L^2(\Omega)\times \ran(\grad_0))\Big).
  \end{multline*}
  For this, however, it is sufficient to prove that 
  \[
     \s Y\coloneqq \dom(\pi\grad_0)\times\dom(\dive\pi^*) \hookrightarrow\hookrightarrow L^2(\Omega)\times \ran(\grad_0),
  \]
  where the domains on the left-hand side are endowed with the respective graph norms. But, by the Theorem of Rellich-Kondrachov, we get $\dom(\pi\grad_0) \hookrightarrow\hookrightarrow L^2(\Omega)$. Hence, with $S=\pi\grad_0$, the assertion follows from Theorem \ref{t:cc}.
\end{Proof}

\begin{Remark} The convergence of solution operators $\s S_{\textnormal{WAVE}}(a_\eps)$ proved in Theorem \ref{t:wpwcdw} implies the convergence of the solution operators in the weak operator topology of $L_\nu^2(\mathbb{R};L^2(\Omega)\times \ran(\grad_0))$ for some $\nu\in \mathbb{R}$. The latter in turn is equivalent to the notion of $G$-convergence introduced by Spagnolo in the first place, see also \cite[p 74]{Gcon1}. The novelty element of Theorem \ref{t:wpwcdw} lies in the fact that the solution operators converge in the weak operator topology of evolutionary mappings. Thus, one might phrase the result of Theorem \ref{t:wpwcdw} as ``evolutionary $G$-convergence'' of the solution operators of the wave equation.
\end{Remark}

\section{Comments} 

We want to give a brief account on the study of the continuous dependence of solutions to partial differential equations on the coefficients. There are some results for particular equations or with both stronger and weaker topologies the coefficients are considered in. The focus in the available literature is on non-linear equations. In \cite{Yaman2011} a particular non-linear equation is considered and the continuous dependence of the solution on some scalar factors is addressed. The so-called Brinkman--Forchheimer equation is discussed with regards to continuous dependence on some bounded functions under the sup-norm in \cite{cCelebi2006,Franchi2003,Tu2007,Payne1999,Liu2009}. The local sup-norm has been considered in \cite{Gianni2009}, where the continuous dependence on the (non-linear) constitutive relations for particular equations of fluid flow in porous media is discussed. A weak topology for the coefficients is considered in \cite{Kim2009}. However, the partial differential equations considered are of a specific form and the underlying spatial domain is the real line. Dealing with time dependent coefficients in a boundary value problem of parabolic type, the author of \cite{Penning1991} shows continuous dependence of the associated evolution families on the coefficients. In \cite{Penning1991}, the coefficients are certain functions considered with the $C^1$-norm. The author of \cite{Tudor1976} studies the continuous dependence of diffusion processes under the $C^0$-norm of the coefficients. Also with regards to strong topologies, the authors of \cite{Kunze2011,Kunze2012} studied continuous dependence results for a class of stochastic partial differential equations. All in all, the results developed in this exposition merely complement the research on the continuous dependence on the coefficients. A main reason being that the class focused on in this exposition is rather general and, thus, applies to many particular settings. In turn, the advantage of being applicable to many equations at the same time inherits the drawback of being certainly not optimal, if one restricts oneself to a specific class of both equations and coefficients. A particular instant for this fact can be found in \cite{Cherednichenko2015}.

With regards to homogenization theory, we shall particularly refer to \cite{TarIntro,CioDon,BenLiPap}, where the continuous dependence of the coefficients has been addressed in the particular situation of homogenization problems. For the concept of $H$-convergence generalizing $G$-convergence we also refer to \cite{Murat1978,Mur1}.

The present exposition does not treat the case of Maxwell's equations with regards to the weak operator topology. We refer to \cite{Waurick2012_HowFar} instead for the rather involved treatment. Note that Maxwell's equations have the particular property that memory effects (as in the case of ordinary differential equations) are likely to occur due to the homogenization procedure. This has also been observed in \cite{BarStr,Well}. The reason being that -- in contrast to the example treated in Section \ref{s:homaw} -- both the spatial derivative operators ($\curl$ and $\curl_0$) have an infinite-dimensional nullspace. Hence, a projection technique similar to the one in Section \ref{s:homaw} applied to Maxwell's equation leads to a system of a partial differential equation and an ordinary differential equation with \emph{infinite-dimensional} state space, see \cite{Waurick2012_HowFar} for more details.

There is still plenty of research to be addressed in the future: the weak operator topology and Maxwell's equations have only been dealt with for time-shift invariant coefficients in \cite{Waurick2012_HowFar}. The results concerning the norm and the strong operator topology being entirely new, there is certainly room for optimizing the respective results. Regarding numerical treatments of evolutionary equations discussed here, a natural way of complementing the present results is to ask for continuity in the unbounded operator $\s A$ in the solution operator induced by $(\partial_t\s M+\s A)u=f$: Does for instance strong resolvent convergence of a sequence $(\s A_n)_n$ imply convergence of the corresponding solution operators $\s S_n = (\partial_t\s M+\s A_n)^{-1}$? In fact, this would give way for the study of Galerkin approximations and hence form a part of numerical analysis for evolutionary equations.

Apart from Theorem \ref{t:cdpden}, all results on the continuity of the solution operator in the coefficients were non-quantitative. Future research may thus be concerned with convergence rate or moduli of continuity of the mappings considered.

 \bibliographystyle{alpha}
\newcommand{\etalchar}[1]{$^{#1}$}

%
%
%
%
\end{document}